\documentclass[11pt]{amsart}

\usepackage{marginnote}
\usepackage{anysize} 
\usepackage{comment}
\marginsize{2cm}{2cm}{2.3cm}{2.3cm}
\usepackage{listings}
\usepackage{xcolor}
\usepackage{amsmath}
\usepackage{mathtools}
\usepackage{booktabs}
\usepackage[all]{xy}
\usepackage[utf8]{inputenc}
\usepackage{varioref}
\usepackage{upgreek}
\usepackage{amsfonts}
\usepackage{amssymb}
\usepackage{bbm}
\usepackage{esint}
\usepackage{graphicx}
\usepackage{subcaption}
\usepackage{tikz}
\usepackage{empheq}
\usepackage{enumitem}
\usepackage{tikz-cd}
\usepackage{todonotes}
\usetikzlibrary{matrix,arrows,decorations.pathmorphing}
\usepackage{mathrsfs}
\usepackage[hypertexnames=false,backref=page,pdftex,
 	pdfpagemode=UseNone,
 	breaklinks=true,
 	extension=pdf,
 	colorlinks=true,
 	linkcolor=blue,
 	citecolor=red,
 	urlcolor=blue,
 ]{hyperref}

\newcommand{\C}{{\mathbb C}}

\newcommand{\tr}{\operatorname{tr}}
\newcommand{\Res}{\operatorname{Res}}
\newcommand{\Hom}{\operatorname{Hom}}

\newcommand{\diag}{\operatorname{diag}}

\newcommand{\sdim}{\operatorname{sdim}}
\newcommand{\fsdim}{\operatorname{sdim}_{\mathrm f}}
\newcommand{\str}{\operatorname{str}}
\newcommand{\DS}{\operatorname{DS}}

\newcommand{\id}{\operatorname{id}}

\renewcommand{\Im}{\operatorname{im}}
\newcommand{\Ker}{\operatorname{ker}}

\newcommand{\End}{\operatorname{End}}
\newcommand{\ad}{\operatorname{ad}}
\newcommand{\ch}{\operatorname{ch}}
\newcommand{\fch}{\operatorname{ch}}

\def\im{{\operatorname{im}}}

\def\calH{{\mathscr{H}}}

\def\defect{{\mathrm{def}}}

\def\calZ{{\mathcal{Z}}}
\def\sp{{\mathfrak{sp}}}

\def\calG{{\mathcal{G}}}

\def\dd{{\mathfrak{d}}}
\def\so{{\mathfrak{so}}}

\def\sch{{\operatorname{sch}}}

\def\gr{{\mathrm{gr}}}

\def\CC{{\mathbb C}}
\def\RR{{\mathbb R}}

\def\su{{\mathfrak{su}}} 

\def\psl{{\mathfrak{psl}}}
\def\rk{{\operatorname{rk}}}

\def\deg{{\operatorname{deg}}}
\def\sl{{\mathfrak{sl}}}

\def\HH{{\mathcal{H}}}

\def\osp{{\mathfrak{osp}}}

\def\hh{{\mathfrak h}}

\def\pp{{\mathfrak p}}
\def\uu{{\mathfrak u}}

\def\kk{{\mathfrak k}}

\def\ll{{\mathfrak l}}
\def\YY{{\mathcal Y}}

\def\qq{{\mathfrak q}}
\def\aa{{\mathfrak a}}

\def\bb{{\mathfrak b}}
\def\even{{\mathfrak{g}_{\bar{0}}}}
\def\odd{{\mathfrak{g}_{\bar{1}}}}
\def\tt{{\mathfrak t}}

\def\nn {{\mathfrak{n}}}
\def\gg{{\mathfrak{g}}}

\def\ZZ{{\mathbb Z}}

\def\gsmod{{\mathfrak{g}}\textbf{-smod}}

\def\at{{\operatorname{at}}}
\def\slmn{{\mathfrak{sl}(m\vert n)}}

\def\glmn{{\mathfrak{gl}(m\vert n)}} 
\def\UE{{\mathfrak U}}
 
\def\calC{{\mathcal C}}

\newcommand{\gl}{{\mathfrak{gl}}}

\newcommand{\Lie}{\operatorname{Lie}}
\newcommand{\Ad}{\operatorname{Ad}}
\newcommand{\hol}{\operatorname{hol}}

\newcommand{\abs}[1]{\left|{#1}\right|}
\newcommand{\norm}[1]{\lVert#1\rVert}

\makeatletter
\newcommand*{\da@rightarrow}{\mathchar"0\hexnumber@\symAMSa 4B }
\newcommand*{\da@leftarrow}{\mathchar"0\hexnumber@\symAMSa 4C }
\newcommand*{\xdashrightarrow}[2][]{%
  \mathrel{%
    \mathpalette{\da@xarrow{#1}{#2}{}\da@rightarrow{\,}{}}{}%
  }%
}
\newcommand{\xdashleftarrow}[2][]{%
  \mathrel{%
    \mathpalette{\da@xarrow{#1}{#2}\da@leftarrow{}{}{\,}}{}%
  }%
}
\newcommand*{\da@xarrow}[7]{%
  \sbox0{$\ifx#7\scriptstyle\scriptscriptstyle\else\scriptstyle\fi#5#1#6\m@th$}%
  \sbox2{$\ifx#7\scriptstyle\scriptscriptstyle\else\scriptstyle\fi#5#2#6\m@th$}%
  \sbox4{$#7\dabar@\m@th$}%
  \dimen@=\wd0 %
  \ifdim\wd2 >\dimen@
    \dimen@=\wd2 %   
  \fi
  \count@=2 %
  \def\da@bars{\dabar@\dabar@}%
  \@whiledim\count@\wd4<\dimen@\do{%
    \advance\count@\@ne
    \expandafter\def\expandafter\da@bars\expandafter{%
      \da@bars
      \dabar@ 
    }%
  }%  
  \mathrel{#3}%
  \mathrel{%   
    \mathop{\da@bars}\limits
    \ifx\\#1\\%
    \else
      _{\copy0}%
    \fi
    \ifx\\#2\\%
    \else
      ^{\copy2}%
    \fi
  }%   
  \mathrel{#4}%
}
\makeatother

\makeatletter
\newsavebox\myboxA
\newsavebox\myboxB
\newlength\mylenA

\newcommand*\xtilde[2][0.8]{%
    \sbox{\myboxA}{$\m@th#2$}%
    \setbox\myboxB\null% Phantom box
    \ht\myboxB=\ht\myboxA%
    \dp\myboxB=\dp\myboxA%
    \wd\myboxB=#1\wd\myboxA% Scale phantom
    \sbox\myboxB{$\m@th\widetilde{\copy\myboxB}$}%  Overlined phantom
    \setlength\mylenA{\the\wd\myboxA}%   calc width diff
    \addtolength\mylenA{-\the\wd\myboxB}%
    \ifdim\wd\myboxB<\wd\myboxA%
       \rlap{\hskip 0.5\mylenA\usebox\myboxB}{\usebox\myboxA}%
    \else
        \hskip -0.5\mylenA\rlap{\usebox\myboxA}{\hskip 0.5\mylenA\usebox\myboxB}%
    \fi}
\newbox\usefulbox

\def\getslant #1{\strip@pt\fontdimen1 #1}

\def\xxtilde #1{\mathchoice
 {{\setbox\usefulbox=\hbox{$\m@th\displaystyle #1$}%
    \dimen@ \getslant\the\textfont\symletters \ht\usefulbox
    \divide\dimen@ \tw@ 
    \kern\dimen@ 
    \xtilde{\kern-\dimen@ \box\usefulbox\kern\dimen@ }\kern-\dimen@ }}
 {{\setbox\usefulbox=\hbox{$\m@th\textstyle #1$}%
    \dimen@ \getslant\the\textfont\symletters \ht\usefulbox
    \divide\dimen@ \tw@ 
    \kern\dimen@ 
    \xtilde{\kern-\dimen@ \box\usefulbox\kern\dimen@ }\kern-\dimen@ }}
 {{\setbox\usefulbox=\hbox{$\m@th\scriptstyle #1$}%
    \dimen@ \getslant\the\scriptfont\symletters \ht\usefulbox
    \divide\dimen@ \tw@ 
    \kern\dimen@ 
    \xtilde{\kern-\dimen@ \box\usefulbox\kern\dimen@ }\kern-\dimen@ }}
 {{\setbox\usefulbox=\hbox{$\m@th\scriptscriptstyle #1$}%
    \dimen@ \getslant\the\scriptscriptfont\symletters \ht\usefulbox
    \divide\dimen@ \tw@ 
    \kern\dimen@ 
    \xtilde{\kern-\dimen@ \box\usefulbox\kern\dimen@ }\kern-\dimen@ }}%
 {}}

\newcommand*\xoverline[2][0.75]{%
    \sbox{\myboxA}{$\m@th#2$}%
    \setbox\myboxB\null% Phantom box
    \ht\myboxB=\ht\myboxA%
    \dp\myboxB=\dp\myboxA%
    \wd\myboxB=#1\wd\myboxA% Scale phantom
    \sbox\myboxB{$\m@th\overline{\copy\myboxB}$}%  Overlined phantom
    \setlength\mylenA{\the\wd\myboxA}%   calc width diff
    \addtolength\mylenA{-\the\wd\myboxB}%
    \ifdim\wd\myboxB<\wd\myboxA%
       \rlap{\hskip 0.5\mylenA\usebox\myboxB}{\usebox\myboxA}%
    \else
        \hskip -0.5\mylenA\rlap{\usebox\myboxA}{\hskip 0.5\mylenA\usebox\myboxB}%
    \fi}

\def\xxoverline #1{\mathchoice
 {{\setbox\usefulbox=\hbox{$\m@th\displaystyle #1$}%
    \dimen@ \getslant\the\textfont\symletters \ht\usefulbox
    \divide\dimen@ \tw@ 
    \kern\dimen@ 
    \overline{\kern-\dimen@ \box\usefulbox\kern\dimen@ }\kern-\dimen@ }}
 {{\setbox\usefulbox=\hbox{$\m@th\textstyle #1$}%
    \dimen@ \getslant\the\textfont\symletters \ht\usefulbox
    \divide\dimen@ \tw@ 
    \kern\dimen@ 
    \xoverline{\kern-\dimen@ \box\usefulbox\kern\dimen@ }\kern-\dimen@ }}
 {{\setbox\usefulbox=\hbox{$\m@th\scriptstyle #1$}%
    \dimen@ \getslant\the\scriptfont\symletters \ht\usefulbox
    \divide\dimen@ \tw@ 
    \kern\dimen@ 
    \xoverline{\kern-\dimen@ \box\usefulbox\kern\dimen@ }\kern-\dimen@ }}
 {{\setbox\usefulbox=\hbox{$\m@th\scriptscriptstyle #1$}%
    \dimen@ \getslant\the\scriptscriptfont\symletters \ht\usefulbox
    \divide\dimen@ \tw@ 
    \kern\dimen@ 
    \xoverline{\kern-\dimen@ \box\usefulbox\kern\dimen@ }\kern-\dimen@ }}%
 {}}
\makeatother

\makeatletter
\newcommand{\mylabel}[2]{#2\def\@currentlabel{#2}\label{#1}}
\makeatother

\makeatletter
\newcommand{\Mac}{}
\DeclareRobustCommand{\Mac}{%
  M%
  \raisebox{\dimexpr\fontcharht\font`M-\height}{%
    \check@mathfonts\fontsize{\sf@size}{0}\selectfont
    c%
  }%
}
\makeatother

\newtheoremstyle{citing}% name
  {}%      Space above, empty = `usual value'
  {}%      Space below
  {\itshape}% Body font
  {}%         Indent amount (empty = no indent, \parindent = para indent)
  {\bfseries}% Thm head font
  {\textbf{.}}%        Punctuation after thm head
  {.5em}%     Space after thm head: " " = normal interword space;
  {\thmnote{#3}}% Thm head spec

\theoremstyle{plain}
\newtheorem{theorem}{Theorem}

\newtheorem{lemma}[theorem]{Lemma}
\newtheorem{corollary}[theorem]{Corollary}

\newtheorem{bigthm}{Theorem}

\newtheorem{proposition}[theorem]{Proposition}

\theoremstyle{remark}
\newtheorem{example}[theorem]{Example}

\theoremstyle{definition}

\newtheorem{definition}[theorem]{Definition}

\numberwithin{equation}{section}

\theoremstyle{remark}
\newtheorem{remark}[theorem]{Remark}

{\theoremstyle{citing}
}

{\theoremstyle{definition}
}

\title[]{Beyond Finite Dimensions:
Duflo--Serganova Cohomology,
the Formal Superdimension, and a Kac--Wakimoto Theorem for Unitarizable Supermodules}

\author{Steffen Schmidt}
\address{Center for Quantum Mathematics, 
University of Southern Denmark, DK-5230 Odense, Denmark}
\email{stschmidt@imada.sdu.dk}

\let\origmaketitle\maketitle
\def\maketitle{
  \begingroup
  \def\uppercasenonmath##1{} % this disables uppercasing title
  \let\MakeUppercase\relax % this disables uppercasing authors
  \origmaketitle
  \endgroup
}

\begin{document}
\thispagestyle{empty}

\begin{abstract}
The main results of this manuscript are fourfold. First, we introduce relative holomorphic discrete series supermodules over basic classical Lie superalgebras $\mathfrak{g}$ and $\mathfrak{sl}(n\vert n)$. These are unitarizable highest weight supermodules whose $\mathfrak{g}_{\bar{0}}$-constituents are Harish-Chandra modules of relative holomorphic discrete series representations. We show that these are exactly the unitarizable highest weight supermodules whose globalizations have supermatrix coefficients square-integrable modulo the center.

 Second, we define the formal superdimension of a relative holomorphic discrete series supermodule as the parity-weighted sum of the formal dimensions of its $\even$-constituents. We show that it equals, up to sign, the value at zero of an $L$-packet supercharacter obtained by averaging over the Weyl group.

Third, we prove a Kac--Wakimoto theorem for formal superdimension: it is nonzero if and only if the highest weight is maximally atypical. For integral highest weights, it equals the superdimension of a finite-dimensional simple supermodule with the same degree of atypicality. We realize these supermodules through derived Zuckerman functors. For nonintegral highest weights, it is computed by generalized Duflo--Serganova cohomology.

Finally, to establish the Kac--Wakimoto conjecture, we develop the theory of generalized Duflo--Serganova cohomology for unitarizable supermodules. We prove that the cohomology of a unitarizable simple supermodule is unitarizable and either zero or simple.
\end{abstract}

\maketitle

\setlength{\parindent}{1em}
\setcounter{tocdepth}{2}
\tableofcontents

\section{Introduction}
\noindent
Let $\gg$ be a basic classical Lie superalgebra or $\gl(m\vert n)$. The \emph{superdimension} of a finite-dimensional $\gg$-supermodule $V=V_{\bar0}\oplus V_{\bar1}$ is
\begin{equation*}
\sdim V=\dim V_{\bar0}-\dim V_{\bar1}.
\end{equation*}
Unlike the dimension of $V$, the superdimension may vanish for nonzero $V$.

For a fixed positive root system, every finite-dimensional simple $\gg$-supermodule is, up to parity reversal, a highest weight supermodule $L(\Lambda)$. The structure of $L(\Lambda)$ depends strongly on the atypicality of $\Lambda$~\cite{Kac,Kac_classification}. The weight $\Lambda$ is \emph{typical} if $\Lambda+\rho$ is not orthogonal to any isotropic root, and \emph{atypical} otherwise. Here $\rho$ is the Weyl vector, and orthogonality is defined by the invariant bilinear form. The \emph{degree of atypicality} $\at(\Lambda)$ is the largest number of linearly independent, mutually orthogonal isotropic roots orthogonal to $\Lambda+\rho$. This number is independent of the choice of positive roots. Typical supermodules have a character theory similar to that of reductive Lie algebras, and their superdimension vanishes whenever $\gg$ has isotropic roots. The character theory of atypical supermodules is more complicated~\cite{Serganova1996,Brundan2003,GrusonSerganova2010}. 

Kac and Wakimoto~\cite{kac1994integrable} conjectured that $\sdim L(\Lambda)\neq0$ if and only if $\Lambda$ is maximally atypical. The conjecture was proved by Serganova~\cite{serganova2011superdimension}
for all basic classical Lie superalgebras. The main tool in Serganova's proof was the Duflo--Serganova functor~\cite{duflo2005associated}. For an odd element $x$ with $[x,x]=0$, it assigns to $L(\Lambda)$ the supermodule $\DS_x(L(\Lambda))=\ker x/\operatorname{im}x$ over the reduced Lie superalgebra $\gg_x$. The functor is symmetric monoidal, preserves superdimension, and lowers atypicality by the rank of $x$. If $\Lambda$ is not maximally atypical, then $\DS_x(L(\Lambda))=0$ for $x$ of maximal rank, so $\sdim L(\Lambda)=0$, while nonvanishing for maximally atypical modules requires a separate argument.

The Duflo--Serganova functor has since become a central tool in the
representation theory of Lie superalgebras~\cite{gorelik2022duflo}. Among
other things, it allows explicit computations of superdimensions for
$\gl(m\vert n)$ and the orthosymplectic Lie
superalgebras~\cite{heidersdorf2021cohomological,GorelikHeidersdorf2022semisimplicity,GorelikHeidersdorfCharacter}.

All of this concerns finite-dimensional supermodules. Unitarizable
supermodules over noncompact real forms, however, are in general
infinite-dimensional, and their superdimension is not defined. The present work extends the notion of superdimension to relative holomorphic discrete series supermodules by introducing a formal superdimension. We also develop generalized Duflo--Serganova cohomology for unitarizable supermodules and use it to prove an analogue of the Kac--Wakimoto theorem. To our knowledge, this is
the first extension of the Kac--Wakimoto theorem. 

In particular, our results apply to the complex superconformal
algebras $\osp(n\vert4)$ in three dimensions, $\sl(4\vert n)$ in four dimensions, $F(4)$ in five
dimensions, and $\osp(8\vert2n)$ in six dimensions. Our formal
superdimension provides a rigorous counterpart to the regularized supertrace
underlying the superconformal index~\cite{romelsberger2006counting,
kinney2007index,SchmidtWalcher}.

\subsection{Results}
\label{subsec::results}
Throughout, $\gg=\even\oplus\odd$ is a basic classical Lie superalgebra with $\odd\neq0$, or $\sl(n\vert n)$.\footnote{The Lie superalgebra $\sl(n\vert n)$ is not simple, and some statements in this introduction require modification in this case. These are carried out in Section~\ref{sec::equal_rank_type_A}.} We study \emph{relative holomorphic discrete series supermodules}: simple unitarizable highest weight supermodules $L(\Lambda)$ whose simple even constituents belong to the relative holomorphic discrete series of $\gg_{\bar0,\RR}$. Each constituent has a formal dimension in the sense of Harish-Chandra, which allows us to define a formal superdimension. When $\gg_{\bar0,\RR}$ is compact, every simple unitarizable supermodule belongs to this class, and its formal superdimension agrees with its ordinary superdimension. In general, the class contains both typical and maximally atypical supermodules, providing a setting for a Kac--Wakimoto theorem beyond finite dimensions.

Four questions guide this work.
\begin{itemize}
\item Can relative holomorphic discrete series supermodules be characterized by square-integrability of their supermatrix coefficients modulo the center (Theorem~\ref{thmA})?
\item How can Duflo--Serganova cohomology be extended to unitarizable supermodules, and what properties does it satisfy (Theorem~\ref{thmB})?
\item How can superdimension be extended to relative holomorphic discrete series supermodules, and can one compute it from their supercharacters (Theorem~\ref{thmC})?
\item Which relative holomorphic discrete series supermodules have nonzero formal superdimension (Theorem~\ref{thmD})?
\end{itemize}

Theorems~\ref{thmA}--\ref{thmC} are independent of one another, whereas the proof of our main result, Theorem~\ref{thmD}, uses Theorems~\ref{thmB} and~\ref{thmC}.

\subsubsection*{Setting}
Let $\omega$ be the conjugate-linear anti-involution of $\gg$ defining a real form $\gg_{\RR}$. A $\gg$-supermodule $\HH$ is \emph{$\omega$-unitarizable} if it carries a positive definite Hermitian form $\langle\cdot,\cdot\rangle$ such that $\langle Xv,w\rangle=\langle v,\omega(X)w\rangle$ for all $X\in\gg$. Every nontrivial simple unitarizable supermodule is of highest or lowest weight type (Theorem~\ref{theorem::unitary-highest-lowest-weight}), and we restrict to the highest weight case. The lowest weight case is entirely analogous. The even part $\gg_{\bar0,\RR}$ is reductive of Hermitian type. Let $\kk$ be the complexification of a maximal compact subalgebra, $\hh\subseteq\kk$ a compactly embedded Cartan subalgebra, and $\Delta_{\bar0}^+=\Delta_c^+\sqcup\Delta_n^+$ a $\kk$-adapted positive system of compact and noncompact even roots, with half-sum $\rho_{\bar0}$. We denote by $W$ and $W_{c}$ the Weyl groups of $\even$ and $\kk$, respectively.

The restriction of a simple unitarizable highest weight supermodule $L(\Lambda)$ to $\even$ is a finite direct sum
\begin{equation}\label{eq::intro_branching}
L(\Lambda)\vert_{\even}
\cong
\bigoplus_\mu L_0(\mu)\otimes E_\mu,
\end{equation}
where the $L_0(\mu)$ are simple unitarizable highest weight $\even$-modules and the $E_\mu$ are finite-dimensional multiplicity super vector spaces. By Harish-Chandra~\cite{HC_VI}, $L_0(\mu)$ belongs to the relative holomorphic discrete series, i.e., the matrix coefficients of its unitary globalization are square-integrable modulo the center, if and only if
\begin{equation}
2\frac{(\mu+\rho_{\bar0},\beta)}{(\beta,\beta)}<0
\qquad\text{for all }\beta\in\Delta_n^+.
\end{equation}
Among unitarizable highest weight $\even$-modules, the relative holomorphic discrete series modules $L_0(\mu)$ are distinguished by the existence of a positive formal dimension, which replaces ordinary dimension in the Schur orthogonality relations. 
We are interested in supermodules $L(\Lambda)$ whose even constituents all belong to the relative holomorphic discrete series. We call such a supermodule a \emph{relative holomorphic discrete series supermodule}.

\subsubsection*{(A) Relative holomorphic discrete series supermodules and square-integrability}

Our first result gives a characterization of these supermodules in terms of the associated unitary representations of the Lie supergroup. For reductive Lie groups, discrete series representations are characterized by the square-integrability of their matrix coefficients modulo the center. We establish the corresponding result for Lie supergroups: a simple unitarizable highest weight supermodule belongs to the relative holomorphic discrete series if and only if the supermatrix coefficients of its unitary globalization are square-integrable modulo the center of the Lie supergroup.

Let $\widetilde{\calG}_{\RR}=(\widetilde{G}_{\RR},\gg_{\RR})$ be the connected, simply connected Lie supergroup associated with $\gg_{\RR}$. By~\cite{Carmeli_Fioresi_Varadarajan_Unitary}, every simple unitarizable highest weight supermodule arises as the Harish-Chandra supermodule $\calH^{[\kk]}$ of an irreducible unitary representation $\Sigma=(\pi,\tau)$ of $\widetilde{\calG}_{\RR}$ on a Hilbert superspace $\calH$, unique up to unitary equivalence. Following~\cite{NeebSalmasian}, its \emph{supermatrix coefficients} are the functions
\begin{equation*}
\varphi_{v,w}(g,D)=\langle\pi(g)\tau(D)v,w\rangle,
\qquad
g\in\widetilde{G}_{\RR},\quad D\in\UE(\gg),
\end{equation*}
where $v,w\in\calH^\infty$ are smooth vectors. When $\gg$ is simple, the center of $\gg_{\RR}$ is trivial, so the center $\mathcal{Z}$ of $\widetilde{\calG}_{\RR}$ is discrete and even. Its underlying group $C$ acts on $\calH$ through a unitary character $\chi_C$.

\begin{bigthm}[Theorem~\ref{thm::super_embedding}]\label{thmA}
Let $\Sigma=(\pi,\tau)$ be an irreducible unitary representation of $\widetilde{\calG}_{\RR}$ on $\calH$. Then the following statements are equivalent.
\begin{enumerate}
\item[(a)] $\calH^{[\kk]}$ is a relative holomorphic discrete series supermodule.
\item[(b)] Every supermatrix coefficient $\varphi_{v,w}$, with homogeneous $v,w\in\calH^\infty$ of the same parity, is square-integrable modulo $\calZ$.
\item[(c)] 
For some nonzero homogeneous $u\in\calH^{[\kk]}$, the coefficient map $\Psi_u(v)(g,D)=\langle\pi(g)\tau(D)v,u\rangle$ extends to a bounded injective map $\calH\to L^2_+(\widetilde{\calG}_{\RR}/\calZ;\chi_C)$ with closed image and intertwines $\pi$ with right translations.

\end{enumerate}
\end{bigthm}

Theorem~\ref{thmA} is the super analogue of the $L^2$-realization of the relative holomorphic discrete series due to Harish-Chandra~\cite{HC_VI} and Neeb~\cite{NeebSquareIntegrable} (see also Theorem~\ref{thm::embedding}). 

\subsubsection*{(B) Generalized Duflo--Serganova cohomology and unitarity}
The Duflo--Serganova functor is fundamental to the finite-dimensional representation theory of Lie superalgebras. However, it cannot be used to study unitarizable supermodules, as any $x\in\odd$ satisfying $[x,x]=0$ and $\omega(x)=ix$ acts trivially on every unitarizable supermodule. Since every nontrivial supermodule over a simple $\gg$ is faithful, the self-commuting variety has no nonzero real points. For $\sl(n\vert n)$, the same conclusion follows from a direct matrix computation (Lemma~\ref{lemm::equal_rank_real_odd}).

We therefore use the generalized Duflo--Serganova functor of~\cite{Weissauer2010,heidersdorf2021cohomological, GorelikSherman2024,Sherman2025}. Let $\YY^{\hom}$ be the set of $x\in\odd$ for which $x^2=\frac12[x,x]$ acts semisimply on $\gg$, and set $c\coloneqq x^2$. For a $\gg$-supermodule $M$ on which $c$ acts semisimply, the space $M^c\coloneqq\ker(c)$ is $x$-stable with $x^2=0$ on it, and
\begin{equation}
\DS_x(M)\coloneqq\frac{\ker(x:M^c\to M^c)}{\im(x:M^c\to M^c)}
\end{equation}
is a supermodule over $\gg_x\coloneqq H(\ker\ad_c,\ad_x)$. For $x^2=0$, this is the original functor. The \emph{real locus} $\YY^{\hom}_\omega\coloneqq\{x\in\YY^{\hom}:\omega(x)=ix\}$ contains nonzero elements whenever $\gg$ has an isotropic odd root.

\begin{bigthm}\label{thmB}
Let $x\in\YY^{\hom}$ have rank $k$, and set $c\coloneqq x^{2}$.
\begin{enumerate}
\item[(a)] Let $M$ be a $\gg$-supermodule with infinitesimal character $\chi$, and assume that $x^2$ acts semisimply on $M$. Then every simple subquotient of $\DS_x(M)$ has an infinitesimal character $\xi$ satisfying
\[
\eta_x^\ast(\xi)=\chi,
\qquad
\at(\xi)=\at(\chi)-k,
\]
where $\eta_x^\ast$ is induced by the natural homomorphism $\eta_x\colon\mathfrak Z(\gg)\to\mathfrak Z(\gg_x)$. In particular, $\DS_x(M)=0$ if $k>\at(\chi)$.
\item[(b)] If $x\in\YY^{\hom}_\omega$, then $\omega$ induces a conjugate-linear anti-involution $\omega_x$ of $\gg_x$. The corresponding real form is determined explicitly. For example, $\su(p,q\vert0,n)$ reduces to $\su(p-a,q-b\vert0,n-k)$, where $0\leq a\leq p$, $0\leq b\leq q$, and $a+b=k$.
\item[(c)] Let $x\in\YY^{\hom}_\omega$, and let $\HH$ be a simple $\omega$-unitarizable highest weight $\gg$-supermodule on which $c$ acts semisimply. Then
\[
\DS_x(\HH)\cong\HH^{c},
\]
and $\DS_x(\HH)$ is either zero or a simple $\omega_x$-unitarizable $\gg_x$-supermodule.
\end{enumerate}
\end{bigthm}
Parts~(a)--(c) are proved in Theorem~\ref{thm::generalized_DS_character}, Proposition~\ref{prop::unitary_real_forms_and_DS}, and Proposition~\ref{prop::form_DS_x_generalized} together with Theorem~\ref{thm::DS_unitary}, respectively. Part~(a) extends~\cite[Theorem~2.1]{serganova2011superdimension} to the generalized DS cohomology. Parts~(b) and~(c) are specific to unitarity. Positivity of the Hermitian form forces $x$ to vanish on $\HH^c$, which identifies the cohomology with an eigenspace and endows it with a positive definite form. Simplicity is a further consequence of unitarity. It contrasts with the finite-dimensional theory, where the Duflo--Serganova cohomology of a simple supermodule is in general a direct sum of several simple supermodules~\cite{heidersdorf2021cohomological}. Thus generalized Duflo--Serganova cohomology maps simple unitarizable highest weight supermodules over $\gg$ to zero or to simple unitarizable supermodules over the smaller Lie superalgebra $\gg_x$ with the induced real form, and it can be iterated.

Moreover, in contrast to the finite-dimensional case, generalized Duflo--Serganova cohomology does not preserve the numerical invariant introduced next (Section~\ref{subsec::formal_superdimension}). It nevertheless controls that invariant through characters, and this is how Theorem~\ref{thmB} enters the proof of Theorem~\ref{thmD}.

\subsubsection*{(C) The formal superdimension} Let $L(\Lambda)$ be a relative holomorphic discrete series supermodule with
decomposition~\eqref{eq::intro_branching}. For each $\mu$, let $\pi_\mu$
denote the relative holomorphic discrete series representation of
$\widetilde G_{\RR}$ whose Harish-Chandra module is $L_0(\mu)$, and let
$d(\pi_\mu)>0$ be its formal dimension. This is the constant in the Schur
orthogonality relations
\begin{equation}\label{eq::intro_schur_orthogonality}
\int_{\widetilde G_{\RR}/Z}|\langle\pi_\mu(g)v,w\rangle|^2\,d(gZ)
=d(\pi_\mu)^{-1}\|v\|^2\|w\|^2,
\end{equation}
where $Z$ is the center of $\widetilde G_{\RR}$ and the Haar measure is
normalized as in Theorem~\ref{thm::formal_dimension}. Explicitly,
\begin{equation}\label{eq::intro_formal_dimension}
d(\pi_\mu)
=\prod_{\alpha\in\Delta_c^+}
\frac{(\mu+\rho_c,\alpha)}{(\rho_c,\alpha)}
\prod_{\beta\in\Delta_n^+}
\frac{|(\mu+\rho_{\bar0},\beta)|}{|(\rho_{\bar0},\beta)|},
\end{equation}
where $\Delta_c^+$ and $\Delta_n^+$ denote the positive compact and
noncompact even roots. The \emph{formal superdimension} of $L(\Lambda)$ is
the parity-weighted sum of these formal dimensions:
\begin{equation}\label{eq::intro_fsdim}
\fsdim L(\Lambda)\coloneqq\sum_\mu\sdim(E_\mu)\,d(\pi_\mu).
\end{equation}
 It is independent of the chosen decomposition, changes sign under parity reversal, is additive, and satisfies $\fsdim(\HH\otimes F)=\fsdim(\HH)\sdim(F)$ for finite-dimensional $F$ whenever the left-hand side is defined. If $\gg_{\bar0,\RR}$ is compact, it equals the ordinary superdimension. Unlike the superdimension, its values cannot all be made integral by any normalization of the Haar measure (Section~\ref{subsec::formal_superdimension}).

To relate the formal superdimension to characters, we average the
supercharacter over the full even Weyl group, as in the Weyl-group
description of discrete series $L$-packets. The \emph{$L$-packet
supercharacter} of $L(\Lambda)$ is
\begin{equation}\label{eq::intro_L_packet}
\Theta_\Lambda
\coloneqq
\frac{1}{|W_c|}\sum_{w\in W}w\bigl(\sch L(\Lambda)\bigr).
\end{equation}

\begin{bigthm}[Theorem~\ref{thm::L_packet_and_sdim}]\label{thmC}
Let $L(\Lambda)$ be a relative holomorphic discrete series supermodule. Then the $L$-packet supercharacter
$\Theta_\Lambda$ extends holomorphically to a neighborhood of $0\in\hh$, and
\begin{equation*}
\Theta_\Lambda(0)
=
(-1)^{|\Delta_n^+|}\fsdim L(\Lambda).
\end{equation*}
\end{bigthm}

For the holomorphic discrete series, Harish-Chandra expressed the normalized formal dimension as the absolute value of Weyl's dimension polynomial evaluated at the Harish-Chandra parameter. Neeb extended this formula to the relative holomorphic discrete series modulo a possibly noncompact center~\cite[Theorem~3.17]{NeebSquareIntegrable}. Theorem~\ref{thmC} establishes the corresponding formula for supermodules. It also expresses $\fsdim L(\Lambda)$ as the value at zero of an averaged supercharacter. This description permits the use of Duflo--Serganova cohomology and derived Zuckerman functors, which is central to proving the Kac--Wakimoto conjecture.

\subsubsection*{(D) The Kac--Wakimoto conjecture for the formal superdimension}
The examples in Section~\ref{subsec::formal_superdimension} suggest that formal superdimension plays the same role for relative holomorphic discrete series supermodules as ordinary superdimension does in the finite-dimensional case. In particular, its nonvanishing characterizes maximal atypicality. Our main result establishes this analogue of the Kac--Wakimoto conjecture.

Section~\ref{subsec::Zuckerman_functors} introduces the Zuckerman functor $\Gamma$ and its right derived functors $R^j\Gamma$.

\begin{bigthm}[Theorem~\ref{thm::KW_formal_superdimension}]\label{thmD}
Let $L(\Lambda)$ be a relative holomorphic discrete series supermodule. Then
\begin{equation*}
\fsdim L(\Lambda)\neq0
\quad\Longleftrightarrow\quad
\at(\Lambda)=\defect(\gg).
\end{equation*}
Moreover, if $\Lambda$ is integral, then $R^{|\Delta_n^+|}\Gamma(L(\Lambda))$ is a finite-dimensional simple $\gg$-supermodule with the same infinitesimal character as $L(\Lambda)$, and
\begin{equation*}
\fsdim L(\Lambda)=\sdim R^{|\Delta_n^+|}\Gamma(L(\Lambda)).
\end{equation*}
\end{bigthm}

The proof of Theorem~\ref{thmD} has two parts. Vanishing comes from Theorem~\ref{thmB}: if $\Lambda$ is not maximally atypical, the generalized Duflo--Serganova cohomology of maximal rank is zero, and by unitarity the $L$-packet supercharacter then vanishes at $0$. Theorem~\ref{thmC} gives $\fsdim L(\Lambda)=0$.

For nonvanishing, integral and nonintegral highest weights behave differently. When $\Lambda$ is integral, derived Zuckerman functors turn $L(\Lambda)$ into a finite-dimensional simple supermodule with the same infinitesimal character, whose supercharacter is, up to sign, the $L$-packet supercharacter of $L(\Lambda)$. Serganova's theorem then applies (Section~\ref{subsec::KW_conjecture_integral_case}). When $\Lambda$ is nonintegral, no such finite-dimensional counterpart exists. Instead, generalized Duflo--Serganova cohomology reduces the computation to a single cohomology module, which is finite-dimensional, a relative holomorphic discrete series module, or a supermodule over $\osp(1\vert2s)$. In each case, its contribution at $0$ is nonzero (Section~\ref{subsec::KW_conjecture_nonintegral_case}).

For $\sl(n\vert n)$ and $\psl(n\vert n)$, atypicality is defined through $\gl(n\vert n)$, and the same results hold. Here only integral highest weights occur (Section~\ref{sec::equal_rank_type_A}).

\subsection{Leitfaden}
In Sections~\ref{sec::basic_classical_Lie_superalgebras}--\ref{sec::relative_holomorphic_discrete_series}, we recall the necessary background on basic classical Lie superalgebras, real forms, unitarity, and the relative holomorphic discrete series. Section~\ref{sec::relative_HDS_supermodules} introduces relative holomorphic discrete series supermodules and characterizes them by the condition that the supermatrix coefficients of their globalizations are square-integrable modulo the center (Theorem~\ref{thmA}). Section~\ref{sec::DS_functor_and_unitarity} develops generalized Duflo--Serganova cohomology for unitarizable supermodules (Theorem~\ref{thmB}). Section~\ref{sec::formal_superdimension} introduces the formal superdimension, gives examples, and proves Theorems~\ref{thmC} and~\ref{thmD}. Finally, in Section~\ref{sec::equal_rank_type_A}, we treat the cases $\sl(n\vert n)$ and $\psl(n\vert n)$ separately. 

\subsection{Conventions}
\label{subsec::conventions_intro}

Unless otherwise stated, we work over the complex numbers $\CC$. We write $\ZZ_{\geq0}$ and $\ZZ_{>0}$ for the nonnegative and positive integers, and $\ZZ_2=\{\bar0,\bar1\}$. A superspace is a $\ZZ_2$-graded vector space $V=V_{\bar0}\oplus V_{\bar1}$. The parity of a homogeneous element $v$ is $p(v)$, and identities involving parities are stated for homogeneous elements and extended linearly.

Supermodules are left supermodules, submodules are graded, and morphisms are even. Accordingly, $\operatorname{Hom}_{\gg}(M,N)$ denotes the space of even intertwiners and $\underline{\operatorname{Hom}}_{\gg}(M,N)$ the full Hom superspace. The parity reversal $\Pi M$ carries the action of $M$ and the grading $(\Pi M)_{\bar a}=M_{\bar a+\bar1}$.

Invariant bilinear forms are written $(\cdot,\cdot)$ and Hermitian forms $\langle\cdot,\cdot\rangle$. Hermitian forms are linear in the first variable, conjugate-linear in the second, and even, so that $V_{\bar0}$ and $V_{\bar1}$ are orthogonal. The adjoint of an operator $T$ is denoted by $T^\ast$.

A conjugate-linear anti-involution of $\gg$ is a parity-preserving conjugate-linear map $\omega$ with $\omega^2=\operatorname{id}$ and $\omega([X,Y])=[\omega(Y),\omega(X)]$. It extends uniquely to a conjugate-linear antiautomorphism of $\UE(\gg)$, again denoted by $\omega$.

\subsubsection*{Acknowledgements}
The author is very grateful to Rainer Weissauer, who not only provided the underlying idea for this manuscript but also contributed to its development through numerous conversations and helpful suggestions. The author further thanks Johannes Walcher for many discussions and helpful suggestions during the writing of the manuscript.

\subsubsection*{AI disclosure}
ChatGPT and Claude were solely used for language editing. 

\section{Basic classical Lie superalgebras}\label{sec::basic_classical_Lie_superalgebras}
\noindent
We briefly review the structure theory of basic classical Lie superalgebras, with particular emphasis on type $A$, which provides the main examples. We then discuss the real forms relevant to unitarity and their maximal compact subalgebras. Finally, we recall highest weight supermodules and define atypicality.

\subsection{Basic classical Lie superalgebras}

A finite-dimensional complex Lie superalgebra $\gg=\even\oplus\odd$ is called \emph{basic classical} if $\gg$ is simple, $\even$ is reductive, and $\gg$ admits a nondegenerate even invariant supersymmetric bilinear form $(\cdot,\cdot)$. In particular,
$
(\even,\odd)=0
$, $(x,y)=(-1)^{p(x)p(y)}(y,x)$,
and
$
([x,y],z)=(x,[y,z])
$
for all homogeneous $x,y,z\in\gg$.

The basic classical Lie superalgebras are the simple complex Lie algebras, regarded as purely even Lie superalgebras, together with the classical Lie superalgebras
\begin{equation}
\begin{aligned}
&A(m\vert n)=\sl(m+1\vert n+1), && m,n\geq0,\quad m\neq n,\\
&A(n\vert n)=\psl(n+1\vert n+1), && n\geq1,\\
&B(m\vert n)=\osp(2m+1\vert2n), && m\geq0,\quad n\geq1,\\
&C(n)=\osp(2\vert2n-2), && n\geq3,\\
&D(m\vert n)=\osp(2m\vert2n), && m\geq2,\quad n\geq1,
\end{aligned}
\end{equation}
and the exceptional Lie superalgebras
\begin{equation}
D(2,1;\upalpha),\qquad F(4),\qquad G(3),
\end{equation}
where $\upalpha\in\CC\setminus\{0,-1\}$; see~\cite{Kac}.  

Throughout, $\gg$ denotes a basic classical Lie superalgebra with $\odd\neq0$, excluding $A(n\vert n)$, which is treated separately in Section~\ref{sec::equal_rank_type_A}. We fix a nondegenerate even invariant supersymmetric bilinear form $(\cdot,\cdot)$ on $\gg$. For $D(2,1;\upalpha)$, we restrict to $\upalpha\in\RR\setminus\{0,-1\}$, since only these parameter values admit nontrivial unitarizable supermodules~\cite{neeb2011lie}.

The super Killing form
\begin{equation}
B(x,y)=\str(\ad_x\circ\ad_y),
\qquad x,y\in\gg,
\label{eq::super_Killing_form}
\end{equation}
is nondegenerate for simple complex Lie algebras and for $A(m\vert n)$ with $m\neq n$, $B(m\vert n)$, $C(n)$, $D(m\vert n)$ with $m\neq n+1$, $F(4)$, and $G(3)$. It vanishes identically for $D(n+1\vert n)$ and $D(2,1;\upalpha)$. Throughout, $(\cdot,\cdot)$ denotes a nondegenerate invariant even supersymmetric bilinear form, chosen as a normalized Killing form, a supertrace form, or the standard invariant form of~\cite[Section~5.4]{musson2012lie}. The normalization is specified in each case.

Fix a Cartan subalgebra $\hh\subseteq\even$. We assume that $\hh$ is contained in the complexification of a maximal compact subalgebra of the relevant real form used to study unitarity, as described in Section~\ref{subsec::real_forms}. Let
$
\Delta=\Delta_{\bar0}\sqcup\Delta_{\bar1}
$
be the corresponding root system. Then
\begin{equation}
\gg
=
\hh\oplus\bigoplus_{\alpha\in\Delta}\gg^\alpha,
\end{equation}
where
$
\gg^\alpha
=
\{x\in\gg:[h,x]=\alpha(h)x
\text{ for all }h\in\hh\}.
$ The set $\Delta_{\bar0}$ is the root system of the complex reductive Lie algebra $\even$. In addition, the Weyl group $W$ of $\gg$ is the Weyl group of the Lie algebra $\even$. We record the root space properties used below.

\begin{proposition}[{\cite[Proposition~2.5.5]{Kac}}]
\label{prop::structure_theory_BCLSA}
The following statements hold.
\begin{enumerate}
    \item[(a)]
    If $\alpha\in\Delta$, then $-\alpha\in\Delta$. In particular,
    $\Delta=-\Delta$, $\Delta_{\bar0}=-\Delta_{\bar0}$, and
    $\Delta_{\bar1}=-\Delta_{\bar1}$.

    \item[(b)]
     Each root space $\gg^\alpha$ is one-dimensional. More precisely,
    $\bigl(\dim(\gg^\alpha)_{\bar0}\,\vert\,\dim(\gg^\alpha)_{\bar1}\bigr)=(1\vert0)$ for $\alpha\in\Delta_{\bar0}$ and
    $\bigl(\dim(\gg^\alpha)_{\bar0}\,\vert\,\dim(\gg^\alpha)_{\bar1}\bigr)=(0\vert1)$ for $\alpha\in\Delta_{\bar1}$.

    \item[(c)]
    If $\alpha,\beta\in\Delta$ and $\alpha+\beta\in\Delta$, then
    $
    [\gg^\alpha,\gg^\beta]\subseteq\gg^{\alpha+\beta}.
    $
    If $\alpha+\beta\notin\Delta\cup\{0\}$, then
    $[\gg^\alpha,\gg^\beta]=0$. Moreover,
    $[\gg^\alpha,\gg^{-\alpha}]$ is one-dimensional and contained in
    $\hh$.

    \item[(d)]
    The restriction of $(\cdot,\cdot)$ to $\hh$ is nondegenerate, and
    $
    (\gg^\alpha,\gg^\beta)=0
    $
    unless $\alpha+\beta=0$.

    \item[(e)]
    For each $\alpha\in\Delta$, let $h_\alpha\in\hh$ be defined by
    $
    (h_\alpha,h)=\alpha(h)
    $
    for all $h\in\hh$. If $0\neq e_\alpha\in\gg^\alpha$ and
    $0\neq e_{-\alpha}\in\gg^{-\alpha}$, then
    $
    [e_\alpha,e_{-\alpha}]
    =
    (e_\alpha,e_{-\alpha})h_\alpha.
    $

    \item[(f)]
    The form on $\hh^{\ast}$ defined by
    $
    (\lambda,\mu)\coloneqq (h_\lambda,h_\mu)
    $
    is nondegenerate and $W$-invariant.
\end{enumerate}
\end{proposition}

We fix a positive system $\Delta^+$, specified below, and write $\Delta^+=\Delta_{\bar0}^+\sqcup\Delta_{\bar1}^+$ for its even and odd roots. The corresponding Weyl
vector is
\begin{equation}
\rho=\rho_{\bar0}-\rho_{\bar1},
\qquad
\rho_{\bar0}
=\frac12\sum_{\alpha\in\Delta_{\bar0}^+}\alpha,
\qquad
\rho_{\bar1}
=\frac12\sum_{\alpha\in\Delta_{\bar1}^+}
\dim(\gg^\alpha)\alpha.
\end{equation}
The associated \emph{simple system} $\Uppi$ consists of the
positive roots that cannot be written as a sum of two positive
roots. Every positive root is a nonnegative integral combination
of roots in $\Uppi$. We set
\begin{equation}
\nn^\pm
=\bigoplus_{\alpha\in\Delta^+}\gg^{\pm\alpha}.
\end{equation}
The root-space decomposition then gives the
\emph{triangular decomposition}
$\gg=\nn^-\oplus\hh\oplus\nn^+$, and
$\mathfrak b=\hh\oplus\nn^+$ is the corresponding
\emph{Borel subalgebra}.

Unlike for ordinary reductive Lie algebras, $W$ need not act transitively on the simple systems of $\gg$. Let $\Uppi$ be a simple system and $\theta\in\Uppi\cap\Delta_{\bar1}$ an isotropic root. The \emph{odd reflection} at $\theta$ produces the simple system $\Uppi_\theta\coloneqq r_\theta(\Uppi)$, where
\begin{equation}
r_\theta(\alpha)=
\begin{cases}
-\theta,&\alpha=\theta,\\
\alpha+\theta,&\alpha\neq\theta,\ (\alpha,\theta)\neq0,\\
\alpha,&\alpha\neq\theta,\ (\alpha,\theta)=0.
\end{cases}
\label{eq::odd-reflection}
\end{equation}
The corresponding positive system is $\Delta_\theta^+=(\Delta^+\setminus\{\theta\})\cup\{-\theta\}$, so the positive even roots remain unchanged. Conversely, any two simple systems with the same positive even roots are related by a sequence of odd reflections; see~\cite{Serganova_reps}.

In contrast to simple Lie algebras, $\gg$ can have isotropic roots, that is, roots $\alpha\in\Delta$ with $(\alpha,\alpha)=0$. Such roots are necessarily odd. We denote the set of isotropic odd roots by
$
\overline{\Delta}_{\bar1}.
$
The defect $\defect(\gg)$ is the maximal number of linearly independent, pairwise orthogonal roots in $\overline{\Delta}_{\bar1}$. 
Its values are
\begin{equation}
\renewcommand{\arraystretch}{1.5}
\begin{array}{c|ccccccc}
\gg
& A(m\vert n),\ m\neq n
& B(m\vert n)
& C(n)
& D(m\vert n)
& D(2,1;\upalpha)
& F(4)
& G(3)\\ \hline
\defect(\gg)
& \min\{m+1,n+1\}
& \min\{m,n\}
& 1
& \min\{m,n\}
& 1
& 1
& 1.
\end{array}
\end{equation}
In particular, $B(0\vert n)=\osp(1\vert2n)$ has defect zero. A simple complex Lie algebra, regarded as a purely even Lie superalgebra, also has defect zero.

Basic classical Lie superalgebras are divided into type~I and type~II. A Lie superalgebra $\gg$ is of type~I if it admits a $\ZZ$-grading
$
\gg=\gg_{-1}\oplus\gg_0\oplus\gg_1
$
with $\gg_0=\even$ and simple $\even$-modules $\gg_{\pm1}$. It is of type~II if $\even$ is semisimple and $\odd$ is a simple $\even$-module. The Lie superalgebras $\gl(m\vert n)$, $A(m\vert n)$, and $C(n)$ are of type~I, whereas $B(m\vert n)$, $D(m\vert n)$, $D(2,1;\upalpha)$, $F(4)$, and $G(3)$ are of type~II. 
Since our examples concern $A(m\vert n)$, we describe its structure in more detail. The root systems of the remaining basic classical Lie superalgebras are collected in Appendix~\ref{appendix::structure_theory} and will be needed only once.

\subsubsection{Basic classical Lie superalgebras of type \texorpdfstring{$A$}{A}}\label{subsec::type_A_structure}

Let $m,n\in\ZZ_{\geq1}$. Denote by $M_n(\CC)$ the space of $n\times n$ complex matrices and by $M_{m,n}(\CC)$ the space of $m\times n$ complex matrices. The \emph{general linear Lie superalgebra} $\glmn$ consists of the block matrices
\begin{equation}\label{eq::matrix_form}
X=\left(\begin{array}{@{}c|c@{}}
A&B\\ \hline
C&D
\end{array}\right),
\end{equation}
where $A\in M_m(\CC)$, $B\in M_{m,n}(\CC)$, $C\in M_{n,m}(\CC)$, and $D\in M_n(\CC)$. Its even part consists of the block-diagonal matrices, while its odd part consists of the off-diagonal matrices. The bracket is the supercommutator.

The \emph{supertrace} of $X\in\glmn$ is
$
\str(X)=\tr(A)-\tr(D).
$
The \emph{special linear Lie superalgebra} is the codimension-one ideal
$
\slmn=\{X\in\glmn:\str(X)=0\}.
$
Its even part is
$
\slmn_{\bar0}\cong\sl(m)\oplus\sl(n)\oplus\CC E_{m\vert n},
$
where $E_{m\vert n}=\diag(nI_m,mI_n)$ and $I_m,I_n$ are the corresponding identity matrices. If $m=n$, then $E_{n\vert n}$ is a scalar multiple of the identity matrix and spans the center of $\sl(n\vert n)$.

We assume throughout that $m+n>2$. If $m\neq n$, then $\slmn$ is simple and
$
\slmn=A(m-1\vert n-1).
$
If $m=n$, then
$
\mathfrak{psl}(n\vert n)\coloneqq\sl(n\vert n)/\CC E_{n\vert n}
$
is simple and equals $A(n-1\vert n-1)$. For the rest of this section, let $\gg\coloneqq A(m-1\vert n-1)$ with $m\neq n$.

Let $\dd\coloneqq\{H=\diag(h_1,\ldots,h_{m+n})\}$ be the abelian Lie subalgebra of diagonal matrices in $\glmn$. We denote by $\hh$ the image in $\gg$ of the diagonal matrices in $\dd$ with vanishing supertrace. Then $\hh$ is a Cartan subalgebra of $\gg$.

Let $(\epsilon_1,\ldots,\epsilon_m,\delta_1,\ldots,\delta_n)$ be the standard basis of $\dd^{\ast}$, defined by $\epsilon_i(H)=h_i$ and $\delta_a(H)=h_{m+a}$. We use this notation also to describe weights and roots of $\gg$. We represent a weight $\lambda\in\hh^{\ast}$ by a tuple $(\lambda_1,\ldots,\lambda_m\vert\lambda'_1,\ldots,\lambda'_n)$ such that
$
\lambda=\sum_{i=1}^m\lambda_i\epsilon_i+\sum_{a=1}^n\lambda'_a\delta_a.
$
Since $\sum_{i=1}^m\epsilon_i-\sum_{a=1}^n\delta_a$ vanishes on the diagonal matrices of vanishing supertrace, this representation is not unique. In particular, adding a multiple of $(1,\ldots,1\vert-1,\ldots,-1)$ does not change the corresponding weight.

With respect to $\hh$, the root system of $\gg$ is
\begin{equation}
\begin{aligned}
\Delta_{\bar0}
&=
\{\pm(\epsilon_i-\epsilon_j):1\leq i< j\leq m\}
\sqcup
\{\pm(\delta_a-\delta_b):1\leq a< b\leq n\},\\
\Delta_{\bar1}
&=
\{\pm(\epsilon_i-\delta_a):1\leq i\leq m,\ 1\leq a\leq n\}.
\end{aligned}
\end{equation}
The Weyl group is $W\cong S_m\times S_n$, where the two factors permute the $\epsilon_i$ and the $\delta_a$, respectively.

Up to a nonzero scalar, the invariant form on the root space is the restriction of the form determined by
\begin{equation}
(\epsilon_i,\epsilon_j)=\delta_{ij},
\qquad
(\delta_a,\delta_b)=-\delta_{ab},
\qquad
(\epsilon_i,\delta_a)=0.
\end{equation}
In particular, every odd root is isotropic, so
$
\overline{\Delta}_{\bar1}=\Delta_{\bar1}.
$

The set $\Delta_{\bar0}$ is the root system of the reductive Lie algebra $\even$. We fix the standard positive system
\begin{equation}
\Delta_{\bar0}^+
=
\{\epsilon_i-\epsilon_k:1\leq i<k\leq m\}
\sqcup
\{\delta_a-\delta_b:1\leq a<b\leq n\}.
\end{equation}
Note that all positive systems $\Delta_{\bar{0}}^{+}$ are conjugate under $W$. As discussed in~\cite{Schmidt,Schmidt_unitary_branching}, two choices of positive odd roots are relevant for unitarity, namely the standard and the nonstandard choice. The standard choice of positive odd roots is
\begin{equation}
\Delta_{\bar1,\mathrm{st}}^+
=
\{\epsilon_i-\delta_a:1\leq i\leq m,\ 1\leq a\leq n\}.
\end{equation}

We write
$
\Delta_{\mathrm{st}}^+=\Delta_{\bar0}^+\sqcup\Delta_{\bar1,\mathrm{st}}^+.
$
This is the distinguished positive system. Its positive root spaces are represented by upper triangular matrices, and the corresponding Borel subalgebra is the standard Borel subalgebra of upper triangular matrices.

The nonstandard choice is defined with respect to a partition $m=p+q$, where $p,q>0$. The $pq$-nonstandard choice is
\begin{equation}\label{eq::non_standard_positive}
\Delta_{\bar1,\mathrm{nst}}^+
=
\{\epsilon_i-\delta_a:1\leq i\leq p,\ 1\leq a\leq n\}
\sqcup
\{-\epsilon_i+\delta_a:p<i\leq m,\ 1\leq a\leq n\}.
\end{equation}
We write
$
\Delta_{\mathrm{nst}}^+=\Delta_{\bar0}^+\sqcup\Delta_{\bar1,\mathrm{nst}}^+
$
for the corresponding positive system.

\subsection{Real forms}\label{subsec::real_forms}

A \emph{real form} of $\gg$ is a real Lie subsuperalgebra $\gg_{\RR}$ such that
$
\gg=\gg_{\RR}\oplus i\gg_{\RR}.
$
Equivalently, it is the fixed-point algebra of a parity-preserving conjugate-linear involutive automorphism. We use the equivalent description by a parity-preserving conjugate-linear anti-involution $\omega$, that is, $\omega^2=\operatorname{id}$ and $\omega([X,Y])=[\omega(Y),\omega(X)]$ for homogeneous $X,Y\in\gg$. The associated real form is
\begin{equation}
    \gg^\omega
    =
    \{X\in\gg_{\bar0}:\omega(X)=-X\}
    \oplus
    \{X\in\gg_{\bar1}:\omega(X)=iX\}.
\end{equation}
Indeed, the map $X\mapsto-i^{p(X)}\omega(X)$ is a parity-preserving conjugate-linear involutive automorphism. In particular, $\gg^\omega_{\bar0}$ is an ordinary real form of $\even$.

\begin{remark}\label{rmk::pairs_conjugation}
For a basic classical Lie superalgebra of type~II, the restriction of $\omega$ to $\even$ determines $\omega$ up to sign on $\odd$. Indeed, $\odd$ is a simple $\even$-module. If two anti-involutions $\omega$ and $\omega'$ agree on $\even$, then Schur's lemma gives
$
\omega'|_{\odd}=c\,\omega|_{\odd}
$
for some $c\in\CC^\times$. Since $[\odd,\odd]=\even$, compatibility with the bracket implies $c^2=1$. Hence
$
\omega'|_{\odd}=\pm\omega|_{\odd}.
$
\end{remark}

The real forms of the basic classical Lie superalgebras were classified by Parker~\cite{ParkerRealForms} and Serganova~\cite{serganova1983classification}. We recall the real forms relevant for unitarity in Theorem~\ref{theorem::unitary-highest-lowest-weight}. 

The classical families $\sl(m\vert n)$ with $m\neq n$, and $\osp(m\vert 2n)$ admit standard matrix realizations. The exceptional Lie superalgebras $D(2,1;\upalpha)$, $F(4)$, and $G(3)$ do not admit such standard defining realizations, although they admit faithful finite-dimensional representations.

For the classical families, fix the standard matrix realization of $\gg$. For a matrix $X \in \gg$, let $\overline{X}$ denote its complex conjugate and let $X^T$ denote its transpose. Set
$
X^\dagger\coloneqq\overline{X}^T.
$
Let $H$ be an even invertible Hermitian matrix satisfying $H^\dagger=H$ and
$
H^{-1}X^\dagger H\in\gg
$
for all $X\in\gg$. Then
\begin{equation}
\omega_H(X)\coloneqq H^{-1}X^\dagger H,
\qquad X\in\gg,
\end{equation}
defines a conjugate-linear anti-involution of $\gg$. Moreover, set
\begin{equation}
    I_{p,q}\coloneqq\operatorname{diag}(I_p,-I_q),
    \qquad
    \Omega_n\coloneqq
    \begin{pmatrix}0&I_n\\-I_n&0\end{pmatrix},
    \qquad
    R_m\coloneqq\bigoplus_{j=1}^{m}\Omega_1.
\end{equation}
For $\osp(m\vert2n)$, we use the realization preserving the even supersymmetric form with matrix $\operatorname{diag}(I_m,\Omega_n)$. The relevant classical real forms are listed below.

\begin{equation}\label{table::real_forms_classical}
\renewcommand{\arraystretch}{1.4}
\begin{array}{c|c|c|c}
    \gg
    & \text{real form}
    & H
    & \gg_{\bar0,\RR} \\ \hline
    \sl(m\vert n)
    & \begin{gathered}
        \su(p,q\vert0,n)\\
        \su(p,q\vert n,0)
      \end{gathered}
    & \begin{gathered}
        \operatorname{diag}(I_{p,q},-I_n)\\
        \operatorname{diag}(I_{p,q},I_n)
      \end{gathered}
    & \mathfrak{s}\bigl(\mathfrak{u}(p,q)\oplus\mathfrak{u}(n)\bigr)
    \\ \hline
    \osp(m\vert2n)
    & \begin{gathered}
        \osp(m,0\vert2n;\RR)\\
        \osp(0,m\vert2n;\RR)\\[1ex]
        \osp^{\ast}(m\vert2n,0)\\
        \osp^{\ast}(m\vert0,2n)
      \end{gathered}
    & \begin{gathered}
        \operatorname{diag}(I_m,i\Omega_n)\\
        \operatorname{diag}(-I_m,i\Omega_n)\\[1ex]
        \operatorname{diag}(iR_{m/2},I_{2n})\\
        \operatorname{diag}(iR_{m/2},-I_{2n})
      \end{gathered}
    & \begin{gathered}
        \so(m,\RR)\oplus\mathfrak{sp}(2n,\RR)
        \\[3ex]
        \so^{\ast}(m)\oplus\mathfrak{usp}(2n),
        \quad m\text{ even}
      \end{gathered}
\end{array}
\end{equation}

For the exceptional families, the real forms relevant for unitarity are determined up to isomorphism by their even parts, except for $D(2,1;\upalpha)$, where the position of the noncompact factor among the three simple factors of $\gg_{\bar0}$ must also be specified. Their even parts are
\begin{equation}\label{table::real_forms_exceptional}
\renewcommand{\arraystretch}{1.4}
\begin{array}{c|c}
    \gg  
    & \text{even part of the real form}\\ \hline
    F(4)
    & \begin{gathered}
        \so(7,\RR)\oplus\sl(2,\RR)
        \\
        \so(2,5)\oplus\su(2)
    \end{gathered}\\ \hline
    G(3)
    & \gg_{2,\mathrm{compact}}\oplus\sl(2,\RR)
      \\ \hline
    D(2,1;\upalpha)
    & \sl(2,\RR)\oplus
      \su(2)\oplus\su(2).
\end{array}
\end{equation}
Here $\gg_{2,\mathrm{compact}}$ denotes the compact real form of $\gg_2$. For $D(2,1;\upalpha)$, we assume $\upalpha\in\RR\setminus\{0,-1\}$. Its even part is the direct sum of three copies of $\sl(2,\mathbb C)$. Once these factors are labeled, the real form is specified by choosing which factor is $\sl(2,\RR)$. The standard parameter equivalences permute the three factors.

All these even real forms are reductive Lie algebras of Hermitian type. More precisely, for every noncompact simple factor $\mathfrak{s}_{\RR}$ with Cartan decomposition
$
\mathfrak{s}_{\RR}=\kk_{\RR}\oplus\mathfrak{p}_{\RR},
$
the maximal compact subalgebra $\kk_{\RR}$ has nonzero center. Equivalently, the associated Riemannian symmetric space carries an invariant complex structure. We use this fact throughout without further mention.

\subsection{Maximal compact subalgebras}\label{subsec::maximal_compact_subalgebra}

For a basic classical Lie superalgebra $\gg$, the even subalgebra $\even$ is reductive. Hence every real form $\gg_{\bar0,\RR}$ of $\even$ is a real reductive Lie algebra. Such a real form admits a \emph{Cartan involution}, that is, an involutive automorphism
$
\theta\colon\gg_{\bar0,\RR}\to\gg_{\bar0,\RR}
$
that acts as the identity on the center $\mathfrak{z}(\gg_{\bar0,\RR})$ and whose restriction to the semisimple ideal $[\gg_{\bar0,\RR},\gg_{\bar0,\RR}]$ is a Cartan involution in the sense of~\cite[Chapter~VI]{Knapp}. Its fixed-point subalgebra
$
\kk_{\RR}\coloneqq\{X\in\gg_{\bar0,\RR}:\theta(X)=X\}
$
is compactly embedded in $\gg_{\bar0,\RR}$. Thus, the subgroup of $\operatorname{Aut}(\gg_{\bar0,\RR})$ generated by $\exp(\ad(\kk_{\RR}))$ has compact closure. The subalgebra $\kk_{\RR}$ is called a \emph{maximal compact subalgebra} of $\gg_{\bar0,\RR}$.

Set
$
\pp_{\RR}\coloneqq\{X\in\gg_{\bar0,\RR}:\theta(X)=-X\}.
$
Then
$
\gg_{\bar0,\RR}=\kk_{\RR}\oplus\pp_{\RR}
$
is the corresponding Cartan decomposition. Let $\kk$ and $\pp$ denote the complexifications of $\kk_{\RR}$ and $\pp_{\RR}$, respectively. Complexification gives
\begin{equation}
\even=\kk\oplus\pp.
\end{equation}
We call $\kk$ the \emph{complexified maximal compact subalgebra} of $\even$. Since $\kk$ and $\pp$ are the $+1$- and $-1$-eigenspaces of $\theta$, respectively,
\begin{equation}
[\kk,\kk]\subseteq\kk,\qquad
[\kk,\pp]\subseteq\pp,\qquad
[\pp,\pp]\subseteq\kk.
\end{equation}

We fix a compactly embedded Cartan subalgebra $\hh_{\RR}\subseteq\kk_{\RR}$ of $\gg_{\bar0,\RR}$ and denote its complexification by $\hh$. By abuse of terminology, we also call $\hh$ compactly embedded. Then $\hh$ is a Cartan subalgebra of $\even$, and hence also of $\gg$. All root decompositions of $\gg$ are taken with respect to this choice of $\hh$. With respect to this choice, the real forms of the previous section are consistent in the sense of~\cite{jakobsen1994full}.

\begin{lemma}\label{lemm::consistent_omega}
For every real form in~\eqref{table::real_forms_classical} or~\eqref{table::real_forms_exceptional}, the associated conjugate-linear anti-involution $\omega$ satisfies $\omega(\gg^\alpha)=\gg^{-\alpha}$ for all $\alpha\in\Delta$.
\end{lemma}

\begin{proof}

We first show that $\alpha(\hh_{\RR})\subseteq i\RR$ for every $\alpha\in\Delta$. For even roots, this follows from the compact embeddedness of $\hh_{\RR}$. Suppose first that $\even$ is semisimple. Since $\odd$ is a finite-dimensional $\even$-module, every odd root is an integral weight of $\even$ and hence a rational linear combination of even roots. The claim therefore follows. The remaining cases are $\sl(m\vert n)$ and $\osp(2\vert2n)$. For $\sl(m\vert n)$, elements of $\hh_{\RR}$ are diagonal with entries $it_1,\ldots,it_m$ and $is_1,\ldots,is_n$, where $t_i,s_a\in\RR$. Thus every odd root $\epsilon_i-\delta_a$ takes the value $i(t_i-s_a)$. For $\osp(2\vert2n)$, the compact Cartan subalgebra acts on the natural orthogonal and symplectic modules with eigenvalues $\pm it$ and $\pm is_a$, respectively, where $t,s_a\in\RR$. Hence the odd roots $\pm\epsilon_1\pm\delta_a$ take values in $i\RR$ as well.

Let $X\in\gg^\alpha$ and $h\in\hh_{\RR}$. Then
\[ [h,\omega(X)]=\omega([h,X])=\overline{\alpha(h)}\omega(X)
=-\alpha(h)\omega(X),\]
where we used that $\omega$ is conjugate-linear and $\omega(h)=-h$ for all $h \in \hh_{\RR}$. We conclude $\omega(\gg^\alpha)\subseteq\gg^{-\alpha}$. Since $\omega^{2}=\id_{\gg}$, the statement follows.
\end{proof}

Since $\hh\subset \kk$ by assumption, we can define the compact root set $\Delta_c\coloneqq\Delta(\kk,\hh)$, which is a subset of $\Delta_{\bar{0}}$. Moreover, we have a decomposition $\Delta_{\bar{0}}=\Delta_{c}\sqcup \Delta_{n}$, where $\Delta_n\coloneqq\Delta_{\bar0}\setminus\Delta_c$ are the noncompact even roots. The Weyl group of $\Delta_{c}$ is denoted by $W_{c}$. It is naturally a subgroup of $W$.

A positive system $\Delta_{\bar0}^{+}\subseteq\Delta_{\bar0}$ is called $\kk$-adapted if $\Delta_n^+\coloneqq\Delta_n\cap\Delta_{\bar0}^+$ is $W_c$-invariant. This condition ensures that the positive noncompact root spaces form a $\kk$-stable subspace and that the corresponding Weyl vector $\rho_n$ is $W_c$-invariant. We fix such a positive system and set $\Delta_c^+\coloneqq\Delta_c\cap\Delta_{\bar0}^+$. The corresponding Weyl vectors are
\begin{equation*}
\rho_c\coloneqq\frac12\sum_{\alpha\in\Delta_c^+}\alpha,
\qquad
\rho_n\coloneqq\frac12\sum_{\alpha\in\Delta_n^+}\alpha.
\end{equation*}
Thus $\rho_{\bar0}=\rho_c+\rho_n$.

\subsection{Highest weight supermodules and atypicality}

Highest and lowest weight supermodules are the relevant simple modules for unitarity. Indeed, every simple unitarizable $\gg$-supermodule is either a highest or a lowest weight supermodule, as recalled in Section~\ref{subsec::unitary_representations_lie_supergroups}. We restrict to the highest weight case and recall the corresponding notion of atypicality.

Fix a positive system $\Delta^+$. Let $\nn^\pm$ be the sums of the positive and negative root spaces, that is, $\nn^{\pm} \coloneqq \sum_{\alpha \in \Delta^{+}}\gg^{\pm \alpha}$. We then have the triangular decomposition $\gg=\nn^{-}\oplus \hh \oplus \nn^{+}$. Moreover, we define the associated Borel subalgebra by $\bb \coloneqq \hh \oplus \nn^{+}$.

\begin{definition}
A $\gg$-supermodule $M$ is a \emph{highest weight supermodule} of highest weight $\Lambda\in\hh^{\ast}$ if there exists a nonzero homogeneous vector $v_\Lambda\in M$ such that $Xv_\Lambda=0$ for all $X\in\nn^+$, $hv_\Lambda=\Lambda(h)v_\Lambda$ for all $h\in\hh$, and $M=\UE(\gg)v_\Lambda$. The vector $v_\Lambda$ is called a \emph{highest weight vector}. 
\end{definition}

\begin{remark}
The parity of a highest weight vector is a convention. If $M$ is a highest weight supermodule with highest weight vector $v_\Lambda$, then $\Pi M$ is a highest weight supermodule of the same highest weight, with highest weight vector of opposite parity. Thus the highest weight alone determines a simple highest weight supermodule only up to parity change. We therefore choose $v_\Lambda$ to be even unless stated otherwise.
\end{remark}

We record the following standard properties of highest weight supermodules.

\begin{lemma}[{\cite[Section~8]{musson2012lie}}]\label{lemm::properties_HWM}
Let $M$ be a highest weight $\gg$-supermodule of highest weight $\Lambda$. Write $\Delta_{\bar0}^+=\{\alpha_1,\ldots,\alpha_k\}$ and $\Delta_{\bar1}^+=\{\beta_1,\ldots,\beta_l\}$, and choose nonzero root vectors $X_i\in\gg^{-\alpha_i}$ and $Y_j\in\gg^{-\beta_j}$. Then:
\begin{enumerate}
    \item[(a)] If $\gg\neq A(1\vert 1)$, the supermodule $M$ is spanned by
    \[
    X_1^{r_1}\cdots X_k^{r_k}Y_1^{s_1}\cdots Y_l^{s_l}v_\Lambda,
    \qquad
    r_i\in\ZZ_{\geq0},\quad s_j\in\{0,1\}.
    \]
    Such a vector has weight
    $
    \Lambda-\sum_{i=1}^k r_i\alpha_i-\sum_{j=1}^l s_j\beta_j
    $
    and parity
    $
    p(v_\Lambda)+\sum_{j=1}^l s_j
    $
    modulo $2$.

   \item[(b)] The supermodule $M$ is a weight supermodule,
$
M=\bigoplus_{\mu\in\hh^{\ast}}M^\mu,
$
where
$
M^\mu=\{v\in M:hv=\mu(h)v\text{ for all }h\in\hh\}.
$
Each weight space is finite-dimensional and $\dim M^\Lambda=1$.

    \item[(c)] Every nonzero quotient of $M$ is a highest weight supermodule of highest weight $\Lambda$. Moreover, $M$ has a unique maximal proper submodule and hence a unique simple quotient.

    \item[(d)] The supermodule $M$ is a quotient of the Verma supermodule
    $
    M(\Lambda)=\UE(\gg)\otimes_{\UE(\bb)}\CC_\Lambda,
    $
    where $\bb=\hh\oplus\nn^+$ and $\CC_\Lambda$ is the one-dimensional $\bb$-supermodule on which $\nn^+$ acts trivially and $\hh$ acts by $\Lambda$.

    \item[(e)] Two simple highest weight $\gg$-supermodules with the same highest weight and the same parity of their highest weight vectors are isomorphic.
\end{enumerate}
\end{lemma}

For each $\Lambda\in\hh^{\ast}$, let $L(\Lambda)$ denote the simple highest weight $\gg$-supermodule of highest weight $\Lambda$. If the choice of positive system is relevant, we write $L(\Lambda;\Delta^+)$. If the Lie superalgebra $\ll$ has to be specified, we write $L_{\ll}(\Lambda)$. Different positive systems can assign different highest weights to the same simple $\gg$-supermodule. Recall that we fix $\Delta_{\bar{0}}^{+}$ throughout. Two positive systems with the same even part are related by odd reflections (\emph{cf.}~\cite{Serganova_reps}), and we have the following result.

\begin{lemma}[{\cite[Lemma~1.40]{cheng2012dualities}}]\label{lemm::odd_reflections_HW}
Let $\Delta^+$ and $\widetilde{\Delta}^+$ be two positive systems with the same positive even roots. Then, for every $\Lambda\in\hh^{\ast}$, there exist $\Lambda'\in\hh^{\ast}$ and $\varepsilon\in\ZZ_2$ such that
\begin{equation*}
L(\Lambda;\Delta^+)
\cong
\Pi^\varepsilon L(\Lambda';\widetilde{\Delta}^+).
\end{equation*}
\end{lemma}

We finally discuss atypicality. Let $M$ be a simple highest weight $\gg$-supermodule of highest weight $\Lambda$. By~\cite[Section~8.2.4]{musson2012lie}, the center $\mathfrak{Z}(\gg)$ of $\UE(\gg)$ acts on $M$ by an algebra homomorphism
$
\chi_\Lambda:\mathfrak{Z}(\gg)\to\CC,
$
called the \emph{infinitesimal character} of $M$. Under the Harish-Chandra homomorphism, $\chi_\Lambda$ is determined by the shifted weight $\Lambda+\rho$. Set
\begin{equation}
A_\Lambda
\coloneqq
\{\alpha\in\overline{\Delta}_{\bar1}^+:(\alpha,\alpha)=0,\;(\Lambda+\rho,\alpha)=0\}.
\end{equation}
The weight $\Lambda$ is \emph{typical} if $A_\Lambda=\varnothing$ and \emph{atypical} otherwise. Its \emph{atypicality} is
\begin{equation}
\at(\Lambda)
=
\max\{r:\alpha_1,\ldots,\alpha_r\in A_\Lambda
\text{ are linearly independent and pairwise orthogonal}\}.
\end{equation}
Thus, $\Lambda$ is typical if and only if $\at(\Lambda)=0$. Since
$\at(\Lambda)\leq\defect(\gg)$, we make the following definition.

\begin{definition}
A weight $\Lambda$ is called \emph{maximally atypical} if
$\at(\Lambda)=\defect(\gg)$.
\end{definition}

Unlike for complex semisimple Lie algebras, the infinitesimal character is constant along certain isotropic odd directions. More precisely, let $\alpha_1,\ldots,\alpha_k\in A_\Lambda$ be linearly independent and pairwise orthogonal. For $t_1,\ldots,t_k\in\CC$ and $w\in W$, set
\begin{equation}
\Lambda'
=
w\left(\Lambda+\rho+\sum_{i=1}^k t_i\alpha_i\right)-\rho.
\end{equation}
Then
$
\chi_\Lambda=\chi_{\Lambda'}.
$

A simple highest weight $\gg$-supermodule is called typical or atypical according to its highest weight. Although $A_\Lambda$ depends on the choice of positive system,
the atypicality does not. Indeed, atypicality is determined
by the infinitesimal character
\cite[Lemmas~5.1--5.2]{duflo2005associated}. More precisely, let $\Delta^+$ and $\widetilde{\Delta}^+$ be two positive systems, and let $\Lambda$ and $\widetilde{\Lambda}$ be the corresponding highest weights of the same simple $\gg$-supermodule. Then
$
\at_{\Delta^+}(\Lambda)=\at_{\widetilde{\Delta}^+}(\widetilde{\Lambda}).
$ Thus atypicality is an invariant of the simple $\gg$-supermodule.

\subsection{Standing conventions}\label{subsec::conventions}
\noindent
We recall the following conventions, which remain in force throughout the article, except in Section~\ref{sec::equal_rank_type_A}, and will be used without further mention.

Throughout, $\gg$ is a basic classical Lie superalgebra with $\odd\neq0$, excluding type $A(n\vert n)$. For $\gg=\sl(m\vert n)$, we assume $m\neq n$ and $m+n>2$, and for $\gg=D(2,1;\upalpha)$, we assume $\upalpha\in\RR\setminus\{0,-1\}$. The defect-zero case $\gg=\osp(1\vert2n)$ is included.

We fix a conjugate-linear anti-involution $\omega$ of $\gg$ whose associated real form $\gg_{\RR}$ is one of the real forms in~\eqref{table::real_forms_classical} and~\eqref{table::real_forms_exceptional}. For the classical families, $\omega=\omega_H$ with $H$ as in~\eqref{table::real_forms_classical}. In particular, $\gg_{\bar0,\RR}$ is reductive of Hermitian type and has at most one noncompact simple ideal.

Let $\gg_{\bar0,\RR}=\kk_{\RR}\oplus\pp_{\RR}$ be the Cartan decomposition, and write $\kk$ and $\pp$ for the complexifications. We fix a Cartan subalgebra $\hh_{\RR}$ of $\gg_{\bar0,\RR}$ contained in $\kk_{\RR}$. In particular, $\hh_{\RR}$ is compactly embedded. We set $\hh\coloneqq\hh_{\RR}\otimes_{\RR}\CC=\hh_{\RR}\oplus i\hh_{\RR}$ and take all root decompositions with respect to $\hh$. Then
\begin{equation}
\omega|_{\hh_{\RR}}=-\mathrm{id},
\qquad
\omega|_{i\hh_{\RR}}=\mathrm{id},
\qquad
\omega(\gg^\alpha)=\gg^{-\alpha}
\quad(\alpha\in\Delta),
\end{equation}
where the last identity is Lemma~\ref{lemm::consistent_omega}. Moreover, we write $\Delta_c=\Delta(\kk,\hh)$ for the compact roots, $\Delta_n=\Delta_{\bar0}\setminus\Delta_c$ for the noncompact even roots, $W=W(\even,\hh)$ for the Weyl group, and $W_c=W(\kk,\hh)$ for the compact Weyl group.

For the chosen real form and Cartan subalgebra $\hh$, define
\begin{equation*}
i\hh_{\RR}^{\ast}
\coloneqq
\{\lambda\in\hh^{\ast}:\lambda(\hh_{\RR})\subseteq i\RR\}
=
\{\lambda\in\hh^{\ast}:\lambda(i\hh_{\RR})\subseteq\RR\}.
\end{equation*}
Then, since $\hh_{\RR}$ is compactly embedded, every root lies in $i\hh_{\RR}^{\ast}$. Moreover, the bilinear form $(\cdot,\cdot)$ is real on $i\hh_{\RR}$ and on $i\hh_{\RR}^{\ast}$, and $h_\alpha\in i\hh_{\RR}$ for every $\alpha\in\Delta$. Thus $(\lambda,\alpha)$ is real for all $\lambda\in i\hh_{\RR}^{\ast}$ and $\alpha\in\Delta$. In what follows, inequalities of the form $\alpha(H)>0$ are understood for $H\in i\hh_{\RR}$.

We fix a positive system $\Delta^+=\Delta_{\bar0}^+\sqcup\Delta_{\bar1}^+$ whose even part
\begin{equation}
\Delta_{\bar0}^+=\Delta_c^+\sqcup\Delta_n^+,
\qquad
\Delta_c^+\coloneqq\Delta_c\cap\Delta_{\bar0}^+,
\qquad
\Delta_n^+\coloneqq\Delta_n\cap\Delta_{\bar0}^+,
\end{equation}
is $\kk$-adapted. The even part $\Delta_{\bar0}^+$ and, unless stated otherwise, the odd part $\Delta_{\bar1}^+$ are specified in Section~\ref{subsec::type_A_structure} for $\sl(m\vert n)$ and in Appendix~\ref{appendix::structure_theory} for the other basic classical Lie superalgebras. The Weyl vectors $\rho$, $\rho_{\bar0}$, $\rho_{\bar1}$, $\rho_c$, and $\rho_n$ refer to these choices. Highest weight supermodules are taken with respect to $\Delta^+$, and highest weight vectors are chosen to be even.

\section{Unitarity for Lie superalgebras and Lie supergroups}\label{sec::unitarity}
\noindent
We introduce unitarizable $\gg$-supermodules and unitary representations of the associated Lie supergroups. We then discuss their properties and their relation through Harish-Chandra supermodules. The standard references are~\cite{jakobsen1994full, neeb2011lie, Carmeli_Fioresi_Varadarajan_Unitary}.

\subsection{Super Hilbert spaces}

Let $V=V_{\bar0}\oplus V_{\bar1}$ be a complex super vector space. A \emph{super Hermitian form} on $V$ is a sesquilinear form $\psi:V\times V\to\CC$, linear in the first variable and conjugate-linear in the second, such that
\begin{equation}
\psi(v,w)=(-1)^{p(v)p(w)}\overline{\psi(w,v)}
\end{equation}
for all homogeneous $v,w\in V$, and $\psi(V_{\bar0},V_{\bar1})=0$.

The restrictions
$
\psi_{\bar0}\coloneqq\psi|_{V_{\bar0}\times V_{\bar0}}
$
and
$
\psi_{\bar1}\coloneqq-i\psi|_{V_{\bar1}\times V_{\bar1}}
$
are ordinary Hermitian forms. Thus, $\psi=\psi_{\bar0}+i\psi_{\bar1}$. We call $\psi$ \emph{nondegenerate} if both forms are nondegenerate, and \emph{super positive definite} if both are positive definite. Equivalently,
$
\psi(v,v)>0
$
for all $0\neq v\in V_{\bar0}$ and
$
-i\psi(v,v)>0
$
for all $0\neq v\in V_{\bar1}$.

Associated with $\psi$ is the ordinary Hermitian form
\begin{equation}
\langle v,w\rangle
=
\begin{cases}
\psi(v,w),&v,w\in V_{\bar0},\\
-i\psi(v,w),&v,w\in V_{\bar1},\\
0,&p(v)\neq p(w).
\end{cases}
\end{equation}
The form $\langle\cdot,\cdot\rangle$ is nondegenerate, respectively positive definite, if and only if $\psi$ is nondegenerate, respectively super positive definite. Thus, super positive definite super Hermitian forms on $V$ correspond to positive definite Hermitian forms for which $V_{\bar0}$ and $V_{\bar1}$ are orthogonal.

\begin{definition}
A \emph{Hermitian superspace} is a complex super vector space $V=V_{\bar0}\oplus V_{\bar1}$ equipped with a super positive definite super Hermitian form $\psi$, or equivalently, with a positive definite Hermitian form $\langle\cdot,\cdot\rangle$ for which $V_{\bar0}$ and $V_{\bar1}$ are orthogonal.
\end{definition}

Throughout, we use the ordinary Hermitian form and regard a Hermitian superspace as a pair $(V,\langle\cdot,\cdot\rangle)$. The associated norm is
$
\|v\|\coloneqq\sqrt{\langle v,v\rangle}.
$

\begin{definition}
A \emph{pre-Hilbert superspace} is a Hermitian superspace $(\HH,\langle\cdot,\cdot\rangle)$. It is called a \emph{Hilbert superspace} if $\HH$ is complete with respect to the induced norm.
\end{definition}

We use $\HH$ for pre-Hilbert superspaces and $\calH$ for Hilbert superspaces to distinguish the two structures.

\subsection{Unitarizable supermodules}

Unitarizability is defined with respect to a real form of $\gg$, or equivalently, a conjugate-linear anti-involution $\omega:\gg\to\gg$. Recall that the corresponding real form is
$
\gg^\omega\coloneqq\{x\in\gg:\omega(x)=-(-i)^{p(x)}x\},
$
where $x$ is homogeneous. Thus, $\omega(x)=-x$ for $x\in\gg^\omega_{\bar0}$ and $\omega(x)=ix$ for $x\in\gg^\omega_{\bar1}$. We further assume that $\gg^\omega_{\bar0}$ contains the fixed compactly embedded Cartan subalgebra $\hh_{\RR}$, whose complexification is $\hh$. In particular, $\omega(\gg^{\alpha})=\gg^{-\alpha}$ for all $\alpha \in \Delta(\gg,\hh)$.  

\begin{definition}
An $\omega$-unitarizable $\gg$-supermodule is a $\gg$-supermodule $\HH$ equipped with a pre-Hilbert superspace structure $(\HH,\langle\cdot,\cdot\rangle)$ such that
\begin{equation*}
\langle xv,w\rangle=\langle v,\omega(x)w\rangle
\end{equation*}
for all $x\in\gg$ and $v,w\in\HH$. The Hermitian form $\langle\cdot,\cdot\rangle$ is called $\omega$-contravariant.
\end{definition}

\begin{remark}
Let $\HH$ be an $\omega$-unitarizable $\gg$-supermodule. For $x\in\gg^\omega_{\bar0}$, one has
$
\langle xv,w\rangle=-\langle v,xw\rangle.
$
Thus, the even part of the real form acts by skew-Hermitian operators. For $x\in\gg^\omega_{\bar1}$, one has
$
\langle xv,w\rangle=-i\langle v,xw\rangle.
$
Equivalently, $e^{i\tfrac{\pi}{4}}x$ acts by a Hermitian operator.
\end{remark}

The anti-involution $\omega$ extends uniquely to a conjugate-linear anti-involution of $\UE(\gg)$, denoted by the same symbol. In particular,
$
\omega(XY)=\omega(Y)\omega(X)
$
for all $X,Y\in\UE(\gg)$. In particular,
$
\langle Xv,w\rangle=\langle v,\omega(X)w\rangle
$
for all $X\in\UE(\gg)$ and $v,w\in\HH$.

\begin{lemma}\label{lemm::properties_unitarizable_supermodules}
Let $\HH$ be an $\omega$-unitarizable $\gg$-supermodule.
\begin{enumerate}
    \item[(a)] If $\ll\subseteq\gg$ is an $\omega$-stable Lie subsuperalgebra, then $\HH$ is an $\omega\vert_{\ll}$-unitarizable $\ll$-supermodule.

    \item[(b)] If $M\subseteq\HH$ is a $\gg$-submodule, then its orthogonal complement $M^\perp$ is also a $\gg$-submodule.

    \item[(c)] If $\HH$ is finite-dimensional, then $\HH$ is completely reducible.
\end{enumerate}
\end{lemma}

\begin{proof}
Part~(a) follows directly from the definition. For~(b), let $v\in M^\perp$, $w\in M$, and $x\in\gg$. Since $\omega(x)w\in M$, we have
$
\langle xv,w\rangle=\langle v,\omega(x)w\rangle=0.
$
Hence $xv\in M^\perp$. Finally, (c) follows from~(b) by induction on $\dim\HH$.
\end{proof}

Unitarity imposes strong restrictions on simple $\gg$-supermodules. In particular, simple unitarizable supermodules associated with the real forms considered above are of highest or lowest weight type.

\begin{theorem}[{\cite[Section~7]{neeb2011lie}, \cite{jakobsen1994full}}]
\label{theorem::unitary-highest-lowest-weight}
Let $\gg$ be a basic classical Lie superalgebra with nonzero odd part. Then the following statements hold.
\begin{enumerate}
\item[(a)]
Every nontrivial simple unitarizable supermodule associated with $\su(p,q\vert0,n)$, $\osp(0,m\vert2n;\RR)$, or $\osp^{\ast}(2r\vert0,2n)$ is of highest weight type. Those associated with $\su(p,q\vert n,0)$, $\osp(m,0\vert2n;\RR)$, or $\osp^{\ast}(2r\vert2n,0)$ are of lowest weight type. In type $A$, we assume $pq\neq0$. If $pq=0$, every simple unitarizable supermodule is finite-dimensional and has both a highest and a lowest weight.

\item[(b)]
Let $\gg$ be $F(4)$, $G(3)$, or $D(2,1;\upalpha)$ with $\upalpha\in\RR\setminus\{0,-1\}$, and fix an even real form $\gg_{\bar0,\RR}\supseteq\hh_{\RR}$ from~\eqref{table::real_forms_exceptional}, including the position of the noncompact factor for $D(2,1;\upalpha)$. Consider the associated pair of conjugate-linear anti-involutions $\{\omega,(-1)^{p(\cdot)}\omega\}$ of Remark~\ref{rmk::pairs_conjugation}. One member, denoted $\omega^{\mathrm{hw}}$, has the property that every nontrivial simple $\omega^{\mathrm{hw}}$-unitarizable supermodule is of highest weight type. For the other member, $\omega^{\mathrm{lw}}\coloneqq(-1)^{p(\cdot)}\omega^{\mathrm{hw}}$, every nontrivial simple $\omega^{\mathrm{lw}}$-unitarizable supermodule is of lowest weight type.
\item[(c)]
Every real form admitting a nontrivial simple unitarizable supermodule is isomorphic to one of the real forms in~\eqref{table::real_forms_classical} or~\eqref{table::real_forms_exceptional}.

\item[(d)]
Up to parity reversal, the simple unitarizable highest and lowest weight supermodules above, together with the trivial supermodules, exhaust all simple unitarizable supermodules.
\end{enumerate}
\end{theorem}

\begin{proof}
We exclude the remaining real forms of $G(3)$ and
$D(2,1;\upalpha)$, which are omitted from
\cite[Theorem~6.2.1]{neeb2011lie}.
The case of $F(4)$ was treated in
\cite{Schmidt_unitary_branching}.

Following~\cite{neeb2011lie}, we first show that every noncompact
simple ideal of the even real algebra $\gg_{\bar{0},\RR}$ must be of
Hermitian type. The statement follows by~\cite[Sections~5.3 and~6.1]{neeb2011lie}, if we show that the closed convex cone
\begin{equation*}
C\coloneqq\overline{\operatorname{Cone}\{[x,x]:x\in\gg_{\bar1,\RR}\}}\subseteq\gg_{\bar0,\RR}
\end{equation*}
is pointed, generating, and invariant. 

The cone is clearly invariant under the adjoint action of the connected even real Lie group with Lie algebra $\gg_{\bar0,\RR}$ and generating, since $[\odd,\odd]=\even$, as $G(3)$ and $D(2,1;\upalpha)$ are basic classical.

It remains to prove that $C$ is pointed, that is, $C\cap(-C)=\{0\}$. By~\cite{jakobsen1994full}, there exists a nontrivial simple unitarizable $\gg$-supermodule $(\tau,\HH)$. Fix such a supermodule. Since $\gg$ is simple, $\tau$ is faithful. Assume $[x,x] \in C\cap (-C)$ for some $x \in \odd$. Since $x \in \odd$, the operator $e^{i\pi/4}\tau(x)$ is symmetric, and therefore $i\tau([x,x])=2(e^{i\pi/4}\tau(x))^2$ is positive semidefinite. We conclude that both $i\tau([x,x])$ and $-i\tau([x,x])$ are positive semidefinite. This is only true if $[x,x]=0$, which shows that $C$ is pointed.

This excludes all non-Hermitian real forms. By Serganova's classification~\cite{serganova1983classification}, the only remaining real form of $G(3)$ has even part $\mathfrak{g}_{2,\mathrm{compact}}\oplus\sl(2,\RR)$. 

The Hermitian real forms of $D(2,1;\upalpha)$ not listed in~\eqref{table::real_forms_exceptional} have even parts $\sl(2,\RR)^{\oplus3}$. In this case, the odd part is $V_1\otimes_{\RR}V_2\otimes_{\RR}V_3$, where each $V_i$ is a two-dimensional real symplectic module. Choose nonzero vectors $v_i\in V_i$ and set $u=v_1\otimes v_2\otimes v_3$. The defining bracket gives $[u,u]=0$, since each summand contains an alternating pairing of a vector with itself. Thus $T=e^{i\pi/4}\tau(u)$ is symmetric and satisfies $T^2=\frac{i}{2}\tau([u,u])=0$. Consequently,
\begin{equation*}
\langle Tv,Tv\rangle=\langle v,T^2v\rangle=0,
\qquad v\in\HH.
\end{equation*}
Hence $\tau(u)=0$, contradicting faithfulness. 
\end{proof}

Unitarizable highest and lowest weight supermodules are related by a duality~\cite{Schmidt_unitary_branching}. More precisely, $\omega$-unitarizable highest weight $\gg$-supermodules of highest weight $\Lambda$ correspond bijectively to $\omega'$-unitarizable lowest weight supermodules of lowest weight $-\Lambda$, where $\omega'=(-1)^{p(\cdot)}\omega$. We therefore restrict our attention to highest weight supermodules and write $L(\Lambda)$ for the simple highest weight $\gg$-supermodule of highest weight $\Lambda$. When $\omega$ is clear from the context, we simply speak of unitarizable highest weight supermodules.

Of particular importance is the restriction to the $\omega$-stable Lie subalgebra $\even$. We denote the simple highest weight $\even$-module of highest weight $\mu\in\hh^{\ast}$ by $L_0(\mu)$.

\begin{lemma}\label{lemm::restriction_even}
Let $L(\Lambda)$ be an $\omega$-unitarizable highest weight $\gg$-supermodule. Then its restriction to $\even$ is a finite direct sum of unitarizable highest weight $\even$-modules.
\end{lemma}

\begin{proof} 
The Lie subalgebra $\even\subset\gg$ is $\omega$-stable. By Lemma~\ref{lemm::properties_unitarizable_supermodules}, $L(\Lambda)$ is completely reducible as a $\even$-module upon forgetting parity. Every simple constituent is a highest weight module, since the weights of $L(\Lambda)$ are bounded above by $\Lambda$ (Lemma~\ref{lemm::properties_HWM}). Finally, the PBW theorem implies that $\UE(\gg)$ is finitely generated over $\UE(\even)$, so only finitely many simple constituents occur.
\end{proof}

\subsection{Unitary representations of Lie supergroups}
\label{subsec::unitary_representations_lie_supergroups}

Having established unitarity for Lie superalgebras, we now turn to unitarity for real Lie supergroups. By a real Lie supergroup $\calG_{\RR}$, we mean a real super Harish-Chandra pair $(G_{\RR},\gg_{\RR})$, as defined below. We restrict to finite-dimensional real Lie superalgebras whose complexifications are basic classical.

\subsubsection{Real super Harish-Chandra pairs}

A \emph{real super Harish-Chandra pair} is a pair
$
\calG_{\RR}=(G_{\RR},\gg_{\RR}),
$
where $G_{\RR}$ is a real Lie group and $\gg_{\RR}=\gg_{\bar0,\RR}\oplus\gg_{\bar1,\RR}$ is a real Lie superalgebra, together with a smooth action
$
\Ad:G_{\RR}\to\operatorname{Aut}(\gg_{\RR})
$
such that $\Lie(G_{\RR})=\gg_{\bar0,\RR}$, the restriction of $\Ad$ to $\gg_{\bar0,\RR}$ is the usual adjoint action of $G_{\RR}$, and
\begin{equation}
\left.\frac{\mathrm d}{\mathrm dt}\right|_{t=0}
\Ad(\exp(tx))y
=
[x,y]
\end{equation}
for all $x\in\gg_{\bar0,\RR}$ and $y\in\gg_{\RR}$. We assume throughout that $G_{\RR}$ is connected.

A morphism of real super Harish-Chandra pairs
$
(\phi,\varphi):(G_{\RR},\gg_{\RR})\to(G'_{\RR},\gg'_{\RR})
$
consists of a Lie group homomorphism $\phi:G_{\RR}\to G'_{\RR}$ and an even Lie superalgebra homomorphism $\varphi:\gg_{\RR}\to\gg'_{\RR}$ satisfying
$
\varphi|_{\gg_{\bar0,\RR}}=\mathrm d\phi
$
and
\begin{equation}
\varphi(\Ad(g)x)=\Ad(\phi(g))\varphi(x)
\end{equation}
for all $g\in G_{\RR}$ and $x\in\gg_{\RR}$.

We assume that the complexification of $\gg_{\RR}$ is $\gg$, that is,
$
\gg_{\RR}\otimes_{\RR}\CC\cong\gg.
$
Throughout, we consider ordinary real forms in the sense of~\cite{FioresiGavariniRealForms}, rather than graded real forms. The real form $\gg_{\RR}$ determines a conjugate-linear anti-involution $\omega$ of $\gg$.

The adjoint action of $G_{\RR}$ extends complex-linearly to $\gg$ and satisfies
\begin{equation}
\omega(\Ad(g)x)=\Ad(g)\omega(x)
\end{equation}
for all $g\in G_{\RR}$ and $x\in\gg$. On the group, the corresponding anti-involution is inversion,
$
\omega_G(g)=g^{-1}.
$
Its differential agrees with $\omega$ on $\gg_{\bar0,\RR}$.

\subsubsection{Unitary representations of \texorpdfstring{$\calG_{\RR}$}{G}}

Let $\calG_{\RR}=(G_{\RR},\gg_{\RR})$ be a real Lie supergroup, and let $\gg$ be the complexification of $\gg_{\RR}$. Recall that a Hilbert superspace is a $\ZZ_2$-graded complex Hilbert space
$
\calH=\calH_{\bar0}\oplus\calH_{\bar1},
$
where $\calH_{\bar0}$ and $\calH_{\bar1}$ are closed orthogonal subspaces. 

For a bounded operator $T$ on $\calH$, let $T^{\ast}$ denote its Hilbert space adjoint, characterized by
$
\langle Tv,w\rangle=\langle v,T^{\ast}w\rangle.
$
The operator $T$ is \emph{unitary} if $T^{\ast}=T^{-1}$. The unitary operators form the group $U(\calH)$, equipped with the strong operator topology, which agrees with the weak operator topology on $U(\calH)$.

A \emph{unitary representation} of $G_{\RR}$ on $\calH$ is a continuous group homomorphism
$
\pi:G_{\RR}\to U(\calH)
$
such that $\pi(g)$ is even for every $g\in G_{\RR}$.

Let $(\pi,\calH)$ be such a representation. Its space of smooth vectors is the dense graded subspace
$
\calH^\infty
=
\{v\in\calH:g\mapsto\pi(g)v\text{ is smooth}\}.
$
The derived representation of $\gg_{\bar0,\RR}=\Lie(G_{\RR})$ is defined on $\calH^\infty$ by
\begin{equation}
d\pi(X)v
=
\left.
\frac{\mathrm d}{\mathrm dt}
\right|_{t=0}
\pi(\exp(tX))v,
\qquad
X\in\gg_{\bar0,\RR},\quad v\in\calH^\infty.
\end{equation}

We extend this notion to real Lie supergroups.

\begin{definition}
A \emph{unitary representation} of the real Lie supergroup $\calG_{\RR}=(G_{\RR},\gg_{\RR})$ on a Hilbert superspace $\calH$ is a pair $\Sigma=(\pi,\tau)$, where $\pi:G_{\RR}\to U(\calH)$ is a continuous unitary representation by even operators and
$
\tau:\gg_{\RR}\to\End_{\CC}(\calH^\infty)
$
is a parity-preserving real-linear representation of the Lie superalgebra $\gg_{\RR}$ satisfying
\begin{equation*}
\tau|_{\gg_{\bar0,\RR}}=d\pi,
\qquad
\tau(\Ad(g)X)=\pi(g)\tau(X)\pi(g)^{-1}
\end{equation*}
for all $g\in G_{\RR}$ and $X\in\gg_{\RR}$. Moreover, the operators satisfy
\begin{equation*}
\langle\tau(X)v,w\rangle
=
\begin{cases}
-\langle v,\tau(X)w\rangle,
& X\in\gg_{\bar0,\RR},\\
-i\langle v,\tau(X)w\rangle,
& X\in\gg_{\bar1,\RR},
\end{cases}
\end{equation*}
for all $v,w\in\calH^\infty$.
\end{definition}

\begin{remark}
The condition for $X\in\gg_{\bar0,\RR}$ follows from $\tau(X)=d\pi(X)$ and the unitarity of $\pi$. For $X\in\gg_{\bar1,\RR}$, it is equivalent to requiring that $e^{i\tfrac{\pi}{4}}\tau(X)$ be symmetric on $\calH^\infty$.
\end{remark}

The representation $\tau$ extends uniquely by complex linearity to $\gg$ and, by the universal property of the universal enveloping superalgebra, to a unital algebra representation
$
\tau:\UE(\gg)\to\End_{\CC}(\calH^\infty).
$
We use the same notation for all these representations. 

Moreover, the conjugate-linear anti-involution $\omega$ determined by $\gg_{\RR}$ extends uniquely to $\UE(\gg)$ by
$
\omega(XY)=\omega(Y)\omega(X)
$ for all $X,Y\in \UE(\gg)$. Unitarity is then equivalent to
\begin{equation}
\langle\tau(D)v,w\rangle
=
\langle v,\tau(\omega(D))w\rangle,
\qquad
D\in\UE(\gg),\quad v,w\in\calH^\infty.
\end{equation}
In particular, $\calH^\infty$ is an $\omega$-unitarizable $\gg$-supermodule.

Every unitary representation $\calH$ of a real Lie supergroup has an algebraic counterpart, called its associated Harish-Chandra supermodule.

\subsection{Harish-Chandra supermodules}\label{subsec::HC_supermodules}
Let $\calG_{\RR}=(G_{\RR},\gg_{\RR})$ be a real Lie supergroup with $G_{\RR}$ connected and $\gg_{\bar0,\RR}$ reductive. Let $\kk_{\RR}\subseteq\gg_{\bar0,\RR}$ be a maximal compact subalgebra, and let $\kk$ denote its complexification, as in Section~\ref{subsec::maximal_compact_subalgebra}.

A $(\gg,\kk)$-supermodule is a $\gg$-supermodule $V$ whose restriction to $\kk$ is locally finite and semisimple. Local finiteness means that $\UE(\kk)v$ is finite-dimensional for every $v\in V$, while semisimplicity means that $V$ is a direct sum of simple $\kk$-modules. It is called a \emph{Harish-Chandra supermodule} if it is finitely generated over $\UE(\gg)$~\cite{Carmeli_Fioresi_Varadarajan_HW,Alldridge}. We describe how irreducible unitary representations of $\calG_{\RR}$ yield Harish-Chandra supermodules.

Let $\Sigma=(\pi,\tau)$ be a unitary representation of $\calG_{\RR}$ on $\calH$. A vector $v\in\calH^\infty$ is called \emph{$\kk$-finite} if the $\kk$-submodule $\UE(\kk)v$ is finite-dimensional. We denote the space of smooth $\kk$-finite vectors by
\begin{equation}
\calH^{[\kk]}
\coloneqq
\left\{
v\in\calH^\infty:
\dim_{\CC}\UE(\kk)v<\infty
\right\}.
\end{equation}

\begin{lemma}
Let $\Sigma=(\pi,\tau)$ be a unitary representation of $\calG_{\RR}$ on $\calH$. Then $\calH^{[\kk]}$ is a $(\gg,\kk)$-supermodule.
\end{lemma}

\begin{proof} By definition, $\calH^{[\kk]}$ is locally finite and semisimple as a $\kk$-module. It remains to prove that $\calH^{[\kk]}$ is also a $\gg$-supermodule. 

Let $x\in\gg$ and $v\in\calH^{[\kk]}$. We show that $\tau(x)v$ is $\kk$-finite. For $y\in\kk$, $z\in\gg$, and $w\in\UE(\kk)v$, we have
\[
\tau(y)\tau(z)w=\tau(z)\tau(y)w+\tau([y,z])w\in\tau(\gg)\UE(\kk)v.
\]
Thus $\tau(\gg)\UE(\kk)v$ is a $\kk$-submodule. It is finite-dimensional since $v$ is $\kk$-finite and $\gg$ is finite-dimensional. In particular, $\tau(x)v$ is $\kk$-finite.
\end{proof}

For an irreducible unitary representation, $\calH^{[\kk]}$ is a Harish-Chandra supermodule and is dense in $\calH$ by the classical density theorem for $\kk$-finite vectors (\emph{cf.}~\cite{Alldridge}). We collect the corresponding results.

\begin{proposition}
Let $\Sigma=(\pi,\tau)$ be a unitary representation of $\calG_{\RR}$ on $\calH$. Then:
\begin{enumerate}
    \item[(a)] If $\Sigma$ is irreducible, then $\calH^{[\kk]}$ is a dense, simple Harish-Chandra supermodule.

    \item[(b)] If $\calH^{[\kk]}$ is dense in $\calH$, then $\Sigma$ is irreducible if and only if $\calH^{[\kk]}$ is simple.

    \item[(c)] Two irreducible unitary representations
    $\Sigma_1=(\pi_1,\tau_1)$ on $\calH_1$ and
    $\Sigma_2=(\pi_2,\tau_2)$ on $\calH_2$
    are unitarily equivalent if and only if their Harish-Chandra supermodules $\calH_1^{[\kk]}$ and $\calH_2^{[\kk]}$ are isomorphic as $\gg$-supermodules.
\end{enumerate}
\end{proposition}

\begin{proof}
Recall that we assume throughout that $G_{\RR}$ is connected, $\gg$ is basic classical, and $\odd\neq0$. 

(a) The assertion is immediate for the trivial representation. Otherwise, the kernel of the complexified derived representation is a graded ideal of $\gg$. 
By simplicity and connectedness, the complexified derived representation is faithful. Hence $\calG_{\RR}$ is $\star$-reduced in the sense of~\cite{neeb2011lie}, that is, its unitary representations separate the elements of $\gg_{\RR}$. Moreover,
$
[\gg_{\bar1,\RR},\gg_{\bar1,\RR}]=\gg_{\bar0,\RR},
$
since $\odd+[\odd,\odd]$ is a nonzero graded ideal. 
Let $C\subseteq\gg_{\bar0,\RR}$ be the closed convex invariant cone generated by $[x,x]$ for $x\in\gg_{\bar1,\RR}$. Since $C$ is pointed and generating and $\hh_{\RR}$ is compactly embedded, \cite[Lemma~5.1.3]{neeb2011lie} implies that $\hh_{\RR}\cap\operatorname{Int}(C)$ contains a regular element. Thus~\cite[Theorem~7.3.2]{neeb2011lie} applies to our fixed Cartan subalgebra and shows that
\[
\calH^{[\hh]}
=
\{v\in\calH^\infty:
\dim_{\CC}\UE(\hh)v<\infty\}
\]
is a dense, simple, unitarizable $\gg$-supermodule.

By Theorem~\ref{theorem::unitary-highest-lowest-weight}, $\calH^{[\hh]}$ is a highest or lowest weight supermodule. Its restriction to $\even$ is a finite direct sum of simple unitarizable highest or lowest weight modules. Moreover, each summand is $\kk$-locally finite and has finite-dimensional $\kk$-isotypic components. Hence the same holds for $\calH^{[\hh]}$, and consequently, $\calH^{[\hh]} \subseteq\calH^{[\kk]}$. The reverse inclusion follows from $\hh\subseteq\kk$. Altogether, we proved that
$
\calH^{[\kk]}=\calH^{[\hh]}
$
is dense, simple, and $\kk$-semisimple. Since every nonzero homogeneous vector generates $\calH^{[\hh]}$, it is also finitely generated over $\UE(\gg)$. Hence $\calH^{[\kk]}$ is a dense, simple Harish-Chandra supermodule.

(b) Assume $\calH^{[\kk]}$ is dense in $\calH$. If $\Sigma$ is irreducible, then we already proved in (a) that $\calH^{[\kk]}$ is simple. Conversely, suppose that $\calH^{[\kk]}$ is dense and simple. Let $\mathcal{W}\subseteq\calH$ be a closed $\ZZ_2$-graded invariant subspace, and let $P$ be the orthogonal projection onto $\mathcal{W}$. Since $\Sigma$ is unitary and $\mathcal{W}$ is invariant, its orthogonal complement is also a $\ZZ_2$-graded invariant subspace. In particular, $P$ intertwines $\Sigma$ and preserves $\calH^{[\kk]}$. The restriction of $P$ to $\calH^{[\kk]}$ is an even $\gg$-endomorphism satisfying $P^2=P$. Its kernel and image are therefore $\ZZ_2$-graded invariant submodules. Since $\calH^{[\kk]}$ is simple, $P=0$ or $P=\id_{\calH^{[\kk]}}$. Since $P$ is bounded and $\calH^{[\kk]}$ is dense, this equality extends to $\calH$. Hence $\mathcal{W}=0$ or $\mathcal{W}=\calH$, proving that $\Sigma$ is irreducible.

(c) A unitary equivalence restricts to an isomorphism of the $\kk$-finite smooth supermodules. Conversely, such an isomorphism identifies their highest weight spaces for the common positive system supplied by~\cite[Theorem~7.3.2]{neeb2011lie}. The uniqueness assertion in part~(iv) of that theorem gives a unitary equivalence of the original representations.
\end{proof}

Let $\widetilde{G}_{\RR}$ be the connected, simply connected Lie group with Lie algebra $\gg_{\bar0,\RR}$. The adjoint action of $\gg_{\bar0,\RR}$ on $\gg_{\RR}$ integrates uniquely to an action of $\widetilde{G}_{\RR}$. We call the resulting super Harish-Chandra pair
$
\widetilde{\calG}_{\RR}=(\widetilde{G}_{\RR},\gg_{\RR})
$
the \emph{simply connected Lie supergroup associated with $\gg_{\RR}$}.

The following theorem identifies simple unitarizable highest weight $\gg$-supermodules with Harish-Chandra supermodules of unitary representations.

\begin{theorem}[{\cite[Theorems~2 and~5]{Carmeli_Fioresi_Varadarajan_Unitary}}]
\label{thm::HW_supergroup_HC}
Let $L(\Lambda)$ be a simple highest weight $\gg$-supermodule. Then $L(\Lambda)$ is $\omega$-unitarizable if and only if there exists an irreducible unitary highest weight representation $\Sigma_{\Lambda}=(\pi_{\Lambda},\tau_{\Lambda})$ of $\widetilde{\calG}_{\RR}$ on $\calH_{\Lambda}$ such that
\begin{equation}
\calH_{\Lambda}^{[\kk]}
\cong
L(\Lambda)
\end{equation}
as $\gg$-supermodules. The representation $\Sigma_{\Lambda}$ is unique up to unitary equivalence.
\end{theorem}

    \begin{remark}
\label{rem::even_globalization}
Theorem~\ref{thm::HW_supergroup_HC} is the super analogue of a classical result for the
even part, which we use repeatedly. Let $L_{0}(\Lambda)$ be an $\omega|_{\even}$-unitarizable
highest weight $\even$-module with respect to $\Delta_{\bar0}^{+}$. Then there exists an
irreducible unitary representation $(\pi_{\Lambda},\calH_{\Lambda})$ of $\widetilde{G}_{\RR}$,
unique up to unitary equivalence, such that
\[
\calH_{\Lambda}^{[\kk]}\cong L_{0}(\Lambda)
\]
as $\even$-modules; see~\cite[Theorem~X.2.8]{NeebHW} or~\cite[Theorem~6]{Carmeli_Fioresi_Varadarajan_Unitary}. 
\end{remark}

\begin{remark}
Let $G_{\RR}$ be a connected Lie group with Lie algebra $\gg_{\bar0,\RR}$, and let $p_{\RR}:\widetilde G_{\RR}\to G_{\RR}$ be the universal covering homomorphism. The representation $\Sigma_{\Lambda}$ descends to a unitary representation of $\calG_{\RR}=(G_{\RR},\gg_{\RR})$ if and only if $\ker p_{\RR}$ acts trivially on $\calH_{\Lambda}$.
\end{remark}

Motivated by Theorem~\ref{thm::HW_supergroup_HC}, we make the following definition.

\begin{definition}
An irreducible unitary representation $\Sigma=(\pi,\tau)$ of $\calG_{\RR}$ on $\calH$ is called a \emph{highest weight representation} if its associated Harish-Chandra supermodule $\calH^{[\kk]}$ is a highest weight $\gg$-supermodule.
\end{definition}

The finite branching of unitarizable highest weight supermodules has the following consequence for unitary representations.

\begin{corollary}\label{cor::finite_decomposition_group_level}
Let $\Sigma=(\pi,\tau)$ be an irreducible unitary highest weight representation of $\calG_{\RR}$ on $\calH$. Then its restriction to $G_{\RR}$ is a finite direct sum of irreducible unitary representations.
\end{corollary}

\begin{proof}
By Lemma~\ref{lemm::restriction_even}, $\calH^{[\kk]}=V_1\oplus\cdots\oplus V_r$ with mutually orthogonal simple unitarizable highest weight $\even$-modules $V_i$, since the complements constructed in its proof are orthogonal. Let $\calH_i$ be the closure of $V_i$. As $\calH^{[\kk]}$ is dense, $\calH=\calH_1\oplus\cdots\oplus\calH_r$.

For $v\in V_i$ and $w\in V_j$ with $i\neq j$, the function $g\mapsto\langle\pi(g)v,w\rangle$ is analytic, the vectors of $\calH^{[\kk]}$ being analytic, and all its derivatives at the identity vanish because $\UE(\even)V_i=V_i\perp V_j$. Since $G_{\RR}$ is connected, it vanishes identically. Hence each $\calH_i$ is $G_{\RR}$-invariant, and $\calH_i^{[\kk]}=\calH^{[\kk]}\cap\calH_i=V_i$.

The orthogonal projection onto a closed $G_{\RR}$-invariant subspace of $\calH_i$ commutes with $\pi$ and restricts to an idempotent endomorphism of the simple $\even$-module $V_i$. It is therefore $0$ or the identity on $V_i$, and by density on $\calH_i$.
\end{proof}

\section{Relative holomorphic discrete series modules}
\label{sec::relative_holomorphic_discrete_series}
\noindent
Fix a real form $\gg_{\RR}=\gg_{\bar0,\RR}\oplus\gg_{\bar1,\RR}$ of a basic classical Lie superalgebra $\gg$, as described in Section~\ref{subsec::real_forms}, with associated conjugate-linear anti-involution $\omega$. Let $\widetilde{G}_{\RR}$ be the connected, simply connected Lie group with Lie algebra $\gg_{\bar0,\RR}$, and let $\widetilde{\calG}_{\RR}=(\widetilde{G}_{\RR},\gg_{\RR})$ be the associated simply connected Lie supergroup.

By Theorem~\ref{thm::HW_supergroup_HC}, every simple $\omega$-unitarizable highest weight $\gg$-supermodule $L(\Lambda)$ is the Harish-Chandra supermodule of an irreducible unitary representation $\Sigma$ of $\widetilde{\calG}_{\RR}$. By Corollary~\ref{cor::finite_decomposition_group_level}, the restriction of $\Sigma$ to $\widetilde{G}_{\RR}$ is a finite direct sum of irreducible unitary representations whose Harish-Chandra modules are unitarizable highest weight $\even$-modules. We introduce the relative holomorphic discrete series among these representations and introduce a numerical invariant, called the \emph{formal dimension or Harish-Chandra degree}, which serves as an analogue of dimension. We follow~\cite{NeebSquareIntegrable}.

\subsection{Relative holomorphic discrete series}
\label{subsec::relative_holomorphic_discrete_series}

Let $L_{0}(\Lambda)$ be a unitarizable highest weight $\even$-module. By unitary globalization (Remark~\ref{rem::even_globalization}), there exists an irreducible unitary representation $(\pi_{\Lambda},\calH_{\Lambda})$ of $\widetilde{G}_{\RR}$ such that
\begin{equation}
\calH_{\Lambda}^{[\kk]}
\cong
L_{0}(\Lambda)
\end{equation}
as $\even$-modules. The \emph{relative holomorphic discrete series} consists of those irreducible unitary highest weight representations $(\pi_{\Lambda},\calH_{\Lambda})$ whose matrix coefficients are square-integrable modulo the center of $\widetilde{G}_{\RR}$.

Let $Z\coloneqq Z(\widetilde{G}_{\RR})$ be the center of $\widetilde{G}_{\RR}$. By Schur's lemma, there exists a continuous unitary character
$
\chi_{\pi_{\Lambda}}:Z\to U(1)
$
such that
$
\pi_{\Lambda}(z)=\chi_{\pi_{\Lambda}}(z)\id_{\calH_{\Lambda}}
$
for all $z\in Z$. We call $\chi_{\pi_{\Lambda}}$ the \emph{central character} of $\pi_{\Lambda}$.

Set $X\coloneqq Z\backslash\widetilde{G}_{\RR}=\widetilde{G}_{\RR}/Z$. The central character defines a homogeneous line bundle
$
E_{\pi_{\Lambda}}\coloneqq\widetilde{G}_{\RR}\times_Z\CC
$
over $X$, where $Z$ acts by
$
z\cdot(g,\zeta)\coloneqq(zg,\chi_{\pi_{\Lambda}}(z)\zeta).
$
Writing $[g,\zeta]$ for the corresponding orbit, the projection is $[g,\zeta]\mapsto Zg$, and the right action of $\widetilde{G}_{\RR}$ is given by
$
[g,\zeta]\cdot x\coloneqq[gx,\zeta].
$
Since $\chi_{\pi_{\Lambda}}$ is unitary, the standard Hermitian metric on $\CC$ induces a $\widetilde{G}_{\RR}$-invariant Hermitian metric on $E_{\pi_{\Lambda}}$.

Let $\calC(\widetilde{G}_{\RR},\CC)$ denote the space of continuous complex-valued functions on $\widetilde{G}_{\RR}$. The continuous sections of $E_{\pi_{\Lambda}}$ identify with
\begin{equation}
\Gamma_{\widetilde{G}_{\RR}}(E_{\pi_{\Lambda}})
\coloneqq
\left\{
f\in\calC(\widetilde{G}_{\RR},\CC):
f(zg)=\chi_{\pi_{\Lambda}}(z)f(g)
\text{ for all }z\in Z,\ g\in\widetilde{G}_{\RR}
\right\}.
\end{equation}
Under this identification, $f$ corresponds to the section
$
\sigma_f(Zg)\coloneqq[g,f(g)],
$
and $\widetilde{G}_{\RR}$ acts by right translation,
$
(x\cdot f)(g)\coloneqq f(gx).
$

The groups $\widetilde{G}_{\RR}$ and $X$ are unimodular, as discussed in~\cite[Remark~3.5 and proof of Proposition~2.2]{NeebSquareIntegrable}. Fix a Haar measure $\mathrm d\mu_X$ on $X$. We denote by $\Gamma_{\widetilde{G}_{\RR}}^{2}(E_{\pi_{\Lambda}})$ the Hilbert space of measurable square-integrable sections of $E_{\pi_{\Lambda}}$, identified up to equality almost everywhere. Equivalently,
\begin{equation}
\Gamma_{\widetilde{G}_{\RR}}^{2}(E_{\pi_{\Lambda}})
=
\left\{
f:\widetilde{G}_{\RR}\to\CC\text{ measurable}:
f(zg)=\chi_{\pi_{\Lambda}}(z)f(g),\
\int_X |f(g)|^2\,\mathrm d\mu_X(Zg)<\infty
\right\},
\end{equation}
where functions are identified if they agree almost everywhere. Since $\chi_{\pi_{\Lambda}}$ is unitary, $|f(g)|$ depends only on $Zg\in X$. The Hermitian form is
\begin{equation}
\langle f,h\rangle
\coloneqq
\int_X f(g)\overline{h(g)}\,\mathrm d\mu_X(Zg).
\end{equation}
The right regular action of $\widetilde{G}_{\RR}$ on $\Gamma_{\widetilde{G}_{\RR}}^{2}(E_{\pi_{\Lambda}})$ is unitary.

\begin{definition}[{\cite[Definition~2.1]{NeebSquareIntegrable}}]
An irreducible unitary highest weight representation $(\pi_{\Lambda},\calH_{\Lambda})$ of $\widetilde{G}_{\RR}$ belongs to the \emph{relative holomorphic discrete series} if there exist nonzero vectors $v,w\in\calH_{\Lambda}$ such that
\begin{equation}
\int_X
\abs{\langle\pi_{\Lambda}(g)v,w\rangle}^{2}
\,\mathrm d\mu_X(Zg)
<\infty.
\end{equation}
\end{definition}

\begin{remark}
The adjective \emph{holomorphic} reflects the fact that these representations admit realizations by holomorphic functions, as established in~\cite{HC_VI}.
\end{remark}

For $v\in\calH_{\Lambda}$, define
\begin{equation}
\Psi_v:\calH_{\Lambda}\longrightarrow
\Gamma_{\widetilde{G}_{\RR}}(E_{\pi_{\Lambda}}),
\qquad
\Psi_v(u)(g)\coloneqq
\langle u,\pi_{\Lambda}(g^{-1})v\rangle.
\end{equation}
Since the Hermitian form is conjugate-linear in the second variable and $\chi_{\pi_{\Lambda}}$ is unitary,
$
\Psi_v(u)(zg)=\chi_{\pi_{\Lambda}}(z)\Psi_v(u)(g)
$
for all $z\in Z$. Thus $\Psi_v(u)$ defines a continuous section of $E_{\pi_{\Lambda}}$. Moreover,
$
\Psi_v(\pi_{\Lambda}(x)u)(g)=\Psi_v(u)(gx)
$
for all $x,g\in\widetilde{G}_{\RR}$. Hence $\Psi_v$ is $\widetilde{G}_{\RR}$-equivariant.

\begin{theorem}[{\cite[Theorem~1.2 and Propositions~2.2--2.3]{NeebSquareIntegrable}}]
\label{thm::embedding}
Suppose that $(\pi_{\Lambda},\calH_{\Lambda})$ belongs to the relative holomorphic discrete series. Then there exists a constant $d(\pi_{\Lambda})>0$, called the \emph{formal dimension or Harish-Chandra degree} of $\pi_{\Lambda}$, such that for every unit vector $v\in\calH_{\Lambda}$,
\begin{equation*}
\sqrt{d(\pi_{\Lambda})}\Psi_v:
\calH_{\Lambda}
\longrightarrow
\Gamma_{\widetilde{G}_{\RR}}^{2}(E_{\pi_{\Lambda}})
\end{equation*}
is an isometric $\widetilde{G}_{\RR}$-equivariant embedding. In particular, $\pi_{\Lambda}$ occurs as a unitary direct summand of $\Gamma_{\widetilde{G}_{\RR}}^{2}(E_{\pi_{\Lambda}})$. The formal dimension depends on the chosen Haar measure on $X$.
\end{theorem}

\begin{remark}
Let $G_{\RR}'$ be another connected Lie group with Lie algebra $\gg_{\bar0,\RR}$, and let
$
p_{\RR}:\widetilde{G}_{\RR}\to G_{\RR}'
$
be the universal covering homomorphism. Suppose that $\pi_{\Lambda}$ descends to an irreducible unitary representation $\pi_{\Lambda}'$ of $G_{\RR}'$. By~\cite[Remark~3.6]{NeebSquareIntegrable}, square-integrability modulo the center is invariant under coverings. Thus $\pi_{\Lambda}$ belongs to the relative holomorphic discrete series of $\widetilde{G}_{\RR}$ if and only if $\pi_{\Lambda}'$ belongs to the relative holomorphic discrete series of $G_{\RR}'$. If $Z(G_{\RR}')$ is compact, this is equivalent to square-integrability on $G_{\RR}'$. In this case, $\pi_{\Lambda}'$ belongs to the holomorphic discrete series of $G_{\RR}'$.
\end{remark}

Harish-Chandra's calculation gives an explicit criterion for membership in the relative holomorphic discrete series. Let $v_{\Lambda}$ be a unit highest weight vector. The criterion characterizes the convergence of
\begin{equation}
\int_X
\abs{\langle\pi_{\Lambda}(g)v_{\Lambda},v_{\Lambda}\rangle}^{2}
\,\mathrm{d}\mu_X(Zg)
\end{equation}
in terms of the highest weight $\Lambda$ and the positive noncompact roots.

\begin{theorem}[{\cite[Lemmas~27 and~29]{HC_VI}, \cite[Theorem~3.10]{NeebSquareIntegrable}}]
\label{thm::HC_condition}
Let $(\pi_{\Lambda},\calH_{\Lambda})$ be an irreducible unitary highest weight representation of $\widetilde{G}_{\RR}$ with highest weight $\Lambda$. Then $(\pi_{\Lambda},\calH_{\Lambda})$ belongs to the relative holomorphic discrete series if and only if
\begin{equation*}
2\frac{(\Lambda+\rho_{\bar0},\beta)}{(\beta,\beta)}<0
\end{equation*}
for every $\beta\in\Delta_n^+$.
\end{theorem}

\begin{remark}
The inequalities in Theorem~\ref{thm::HC_condition} are known as the \emph{Harish-Chandra condition}.
\end{remark}

By Remark~\ref{rem::even_globalization}, every unitarizable highest weight $\even$-module $L_{0}(\Lambda)$ admits an irreducible unitary globalization $\calH_{\Lambda}$ to $\widetilde{G}_{\RR}$, such that $\calH_{\Lambda}^{[\kk]} \cong L_{0}(\Lambda)$. This yields the following definition. 

\begin{definition} A unitarizable highest weight $\even$-module $L_{0}(\Lambda)$ is called a \emph{relative holomorphic discrete series module} if it is isomorphic to the Harish-Chandra module of a relative holomorphic discrete series representation of $\widetilde{G}_{\RR}$.
\end{definition}

\subsection{Harish-Chandra parametrization}
\label{subsec::HC_parameterization}
We parametrize the relative holomorphic discrete series modules. To do so, we regard $i\hh_{\RR}^{\ast}$ as a real form of $\hh^{\ast}$. 

The set of \emph{compactly integral weights} is
\begin{equation}
P_{c}
\coloneqq
\left\{
\lambda\in i\hh_{\RR}^{\ast}:
2\frac{(\lambda,\alpha)}{(\alpha,\alpha)}\in\ZZ
\text{ for all }\alpha\in\Delta_{c}
\right\}.
\end{equation}
Since $\widetilde{G}_{\RR}$ is simply connected, these are precisely the integrality conditions associated with the compact roots.

Fix the $\kk$-adapted positive system
$
\Delta_{\bar0}^{+}=\Delta_{c}^{+}\sqcup\Delta_{n}^{+}
$
introduced above. Define
\begin{equation}
C_{\hol}
\coloneqq
\left\{
\lambda\in i\hh_{\RR}^{\ast}:
\frac{(\lambda,\alpha)}{(\alpha,\alpha)}>0
\text{ for all }\alpha\in\Delta_{c}^{+},
\quad
\frac{(\lambda,\beta)}{(\beta,\beta)}<0
\text{ for all }\beta\in\Delta_{n}^{+}
\right\}.
\end{equation}
Since $\Delta_{\bar0}^{+}$ is $\kk$-adapted, these inequalities define a Weyl chamber. We call the chambers in the $W_{c}$-orbit of $C_{\hol}$ \emph{holomorphic Weyl chambers}.

Set
\begin{equation}
\mathcal{D}_{\hol}
\coloneqq
P_{c}\cap\bigcup_{w\in W_{c}}wC_{\hol}.
\end{equation}
Every weight in $\mathcal{D}_{\hol}$ is regular, that is, $(\Lambda,\alpha)\neq 0$ for all $\alpha \in \Delta_{\bar{0}}$ if $\Lambda \in \mathcal{D}_{\hol}$. The isomorphism classes of relative holomorphic discrete series representations of $\widetilde{G}_{\RR}$ are parametrized by
\begin{equation}
\mathcal{D}_{\hol}/W_{c}.
\end{equation}
More precisely, each $W_{c}$-orbit has a unique representative $\lambda\in C_{\hol}$. The corresponding representation has highest weight
$
\Lambda=\lambda-\rho_{\bar0}.
$
Equivalently, a highest weight $\Lambda$ belongs to the relative holomorphic discrete series if and only if it is dominant integral with respect to $\Delta_{c}^{+}$ and satisfies the Harish-Chandra condition
$
2\frac{(\Lambda+\rho_{\bar0},\beta)}{(\beta,\beta)}<0
$
for every $\beta\in\Delta_{n}^{+}$.

\subsection{Formal dimension}
\label{subsec::formal_dimension}

The positive constant $d(\pi_{\Lambda})$ in Theorem~\ref{thm::embedding} is called the \emph{Harish-Chandra degree}, or \emph{formal dimension}, of the relative holomorphic discrete series representation $(\pi_{\Lambda},\calH_{\Lambda})$. It is characterized by the following orthogonality relations.

\begin{proposition}[{\cite[Proposition~2.3]{NeebSquareIntegrable}}]
\label{prop::orthogonality_relations}
Let $(\pi_i,\calH_i)$, $i=1,2$, be relative holomorphic discrete series representations of $\widetilde{G}_{\RR}$ with the same central character. Then
\[
\int_X
\langle\pi_1(g)v,w\rangle
\overline{\langle\pi_2(g)v',w'\rangle}
\,\mathrm{d}\mu_X(Zg)
=
\begin{cases}
0, & \pi_1\ncong\pi_2,\\ d(\pi_\Lambda)^{-1}\langle w',w\rangle\langle v,v'\rangle,
& \pi_1=\pi_2=\pi_\Lambda,
\end{cases}
\]
for $v,w\in\calH_1$ and $v',w'\in\calH_2$.
\end{proposition}

Setting $v'=v$ and $w'=w$ yields
\begin{equation}
\int_X
\abs{\langle\pi_\Lambda(g)v,w\rangle}^{2}
\,\mathrm{d}\mu_X(Zg)
=
\frac{\norm{v}^{2}\norm{w}^{2}}{d(\pi_\Lambda)}.
\end{equation}
Replacing $\mathrm{d}\mu_X$ by $c\,\mathrm{d}\mu_X$, $c>0$, replaces
$d(\pi_\Lambda)$ by $c^{-1}d(\pi_\Lambda)$. For a suitable normalization
of $\mathrm{d}\mu_X$, the formal dimension is given by the following
product formula.

\begin{theorem}[{\cite[Theorem~3.17]{NeebSquareIntegrable}}]
\label{thm::formal_dimension}
There exists a normalization of the Haar measure $\mathrm{d}\mu_X$ on
$X=\widetilde{G}_{\RR}/Z$, depending only on $\widetilde{G}_{\RR}$, such that
for every relative holomorphic discrete series representation
$(\pi_{\Lambda},\calH_{\Lambda})$ of $\widetilde{G}_{\RR}$,
\begin{equation*}
d(\pi_{\Lambda})
=
\prod_{\alpha\in\Delta_{c}^{+}}
\frac{(\Lambda+\rho_{c},\alpha)}
     {(\rho_{c},\alpha)}
\prod_{\beta\in\Delta_{n}^{+}}
\frac{\abs{(\Lambda+\rho_{\bar0},\beta)}}
     {\abs{(\rho_{\bar0},\beta)}}.
\end{equation*}
\end{theorem}

\begin{remark}
The normalization in Theorem~\ref{thm::formal_dimension} does not depend on $\Lambda$. Since replacing $\mathrm{d}\mu_X$ by $c\,\mathrm{d}\mu_X$, $c>0$, multiplies every formal dimension by $c^{-1}$, the product formula holds for an arbitrary normalization up to a positive factor independent of $\Lambda$. In particular, the ratios $d(\pi_{\Lambda})/d(\pi_{\Lambda'})$ do not depend on the normalization.
\end{remark}

\begin{remark}
The irreducible $\kk$-module $F_{\Lambda}=\UE(\kk)v_{\Lambda}\subset L_{0}(\Lambda)$ generated by the highest weight vector $v_{\Lambda}$ of $L_0(\Lambda)$ has highest weight $\Lambda$ and dimension
\begin{equation*}
\dim F_{\Lambda} = \prod_{\alpha\in\Delta_c^+}
\frac{(\Lambda+\rho_c,\alpha)}{(\rho_c,\alpha)}.\end{equation*}
By Theorem~\ref{thm::formal_dimension}, $d(\pi_{\Lambda})$ is the product of $\dim F_{\Lambda}$ and
$\prod_{\beta\in\Delta_n^+}\abs{(\Lambda+\rho_{\bar0},\beta)}/\abs{(\rho_{\bar0},\beta)}$.
\end{remark}

\section{Relative holomorphic discrete series supermodules}\label{sec::relative_HDS_supermodules}
\noindent
For a reductive Lie group, the relative holomorphic discrete series is characterized by the square-integrability of matrix coefficients modulo the center. We prove the analogous characterization for relative holomorphic discrete series supermodules and obtain an $L^2$-realization of their unitary globalizations.

Throughout this section, $\widetilde{\calG}_{\RR}=(\widetilde{G}_{\RR},\gg_{\RR})$ is the simply connected Lie supergroup associated with $\gg_{\RR}$, and $\Sigma=(\pi,\tau)$ is an irreducible unitary representation of $\widetilde{\calG}_{\RR}$ on a Hilbert superspace $\calH$, with simple Harish-Chandra supermodule $\calH^{[\kk]}$ (Section~\ref{subsec::HC_supermodules}).

Recall that a unitarizable highest weight $\gg$-supermodule has only finitely many simple $\even$-constituents.

\begin{definition}
A unitarizable highest weight $\gg$-supermodule $\HH$ is a \emph{relative holomorphic discrete series supermodule} if each simple $\even$-constituent of $\HH$ is the Harish-Chandra module of a relative holomorphic discrete series representation of $\widetilde{G}_{\RR}$.
\end{definition}

The main result of this section, Theorem~\ref{thm::super_embedding}, expresses this condition through $\Sigma$. The supermodule $\calH^{[\kk]}$ is a relative holomorphic discrete series supermodule if and only if the supermatrix coefficients $\varphi_{v,w}$, for $v,w\in\calH^\infty$ homogeneous of the same parity, are square-integrable modulo the center of $\widetilde{\calG}_{\RR}$. 
We therefore first introduce the notion of supermatrix coefficients.

\subsection{Supermatrix coefficients}

Following~\cite[Section~6.2]{NeebSalmasian}, we equip $\widetilde{G}_{\RR}\times\UE(\gg)$ with the involutive monoid structure
\begin{equation}
(g_1,D_1)(g_2,D_2)
=
(g_1g_2,(g_2^{-1}\cdot D_1)D_2),
\qquad
(g,D)^\ast
=
(g^{-1},g\cdot \omega(D)),
\end{equation}
where $g\cdot D$ denotes the adjoint action of $\widetilde{G}_{\RR}$ on $\UE(\gg)$. Here $D\mapsto \omega(D)$ is the conjugate-linear anti-involution of $\UE(\gg)$ extending $\omega$, characterized by
\begin{equation}
\omega(D_1D_2)=\omega(D_2)\omega(D_1),
\qquad
\omega(1)=1,
\end{equation}
for $D_1,D_2\in\UE(\gg)$. In particular, our convention gives $\omega(x)=-x$ for $x\in\gg_{\bar0,\RR}$ and $\omega(x)=ix$ for $x\in\gg_{\bar1,\RR}$. Our convention differs from that of~\cite{NeebSalmasian} by the sign of the anti-involution on the odd part. The two conventions are related by complex conjugation of representations: if $\tau(x)^\ast=i\tau(x)$, then $\overline{\tau(x)}^\ast=-i\overline{\tau(x)}$ for $x\in\gg_{\bar1,\RR}$.

A \emph{smooth superfunction} on $\widetilde{\calG}_{\RR}$ is a $\CC$-linear map
$f:\UE(\gg)\to C^\infty(\widetilde{G}_{\RR},\CC)$
satisfying the conditions of~\cite[Section~6.3]{NeebSalmasian}. In particular,
\[
f(xD)=L_xf(D),
\qquad
x\in\even,\quad D\in\UE(\gg),
\]
where $L_x$ is the left-invariant differential operator on $\widetilde{G}_{\RR}$ associated with $x$. For $x\in\gg_{\bar0,\RR}$, it is given by
$
(L_xh)(g)=\left.\frac{\mathrm d}{\mathrm dt}\right|_{t=0}h(g\exp(tx)),
$
and extends $\CC$-linearly to $\even$. Moreover, for every $n\geq0$, the map
\[
(x_1,\ldots,x_n,g)
\longmapsto
f(x_1\cdots x_n)(g)
\]
from $\gg^n\times\widetilde{G}_{\RR}$ to $\CC$ is smooth.

We denote the space of smooth superfunctions by $C^\infty(\widetilde{\calG}_{\RR})$. For $f\in C^\infty(\widetilde{\calG}_{\RR})$, define
$
\check f(g,D)\coloneqq f(D)(g).
$
This identifies smooth superfunctions with functions on $\widetilde{G}_{\RR}\times\UE(\gg)$ that are $\CC$-linear in the second variable and satisfy the preceding equivariance and smoothness conditions.

\begin{definition}
For $v,w\in\calH^\infty$, the \emph{supermatrix coefficient} of $\Sigma$ associated with $v,w$ is the function
\begin{equation*}
\varphi_{v,w}:\widetilde{G}_{\RR}\times\UE(\gg)\longrightarrow\CC,
\qquad
\varphi_{v,w}(g,D)
\coloneqq
\langle\pi(g)\tau(D)v,w\rangle.
\end{equation*}
\end{definition}

By~\cite[Definition~6.8]{NeebSalmasian}, the supermatrix coefficients satisfy
$
\varphi_{v,w}(s)=\overline{\varphi_{w,v}(s^\ast)}
$
for every $s\in\widetilde{G}_{\RR}\times\UE(\gg)$. In particular, $\varphi_{v,v}$ is positive definite.

The argument of~\cite[Proposition~6.9]{NeebSalmasian} extends to homogeneous vectors $v,w\in\calH^\infty$ of arbitrary parity. Indeed, the defining equivariance condition holds, and for every $n\geq0$, the map
\begin{equation}
(x_1,\ldots,x_n,g)\longmapsto
\langle\pi(g)\tau(x_1\cdots x_n)v,w\rangle
\end{equation}
is smooth. Consequently, there exists a smooth superfunction $f_{v,w}\in C^\infty(\widetilde{\calG}_{\RR})$ of parity $p(v)+p(w)$ such that $\check f_{v,w}=\varphi_{v,w}$. By linearity, this extends to arbitrary $v,w\in\calH^\infty$.

\subsection{A Hilbert superspace of equivariant superfunctions}

We distinguish the center of the Lie supergroup from that of its
underlying even group. Let $\calZ$ denote the center of
$\widetilde{\calG}_{\RR}$, and set
\begin{equation}
Z\coloneqq Z(\widetilde{G}_{\RR}),
\qquad
C\coloneqq
\{z\in Z:\operatorname{Ad}(z)|_{\gg}=\id_{\gg}\}.
\end{equation}

Then $C$ is the underlying group of $\calZ$. The Lie superalgebra of $\calZ$ is the center of $\gg_{\RR}$, which is zero since $\gg$ is simple. Hence $\calZ=(C,0)$ is discrete and purely even. The center $Z$ of $\widetilde{G}_{\RR}$, by contrast, need not be discrete, and its elements may act nontrivially on $\odd$ and on $\UE(\gg)$. In this section, square-integrability modulo the center always refers to $\calZ$.

Let $c\in C$. Then $\Ad(c)$ is trivial on $\gg_{\RR}$, and the even operator $\pi(c)$ commutes with $\pi(\widetilde{G}_{\RR})$ and $\tau(\gg)$. As $\Sigma$ is irreducible, Schur's lemma gives a unitary character $\chi_C\colon C\to U(1)$ with
\begin{equation}
\pi(c)=\chi_C(c)\id_{\calH},
\qquad c\in C.
\end{equation}
Hence
\begin{equation}\label{eq::C_equivariance}
\varphi_{v,w}(gc,D)=\chi_C(c)\,\varphi_{v,w}(g,D),
\qquad
g\in\widetilde{G}_{\RR},\quad c\in C,\quad D\in\UE(\gg).
\end{equation}
The character $\chi_C$ is the same for all $v,w$, whereas $Z$ may act on different even constituents of $\calH$ by different characters.

\begin{definition}
Let $C^\infty_{\chi_C}(\widetilde{\calG}_{\RR})$ be the space of all $f\in C^\infty(\widetilde{\calG}_{\RR})$ with $\check f(gc,D)=\chi_C(c)\check f(g,D)$ for $g\in\widetilde{G}_{\RR}$, $c\in C$, and $D\in\UE(\gg)$.
\end{definition}

Since $C$ is central in $\widetilde{G}_{\RR}$, left and right $C$-equivariance agree. By~\eqref{eq::C_equivariance}, every $f_{v,w}$ with $v,w\in\calH^\infty$ lies in $C^\infty_{\chi_C}(\widetilde{\calG}_{\RR})$. In general, such superfunctions do not descend to $\widetilde{\calG}_{\RR}/\calZ$ but are sections of the line bundle associated with $\chi_C$.

Fix a homogeneous basis $(X_{a})_{a\in A}$ of $\bigwedge\gg_{\bar1}$, and let $\beta\colon\bigwedge\gg_{\bar1}\to\UE(\gg)$ be the PBW antisymmetrization map. Since $\UE(\gg)=\UE(\even)\,\beta\bigl(\bigwedge\gg_{\bar1}\bigr)$ and $f(xD)=L_xf(D)$ for $x\in\even$, a smooth superfunction $f$ is determined by the smooth functions
\begin{equation}\label{eq::coefficient_functions}
f_a(g)\coloneqq\check f(g,\beta(X_{a})),
\qquad a\in A.
\end{equation}
Let
\begin{equation}
E_{\chi_C}
\coloneqq
\widetilde{G}_{\RR}\times_C\CC
\longrightarrow
\widetilde{G}_{\RR}/C,
\qquad
(gc,\chi_C(c)t)\sim(g,t)
\quad(c\in C,\ t\in\CC).
\end{equation}
For $f\in C^\infty_{\chi_C}(\widetilde{\calG}_{\RR})$ and $D\in\UE(\gg)$, one has $\check f(gc,D)=\chi_C(c)\check f(g,D)$, so $gC\mapsto[g,\check f(g,D)]$ is a well-defined section of $E_{\chi_C}$; in particular, so is each $f_a$. As $\chi_C$ is unitary, the standard Hermitian metric on $\CC$ induces one on $E_{\chi_C}$, and $|\check f(g,D)|^2$ is a function on $\widetilde{G}_{\RR}/C$.

Let $\mathscr D_{\chi_C}(\widetilde{\calG}_{\RR})$ be the subspace of those $f\in C^\infty_{\chi_C}(\widetilde{\calG}_{\RR})$ for which the supports of all $f_a$ have compact image in $\widetilde{G}_{\RR}/C$. This condition is independent of the basis.

Let $\mathcal M\coloneqq\widetilde{\calG}_{\RR}/\calZ$. As $\calZ$ is closed, central, and discrete, $\mathcal M$ is a Lie supergroup with underlying group $\widetilde{G}_{\RR}/C$ and Lie superalgebra $\gg_{\RR}$. The fiber of its Berezinian line bundle at the identity coset is $\operatorname{Ber}(\gg_{\RR}^{\ast})$, and $C$ acts on it trivially, through $\Ad$. By~\cite[Theorem~4.13]{Alldridge_Hilgert}, $\mathcal M$ admits a nonzero $\widetilde{\calG}_{\RR}$-invariant Berezinian section, unique up to a scalar. We fix one and denote it by $\omega_{\mathcal M}$.

For $f,h\in\mathscr D_{\chi_C}(\widetilde{\calG}_{\RR})$,
the product $\overline f h$ is invariant under $\calZ$, since
$
\overline{\chi_C(c)}\chi_C(c)=1$ for all $c\in C$.
It therefore descends to a compactly supported smooth
superfunction on $\mathcal M$. Here the bar denotes conjugation
of complex superfunctions on the real Lie supergroup. We define
the \emph{Berezin pairing} by
\begin{equation}
\label{eq::Berezin_L2_pairing}
(f,h)_{\mathrm{Ber}}
\coloneqq
\int_{\mathcal M}
\overline f h\,\omega_{\mathcal M},
\qquad
f,h\in\mathscr D_{\chi_C}(\widetilde{\calG}_{\RR}).
\end{equation}

Suppose that $\odd\neq0$, and let $q=\dim\odd$. Then it is well known that the Berezin pairing is not positive definite. Indeed, choose a chart of $\mathcal M$ with odd coordinates $\theta_1,\ldots,\theta_q$, take the conjugation with $\overline{fh}=\overline{f}\,\overline{h}$ and $\overline{\theta_j}=\theta_j$, and let $\varphi\neq0$ be a smooth function with compact support in the chart. In a local trivialization of the line bundle associated with $\chi_C$, the superfunction $f=\varphi\theta_1\cdots\theta_q$, extended by zero, is a nonzero element of $\mathscr D_{\chi_C}(\widetilde{\calG}_{\RR})$ with $\overline{f}f=0$. Hence $(f,f)_{\mathrm{Ber}}=0$ for every invariant Berezinian, and the Berezin pairing does not define a Hilbert norm.

We therefore define a positive $L^2$-space through~\eqref{eq::coefficient_functions}. Fix a Haar measure $\mathrm{d}(gC)$ on $\widetilde{G}_{\RR}/C$. For $f\in\mathscr D_{\chi_C}(\widetilde{\calG}_{\RR})$, we define
\begin{equation}
\label{eq::positive_super_L2_right}
\|f\|_2^2
\coloneqq
\sum_{a\in A}
\int_{\widetilde{G}_{\RR}/C}
\left|
\check f\bigl(g,\operatorname{Ad}(g^{-1})\beta(X_{a})\bigr)
\right|^2\,\mathrm{d}(gC).
\end{equation}
Each integrand is well-defined on $\widetilde{G}_{\RR}/C$, since $\operatorname{Ad}(c)=\id$ for $c\in C$ and $\check f$ transforms under $C$ by the unitary character $\chi_C$. Moreover, 
for $f\in\mathscr D_{\chi_C}(\widetilde{\calG}_{\RR})$, each integrand is continuous with compact support on $\widetilde{G}_{\RR}/C$, so its integral is finite.

\begin{lemma}\label{lemm::norm}
The expression~\eqref{eq::positive_super_L2_right} defines a norm on $\mathscr D_{\chi_C}(\widetilde{\calG}_{\RR})$. Different choices of homogeneous basis of $\bigwedge(\odd)$ yield equivalent norms.
\end{lemma}
\begin{proof}
It is enough to prove positive definiteness. If $\|f\|_2=0$, continuity implies
\[
\check f\bigl(g,\operatorname{Ad}(g^{-1})\beta(X_a)\bigr)=0
\]
for every $g$ and $a\in A$. 
By $\operatorname{Ad}$-equivariance of $\beta$, the elements $\operatorname{Ad}(g^{-1})\beta(X_a)$ form a basis of $\beta(\bigwedge\odd)$ for every $g$. Their vanishing under $\check f(g,\cdot)$ therefore implies
$
\check f(g,\beta(X_a))=0
$
for all $g$ and $a\in A$. However, these functions determine $f$, giving $f=0$.

Finally, let $(Y_a)_{a\in A}$ be another homogeneous basis of $\bigwedge\odd$, with $Y_a=\sum_{b\in A}M_{ab}X_b$ for an invertible matrix $M$. By linearity of $\beta$ and $\check f$, the corresponding coefficient functions satisfy
\[
\check f\bigl(g,\operatorname{Ad}(g^{-1})\beta(Y_a)\bigr)
=
\sum_{b\in A}M_{ab}\check f\bigl(g,\operatorname{Ad}(g^{-1})\beta(X_b)\bigr).
\]
The equivalence of Euclidean norms on $\CC^{|A|}$ gives constants $c_{1},c_{2}>0$, independent of $g$, such that
\[
c_{1}\sum_{a\in A}\left|\check f\bigl(g,\operatorname{Ad}(g^{-1})\beta(X_a)\bigr)\right|^2
\leq
\sum_{a\in A}\left|\check f\bigl(g,\operatorname{Ad}(g^{-1})\beta(Y_a)\bigr)\right|^2
\leq
c_{2}\sum_{a\in A}\left|\check f\bigl(g,\operatorname{Ad}(g^{-1})\beta(X_a)\bigr)\right|^2.
\]
Integrating over $\widetilde{G}_{\RR}/C$ proves the equivalence of the resulting norms.
\end{proof}

We are now ready to define the appropriate space.

\begin{definition}
The \emph{positive $L^2$-space of $\chi_C$-equivariant superfunctions} is the Hilbert completion
\begin{equation*}
L^2_+(\widetilde{\calG}_{\RR}/\calZ;\chi_C)
\coloneqq
\overline{
\mathscr D_{\chi_C}(\widetilde{\calG}_{\RR})
}^{\,\|\cdot\|_2}.
\end{equation*}
\end{definition}

\begin{remark}\label{rmk::right_translations_are_unitary}
\leavevmode
\begin{enumerate}
\item[(a)] 
The natural $\ZZ_{2}$-grading of $\bigwedge\odd$ defines an orthogonal $\mathbb Z_2$-grading on $L^2_+(\widetilde{\calG}_{\RR}/\calZ;\chi_C)$.
 
\item[(b)]
The even right regular action satisfies
$
\check{(R_hf)}(g,D)
=
\check f\bigl(gh,\operatorname{Ad}(h^{-1})D\bigr),
$
and consequently
\[
\check{(R_hf)}
\bigl(g,\operatorname{Ad}(g^{-1})\beta(X_{a})\bigr)
=
\check f\bigl(gh,\operatorname{Ad}((gh)^{-1})\beta(X_{a})\bigr).
\]
The Haar measure on the unimodular group
$\widetilde{G}_{\RR}/C$ is right invariant, so that $\|R_hf\|_2=\|f\|_2$.
Thus right translations extend to a strongly continuous unitary
representation of $\widetilde{G}_{\RR}$ on the completion.
We do not claim that the full Lie supergroup acts unitarily
on this Hilbert superspace.
\item[(c)] One has an isometric identification
\begin{equation*}
L^2_+(\widetilde{\calG}_{\RR}/\calZ;\chi_C)
\cong
\Gamma^2_{\widetilde{G}_{\RR}}(E_{\chi_C})
\otimes\bigwedge\gg_{\bar1}^{\ast},
\end{equation*}
where $\Gamma^2_{\widetilde{G}_{\RR}}(E_{\chi_C})$ is the Hilbert space of square-integrable sections of the associated line bundle over $\widetilde{G}_{\RR}/C$. Indeed, smooth compactly supported sections are dense in this space. Consequently, a smooth supermatrix coefficient $\varphi_{v,w}$ belongs to the positive $L^2$-space precisely when
\begin{equation*}
\sum_{a\in A}
\int_{\widetilde{G}_{\RR}/C}
\left|
\left\langle
\pi(g)\tau\bigl(\operatorname{Ad}(g^{-1})\beta(X_{a})\bigr)v,w
\right\rangle
\right|^2\,\mathrm{d}(gC)
<\infty.
\end{equation*}

\end{enumerate}
\end{remark}

We can now define square-integrability modulo $\calZ$ for the supermatrix coefficients of $\Sigma$.

\begin{definition}
The supermatrix coefficients of $\Sigma$ are called
\emph{square-integrable modulo the center} if
\begin{equation*}
\varphi_{v,w}\in
L^2_+(\widetilde{\calG}_{\RR}/\calZ;\chi_C)
\end{equation*}
for all homogeneous $v,w\in\calH^\infty$ of the same parity.
\end{definition}

\subsection{Relative holomorphic discrete series supermodules and square-integrable supermatrix coefficients}

We retain the notation of the preceding subsection. In particular, $Z$ is the center of $\widetilde{G}_{\RR}$, whereas $\calZ=(C,0)$ is the center of $\widetilde{\calG}_{\RR}$. The latter acts on $\calH$ through the unitary character $\chi_C$.

Let $\{X_{a}\}_{a\in A}$ be the fixed homogeneous basis of $\bigwedge\odd$, and let $\beta$ denote PBW antisymmetrization. The supermatrix coefficient $\varphi_{v,w}$ is square-integrable modulo $\calZ$ precisely when each
\begin{equation}
g\longmapsto
\varphi_{v,w}\bigl(g,\operatorname{Ad}(g^{-1})\beta(X_{a})\bigr)
=
\left\langle
\pi(g)\tau\bigl(\operatorname{Ad}(g^{-1})\beta(X_{a})\bigr)v,w
\right\rangle
\end{equation}
is square-integrable modulo $C$. We first compare this condition with ordinary square-integrability modulo $Z$. The following lemma is central.

\begin{lemma}\label{lemm::Z_mod_C_compact}
    $Z/C$ is compact.
\end{lemma}

\begin{proof}
Since $\even$ is reductive, its real form is the direct sum of
its center and a semisimple Lie algebra. By
\cite[Corollary~1.2.8]{musson2012lie}, the center has dimension at most one.
Thus
$
\widetilde{G}_{\RR}\cong G_{\mathrm{ss}}$ or $
\widetilde{G}_{\RR}\cong G_{\mathrm{ss}}\times(\RR,+)$,
where $G_{\mathrm{ss}}$ is simply connected and semisimple.

By~\cite[Proposition~7.9]{Knapp}, the matrix group
$\operatorname{Ad}(G_{\mathrm{ss}})\subseteq\operatorname{GL}(\gg)$
has finite center. Hence the image of $Z(G_{\mathrm{ss}})$ is finite,
so
$
Z(G_{\mathrm{ss}})/(C\cap G_{\mathrm{ss}})
$
is finite. This proves the assertion when the second factor is absent.

Otherwise, let $T$ generate the one-dimensional real center.
By~\cite[Corollary~1.2.8]{musson2012lie}, $\operatorname{ad}(T)$ vanishes on $\even$ and
acts on the two odd summands $\gg_{1}$ and $\gg_{-1}$ by scalars $c$ and $-c$, where
$c\neq0$ because $\gg$ has trivial center. Since $T$ is central in $\gg_{\bar0,\RR}$, it
lies in $\hh_{\RR}$, and the first step of the proof of Lemma~\ref{lemm::consistent_omega}
shows that $\alpha(\hh_{\RR})\subseteq i\RR$ for every odd root $\alpha$. Hence
$c\in i\RR\setminus\{0\}$, and $\operatorname{ad}(T)$ acts on the two odd summands with
eigenvalues $ia$ and $-ia$, where $a\coloneqq|c|>0$. Therefore
\[
\operatorname{Ad}\!\left(\exp\!\left(\tfrac{2\pi}{a}T\right)\right)
=\id_{\gg}.
\]
Identifying the second factor with $t\mapsto\exp(tT)$, we obtain
a continuous surjection
\[
\frac{Z(G_{\mathrm{ss}})}{C\cap G_{\mathrm{ss}}}
\times
\frac{\RR}{(2\pi/a)\ZZ}
\longrightarrow Z/C.
\]
The source is a finite group times a circle, hence compact.
Thus $Z/C$ is compact.
\end{proof}

\begin{remark}
Lemma~\ref{lemm::Z_mod_C_compact} remains true if we work with $\gg=\sl(n\vert n)$ and real forms $\su(p,q\vert 0,n)$ with $p+q=n$. In this case, the center of $\widetilde G_\RR$ has Lie algebra $i\RR I_{2n}$, and the corresponding one-parameter subgroup lies in $C$. Thus $Z/C$ is a quotient of the finite group $Z(G_{\mathrm{ss}})/(C\cap G_{\mathrm{ss}})$ and is therefore compact. Here, we use the same notation as in the proof of Lemma~\ref{lemm::Z_mod_C_compact}.
\end{remark}

We equip the compact group $Z/C$ with a normalized Haar measure and choose compatible invariant measures on $\widetilde{G}_{\RR}/C$ and $\widetilde{G}_{\RR}/Z$ using Weil's integration formula; see~\cite[Theorem~2.52]{Folland}.

\begin{lemma}\label{lemm::integral_comparison}
Let $\pi$ be an irreducible unitary representation of $\widetilde{G}_{\RR}$. Then
\begin{equation*}
\int_{\widetilde{G}_{\RR}/C}
|\langle\pi(g)v,w\rangle|^2\,\mathrm{d}(gC)
=
\int_{\widetilde{G}_{\RR}/Z}
|\langle\pi(g)v,w\rangle|^2\,\mathrm{d}(gZ)
\end{equation*}
for all vectors $v,w$, with equality in $[0,\infty]$.
\end{lemma}

\begin{proof}
By Schur's lemma, $Z$ acts through a unitary character. Hence
$
|\langle\pi(gz)v,w\rangle|
=
|\langle\pi(g)v,w\rangle|
$
for every $z\in Z$. The result follows from Weil's integration formula and the normalization of Haar measure on $Z/C$, including for infinite integrals.
\end{proof}

As in the reductive Lie group case, we seek to realize $\calH$ as a closed subspace of $L^2_+(\widetilde{\calG}_{\RR}/\calZ;\chi_C)$. Fix a nonzero homogeneous vector $u\in\calH^{[\kk]}$.
For $v\in\calH^\infty$, define $\Psi_u(v)$ by
\begin{equation}
\label{eq::coefficient_transform}
\check{(\Psi_u(v))}(g,D)
=
\langle\pi(g)\tau(D)v,u\rangle.
\end{equation}

Recall the right regular actions:
\begin{equation}
\begin{aligned}
\check{(R_hf)}(g,D)
&=
\check f\bigl(gh,\operatorname{Ad}(h^{-1})D\bigr),
\\
\check{(R_xf)}(g,D)
&=
\check f(g,Dx),
\end{aligned}
\end{equation}
for $h\in\widetilde{G}_{\RR}$ and $x\in\gg$.
These actions preserve
$C^\infty_{\chi_C}(\widetilde{\calG}_{\RR})$.

\begin{lemma}
The map
$
\Psi_u:\calH^\infty
\longrightarrow C^\infty_{\chi_C}(\widetilde{\calG}_{\RR})
$
is complex-linear, has parity $p(u)$, and satisfies
\[
\Psi_u(\pi(h)v)=R_h\Psi_u(v),
\qquad
\Psi_u(\tau(x)v)=R_x\Psi_u(v).
\]
\end{lemma}

\begin{proof}
As the inner product is linear in the first variable, $\Psi_u$ is complex-linear, and $\pi(c)=\chi_C(c)\id_{\calH}$ gives
\[
\check{(\Psi_u(v))}(gc,D)=\chi_C(c)\check{(\Psi_u(v))}(g,D),
\qquad c\in C.
\]
For $h\in\widetilde{G}_{\RR}$ and $x\in\gg$, the identities $\tau(D)\pi(h)=\pi(h)\tau(\operatorname{Ad}(h^{-1})D)$, which follow from the covariance of $\tau$, and $\tau(D)\tau(x)=\tau(Dx)$ yield
\[
\begin{aligned}
\check{(\Psi_u(\pi(h)v))}(g,D)
&=\langle\pi(gh)\tau(\operatorname{Ad}(h^{-1})D)v,u\rangle
=\check{(R_h\Psi_u(v))}(g,D),\\
\check{(\Psi_u(\tau(x)v))}(g,D)
&=\langle\pi(g)\tau(Dx)v,u\rangle
=\check{(R_x\Psi_u(v))}(g,D).
\end{aligned}
\]
If $v$ and $D$ are homogeneous, then $\pi(g)\tau(D)v$ has parity $p(D)+p(v)$. Since $\calH_{\bar0}\perp\calH_{\bar1}$, the coefficient vanishes unless $p(D)=p(u)+p(v)$, so $\Psi_u$ has parity $p(u)$.
\end{proof}

We can now characterize relative holomorphic discrete series supermodules by their supermatrix coefficients. This is the supersymmetric version of Theorem~\ref{thm::embedding}.

\begin{theorem}\label{thm::super_embedding}
Let $\Sigma=(\pi,\tau)$ be an irreducible unitary highest weight representation of $\widetilde{\calG}_{\RR}$ with central character $\chi_C$. The following are equivalent:
\begin{enumerate}
\item[(a)] $\calH^{[\kk]}$ is a relative holomorphic discrete series supermodule;
\item[(b)] $\varphi_{v,w}$ is square-integrable modulo $\calZ$ for all homogeneous $v,w\in\calH^\infty$ of the same parity;
\item[(c)] for some nonzero homogeneous $u\in\calH^{[\kk]}$, $\Psi_u$ extends to an injective bounded operator $\Psi_u\colon\calH\to L^2_+(\widetilde{\calG}_{\RR}/\calZ;\chi_C)$ with closed image, intertwining $\pi$ with right translations.
\end{enumerate}
\end{theorem}

\begin{proof}

Let $\calH=\bigoplus_{j=1}^{N}\calH_j$ be the orthogonal decomposition into irreducible even constituents, each of highest weight type by Lemma~\ref{lemm::restriction_even} and Corollary~\ref{cor::finite_decomposition_group_level}. We write $\pi_j$ for the representation on $\calH_j$ and $P_j$ for the associated projection onto $\calH_j$. 
By the definition of $\varphi_{v,w}$ and $\omega$-contravariance of $\langle\cdot,\cdot\rangle$, we obtain, for $v,w\in\calH^\infty$,
\begin{equation}\label{eq::right_coordinates_identity}
\begin{aligned}
\varphi_{v,w}\bigl(g,\operatorname{Ad}(g^{-1})\beta(X_a)\bigr)
=
\langle\pi(g)v,\tau(\omega(\beta(X_a)))w\rangle=
\sum_{j=1}^{N}\langle\pi(g)P_jv,P_j\tau(\omega(\beta(X_a)))w\rangle.
\end{aligned}
\end{equation}

Assume (a) holds. Then every summand is square-integrable modulo $Z$ by Theorem~\ref{thm::embedding}, hence modulo $C$ by Lemma~\ref{lemm::integral_comparison}. (b) follows from Remark~\ref{rmk::right_translations_are_unitary}(c). 

Assume (b) holds. Without loss of generality, we can take $X_{a_0}=1$ by Lemma~\ref{lemm::norm}. For homogeneous $0\neq v\in\calH_j\cap\calH^\infty$, the $a_0$-coefficient $g\mapsto\langle\pi(g)v,v\rangle$ of $\varphi_{v,v}$ is square-integrable modulo $C$, hence modulo $Z$, and $\pi_j$ lies in the relative holomorphic discrete series. This implies (a).

Next, we show that (a) and (b) imply (c). Let $u\in\calH^{[\kk]}$ be nonzero and homogeneous, and set $w_a\coloneqq\tau(\omega(\beta(X_{a})))u$. 
By~\eqref{eq::coefficient_transform}, \eqref{eq::right_coordinates_identity}, the orthogonality relations, and Lemma~\ref{lemm::integral_comparison}, we obtain, for every $v\in\calH^\infty$,
\begin{equation}
\label{eq::right_coefficient_norm}
\begin{aligned}
\|\Psi_u(v)\|_2^2
&=\sum_{a\in A}\int_{\widetilde{G}_{\RR}/C}
\left|\sum_{j=1}^N\langle\pi(g)P_jv,P_jw_a\rangle\right|^2\,\mathrm{d}(gC)\\
&\leq N\sum_{a\in A}\sum_{j=1}^N
\int_{\widetilde{G}_{\RR}/C}
|\langle\pi(g)P_jv,P_jw_a\rangle|^2\,\mathrm{d}(gC)\\
&=N\sum_{a\in A}\sum_{j=1}^N
d(\pi_j)^{-1}\|P_jv\|^2\|P_jw_a\|^2\\
&\leq b_u\|v\|_{\calH}^2,
\end{aligned}
\end{equation}
where
$
b_u\coloneqq N\sum_{a\in A}\sum_{j=1}^N d(\pi_j)^{-1}\|P_jw_a\|^2.
$
Hence $\Psi_u$ extends to a bounded operator on $\calH$, and both~\eqref{eq::right_coefficient_norm} and the intertwining relation $\Psi_u\circ\pi(h)=R_h\circ\Psi_u$ extend by continuity.

We next prove injectivity. Suppose $\Psi_u(v)=0$. By continuity, the coefficient functions $g\mapsto\langle\pi(g)v,w_a\rangle$ vanish identically. Using
$
\langle\pi(g)v,w_a\rangle=\langle v,\pi(g^{-1})w_a\rangle,
$
differentiation gives
$
\langle v,\tau(\UE(\even))w_a\rangle=0
$
for every $a\in A$.
The vectors $w_a$ generate $\calH^{[\kk]}$ under the action of $\UE(\even)$. Using $\UE(\gg)=\UE(\even)\otimes \bigwedge(\odd)$ and simplicity of the Harish-Chandra supermodule $\calH^{[\kk]}$, we conclude
\[
\sum_{a\in A}\tau(\UE(\even))w_a
=
\tau(\UE(\gg))u
=
\calH^{[\kk]}.
\]
However, $\calH^{[\kk]}$ is dense in $\calH$, which forces $v=0$.

By part (b) of Remark~\ref{rmk::right_translations_are_unitary}, right translations are unitary.
Consequently, $\Psi_u^\ast\Psi_u$ belongs to the commutant of
$\pi$. This commutant is finite-dimensional because $\pi$ is
a finite direct sum of irreducible representations.
Its positive injective element $\Psi_u^\ast\Psi_u$ is therefore
bounded below by a positive scalar. Hence, for some $a_u>0$,
\begin{equation}
\label{eq::coefficient_transform_bounds}
a_u\|v\|_{\calH}^2
\leq
\|\Psi_u(v)\|_2^2
\leq
b_u\|v\|_{\calH}^2,
\qquad v\in\calH.
\end{equation}
This shows that the image of $\Psi_{u}$ is closed and (c) follows. Here, the argument applies to every nonzero homogeneous
$u\in\calH^{[\kk]}$.

Finally, assume (c). Choose a nonzero smooth vector $v\in\calH_j^{\infty}$. By injectivity of $\Psi_u$ and~\eqref{eq::right_coefficient_norm}, at least one matrix coefficient
\[
g\longmapsto\langle\pi(g)v,P_jw_a\rangle
\]
is nonzero and square-integrable modulo $C$. By Lemma~\ref{lemm::integral_comparison}, it is also square-integrable modulo $Z$. Hence $\pi_j$ belongs to the relative holomorphic discrete series. This proves (a).
\end{proof}

\section{The Duflo--Serganova cohomology and unitarity}
\label{sec::DS_functor_and_unitarity}
\noindent
In this section, we first recall Duflo--Serganova (DS) cohomology~\cite{duflo2005associated,GorelikSherman2024}. The construction rests on the following observation. If $x\in\odd$ satisfies $[x,x]=0$, then $x^2=0$, so $x$ defines a differential on every $\gg$-supermodule $M$. Its cohomology is a supermodule over the Lie superalgebra $\gg_x\coloneqq\ker(\ad x)/\im(\ad x)$. On finite-dimensional supermodules, taking DS cohomology preserves superdimension. Moreover, for basic classical $\gg$ and $x\neq0$, the defect decreases. This simplifies the representation theory. 

However, the DS cohomology is not appropriate to treat unitarizable supermodules since no nonzero real odd element satisfies $[x,x]=0$. We therefore use generalized DS cohomology~\cite{GorelikSherman2024,Sherman2025}, which only requires $x^2$ to act semisimply on $M$. We determine the infinitesimal characters and atypicality of generalized DS cohomology (Theorem~\ref{thm::generalized_DS_character}) and the real forms of $\gg_x$ (Proposition~\ref{prop::unitary_real_forms_and_DS}). Our main result, Theorem~\ref{thm::DS_unitary}, shows that the generalized DS cohomology of a simple $\omega$-unitarizable highest weight supermodule is either zero or simple and $\omega_x$-unitarizable.

\subsection{Duflo--Serganova functor}
We begin by recalling the construction of the DS functor and the properties we will need. Following~\cite{duflo2005associated,gorelik2022duflo}, we define:

\begin{definition}
The \emph{self-commuting variety} of $\gg$ is
\[
\YY\coloneqq \{x \in \odd : [x,x]=0\}.
\]
\end{definition}

Each $x\in\YY$ determines a DS functor. Let $G$ be the connected algebraic group associated with $\even$.

Let $S$ be the set of subsets of mutually orthogonal, linearly independent isotropic roots in $\Delta_{\bar{1}}$. The Weyl group $W$ acts naturally on $S$. For
$
A=\{\alpha_{1},\dots,\alpha_{k}\}\in S,
$
choose nonzero $x_{i}\in\gg^{\alpha_{i}}$ and set
$
x=x_{1}+\dots+x_{k}\in\YY.
$
The $G$-orbit of $x$ is independent of the choice of the nonzero root vectors $x_{i}$. Hence we obtain a map
\begin{equation}
\Phi':S\to\YY/G, \qquad \Phi'(A)\coloneqq Gx.
\end{equation}
If $B=w(A)$ for some $w\in W$, then $\Phi'(A)=\Phi'(B)$. Thus, $\Phi'$ induces a map
\begin{equation}
\Phi:S/W\longrightarrow\YY/G.
\end{equation}
By~\cite[Theorem~4.2]{duflo2005associated}, the map $\Phi$ is a bijection for all basic classical Lie superalgebras except $\mathfrak{psl}(n\vert n)$, for which the corresponding orbit classification is formulated using $\mathfrak{gl}(n\vert n)$. Hence the $G$-orbits in $\YY$ are in one-to-one correspondence with the $W$-orbits in $S$. In particular, $\YY$ has finitely many $G$-orbits.

Each orbit has a well-defined \emph{rank}~\cite[Section~4]{duflo2005associated}. Fix $x\in\YY$. There exist $g\in G$ and mutually orthogonal, linearly independent isotropic roots $\alpha_{1},\dotsc,\alpha_{k}$ such that
\begin{equation}
\Ad(g)x=x_{1}+\dotsb+x_{k},
\end{equation}
where $0\neq x_{i}\in\gg^{\alpha_{i}}$. The integer
$
\rk(x)\coloneqq k
$
is independent of these choices and is called the \emph{rank} of $x$. For $\gg=\glmn$ or $\gg=\slmn$, this agrees with the rank of $x$ in the standard representation.

We now associate a subset of $\YY$ with a $\gg$-supermodule $M$. Every $x\in\YY$ defines an odd endomorphism $x_{M}\in\End_{\CC}(M)$ satisfying
$
x_{M}^{2}=0.
$
Set
\begin{equation}
M_{x}\coloneqq \Ker(x_{M})/\Im(x_{M}),
\end{equation}
and define the \emph{associated variety} of $M$ by
\begin{equation}
\YY_{M}\coloneqq \{x\in\YY:M_{x}\neq\{0\}\}\subset\YY.
\end{equation}
For finite-dimensional $M$, the set $\YY_{M}$ is a $G$-invariant Zariski-closed subvariety of $\YY$~\cite{duflo2005associated}. For general highest weight modules, $\YY_{M}$ need not be Zariski closed.

\begin{proposition}
\label{VarietyTypical}
Let $M$ be a highest weight $\gg$-supermodule with typical highest weight $\Lambda$. Then 
\[
\YY_{M}=\{0\}.
\]
\end{proposition}
\begin{proof}
    Since $M$ is a highest weight supermodule with typical highest weight $\Lambda$, it admits the infinitesimal character $\chi_\Lambda$, which has atypicality zero. By~\cite[Theorem~5.3]{duflo2005associated}, $\YY_M\subseteq\{0\}$. Since $M\neq0$, we conclude that $\YY_M=\{0\}$.
\end{proof}

\begin{remark}
For $\gg=\sl(m\vert n)$, the condition $\YY_M=\{0\}$ implies that $\Lambda$ is typical whenever $M$ is a Verma supermodule or a simple highest weight module, by~\cite[Proposition~5.14]{coulembier2017homological}.
\end{remark}

\subsubsection{The Duflo--Serganova functor and its basic properties}
We recall the Duflo--Serganova functor, or DS functor for short, following~\cite{duflo2005associated,serganova2011superdimension}.

First, we associate a Lie superalgebra $\gg_{x}$ with every $x\in\YY$. Fix $x\in\YY$ and set $k\coloneqq \rk(x)$. Since
$
\ad_{x}^{2}=0,
$
we have
$
\Im(\ad_{x})\subset\Ker(\ad_{x}).
$
Moreover, $\Im(\ad_{x})=[x,\gg]$ is an ideal in $\Ker(\ad_{x})$. Hence
\begin{equation}
\gg_{x}\coloneqq \Ker(\ad_{x})/\Im(\ad_{x})
\end{equation}
is a Lie superalgebra. The defect of $\gg_{x}$ is related to the rank of $x$ by~\cite[Section~4]{gorelik2022duflo}
\begin{equation}
\rk(x)=\defect(\gg)-\defect(\gg_{x}).
\end{equation}
For example,
$
\glmn_{x}\cong\mathfrak{gl}(m-k\vert n-k),
$
and, if $m\neq n$,
$
\slmn_{x}\cong\sl(m-k\vert n-k).
$

Suppose that
$
x=x_{1}+\dotsb+x_{k}
$
with $0\neq x_{i}\in\gg^{\alpha_{i}}$, where the roots $\alpha_{1},\dotsc,\alpha_{k}$ are mutually orthogonal and isotropic. Then $\gg_{x}$ has Cartan subalgebra
\begin{equation}
\hh_{x}\coloneqq 
(\Ker(\alpha_{1})\cap\dotsb\cap\Ker(\alpha_{k}))
/
(\hh^{\alpha_{1}}\oplus\dotsb\oplus\hh^{\alpha_{k}}),
\end{equation}
where
$
\hh^{\alpha}=[\gg^{\alpha},\gg^{-\alpha}].
$
The roots of $(\gg_{x},\hh_{x})$ are induced by
\begin{equation}
\Delta_{x}
=
\{\alpha\in\Delta:(\alpha,\alpha_{i})=0
\text{ for }i=1,\dotsc,k\}
\setminus
\{\pm\alpha_{1},\dotsc,\pm\alpha_{k}\}.
\end{equation}
Here each $\alpha\in\Delta_x$ is understood through the functional it induces on $\hh_x$. Indeed, $\alpha$ restricts to $\Ker(\alpha_1)\cap\dotsb\cap\Ker(\alpha_k)$ and vanishes on each $\hh^{\alpha_i}$, since $(\alpha,\alpha_i)=0$. The same convention is made explicit in~\eqref{eq::generalized_gx_roots} for the generalized construction.

Next, we define the DS functor. Let $M$ be a $\gg$-supermodule and recall
$
M_{x}\coloneqq \Ker(x_{M})/\Im(x_{M}).
$
The space $\Ker(x_{M})$ is invariant under $\Ker(\ad_{x})$, and $[x,\gg]$ acts trivially on the quotient $M_{x}$. Thus, $M_{x}$ is naturally a $\gg_{x}$-supermodule. The assignment
\begin{equation}
\DS_{x}(M)\coloneqq M_{x}
\end{equation}
defines a functor from the category of $\gg$-supermodules to the category of $\gg_{x}$-supermodules. This functor is called the \emph{Duflo--Serganova functor}, or the \emph{DS functor}.

We recall its basic properties.

\begin{lemma}[{\cite[Section~2]{gorelik2022duflo}}]\label{DSTensorFunctor}
\begin{enumerate}
\item[(a)] For every $x\in\YY$, the Duflo--Serganova functor
$
\DS_{x}:\gsmod\to\gg_{x}\textbf{-smod}
$
is additive and symmetric monoidal.
\item[(b)] If $M$ is a finite-dimensional $\gg$-supermodule, then $\sdim(M)=\sdim(\DS_{x}(M))$.
\item[(c)] Let $x\in\YY$, and let
$
0\to M'\xrightarrow[]{\psi}M\xrightarrow[]{\phi}M''\to0
$
be an exact sequence of $\gg$-supermodules. Then there is an exact sequence of $\gg_{x}$-supermodules
\[
0\to E\to\DS_{x}(M')\to\DS_{x}(M)\to\DS_{x}(M'')\to\Pi(E)\to0
\]
for some $\gg_{x}$-supermodule $E$. In particular, the DS functor is middle exact and satisfies
$
\DS_{x}(\Pi(M))=\Pi(\DS_{x}(M))
$
for every $\gg$-supermodule $M$.
\end{enumerate}
\end{lemma}

We next restrict to highest weight $\gg$-supermodules. The associated variety gives an immediate vanishing criterion for modules with typical highest weight.

\begin{proposition}
Let $M$ be either a Verma supermodule or a simple highest weight $\gg$-supermodule with typical highest weight. Then
\[
\DS_{x}(M)=\{0\}
\]
for every $x\in\YY\setminus\{0\}$.
\end{proposition}

\begin{proof}
By Proposition~\ref{VarietyTypical}, we have $\YY_{M}=\{0\}$. Since
$
\DS_{x}(M)=M_{x}
$
and
$
\YY_{M}=\{x\in\YY:M_{x}\neq\{0\}\},
$
the assertion follows.
\end{proof}

We now relate the rank of $x$ to the degree of atypicality. We use the construction of~\cite[Section~6]{duflo2005associated}; see also~\cite{serganova2011superdimension}.

Let $\UE(\gg)^{\ad_{x}}$ denote the subsuperalgebra of $\ad_{x}$-invariants in $\UE(\gg)$, and let $I_{x}$ be the left ideal in $\UE(\gg)$ generated by $[x,\gg]$. After choosing a realization of $\gg_{x}$ inside $\Ker(\ad_{x})$, consider
\begin{equation}
\UE(\gg_{x})
\xrightarrow{\iota}
\UE(\gg)^{\ad_{x}}
\xrightarrow{\pi}
\UE(\gg)^{\ad_{x}}/
\bigl(I_{x}\cap\UE(\gg)^{\ad_{x}}\bigr).
\end{equation}
Set $\phi\coloneqq \pi\circ\iota$. By~\cite[Lemma~6.6]{duflo2005associated}, this map
is an isomorphism of super vector spaces. We therefore obtain a surjective map
\begin{equation}
\eta\coloneqq \phi^{-1}\circ\pi:
\UE(\gg)^{\ad_{x}}\longrightarrow\UE(\gg_{x}).
\end{equation}
For every $u\in\UE(\gg)^{\ad_{x}}$, the induced actions on $M_{x}$ satisfy
\begin{equation}
u m=\eta(u)m,
\qquad
m\in M_{x}.
\end{equation}

Let $\mathfrak{Z}(\gg_{x})$ denote the center of $\UE(\gg_{x})$. Since
$
\mathfrak{Z}(\gg)\subset\UE(\gg)^{\ad_{x}}
$
and $\eta$ is a homomorphism of $\gg_{x}$-supermodules for the adjoint action, we have
$
\eta(\mathfrak{Z}(\gg))\subset\mathfrak{Z}(\gg_{x}).
$
Thus, $\eta$ induces the dual map
\begin{equation}
\eta^{\ast}:
\Hom(\mathfrak{Z}(\gg_{x}),\CC)
\longrightarrow
\Hom(\mathfrak{Z}(\gg),\CC), 
\end{equation}
where $\Hom$ denotes the set of unital $\CC$-algebra homomorphisms here. If $\gg\neq\mathfrak{osp}(2l\vert 2n),F(4),D(2,1;\upalpha)$, then
$\eta^{\ast}$ is injective. In the remaining cases, each fiber of
$\eta^{\ast}$ has at most two elements
\cite[Theorem~6.12]{gorelik2022duflo}.

Let $M$ be a highest weight $\gg$-supermodule with infinitesimal character $\chi$. If $m\in\Ker(x_{M})$ and $[m]\in M_{x}$ denotes its class, then for every $z\in\mathfrak{Z}(\gg)$,
\begin{equation}
\chi(z)[m]=z[m]=\eta(z)[m].
\end{equation}
Hence, if $M_{x}$ contains a submodule with infinitesimal character $\xi$, then
$
\eta^{\ast}(\xi)=\chi.
$
Moreover, if $\xi$ is an infinitesimal character of $\gg_{x}$, then~\cite[Corollary~6.9]{duflo2005associated}
\begin{equation}
\at(\eta^{\ast}(\xi))
=
\at(\xi)+\rk(x).
\end{equation}
This gives the following theorem.

\begin{theorem}[{\cite[Theorem~2.1]{serganova2011superdimension}}]\label{DScharacter}
Let $M$ be a simple highest weight $\gg$-supermodule with infinitesimal character $\chi$. Then $M_{x}$ is a direct sum of $\gg_{x}$-supermodules with generalized infinitesimal characters in $(\eta^{\ast})^{-1}(\chi)$. Every infinitesimal character in $(\eta^{\ast})^{-1}(\chi)$ has degree of atypicality
\[
\at(\chi)-\rk(x).
\]
\end{theorem}

The following consequence will be used repeatedly.

\begin{corollary}\label{DSNM}
Let $x\in\YY\setminus\{0\}$ have rank $l$, and let $M$ be a simple highest weight $\gg$-supermodule whose infinitesimal character has degree of atypicality $k$. If $k<l$, then
\[
\DS_{x}(M)=\{0\}.
\]
\end{corollary}

\begin{proof}
If $\DS_{x}(M)\neq\{0\}$, Theorem~\ref{DScharacter} gives an infinitesimal character of $\gg_{x}$ of degree of atypicality $k-l$. Since degrees of atypicality are nonnegative, this is impossible when $k<l$.
\end{proof}

For a unitarizable highest weight $\gg$-supermodule $\HH$ and an element $x\in\YY$, it remains to determine when $\DS_{x}(\HH)$ is unitarizable. The ordinary DS functor is not well adapted to this question because the conjugate-linear anti-involution $\omega$ on $\gg$ does not, in general, induce a conjugate-linear anti-involution on $\gg_{x}$.

There is a more basic obstruction. Suppose that $x\in\YY$ satisfies
$
\omega(x)=\lambda x
$
for some $\lambda\in\CC$ with $|\lambda|=1$. Let $x_{\HH}$ denote the action of $x$ on $\HH$, where $\HH$ is some nontrivial $\omega$-unitarizable highest weight supermodule. Unitarity gives
$
x_{\HH}^{\ast}=\lambda x_{\HH},
$
while $x\in\YY$ implies $x_{\HH}^{2}=0$. Hence, for every $v\in\HH$,
\begin{equation}
\langle x_{\HH}v,x_{\HH}v\rangle
=
\langle v,x_{\HH}^{\ast}x_{\HH}v\rangle 
=
\overline{\lambda}\langle v,x_{\HH}^{2}v\rangle
=
0.
\end{equation}
Thus $x_{\HH}\equiv 0$, since $\langle\cdot,\cdot\rangle$ is positive definite. Consequently,
$
\DS_{x}(\HH)=\HH_{x}=\HH.
$
Since $\gg$ is simple, every nontrivial $\gg$-supermodule is faithful. Hence, no nonzero $x\in\YY$ can satisfy $\omega(x)=ix$, and the self-commuting variety has no nonzero real points.

\subsection{Generalized Duflo--Serganova functors}
\label{subsec::generalized_DS_functor} 
We briefly introduce the generalized DS functor following~\cite{GorelikSherman2024,Sherman2025}. The idea underlying the generalized DS functor first appeared in the work of Weissauer~\cite{Weissauer2010}. This approach was subsequently developed by Heidersdorf and Weissauer~\cite[Sections~5--6]{heidersdorf2021cohomological}, where it appears under the name Dirac functor. These works are the foundations for the construction considered here.

The generalized DS functor assigns a cohomology, analogous to that of the ordinary DS functor, to odd elements $x\in\odd$ whose squares are semisimple elements of $\even$. This generalization allows us to treat unitarizable supermodules.

The \emph{homological variety} of $\gg$ is
\begin{equation}
\mathcal{Y}^{\mathrm{hom}}
\coloneqq 
\left\{
x\in\gg_{\bar1}
: 
\ad_{x^{2}}\text{ acts semisimply on }\gg
\right\},
\end{equation}
where $x^{2}\coloneqq \frac{1}{2}[x,x]$. 
The variety $\mathcal{Y}^{\mathrm{hom}}$ is stable under the adjoint action of $G$ and under scalar multiplication, since $\ad_{(tx)^2}=t^2\ad_{x^2}$ and $\ad_{(g\cdot x)^2}=\operatorname{Ad}(g)\,\ad_{x^2}\operatorname{Ad}(g)^{-1}$. Hence the generalized Duflo--Serganova cohomology depends only on the $G\times\CC^\times$-orbit of $x$, and it suffices to consider one representative of each orbit.

Fix $x\in\mathcal{Y}^{\mathrm{hom}}$ and set $c\coloneqq x^{2}$. Let $M$ be a $\gg$-supermodule on which $c$ acts semisimply, and set $M^{c}\coloneqq \ker(c:M\to M)$. Since $[x,c]=0$, the space $M^{c}$ is $x$-stable, and $x^{2}=0$ on $M^{c}$. Hence one defines
\begin{equation}\label{eq::generalized_DS}
\DS_{x}(M)
\coloneqq 
H(M^{c},x)
=
\frac{\ker(x:M^{c}\to M^{c})}{\im(x:M^{c}\to M^{c})}
=
\frac{\ker x}{xM\cap\ker x}.
\end{equation}
This is the generalized Duflo--Serganova functor; see~\cite[Section~7.1]{GorelikSherman2024} and~\cite[Section~1]{Sherman2025}. If $x^{2}=0$, then $M^{c}=M$, and~\eqref{eq::generalized_DS} reduces to the original Duflo--Serganova functor~\cite{duflo2005associated}.

The corresponding Lie superalgebra is defined by the same cohomological construction. Put $\gg^{c}\coloneqq \ker\ad_{c}$. Since $(\ad_{x})^{2}=\ad_{c}$, the restriction of $\ad_{x}$ to $\gg^{c}$ is square-zero. Hence
\begin{equation}\label{eq::generalized_gx}
\gg_{x}\coloneqq H(\gg^{c},\ad_{x})
=\frac{\gg^{x}}{[x,\gg^{c}]},
\qquad
\gg^{x}\coloneqq \ker\ad_{x}.
\end{equation}

The Lie superalgebra $\gg_{x}$ has a simple root-theoretic description for general $x\in\YY^{\hom}$; see~\cite[Lemma~3.9 and Theorem~3.12]{Sherman2025}. Set $c=x^{2}$. 
By~\cite[Lemma~3.9]{Sherman2025}, we choose a Cartan subalgebra $\tt\subseteq\even$ containing $c$ such that $x$ admits the normal form~\eqref{eq::normal_form}. This Cartan subalgebra need not coincide with the fixed Cartan subalgebra $\hh$.
 We have
\begin{equation}
\gg^{c}
=
\tt\oplus
\bigoplus_{\substack{\alpha\in\Delta(\gg,\tt)\\ \alpha(c)=0}}
\gg^{\alpha},
\end{equation}
and hence $\Delta(\gg^{c},\tt)=\{\alpha\in\Delta(\gg,\tt):\alpha(c)=0\}$.

By~\cite[Lemma~3.9]{Sherman2025}, there exist mutually orthogonal, linearly independent isotropic roots $\alpha_{1},\ldots,\alpha_{k}$ such that
\begin{equation}\label{eq::normal_form}
x
=
u_{\alpha_{1}}+\cdots+u_{\alpha_{k}}
+
z_{1}v_{\alpha_{1}}+\cdots+z_{k}v_{\alpha_{k}},
\end{equation}
where $u_{\alpha_{i}}\in\gg^{\alpha_{i}}\setminus\{0\}$, $v_{\alpha_{i}}\in\gg^{-\alpha_{i}}\setminus\{0\}$, and $z_{i}\in\CC$. The integer $k$ is called the rank of $x$, and we write $\rk(x)=k$. Set $t_{\alpha_{i}}\coloneqq [u_{\alpha_{i}},v_{\alpha_{i}}]\in\tt$. Then $\CC t_{\alpha_{i}}\oplus\CC u_{\alpha_{i}}\oplus\CC v_{\alpha_{i}}$ is isomorphic to $\sl(1\vert1)$, and the corresponding subalgebras commute for different $i$. In particular,
\begin{equation}
c=x^{2}=\sum_{i=1}^{k}z_{i}t_{\alpha_{i}}.
\end{equation}

Set $\tt'\coloneqq \bigoplus_{i=1}^{k}\CC t_{\alpha_{i}}$. Then $\tt'$ is a maximal toral subalgebra of $[x,\gg^{c}]$, and the centralizer of $\tt'$ in $\gg$ is~\cite[Theorem~3.12]{Sherman2025}
\begin{equation}
\mathfrak{c}_{\gg}(\tt')
\cong
\gg_{x}\times\mathfrak{gl}(1\vert1)^{k}.
\end{equation}
Thus $\gg_{x}$ is obtained from the centralizer of $\tt'$ by removing the $k$ copies of $\mathfrak{gl}(1\vert1)$ associated with the roots $\alpha_{1},\ldots,\alpha_{k}$.

More explicitly, set
$
\tt_{A}
\coloneqq 
\ker\alpha_{1}\cap\cdots\cap\ker\alpha_{k}.
$
Since $t_{\alpha_{i}}\in\tt_{A}$, the quotient
$
\tt_{x}
\coloneqq 
\tt_{A}/\tt'
$
is a Cartan subalgebra of $\gg_{x}$.
Its root system is
\begin{equation}\label{eq::generalized_gx_roots}
\Delta_{x}
=
\left\{
\alpha|_{\tt_{x}}
:
\alpha\in\Delta(\gg,\tt)\setminus
\{\pm\alpha_{1},\ldots,\pm\alpha_{k}\},
\quad
(\alpha,\alpha_{i})=0
\text{ for all }i,
\quad
\alpha|_{\tt_{x}}\neq0
\right\}.
\end{equation}

The decomposition of $\mathfrak{c}_{\gg}(\tt')$ provides a noncanonical realization of $\gg_x$ as a Lie subsuperalgebra of $\gg$. Choosing the corresponding Cartan subalgebra gives a complement of $\tt'$ in $\tt_{A}$.

The isomorphism type of $\gg_{x}$ depends only on the isotropic root set
$A=\{\alpha_{1},\ldots,\alpha_{k}\}$ and not on the coefficients $z_{i}$, since
$\mathfrak{c}_{\gg}(\tt')\cong\gg_{x}\times\mathfrak{gl}(1\vert1)^{k}$.
Thus $\gg_{x}$ is isomorphic to the Lie superalgebra $\gg_{u}$ that the Duflo--Serganova construction associates with the
square-zero element $u=u_{\alpha_{1}}+\cdots+u_{\alpha_{k}}$. Consequently, the usual
rank-reduction formulas remain valid for every $x\in\YY^{\hom}$. More precisely, if $x$ has rank $k$, then
\begin{equation}
\begin{array}{c|c}
\gg & \gg_{x}\\
\hline
\mathfrak{gl}(m\vert n)
&
\mathfrak{gl}(m-k\vert n-k)
\\
\sl(m\vert n),\ m\neq n
&
\sl(m-k\vert n-k)
\\
\mathfrak{osp}(m\vert2n)
&
\mathfrak{osp}(m-2k\vert2n-2k).
\end{array}
\end{equation}
Here $0\leq k\leq\defect(\gg)$. These formulas include the types
$B(m\vert n)$, $C(n)$, and $D(m\vert n)$. In particular, for
$\gg=C(n)=\mathfrak{osp}(2\vert2n-2)$ and $k=1$, one obtains
$\gg_{x}\cong\mathfrak{sp}(2n-4)$.

The exceptional basic Lie superalgebras have defect one. Hence a nonzero $x\in\YY^{\hom}$ has rank one, and
\begin{equation}
D(2,1;\upalpha)_{x}\cong\CC,
\qquad
F(4)_{x}\cong\sl(3),
\qquad
G(3)_{x}\cong\sl(2).
\end{equation}

Moreover, $\DS_{x}(M)$ is naturally a $\gg_{x}$-supermodule, and $\DS_{x}$ defines a functor on the category of $\gg$-supermodules on which $x^{2}$ acts semisimply~\cite[Section~7.1]{GorelikSherman2024}. In particular, the construction applies to a highest weight supermodule whenever $x^{2}$ acts semisimply on it.

\begin{proposition}\label{prop::generalized_DS_properties}
Let $x\in\mathcal{Y}^{\mathrm{hom}}$, and let $M,N$ be $\gg$-supermodules on which $x^{2}$ acts semisimply.
\begin{enumerate}
\item[(a)] The functor $\DS_{x}$ is symmetric monoidal and commutes with parity change. Thus $\DS_{x}(M\otimes N)\cong\DS_{x}(M)\otimes\DS_{x}(N)$ and $\DS_{x}(\Pi M)\cong\Pi\DS_{x}(M)$.

\item[(b)] If $M$ is finite-dimensional, then $\sdim\DS_{x}(M)=\sdim M$.

\item[(c)] If $0\to M'\to M\to M''\to0$ is a short exact sequence of $\gg$-supermodules on which $x^{2}$ acts semisimply, then there exists a $\gg_{x}$-supermodule $E$ such that $0\to E\to\DS_{x}(M')\to\DS_{x}(M)\to\DS_{x}(M'')\to\Pi E\to0$ is exact. In particular, $\DS_{x}$ is middle exact, but it is not exact in general.
\end{enumerate}
\end{proposition}

\begin{proof}
Part~(a) follows from~\cite[Section~7.1]{GorelikSherman2024}. For part~(b), decompose $M$ into the eigenspaces of $c=x^{2}$. On every nonzero eigenspace, $x$ is an odd isomorphism, so its superdimension vanishes. Hence $\sdim M=\sdim M^{c}$, and the usual cohomological argument~\cite[Section~2]{gorelik2022duflo} gives $\sdim M^{c}=\sdim\DS_{x}(M)$. For part~(c), the short exact sequence restricts to a short exact sequence on the zero eigenspaces of $x^{2}$. The associated $2$-periodic long exact sequence in $x$-cohomology gives the stated exact sequence; see~\cite[Section~7.1.2]{GorelikSherman2024} and~\cite[Section~2]{gorelik2022duflo}.
\end{proof}

Let $M$ be a $\gg$-supermodule such that $x^{2}$ acts semisimply on $M$ for every $x\in\YY^{\hom}$. We define its \emph{semisimple support} by
\begin{equation}\label{eq::semisimple_support}
X_{M}^{\mathrm{ss}}
\coloneqq 
\left\{
x\in\YY^{\hom}:\DS_{x}(M)\neq0
\right\}.
\end{equation}
This agrees with the usual definition for finite-dimensional supermodules; see~\cite[Section~10.3]{gorelik2022duflo}. The standard properties of semisimple support remain valid in this setting.

\begin{proposition}\label{prop::semisimple_support_properties}
Let $M,N$ be $\gg$-supermodules on which $x^{2}$ acts semisimply for every $x\in\YY^{\hom}$.
\begin{enumerate}
\item[(a)] $X_{M\oplus N}^{\mathrm{ss}}=X_{M}^{\mathrm{ss}}\cup X_{N}^{\mathrm{ss}}$.

\item[(b)] $X_{M\otimes N}^{\mathrm{ss}}=X_{M}^{\mathrm{ss}}\cap X_{N}^{\mathrm{ss}}$.

\item[(c)] $X_{\Pi M}^{\mathrm{ss}}=X_{M}^{\mathrm{ss}}$.

\item[(d)] If $0\to M'\to M\to M''\to0$ is exact and $x^{2}$ acts semisimply on all three modules, then each of $X_{M'}^{\mathrm{ss}}$, $X_{M}^{\mathrm{ss}}$, and $X_{M''}^{\mathrm{ss}}$ is contained in the union of the other two.

\item[(e)] The set $X_{M}^{\mathrm{ss}}$ is stable under multiplication by $\CC^{\times}$.
\end{enumerate}
\end{proposition}

\begin{proof}
Parts~(a) and~(c) follow directly from the definition. Part~(b) follows from the tensor isomorphism $\DS_{x}(M\otimes N)\cong\DS_{x}(M)\otimes\DS_{x}(N)$. Part~(d) follows from the $2$-periodic long exact sequence in $\DS_{x}$-cohomology. For part~(e), one has $\DS_{ax}(M)\cong\DS_{x}(M)$ for every $a\in\CC^{\times}$. 
\end{proof}

If $M$ integrates to a representation of the even group $G$, then $X_{M}^{\mathrm{ss}}$ is also $G$-stable. Hence it is a union of $G\times\CC^{\times}$-orbits in $\YY^{\hom}$. For a general infinite-dimensional highest weight supermodule, this invariance requires such an integrated $G$-action and is not automatic.

Finally, we extend the relation between infinitesimal characters and atypicality to the generalized Duflo--Serganova functor. We retain the basic classical setting and the root-theoretic realization introduced above.

Fix $x\in\YY^{\hom}$, and set $c\coloneqq x^2$. Since $(\ad_x)^2=\ad_c$, the derivation $\ad_x$ preserves $\UE(\gg)^{\ad_c}$ and restricts to a square-zero derivation on this subalgebra. Moreover, every $\ad_x$-invariant element of $\UE(\gg)$ belongs to $\UE(\gg)^{\ad_c}$. Hence
\begin{equation}\label{eq::generalized_DS_enveloping_cohomology}
H(\UE(\gg)^{\ad_c},\ad_x)
=
\frac{\UE(\gg)^{\ad_x}}{[x,\UE(\gg)^{\ad_c}]}.
\end{equation}
The subspace $[x,\UE(\gg)^{\ad_c}]$ is a two-sided ideal in $\UE(\gg)^{\ad_x}$. Thus $H(\UE(\gg)^{\ad_c},\ad_x)$ is an associative superalgebra.

The inclusion $\gg\hookrightarrow\UE(\gg)$ induces a canonical homomorphism
\begin{equation}\label{eq::generalized_DS_PBW}
\jmath_x:
\UE(\gg_x)\longrightarrow H(\UE(\gg)^{\ad_c},\ad_x).
\end{equation}

\begin{lemma}\label{lemm::j_math_is_iso}
The homomorphism $\jmath_x$ is an isomorphism of associative superalgebras.
\end{lemma}

\begin{proof}
We adapt the argument of~\cite[Lemma~6.6]{duflo2005associated}. Set $d\coloneqq\ad_x$ and write $\gg=\gg^c\oplus W$, where $W$ is the sum of the nonzero eigenspaces of $\ad_c$. Since $d^2=\ad_c$, the operator $h=\frac12d\circ\ad_c^{-1}$ on $W$, extended by zero on $\gg^c$ to an odd derivation of $S(\gg)$, satisfies $dh+hd=k\operatorname{id}$ on $S(\gg^c)\otimes S^k(W)$. As $h$ commutes with $\ad_c$, all invariant summands with $k>0$ are contractible. Hence, by the Künneth formula,
\begin{equation*}
H(S(\gg)^{\ad_c},d)\cong H(S(\gg^c),d)\cong S(H(\gg^c,d))=S(\gg_x).
\end{equation*}
Using that PBW supersymmetrization is a filtered isomorphism of complexes, that is, $S(\gg)^{\ad_c}\cong\UE(\gg)^{\ad_c}$, we obtain $\gr H(\UE(\gg)^{\ad_c},d)\cong S(\gg_x)$. As $\gr\UE(\gg_x)\cong S(\gg_x)$, the map $\gr(\jmath_x)$ is the identity, and therefore $\jmath_x$ is an isomorphism.
\end{proof}

As for ordinary DS cohomology, using the quotient homomorphism
$
q_x:\UE(\gg)^{\ad_x}
\longrightarrow H(\UE(\gg)^{\ad_c},\ad_x),
$
we define
\begin{equation}\label{eq::generalized_DS_eta}
\eta_x\coloneqq\jmath_x^{-1}\circ q_x:
\UE(\gg)^{\ad_x}\longrightarrow\UE(\gg_x).
\end{equation}
This is a surjective homomorphism of associative superalgebras by Lemma~\ref{lemm::j_math_is_iso}. Moreover, its restriction to the center of $\UE(\gg)$ takes values in $\mathfrak{Z}(\gg_x)$. Indeed, every central element of $\UE(\gg)$ commutes with $\UE(\gg)^{\ad_x}$, and surjectivity of $\eta_x$ implies that its image commutes with all of $\UE(\gg_x)$.

To verify compatibility with supermodule actions, let $M$ be a $\gg$-supermodule on which $c$ acts semisimply. For homogeneous $a\in\UE(\gg)^{\ad_c}$ and $m\in M^c\cap\ker x$, we have
\begin{equation}
[x,a]m=x(am),
\qquad am\in M^c.
\end{equation}
Hence $[x,\UE(\gg)^{\ad_c}]$ acts trivially on $\DS_x(M)$. Consequently,
\begin{equation}\label{eq::generalized_DS_eta_action}
a[m]=\eta_x(a)[m],
\qquad
a\in\UE(\gg)^{\ad_x},\quad [m]\in\DS_x(M).
\end{equation}
We use the same notation for the induced homomorphism
$
\eta_x:\mathfrak{Z}(\gg)\longrightarrow\mathfrak{Z}(\gg_x)
$
and denote its pullback on infinitesimal characters by
$
\eta_x^{\ast}(\xi)\coloneqq\xi\circ\eta_x.
$

We next compare $\eta_x$ with the ordinary Duflo--Serganova homomorphism. 
For this $x$, choose a Cartan subalgebra $\tt\subseteq\even$ containing $c=x^2$ as in~\cite[Lemma~3.9]{Sherman2025}. With respect to this choice, write
\begin{equation}
x=\sum_{i=1}^{k}
\bigl(u_{\alpha_i}+z_i v_{\alpha_i}\bigr),
\qquad
u\coloneqq\sum_{i=1}^{k}u_{\alpha_i}.
\end{equation}
Then $u^2=0$ and $\rk(u)=\rk(x)=k$.

Choose a common root-theoretic realization $\mathfrak{l}\subseteq\gg$ of $\gg_x$ and $\gg_u$, and let $\tt_x\subseteq\tt$ denote its Cartan subalgebra. We choose $\mathfrak{l}$ as a Lie subsuperalgebra such that
\begin{equation}
\tt_{A}=\tt_x\oplus\tt',
\qquad
\tt'=\bigoplus_{i=1}^{k}\CC t_{\alpha_i}.
\end{equation}
Both $x$ and $u$ centralize $\mathfrak{l}$.

\begin{lemma}\label{lem::generalized_DS_eta_compatibility}
Under the common realization $\gg_x\cong\mathfrak{l}\cong\gg_u$, the homomorphisms on centers satisfy
\begin{equation*}
\eta_x=\eta_u:
\mathfrak{Z}(\gg)\longrightarrow\mathfrak{Z}(\mathfrak{l}).
\end{equation*}
\end{lemma}

\begin{proof}
Choose a Borel subalgebra of $\gg$ for which $\alpha_1,\ldots,\alpha_k$ are simple roots and whose intersection with $\mathfrak{l}$ is a Borel subalgebra of $\mathfrak{l}$. Such a Borel subalgebra exists because the roots $\alpha_1,\ldots,\alpha_k$ are linearly independent, pairwise orthogonal, and isotropic.

Fix $\mu\in\tt_x^{\ast}$. Since $\tt_{A}=\tt_x\oplus\tt'$, we can choose $\lambda\in\tt^{\ast}$ such that
\begin{equation*}
\lambda|_{\tt_x}=\mu,
\qquad
\lambda(t_{\alpha_i})=0,
\quad 1\leq i\leq k.
\end{equation*}
Let $L(\lambda)$ be the corresponding simple highest weight $\gg$-supermodule, with highest weight vector $v_\lambda$.

For each $i$, the vector $v_{\alpha_i}v_\lambda$ is singular because $\alpha_i$ is an isotropic simple root and $\lambda(t_{\alpha_i})=0$. Since $L(\lambda)$ is simple, this vector vanishes. We conclude that
$
L(\lambda)^{\lambda-\alpha_i}=0.
$
The highest weight condition also gives
$
u_{\alpha_i}v_\lambda=0
$
and
$
L(\lambda)^{\lambda+\alpha_i}=0.
$
Thus $xv_\lambda=0$. Moreover, $v_\lambda\notin xL(\lambda)$, since a preimage would have to belong to one of these vanishing weight spaces.

Since $c\in\tt$ acts semisimply on $L(\lambda)$, the vector $v_\lambda$ defines a nonzero class in $\DS_x(L(\lambda))$. This class is a highest weight vector for $\mathfrak{l}$ of weight $\mu$. Hence~\eqref{eq::generalized_DS_eta_action} gives
\begin{equation*}
\chi_\mu^{\mathfrak{l}}(\eta_x(z))
=
\chi_\lambda^\gg(z),
\qquad
z\in\mathfrak{Z}(\gg).
\end{equation*}
Here, $\chi_\mu^{\ll}$ and $\chi_\lambda^{\gg}$ denote the infinitesimal characters of the simple highest weight supermodules of $\ll$ and $\gg$ with highest weights $\mu$ and $\lambda$, respectively. 

The same argument applies to $u$, and therefore
$
\chi_\mu^{\mathfrak{l}}(\eta_x(z))
=
\chi_\mu^{\mathfrak{l}}(\eta_u(z))
$
for every $\mu\in\tt_x^{\ast}$. Thus the Harish-Chandra images of $\eta_x(z)$ and $\eta_u(z)$ agree. Since the Harish-Chandra homomorphism is injective, we conclude that $\eta_x(z)=\eta_u(z)$.
\end{proof}

By Lemma~\ref{lem::generalized_DS_eta_compatibility} and~\cite[Corollary~6.9]{duflo2005associated}, every infinitesimal character $\xi$ of $\gg_x$ satisfies
\begin{equation}\label{eq::generalized_DS_atypicality}
\at\bigl(\eta_x^{\ast}(\xi)\bigr)
=
\at(\xi)+\rk(x).
\end{equation}
Moreover, the description of the fibers in~\cite[Theorem~6.12]{gorelik2022duflo} also applies to $\eta_x^{\ast}$.

Let $V$ be a $\gg_x$-supermodule, and let $\xi$ be an infinitesimal character of $\gg_x$. Set
\begin{equation}
\mathfrak{m}_{\xi}\coloneqq\ker\xi,
\qquad
V_{\xi}\coloneqq
\left\{
v\in V:
\mathfrak{m}_{\xi}^{\,q}v=0
\text{ for some }q\geq1
\right\}.
\end{equation}
We call $V_{\xi}$ the \emph{generalized infinitesimal-character submodule} associated with $\xi$. The center $\mathfrak{Z}(\gg_x)$ acts \emph{locally finitely} on $V$ if
$
\dim_{\CC}\mathfrak{Z}(\gg_x)v<\infty$ for every $v\in V$.

Combining these results gives the following theorem.

\begin{theorem}\label{thm::generalized_DS_character}
Let $x\in\YY^{\hom}$ have rank $k$, and let $M$ be a $\gg$-supermodule with infinitesimal character $\chi$ on which $x^2$ acts semisimply. Then every simple subquotient of $\DS_x(M)$ has an infinitesimal character $\xi$ satisfying
\[
\eta_x^{\ast}(\xi)=\chi,
\qquad
\at(\xi)=\at(\chi)-k.
\]
In particular, $\DS_x(M)=0$ if $k>\at(\chi)$. If $\mathfrak{Z}(\gg_x)$ acts locally finitely on $\DS_x(M)$, then
\[
\DS_x(M)
=
\bigoplus_{\xi\in(\eta_x^{\ast})^{-1}(\chi)}
\DS_x(M)_{\xi}.
\]
\end{theorem}

\begin{proof}
We first note that, since $M$ has infinitesimal character $\chi$, equation~\eqref{eq::generalized_DS_eta_action} gives
\begin{equation}\label{eq::generalized_DS_scalar_action}
\eta_x(z)v=\chi(z)v,
\qquad
z\in\mathfrak{Z}(\gg),\quad v\in\DS_x(M).
\end{equation}

Let $V$ be a simple subquotient of $\DS_x(M)$. Since $V$ is a simple supermodule over the finite-dimensional Lie superalgebra $\gg_x$, it has an infinitesimal character $\xi$. Equation~\eqref{eq::generalized_DS_scalar_action} gives $\xi\circ\eta_x=\chi$, and hence $\at(\xi)=\at(\chi)-k\geq 0$ by~\eqref{eq::generalized_DS_atypicality}.

Finally, suppose that $\mathfrak{Z}(\gg_x)$ acts locally finitely on $\DS_x(M)$, that is, for every $v\in\DS_x(M)$, the space $\{zv:z\in\mathfrak{Z}(\gg_x)\}$ is finite-dimensional. Since $\mathfrak{Z}(\gg_x)$ is commutative, each of these spaces decomposes into simultaneous generalized eigenspaces. These decompositions are compatible with inclusions and therefore give
\begin{equation*}
\DS_x(M)=\bigoplus_\xi\DS_x(M)_\xi.
\end{equation*}
It remains to show that $\eta_x^\ast(\xi)=\chi$ whenever $\DS_x(M)_\xi\neq0$. On each such generalized eigenspace, $\eta_x(z)-\xi(\eta_x(z))$ acts locally nilpotently for every $z\in\mathfrak{Z}(\gg)$. By~\eqref{eq::generalized_DS_scalar_action}, this operator acts as the scalar $\chi(z)-\xi(\eta_x(z))$, which must therefore vanish. This completes the proof.
\end{proof}

\begin{corollary}
    Let $M$ be a $\gg$-supermodule. If $M$ has infinitesimal character $\chi$ and $x^2$ acts semisimply for every $x\in\YY^{\hom}$, then
\begin{equation*}\label{eq::semisimple_support_atypicality_bound}
X_M^{\mathrm{ss}}
\subseteq
\{
x\in\YY^{\hom}:
\rk(x)\leq\at(\chi)\}.
\end{equation*}
\end{corollary} 

\subsection{Generalized Duflo--Serganova functors and unitarity}

Let $\omega$ be a conjugate-linear anti-involution of $\gg$ associated with a real form listed in~\eqref{table::real_forms_classical} or~\eqref{table::real_forms_exceptional}. In contrast to the self-commuting variety, the homological variety $\YY^{\hom}$ can contain nonzero elements compatible with $\omega$. We first establish when $\omega$ descends to $\gg_x$ and then identify a natural real locus of $\YY^{\hom}$. This provides the basis for studying generalized Duflo--Serganova functors and unitarity.

\begin{lemma}\label{lemm::real_form_descend}
Let $x\in\YY^{\hom}$ satisfy $\omega(x)=\lambda x$ for some $\lambda\in\CC^{\times}$ with $\vert \lambda\vert=1$. Then $\omega$ induces a conjugate-linear anti-involution $\omega_x$ on $\gg_x$.
\end{lemma}

\begin{proof}
Since
\[
\omega(c)
=\frac12[\omega(x),\omega(x)]
=\lambda^2c,
\]
the map $\omega$ preserves both $\gg^c$ and $\gg^x$. Moreover, for every homogeneous $y\in\gg^c$,
\[
\omega([x,y])
=[\omega(y),\omega(x)]
=-\lambda(-1)^{p(y)}[x,\omega(y)]
\in[x,\gg^c].
\]
Thus $[x,\gg^c]$ is $\omega$-stable, and
\[
\omega_x\bigl(a+[x,\gg^c]\bigr)
\coloneqq\omega(a)+[x,\gg^c],
\qquad a\in\gg^x,
\]
defines a conjugate-linear anti-involution on $\gg_x$.
\end{proof}

Recall that the real form associated with $\omega$ is the real Lie superalgebra
$\gg^\omega=\gg^\omega_{\bar0}\oplus\gg^\omega_{\bar1}$, where
\begin{equation}
\gg^\omega_{\bar{z}}
\coloneqq
\{x\in\gg_{\bar{z}}:\omega(x)=-(-i)^{p(x)} x\},
\qquad \bar{z}\in \ZZ_{2}.
\end{equation}
 Since $\YY^{\hom}\subseteq\odd$, an element $x\in\YY^{\hom}$ belongs to $\gg^\omega$ precisely when $\omega(x)=ix$. We therefore define the \emph{real locus} of $\YY^{\hom}$ by
\begin{equation}
\YY^{\hom}_\omega
\coloneqq
\YY^{\hom}\cap\gg^\omega
=
\{x\in\YY^{\hom}:\omega(x)=ix\}.
\end{equation}

If $\gg$ has an isotropic odd root, then $\YY^{\hom}_\omega$ contains nonzero elements. Indeed, let $\alpha$ be such a root and choose $0\neq X_\alpha\in\gg^\alpha$. Since $\omega(\gg^\alpha)=\gg^{-\alpha}$, the element
$
x\coloneqq X_\alpha-i\omega(X_\alpha)
$
is nonzero and satisfies $\omega(x)=ix$. Moreover,
$
x^2=-i[X_\alpha,\omega(X_\alpha)]\in\hh
$
is semisimple. Hence $x\in\YY^{\hom}_\omega$.

By our standing assumptions, $\gg$ is basic classical and admits a nontrivial $\omega$-unitarizable supermodule. Consequently, $x\in\YY_\omega^{\hom}$ with $x^2=0$ must vanish. Indeed, unitarity forces $x$ to act trivially, while simplicity of $\gg$ ensures that every nontrivial supermodule is faithful. Henceforth, we consider only $x\neq0$.

We describe the real forms $\gg_x^{\omega_x}$ for $x\in\YY_\omega^{\hom}$. As a first step, we conjugate $c=x^2$ into our fixed compactly embedded Cartan subalgebra $\hh_{\RR}$. Recall that $\hh_{\RR}=\hh\cap\gg_{\bar0,\RR}$.

Let $\operatorname{Int}(\gg_{\bar0,\RR})$ denote the subgroup of $\operatorname{Aut}_{\bar0}(\gg)$ generated by $\exp(\ad_y)$ for $y\in\gg_{\bar0,\RR}$.

\begin{lemma}\label{lemm::conjugation_inner}
Let $x\in\YY_\omega^{\hom}$, and set $c=x^2$. Then there exists $\varphi_0\in\operatorname{Int}(\gg_{\bar0,\RR})$ such that $\varphi_0(c)\in \hh_{\RR}$. Moreover, every element of $\operatorname{Int}(\gg_{\bar0,\RR})$ commutes with $\omega$.
\end{lemma}

\begin{proof}
Since $\omega(x)=ix$, we have $\omega(c)=\omega(x)^2=-c$, so $c\in\gg_{\bar0,\RR}$. Consider the invariant convex cone
\[
C\coloneqq\operatorname{cone}\{[y,y]:y\in\gg_{\bar1,\RR}\}
\subseteq\gg_{\bar0,\RR}.
\]
Since $\gg$ is simple and admits a nontrivial $\omega$-unitarizable supermodule, $C$ is pointed by~\cite[Proposition~6.1.2(a)]{neeb2011lie}. As $2c=[x,x]\in C$, the invariant convex cone generated by the adjoint orbit of $c$ is also pointed.

Since $\gg_{\bar0,\RR}$ is reductive and $c$ is semisimple, $c$ is elliptic, that is, $\ad_c$ has purely imaginary spectrum. Indeed, a semisimple element of a pointed invariant convex cone in a reductive Lie algebra is elliptic~\cite[Proposition~3.10]{Neeb_Oeh},~\cite[Theorem~1.2]{NeebOlshanskii}. An elliptic element generates a compactly embedded subalgebra. It is therefore contained in a maximal compactly embedded subalgebra, which is conjugate to $\kk_{\RR}$, and hence it is conjugate into a compactly embedded Cartan subalgebra. Such Cartan subalgebras are conjugate under the real inner group by~\cite[Proposition~6.61]{Knapp}. Hence there exists $\varphi_0\in\operatorname{Int}(\gg_{\bar0,\RR})$ such that $\varphi_0(c)\in \hh_{\RR}$.

Finally, $\omega$ commutes with $\ad_{y}$ whenever $y \in \gg_{\bar{0},\RR}$: 
\[
\omega([y,x])=[\omega(x),\omega(y)]=[\omega(x),-y]=[y,\omega(x)], \qquad x \in \gg.
\]
By definition, $\omega$ also commutes with $\exp(\ad_y)$, and therefore with every element of $\operatorname{Int}(\gg_{\bar0,\RR})$. 
\end{proof}

We now show that the normal form in~\eqref{eq::normal_form} can be obtained relative to our fixed compactly embedded Cartan subalgebra $\hh$ by an even inner automorphism.

\begin{lemma}\label{lemm::omega_fixed_cartan_normal_form}
Every $x\in\YY_\omega^{\hom}$ is conjugate under $\operatorname{Int}(\gg_{\bar0,\RR})$ to an element in normal form~\eqref{eq::normal_form} relative to $\hh$.
\end{lemma}

\begin{proof}

By Lemma~\ref{lemm::conjugation_inner}, we may assume that $c=x^2\in\hh_{\RR}$. Following~\cite[proof of Lemma~3.9]{Sherman2025}, we consider the decomposition
\[
\gg^c=\mathfrak a\oplus\mathfrak l_1\oplus\cdots\oplus\mathfrak l_s,
\]
where $\mathfrak a\subseteq\hh$ is even and central, and each $\mathfrak l_j$ is simple or isomorphic to $\mathfrak{gl}(d\vert d)$. Each $\mathfrak l_j$ corresponds to a component of the root subsystem $\{\alpha:\alpha(c)=0\}$ and is therefore stable under $\omega$, since $\omega(\gg^\alpha)=\gg^{-\alpha}$.

Using this decomposition and $[x,c]=0$, we obtain $x=\sum_jx_j$, where $x_j\in(\ll_j)_{\bar1}$. By $\omega$-stability of $\ll_j$, we have $\omega(x_j)=ix_j$. Moreover, since $c=x^2$ is central in $\gg^c$, each $x_j^2$ is central in $\ll_j$. If $\ll_j$ is simple, then $x_j^2=0$, which implies $x_j=0$ by the existence of a nontrivial unitarizable supermodule, as above. Thus only factors isomorphic to $\gl(d\vert d)$ contribute.

On these factors, we can describe the components $x_j$ explicitly using a matrix realization of $\omega$. This allows us to conjugate each $x_j$ into normal form relative to $\hh$. 

Fix an $x_j\neq0$, and 
choose an identification $\ll_j\cong\gl(d\vert d)$ under which $\hh\cap\ll_j$ is diagonal and $\omega$ takes the form
\[
\omega(X)=H^{-1}X^\dagger H,\qquad H=\operatorname{diag}(H_0,H_1),
\]
where $H_0,H_1$ are nonsingular real diagonal matrices. Such an identification exists because $\hh_\RR$ is compactly embedded, so that $\omega(\gg^{\alpha})=\gg^{-\alpha}$.

We show that $H_0$ and $H_1$ are definite. 
In standard matrix form, $x_j=\left(\begin{smallmatrix}0&B\\C&0\end{smallmatrix}\right)$.  Using $x_j^2=-i\lambda I\in\hh\cap\ll_j$ for some $\lambda\in\RR^\times$ and $\omega(x_j)=ix_j$, we obtain
\[
C=-iH_1^{-1}B^\dagger H_0,\qquad B^\dagger H_0B=\lambda H_1.
\]

In particular, $H_0$ and $H_1$ are either both definite or both indefinite. Suppose that both are indefinite. We can then choose nonzero isotropic vectors $u$ and $v$ for $H_0$ and $H_1$, respectively, such that
\[
y=\begin{pmatrix}0&uv^\dagger H_1\\-i vu^\dagger H_0&0\end{pmatrix}
\]
satisfies $y\neq0$, $\omega(y)=iy$, and $y^2=0$. This contradicts the fact that $\YY$ has zero real locus.

Thus both forms are definite, and after a diagonal change of basis we may assume $H_0=\varepsilon_0 I$ and $H_1=\varepsilon_1 I$, where $\varepsilon_0,\varepsilon_1\in\{\pm1\}$. This forces $B=tU$ and $C=-i\lambda t^{-1}U^\dagger$ for some $t>0$ and $U\in U(d)$.

The statement now follows by diagonalizing $B$ through an even unitary conjugation, which brings $x_j$ into the required normal form relative to $\hh$.
\end{proof}

We established that every $x\in\YY_\omega^{\hom}$ of rank $k$ is conjugate, by an even inner automorphism commuting with $\omega$, to an element of the form
\begin{equation}\label{eq::decomp_standard}
x=\sum_{i=1}^{k}x_i,
\qquad
x_i\coloneqq u_{\alpha_i}+z_i v_{\alpha_i},
\end{equation}
where $\alpha_1,\ldots,\alpha_k\in\Delta$ are linearly independent, mutually orthogonal isotropic roots, $0\neq u_{\alpha_i}\in\gg^{\alpha_i}$, $0\neq v_{\alpha_i}\in\gg^{-\alpha_i}$, and $z_i\in\CC^\times$. Moreover, using $\omega(\gg^{\alpha})=\gg^{-\alpha}$, every $x_i$ belongs to $\YY_\omega^{\hom}$.

\begin{lemma}\label{lemm::omega_root_pairs}
If $x\in\YY_\omega^{\hom}$ is in normal form
\eqref{eq::decomp_standard}, then
$
x_i= u_{\alpha_i}+z_i v_{\alpha_i}
\in\YY_\omega^{\hom}
$
for every $i$.
\end{lemma}

The constructions of $\YY_\omega^{\hom}$, $\gg_x$, and $\omega_x$ are preserved by automorphisms commuting with $\omega$. By Lemma~\ref{lemm::conjugation_inner}, every element of $\operatorname{Int}(\gg_{\bar0,\RR})$ commutes with $\omega$. Thus, for any $x\in\YY_\omega^{\hom}$ and $\varphi\in\operatorname{Int}(\gg_{\bar0,\RR})$, we have $\varphi(x)\in\YY_\omega^{\hom}$, $\gg_x\cong\gg_{\varphi(x)}$, and $\gg_x^{\omega_x}\cong\gg_{\varphi(x)}^{\omega_{\varphi(x)}}$. To describe $\gg_x^{\omega_x}$, we can therefore assume that $x$ is in normal form~\eqref{eq::decomp_standard}.

We determine the real form $\gg_x^{\omega_x}$. In the proof of Lemma~\ref{lemm::real_form_descend}, we have seen that $\omega_{x}$ is given on $\gg_{x}$ by 
\begin{equation}
\omega_x\bigl(y+[x,\gg^c]\bigr)
\coloneqq
\omega(y)+[x,\gg^c],
\qquad y\in\gg^x.
\end{equation}
As $x$ is in normal form,  it decomposes as $x=\sum_{i=1}^{k}x_i$, where $x_i=u_{\alpha_i}+z_i v_{\alpha_i}$ and the $\alpha_{i}$ are pairwise orthogonal odd isotropic roots. Set
\begin{equation}
\tt'\coloneqq\bigoplus_{i=1}^{k}\CC t_{\alpha_i},
\qquad
\CC t_{\alpha_i}=[\gg^{\alpha_i},\gg^{-\alpha_i}].
\end{equation}
The subalgebra $\tt'$ is $\omega$-stable, since the $x_i^2$ are nonzero multiples of $t_{\alpha_i}$ and $\omega(x_i^2)=-x_i^2$. In particular, $\omega$ preserves the centralizer of $\tt'$ in $\gg$, denoted by $\mathfrak c_{\gg}(\tt')$. By~\cite[Theorem~3.12]{Sherman2025}, we have
\begin{equation}\label{eq::decomp_def_g_x}
\mathfrak c_{\gg}(\tt')
\cong
\gg_x\times\mathfrak{gl}(1\vert1)^k.
\end{equation}
Each \(\mathfrak{gl}(1\vert1)\)-summand corresponds to one isotropic root pair \(\{\alpha_i,-\alpha_i\}\).

We can now describe the induced real structure. Let $\mathfrak s$ be the subalgebra generated by $\gg^{\pm\alpha_i}$ for $1\leq i\leq k$. It is $\omega$-stable, since $\omega(\gg^{\alpha_i})=\gg^{-\alpha_i}$. Under the decomposition~\eqref{eq::decomp_def_g_x}, $\mathfrak s\cong\sl(1\vert1)^{\oplus k}$, with each $\sl(1\vert1)$ contained in the corresponding $\gl(1\vert1)$-factor. The centralizer of $\sl(1\vert1)$ in $\gl(1\vert1)$ is its one-dimensional center, which is contained in $\tt'$. This implies 
\begin{equation}
\mathfrak c_{\mathfrak c_{\gg}(\tt')}(\mathfrak s)/\tt'\cong\gg_x,
\end{equation}
and the isomorphism intertwines the induced anti-involution with $\omega_x$, since the $\gl(1\vert1)$-factors have trivial $\ad_x$-cohomology. Thus $\gg_x^{\omega_x}$ is obtained by removing the root blocks associated with $\alpha_1,\ldots,\alpha_k$ from the corresponding centralizer. The preceding discussion immediately yields the following proposition.

\begin{proposition}\label{prop::unitary_real_forms_and_DS}
Let $x\in\YY_\omega^{\hom}$ have rank $k$. Then the following descriptions hold. 
\begin{enumerate}
\item[(a)] If $\gg=\sl(m\vert n)$ with $m\neq n$ and $\gg^\omega=\su(p,q\vert0,n)$, then
\[
\gg_x^{\omega_x}
\cong\su(p-a,q-b\vert0,n-k),
\]
where $(a,b)$ is the signature of the Hermitian form $\operatorname{diag}(I_{p,q},-I_n)$ on the even part of $x^2\CC^{p+q\vert n}$.

\item[(b)] If $\gg^\omega=\osp(0,m\vert2n;\RR)$ has even part $\so(m,\RR)\oplus\sp(2n,\RR)$, then
\[
\gg_x^{\omega_x}
\cong\osp(0,m-2k\vert2n-2k;\RR),
\qquad
0\leq k\leq\min(\lfloor m/2\rfloor,n).
\]

\item[(c)] If $\gg^\omega=\osp^{\ast}(2r\vert0,2n)$ has even part $\so^{\ast}(2r)\oplus\mathfrak{usp}(2n)$, then
\[
\gg_x^{\omega_x}
\cong\osp^{\ast}(2r-2k\vert0,2n-2k),
\qquad
0\leq k\leq\min(r,n).
\]
In particular, for the two real forms of $C(n)$ and $k=1$, the real forms $\gg_x^{\omega_x}$ are $\sp(2n-4,\RR)$ and $\mathfrak{usp}(2n-4)$, respectively.

\item[(d)] For the exceptional Lie superalgebras, every nonzero $x\in\YY_\omega^{\hom}$ has rank $k=1$. The corresponding real forms are
\[
\begin{array}{c|c|c}
\gg & (\gg^\omega)_{\bar0} & \gg_x^{\omega_x}\\
\hline
F(4)
& \so(7,\RR)\oplus\sl(2,\RR)
& \su(3)
\\[2pt]
F(4)
& \so(2,5)\oplus\su(2)
& \su(1,2)
\\[2pt]
G(3)
& \gg_{2,\mathrm{compact}}\oplus\sl(2,\RR)
& \su(2)
\\[2pt]
D(2,1;\upalpha)
& \sl(2,\RR)\oplus\su(2)\oplus\su(2)
& \mathfrak u(1).
\end{array}
\]
Here $\mathfrak u(1)$ denotes the one-dimensional real toral Lie algebra.
\end{enumerate}
\end{proposition}

Finally, we describe the $\gg_x$-supermodules $\DS_{x}(\HH)$ obtained from $\omega$-unitarizable $\gg$-supermodules $\HH$.

\subsubsection{Generalized Duflo--Serganova cohomology of unitarizable supermodules}
Let $\HH$ be an $\omega$-unitarizable simple $\gg$-supermodule, and let $x\in\YY_\omega^{\hom}$ such that $c=x^2$ acts semisimply on $\HH$. This condition holds whenever $x$ has the normal form~\eqref{eq::decomp_standard}, since $\HH$ is a highest weight supermodule.

\begin{proposition}\label{prop::form_DS_x_generalized} $\DS_{x}(\HH)\cong \HH^{c}$ and $\DS_{x}(\HH)$ is $\omega_{x}$-unitarizable.
\end{proposition}

\begin{proof} Since $\omega(x)=ix$, every
    $v\in\HH^c$ satisfies
    \[
        \langle xv,xv\rangle
        =\langle v,\omega(x)xv\rangle
        =\langle v,icv\rangle
        =0.
    \]
    Positive definiteness therefore implies $xv=0$.
    Consequently,
    \[
        \HH^c=\ker x,
        \qquad
        \DS_x(\HH)=H(\HH^c,x)\cong\HH^c.
    \]

    The restricted form on $\HH^c$ is positive definite and satisfies for all $y \in \gg^{x}$ and $v,w \in \HH^{c}$
    \[
        \langle yv,w\rangle=\langle v,\omega(y)w\rangle.
    \]
    Since $\omega$ induces $\omega_x$, this proves
    $\omega_x$-unitarizability of $\DS_x(\HH)$.
\end{proof}

The next step is to determine the structure of $\DS_x(\HH)$ as a $\gg_x$-supermodule. In particular, we ask whether $\DS_x(\HH)$ is simple or decomposes as a finite direct sum of unitarizable highest weight $\gg_x$-supermodules. The following lemma is central.

\begin{lemma}\label{lemm::ic_positive_semidefinite}
    Let $\HH$ be an $\omega$-unitarizable $\gg$-supermodule
    and $x\in\YY_\omega^{\hom}$. Then $ic$ is a positive semidefinite operator.
\end{lemma}

\begin{proof}
    Since $\omega(x)=ix$, for every $v\in\HH$ we have
    \[
        \langle v,icv\rangle
        =\langle v,\omega(x)xv\rangle
        =\langle xv,xv\rangle\geq0.\qedhere
    \]
\end{proof}

\begin{proposition}\label{prop::DS_composition}
   Let $x=\sum_{i=1}^{k}x_{i}$ be of rank $k$, written as in~\eqref{eq::decomp_standard} as a sum of pairwise commuting rank-one elements $x_i\in\YY^{\hom}_\omega$. Then
   \[
\operatorname{DS}_x(\mathcal H)
\cong
\operatorname{DS}_{x_k}\circ\cdots\circ
\operatorname{DS}_{x_1}(\mathcal H).
\]
\end{proposition}

\begin{proof}
Since $[x_i,x_j]=0$ for $i\neq j$, we have $x^2=\sum_i x_i^2$. Therefore, for every $v\in\HH$,
\[
    \langle xv,xv\rangle
    =\langle v,ix^2v\rangle
    =\sum_{i=1}^k\langle v,ix_i^2v\rangle
    =\sum_{i=1}^k\langle x_iv,x_iv\rangle.
\]
The Hermitian form is positive definite. Hence $xv=0$ if and only if $x_iv=0$ for every $i$. Thus
\[
    \ker x=\bigcap_{i=1}^k\ker x_i.
\]

We now consider the iterated DS cohomology. By Proposition~\ref{prop::form_DS_x_generalized}, each DS cohomology is naturally identified with the kernel of its defining operator. Moreover, $[x_i,x_j]=0$ implies that $x_j$ preserves $\ker x_i$. Hence, after the first $r$ steps, the resulting supermodule is naturally identified with
$
    \bigcap_{i=1}^r\ker x_i.
$
This proves the claim as an isomorphism of $\gg_x$-supermodules. Here, we use
that both sides carry the action induced by $\gg$ on $\bigcap_{i=1}^k\ker x_i$.
\end{proof}

The preceding proposition reduces the analysis to elements $x\in\YY_\omega^{\hom}$ of rank one. By~\cite[Lemma~3.9]{Sherman2025}, there exists a Cartan subalgebra $\tt$ of $\gg$ such that $x=u_\alpha+zv_\alpha$, where $\alpha\in\Delta(\gg,\tt)_{\bar1}$ is isotropic, $z\in\CC^\times$, and $u_\alpha\in\gg^\alpha$, $v_\alpha\in\gg^{-\alpha}$. Thus $c=x^2=z[u_\alpha,v_\alpha]\in\tt$.

By Lemma~\ref{lemm::conjugation_inner}, there exists $\varphi_0\in\operatorname{Int}(\gg_{\bar0,\RR})$ such that $\varphi_0(c)\in \hh_{\RR}$. Hence $\ad_{ic}$ is conjugate to $\ad_{i\varphi_0(c)}$ and is therefore semisimple with real eigenvalues. Consequently, $c$ determines the triangular decomposition
\begin{equation}
\gg=\gg_-\oplus\gg^c\oplus\gg_+,
\qquad
\gg_\pm\coloneqq\bigoplus_{\pm t>0}\gg_t,
\qquad
\gg_t\coloneqq\{a\in\gg:[ic,a]=ta\}.
\end{equation}
Since $[\gg_s,\gg_t]\subseteq\gg_{s+t}$, the subspaces $\gg_\pm$ are Lie subsuperalgebras and $[\gg^c,\gg_\pm]\subseteq\gg_\pm$.

\begin{lemma}\label{lemm::DS_x_simple_if_x_is_of_rank_1}
Let $x\in\YY_\omega^{\hom}$ have rank one, and let $\HH$ be a simple $\omega$-unitarizable highest weight $\gg$-supermodule. Assume $c=x^{2}$ acts semisimply on $\HH$. Then $\DS_x(\HH)$ is either zero or a simple $\omega_x$-unitarizable $\gg_x$-supermodule.
\end{lemma}

\begin{proof}
Assume $\HH^c\neq0$. By Proposition~\ref{prop::form_DS_x_generalized}, we have $\DS_x(\HH)=\HH^c$, which is $\omega_x$-unitarizable.

Since $ic$ acts positively semidefinitely on $\HH$ by Lemma~\ref{lemm::ic_positive_semidefinite}, we have $\gg_-\HH^{c}=0$. For every nonzero homogeneous $v\in \HH^{c}$, simplicity and the PBW theorem give
\[
\HH=\UE(\gg)v
=\UE(\gg_+)\UE(\gg^c)v.
\]
Taking the zero eigenspace of $ic$ yields $\HH^c=\UE(\gg^c)v$. Thus $\HH^c$ is simple as a $\gg^c$-supermodule.

By~\cite[Theorem~3.12]{Sherman2025}, the rank-one centralizer admits a decomposition
\[
\gg^c\cong\gg_x\times\mathfrak{gl}(1\vert1),
\]
where $x=u_\alpha+zv_\alpha$ belongs to the
$\mathfrak{gl}(1\vert1)$-factor. Choose an even element
$h$ in this factor such that
\[
[h,u_\alpha]=u_\alpha,
\qquad
[h,v_\alpha]=-v_\alpha.
\]
Since $x\HH^c=0$ and $h$ preserves $\HH^c$, we have
$[h,x]\HH^c=0$. Consequently,
\[
(u_\alpha+zv_\alpha)\HH^c
=
(u_\alpha-zv_\alpha)\HH^c=0.
\]
Thus $u_\alpha$, $v_\alpha$, and their bracket act
trivially on $\HH^c$.

It follows that $h$ commutes with the $\gg^c$-action on $\HH^c$ and therefore acts by a scalar by Schur's lemma. Hence every $\gg_x$-submodule
of $\HH^c$ is also a $\gg^c$-submodule.
Since $\HH^c$ is simple over $\gg^c$, it is simple
over $\gg_x$.
\end{proof}

Combining Proposition~\ref{prop::DS_composition} with Lemma~\ref{lemm::DS_x_simple_if_x_is_of_rank_1}, we obtain the following theorem.

\begin{theorem}\label{thm::DS_unitary}
    Let $\HH$ be a simple $\omega$-unitarizable highest weight $\gg$-supermodule, and let $x \in \YY_{\omega}^{\hom}$. Assume $c=x^{2}$ acts semisimply on $\HH$. Then $\DS_{x}(\HH)$ is either zero or a simple $\omega_{x}$-unitarizable $\gg_{x}$-supermodule.
\end{theorem}

\begin{proof}

Let $x=\sum_i x_i$ be the normal form of $x$ (see~\eqref{eq::normal_form}). We first prove that each $x_{i}^{2}$ acts semisimply on $\HH$. If $\HH^c=0$, then $\DS_x(\HH)=0$. Assume $\HH^c\neq\{0\}$ and let $0\neq v\in\HH^c$. Then $x_iv=x_i^2v=0$ for all $i$ by positive definiteness of $\langle\cdot,\cdot\rangle$. Indeed,
\[
0=\langle v,icv\rangle=\sum_{i=1}^k\langle x_iv,x_iv\rangle.
\]
Moreover, since $\HH$ is simple, we have $\HH=\UE(\gg)v$, and the action of $x_i^2$ on $\HH$ is induced by the action of $\ad_{x_i^2}$ on $\UE(\gg)$. 
Since $x_i^2\in\hh$ for all $i$, the adjoint actions $\ad_{x_i^2}$ are simultaneously semisimple on $\UE(\gg)$. Consequently, the $x_i^2$ act semisimply on $\HH$.

The statement now follows by successive application of Proposition~\ref{prop::DS_composition} and Lemma~\ref{lemm::DS_x_simple_if_x_is_of_rank_1}, 
since the $x_i$ commute.
\end{proof}

If $x$ is not of maximal rank, then $\gg_x$ is again a basic classical Lie superalgebra with $\odd\neq0$. By Proposition~\ref{prop::unitary_real_forms_and_DS} and Theorem~\ref{theorem::unitary-highest-lowest-weight}, $\DS_x(\HH)$ is a highest weight $\gg_x$-supermodule. We establish the analogous result when $x$ has maximal rank and is in normal form~\eqref{eq::decomp_standard} relative to $\hh$.

We use the root-space realization of $\gg_x$, that is, we identify the root spaces of $\gg_{x}$ with $\gg^{\gamma}$, where $\gamma \in \Delta_{x}$. Then $\nn_x^{\pm}$ is represented by $\bigoplus_{\gamma\in\Delta_x\cap\Delta^{\pm}}\gg^\gamma$, and the triangular decomposition $\gg=\nn^-\oplus\hh\oplus\nn^+$ induces the triangular decomposition $\gg_x=\nn_x^-\oplus\hh_x\oplus\nn_x^+$.

\begin{lemma}\label{lemm::DS_is_HW}
Let $\HH$ be a simple unitarizable highest weight $\gg$-supermodule, and let $x\in\YY_\omega^{\hom}$ be in normal form~\eqref{eq::decomp_standard} relative to $\hh$. Assume $x$ has maximal rank. If $\DS_x(\HH)\neq0$, then $\DS_x(\HH)$ is a highest weight $\gg_x$-supermodule.
\end{lemma}

\begin{proof} 
Since $x$ is in normal form~\eqref{eq::decomp_standard}, $c=x^2\in\hh$ acts semisimply on the highest weight supermodule $\HH$. By assumption, $\DS_x(\HH)\cong\HH^c\neq \{0\}$ and is therefore simple by Theorem~\ref{thm::DS_unitary}.

Let $v\in\HH^c$ be a nonzero weight vector. Under the above identification, $\nn_x^+$ is represented by elements of $\nn^+$ commuting with $x$, and hence also with $c=x^2$. Since these elements preserve $\HH^c$, we obtain
\begin{equation*}
\UE(\nn_x^+)v\subseteq\UE(\nn^+)v\cap\HH^c.
\end{equation*}
Since $\HH$ is a highest weight supermodule, $\UE(\nn^+)v$ is finite-dimensional, and hence so is $\UE(\nn_x^+)v$. Choose a nonzero weight vector $w\in\UE(\nn_x^+)v$ of maximal weight. Then $\nn_x^+w=0$ by maximality, and simplicity implies $\DS_x(\HH)=\UE(\gg_x)w$.
\end{proof}

We end this section by describing $\DS_x(L(\Lambda))$ when $x=\sum_{j=1}^k x_j$ is in normal form~\eqref{eq::normal_form} relative to $\hh$. We first introduce some notation to replace the positive system. Set $c=x^2\in\hh$,
$\widetilde{\nn}^{\pm}\coloneqq(\nn^{\pm}\cap\gg^c)\oplus\gg_{\mp}$, and $\widetilde{\bb}\coloneqq\hh\oplus\widetilde{\nn}^+$. Then
\begin{equation*}
\gg=\widetilde{\nn}^-\oplus\hh\oplus\widetilde{\nn}^+,
\end{equation*}
and $\widetilde{\bb}$ is a Borel subalgebra. 
Note that the positive roots of $\widetilde{\bb}$ are those $\alpha\in\Delta$ with $i\alpha(c)<0$, together with those $\alpha\in\Delta^+$ satisfying $i\alpha(c)=0$.

\begin{lemma}\label{lemm::simple_new_system}
Let $L(\Lambda;\Delta^+)$ be a simple $\omega$-unitarizable highest weight $\gg$-supermodule. Then there exists a weight $\widetilde{\Lambda}\in\hh^*$ such that
\begin{equation*}
L(\Lambda;\Delta^+)\cong L(\widetilde{\Lambda};\widetilde{\Delta}^+),
\end{equation*}
where $\widetilde{\Delta}^+$ is the positive root system corresponding to $\widetilde{\bb}$.
\end{lemma}

\begin{proof}
    Since $L(\Lambda;\Delta^{+})$ is unitarizable, $ic$ acts as a positive semidefinite operator by Lemma~\ref{lemm::ic_positive_semidefinite}. On the other hand, by definition, $\gg_-$ consists of the negative eigenspaces of $\ad_{ic}$. Thus its elements strictly decrease the eigenvalues of $ic$. If $v_\Lambda$ is the highest weight vector of $L(\Lambda;\Delta^+)$, we conclude that $\UE(\gg_-)v_\Lambda$ is finite-dimensional.

Since $\nn^+v_\Lambda=0$, the subalgebra $\nn^+\cap\gg^c$ annihilates $v_\Lambda$. Together with $[\gg^c,\gg_-]\subseteq\gg_-$, this shows that $\UE(\gg_-)v_\Lambda$ is stable under $\nn^+\cap\gg^c$, which also preserves the $ic$-eigenspaces.

Choose a nonzero weight vector $w\in\UE(\gg_-)v_\Lambda$ of minimal $ic$-eigenvalue and maximal weight among vectors with this eigenvalue. Then $\gg_-w=0$ by minimality and $(\nn^+\cap\gg^c)w=0$ by maximality. Thus $\widetilde{\nn}^+w=0$, so $w$ is a highest weight vector with respect to $\widetilde{\bb}$.
\end{proof}

\begin{proposition}\label{prop::DS_highest_weight_general}
Let $L(\Lambda)$ be a simple $\omega$-unitarizable highest weight $\gg$-supermodule with respect to $\Delta^+$, and let $x\in\YY_\omega^{\hom}$ be in normal form~\eqref{eq::decomp_standard} relative to $\hh$. Let $\widetilde{\Lambda}$ be the highest weight of $L(\Lambda)$ with respect to $\widetilde{\bb}$. Then
\begin{equation*}
\DS_x(L(\Lambda))\neq0
\iff \widetilde{\Lambda}(c)=0.
\end{equation*}
In this case,
\begin{equation*}
\DS_x(L(\Lambda))\cong L_{\gg_x}(\widetilde{\Lambda}\vert_{\hh_x}),
\end{equation*}
where the highest weight is taken with respect to $\bb_x=\hh_x\oplus\nn_x^+$.
\end{proposition}
\begin{proof}

Let $w$ be a highest weight vector of weight $\widetilde{\Lambda}$ with respect to $\widetilde{\bb}$. Then $\gg_-w=0$, and by the PBW theorem,
\begin{equation*}
L(\Lambda)=\UE(\gg_+)\UE(\gg^c)w.
\end{equation*}
Since $\gg_+$ strictly increases the eigenvalues of $ic$, the smallest eigenvalue on $L(\Lambda)$ is $i\widetilde{\Lambda}(c)$, with corresponding eigenspace $\UE(\gg^c)w$. As $ic$ is positive semidefinite, Proposition~\ref{prop::form_DS_x_generalized} implies
\begin{equation*}
\DS_x(L(\Lambda))\neq0\iff\widetilde{\Lambda}(c)=0.
\end{equation*}

In this case, $\DS_x(L(\Lambda))\cong L(\Lambda)^c=\UE(\gg^c)w$. As in the proofs of Lemmas~\ref{lemm::DS_is_HW} and~\ref{lemm::simple_new_system}, $\nn_x^+$ is represented by elements of $\nn^+\cap\gg^c$, which annihilate $w$. Hence $w$ is a highest weight vector for $\gg_x$ of weight $\widetilde{\Lambda}\vert_{\hh_x}$. Simplicity of $\DS_x(L(\Lambda))$ by Theorem~\ref{thm::DS_unitary} completes the proof.

\end{proof}

\section{Formal superdimension}\label{sec::formal_superdimension}
\noindent
A central result in the representation theory of Lie superalgebras is Serganova's proof of the Kac--Wakimoto conjecture~\cite{serganova2011superdimension}: a finite-dimensional simple $\gg$-supermodule $L(\Lambda)$ has nonzero superdimension
$\sdim L(\Lambda)=\dim L(\Lambda)_{\bar0}-\dim L(\Lambda)_{\bar1}$
if and only if $\Lambda$ is maximally atypical, that is, $\at(\Lambda)=\defect(\gg)$. In this section, we extend the notion of superdimension and the Kac--Wakimoto theorem to relative holomorphic discrete series supermodules.

The central idea is the following. Relative holomorphic discrete series supermodules are precisely those unitarizable highest weight supermodules which decompose under $\even$ into a finite direct sum of simple unitarizable highest weight modules, which carry an analogous notion of dimension, namely Harish-Chandra's formal dimension. We use this to define the \emph{formal superdimension} $\fsdim$ as the parity-signed sum of formal dimensions of the $\even$-constituents.  

This definition shares common properties with the usual superdimension of finite-dimensional supermodules. For example, it is additive, multiplicative under tensor products with finite-dimensional supermodules whenever it is defined, and changes sign under parity reversal. 

Nevertheless, two important differences exist. First, it is not preserved under generalized DS cohomology, and second, in the infinite-dimensional case, it is not the evaluation of the supercharacter at the identity. We show that $\fsdim$ equals the evaluation of the $L$-packet supercharacter, which is the average of the supercharacter over the Weyl group, normalized by the compact Weyl group. This identification is central to all proofs concerning $\fsdim$, because it reduces the analysis to supercharacters.  

Our main result is the proof of a Kac--Wakimoto conjecture for $\fsdim$. The proof is twofold. We first treat the integral case in Section~\ref{subsec::KW_conjecture_integral_case} and construct, using derived Zuckerman functors, finite-dimensional simple supermodules which have the same $L$-packet supercharacter, and hence the same superdimension. For the nonintegral case, we relate $\fsdim$ to the superdimension of the DS cohomology with respect to an element of maximal rank in Section~\ref{subsec::KW_conjecture_nonintegral_case}. The latter is either a simple finite-dimensional supermodule, a relative holomorphic discrete series module, or a supermodule of $\osp(1\vert 2s)$, and is never zero. We summarize both results in Section~\ref{subsec::KW_main}.

\subsection{Definition, examples, and properties}\label{subsec::formal_superdimension}
Let $L(\Lambda)$ be a relative holomorphic discrete series $\gg$-supermodule. By definition, its restriction to $\even$ decomposes as
\begin{equation}
L(\Lambda)|_{\even} \cong \bigoplus_\mu L_0(\mu)\otimes E_\mu,
\end{equation}
where the $L_0(\mu)$ are pairwise nonisomorphic simple relative holomorphic discrete series $\even$-modules, regarded as purely even, and the $E_\mu$ are finite-dimensional multiplicity superspaces. The sum is finite, and the grading of $E_\mu$ records the parity of each constituent. In particular, $E_{\mu}$ has a well-defined superdimension $\sdim(E_\mu)
\coloneqq
\dim(E_\mu)_{\bar0}-\dim(E_\mu)_{\bar1}.$

By Theorem~\ref{thm::formal_dimension}, each $L_0(\mu)$ has a positive formal dimension
\begin{equation}
d(\pi_\mu)
=
\prod_{\alpha\in\Delta_c^+}
\frac{(\mu+\rho_c,\alpha)}
     {(\rho_c,\alpha)}
\prod_{\beta\in\Delta_n^+}
\frac{\abs{(\mu+\rho_{\bar0},\beta)}}
     {\abs{(\rho_{\bar0},\beta)}},
\end{equation}
with respect to the normalization of the Haar measure on $X\coloneqq \widetilde{G}_{\RR}/Z$ given by Theorem~\ref{thm::formal_dimension}, which we fix from now on.

We use these formal dimensions to define a weighted superdimension, with the parity of each constituent determining its sign.

\begin{definition}
Let $M$ be a $\gg$-supermodule whose restriction to $\even$ admits a finite decomposition
$
M|_{\even}\cong\bigoplus_\mu L_0(\mu)\otimes E_\mu,
$
where the $L_0(\mu)$ are pairwise nonisomorphic relative holomorphic discrete series $\even$-modules and the $E_\mu$ are finite-dimensional multiplicity superspaces. The \emph{formal superdimension} of $M$ is
\begin{equation*}
\fsdim(M)\coloneqq\sum_\mu\sdim(E_\mu)\,d(\pi_\mu).
\end{equation*}
\end{definition}

\begin{remark}
The terminology reflects the analogy with the ordinary superdimension. The
formal dimensions $d(\pi_\mu)$ are defined by Schur orthogonality and reduce
to ordinary dimensions for compact groups with normalized Haar measure. In
particular, if $\Delta_n^+=\varnothing$, the product formula for $d(\pi_\mu)$
becomes Weyl's dimension formula, and $\fsdim(M)$ equals the usual
superdimension $\sdim M$. Following Harish-Chandra's term \emph{formal degree}, one
could also speak of the \emph{formal superdegree}. The term \emph{formal superdimension} is preferred here, both to emphasize the connection with the Kac--Wakimoto theorem and to avoid confusion with the degree of atypicality.
\end{remark}

\begin{remark}\label{rmk::Haar_normalization_sdim} The definition depends on the normalization of the Haar measure on $X\coloneqq \widetilde{G}_{\RR}/Z$ given by Theorem~\ref{thm::formal_dimension}. Replacing $\mathrm{d}\mu_{X}$ by $c\,\mathrm{d}\mu_{X}$ for some $c>0$ multiplies $\fsdim(M)$ by $c^{-1}$. This means that neither the sign of $\fsdim(M)$ nor its vanishing depends on this normalization. 
\end{remark}

\begin{remark}
The definition of $\fsdim$ does not depend on the chosen decomposition. Indeed, the restriction of $M$ to $\even$ is semisimple of finite length, so its simple constituents are uniquely determined up to isomorphism. Using Schur's lemma, we can identify the multiplicity superspaces by
\begin{equation*}
E_\mu\cong\underline{\operatorname{Hom}}_{\even}(L_0(\mu),M),
\qquad
(E_\mu)_{\bar z}
\cong
\operatorname{Hom}_{\even}(L_0(\mu),M_{\bar z}),
\quad \bar z\in\ZZ_2,
\end{equation*}
and both the isomorphism classes of the constituents and the integers $\sdim(E_\mu)$ are intrinsic to $M$. 
\end{remark}

\begin{remark}
Let $M\cong L(\Lambda;\Delta^+)$ be a simple highest weight $\gg$-supermodule, and let $\widetilde{\Delta}^+$ be another positive system with $\Delta^{+}_{\bar{0}}=\widetilde{\Delta}_{\bar{0}}^+$. Then there exist a highest weight $\Lambda'$ and $\varepsilon\in \ZZ_{2}$ (see Lemma~\ref{lemm::odd_reflections_HW}) such that $L(\Lambda;\Delta^+)\cong\Pi^\varepsilon L(\Lambda';\widetilde{\Delta}^+)$. In particular, the formal superdimension changes only by a sign:
\begin{equation*}
\fsdim L(\Lambda;\Delta^+)=(-1)^\varepsilon\fsdim L(\Lambda';\widetilde{\Delta}^+).
\end{equation*}
\end{remark}

The definition applies to relative holomorphic discrete series supermodules $L(\Lambda)$. To compute their formal superdimensions, one needs their decomposition into simple $\even$-modules. This branching problem was solved for arbitrary atypicality in~\cite{Schmidt_unitary_branching}, which provides both recursive and closed formulas. The latter also yields a closed formula for the formal superdimension. All examples below are computed using the recursive formula.

\begin{example}
Let $\gg=\sl(2\vert1)$ with real form $\su(1,1\vert0,1)$ and standard simple system
$\{\epsilon_1-\epsilon_2,\epsilon_2-\delta\}$.
We take all highest weight vectors to be even. The normalization fixed above yields $d(\pi_{(a,b\vert c)})=b-a-1$.

We consider two infinite-dimensional unitarizable supermodules. In both cases, every simple $\even$-constituent is a relative holomorphic discrete series module. Thus both supermodules are relative holomorphic discrete series supermodules.

\begin{enumerate}
\item[(a)] \emph{Maximally atypical example.}
The supermodule $L(-2,0\vert0)$ has highest weight of maximal atypicality and decomposes as
\begin{equation*}
L(-2,0\vert0)|_{\even}
\cong
L_0(-2,0\vert0)\oplus\Pi L_0(-3,0\vert1).
\end{equation*}
Its formal superdimension is
$
\fsdim L(-2,0\vert0)=1-2=-1.
$

\item[(b)] \emph{Typical example.}
The supermodule $L(-3,1\vert0)$ has typical highest weight and decomposes as
\begin{equation*}
\begin{aligned}
L(-3,1\vert0)|_{\even}
\cong{}
L_0(-3,1\vert0)
\oplus\Pi L_0(-4,1\vert1)\oplus\Pi L_0(-3,0\vert1)
\oplus L_0(-4,0\vert2).
\end{aligned}
\end{equation*}
Its formal superdimension is
$
\fsdim L(-3,1\vert0)=3-4-2+3=0.
$
\end{enumerate}
\end{example}

\begin{example}\label{ex::sl23_formal_superdimension}
Let $\gg=\sl(2\vert3)$ with real form $\su(1,1\vert0,3)$ and
nonstandard simple system
$
\{\epsilon_1-\delta_1,\delta_1-\delta_2,
  \delta_2-\delta_3,\delta_3-\epsilon_2\}.
$
The corresponding Weyl vector is $\rho=(-1,1\vert1,0,-1)$.
We take the highest weight vector of each $\gg$-supermodule to be even. Moreover, with the
normalization fixed above, the formal dimension is
\[
d(\pi_{(a,b\vert c_1,c_2,c_3)})
=(b-a-1)
\frac{(c_1-c_2+1)(c_1-c_3+2)(c_2-c_3+1)}{2}.
\]
We consider three infinite-dimensional unitarizable supermodules.
In each case, every simple $\even$-constituent is a relative holomorphic
discrete series module. Thus all three supermodules are relative
holomorphic discrete series supermodules.

\begin{enumerate}
\item[(a)] \emph{Maximally atypical example.}
The supermodule $L(-2,0\vert2,0,0)$ has maximal atypicality and
decomposes as
\[
\begin{aligned}
L(-2,0\vert2,0,0)\vert_{\even}
\cong{}&L_0(-2,0\vert2,0,0)
\oplus\Pi L_0(-2,1\vert1,0,0)\\
&\oplus\Pi L_0(-3,0\vert2,1,0)
\oplus L_0(-3,1\vert1,1,0)\\
&\oplus L_0(-4,0\vert2,1,1)
\oplus\Pi L_0(-4,1\vert1,1,1).
\end{aligned}
\]
Indeed, for $\Lambda=(-2,0\vert2,0,0)$, the two linearly independent,
orthogonal isotropic roots $\epsilon_1-\delta_1$ and
$\delta_3-\epsilon_2$ satisfy
\[
(\Lambda+\rho,\epsilon_1-\delta_1)
=(\Lambda+\rho,\delta_3-\epsilon_2)=0.
\]
Since $\gg$ has defect two, $\Lambda$ is maximally atypical.
Its formal superdimension is
\[
\fsdim L(-2,0\vert2,0,0)=6-6-16+9+9-4=-2.
\]

\item[(b)] \emph{Maximally atypical example.}
The supermodule $L(-3,0\vert3,0,0)$ also has maximal atypicality and
decomposes as
\[
\begin{aligned}
L(-3,0\vert3,0,0)\vert_{\even}
\cong{}&L_0(-3,0\vert3,0,0)
\oplus\Pi L_0(-3,1\vert2,0,0)\\
&\oplus\Pi L_0(-4,0\vert3,1,0)
\oplus L_0(-4,1\vert2,1,0)\\
&\oplus L_0(-5,0\vert3,1,1)
\oplus\Pi L_0(-5,1\vert2,1,1).
\end{aligned}
\]
The same two isotropic roots as for (a) are orthogonal to $\Lambda+\rho$.
Its formal superdimension is
\[
\fsdim L(-3,0\vert3,0,0)=20-18-45+32+24-15=-2.
\]

\item[(c)] \emph{Typical example.}
The supermodule $L(-2,2\vert0,0,0)$ is typical and decomposes as
\[
\begin{aligned}
L(-2,2\vert0,0,0)\vert_{\even}
\cong{}&L_0(-2,2\vert0,0,0)
\oplus L_0(-2,4\vert0,-1,-1)\\
&\oplus L_0(-3,3\vert0,0,0)
\oplus L_0(-3,3\vert1,0,-1)\\
&\oplus L_0(-3,5\vert0,-1,-1)
\oplus L_0(-4,2\vert1,1,0)\\
&\oplus L_0(-4,4\vert0,0,0)
\oplus L_0(-4,4\vert1,0,-1)\\
&\oplus L_0(-5,3\vert1,1,0)
\oplus L_0(-5,5\vert0,0,0)\\
&\oplus\Pi L_0(-2,3\vert0,0,-1)
\oplus\Pi L_0(-2,5\vert-1,-1,-1)\\
&\oplus\Pi L_0(-3,2\vert1,0,0)
\oplus\Pi L_0(-3,4\vert0,0,-1)\\
&\oplus\Pi L_0(-3,4\vert1,-1,-1)
\oplus\Pi L_0(-4,3\vert1,0,0)\\
&\oplus\Pi L_0(-4,3\vert1,1,-1)
\oplus\Pi L_0(-4,5\vert0,0,-1)\\
&\oplus\Pi L_0(-5,2\vert1,1,1)
\oplus\Pi L_0(-5,4\vert1,0,0).
\end{aligned}
\]
Here $\Lambda+\rho=(-3,3\vert1,0,-1)$ is not orthogonal to any
odd root. Its formal superdimension is
\[
\begin{aligned}
\fsdim L(-2,2\vert0,0,0)
={}&(3+15+5+40+21+15+7+56+21+9)\\
&-(12+6+12+18+36+18+36+24+6+24)=0.
\end{aligned}
\]
\end{enumerate}
\end{example}

\begin{example}
We take all highest weight vectors to be even.

\begin{enumerate}
\item[(a)] \emph{The case $\gg=\mathfrak{osp}(1\vert4)$.}
Consider the real form $\mathfrak{osp}(0,1\vert4;\mathbb R)$ with standard simple system $\{\delta_1-\delta_2,\delta_2\}$.  For the $\even\cong\mathfrak{sp}(4)$-modules below, the formal dimension is
\begin{equation*}
d(\pi_{(a,b)})
=
\frac{(a-b+1)(-a-b-3)(-a-2)(-b-1)}{6}.
\end{equation*}

The infinite-dimensional unitarizable supermodule $L(-\frac52,-\frac52)$ decomposes as
\begin{equation*}
\begin{aligned}
L\left(-\tfrac52,-\tfrac52\right)|_{\even}
\cong{}&
L_0\left(-\tfrac52,-\tfrac52\right)
\oplus\Pi L_0\left(-\tfrac52,-\tfrac72\right) \oplus L_0\left(-\tfrac72,-\tfrac72\right).
\end{aligned}
\end{equation*}
Every simple $\even$-constituent is a relative holomorphic discrete series module. Thus $L(-\frac52,-\frac52)$ is a relative holomorphic discrete series supermodule. Its formal superdimension is
\[
\fsdim L\left(-\tfrac52,-\tfrac52\right)
=\tfrac14-\tfrac54+\tfrac52=\tfrac32.
\]

\item[(b)] \emph{The case $\gg=\sl(3\vert1)$.}
Consider the real form $\su(1,2\vert0,1)$ with standard simple system
$\{\epsilon_1-\epsilon_2,\epsilon_2-\epsilon_3,\epsilon_3-\delta\}$. 

For every real $t>2$, the infinite-dimensional unitarizable supermodule $L(-t,0,0\vert0)$ decomposes as
\begin{equation*}
L(-t,0,0\vert0)|_{\even}
\cong
L_0(-t,0,0\vert0)
\oplus\Pi L_0(-t-1,0,0\vert1).
\end{equation*}
Both constituents are relative holomorphic discrete series modules. Thus $L(-t,0,0\vert0)$ is a relative holomorphic discrete series supermodule. Since
$
d(\pi_{(-t,0,0\vert0)})=\frac{(t-1)(t-2)}{2}
$
and
$
d(\pi_{(-t-1,0,0\vert1)})=\frac{t(t-1)}{2},
$
its formal superdimension is
\begin{equation*}
\fsdim L(-t,0,0\vert0)
=
\frac{(t-1)(t-2)}{2}
-
\frac{t(t-1)}{2}
=
1-t.
\end{equation*}
In particular,
$
\fsdim L(-\frac52,0,0\vert0)
=\frac38-\frac{15}{8}=-\frac32.
$ These supermodules are maximally atypical. Indeed, we have
$
(\Lambda+\rho,\epsilon_3-\delta)=0,
$ where $\rho=(\frac12,-\frac12,-\frac32\vert\frac32)$.
\end{enumerate}
\end{example}

\begin{remark}
By Remark~\ref{rmk::Haar_normalization_sdim}, ratios of formal superdimensions are independent of the normalization. To see that these ratios need not be rational, consider $\gg=\sl(3\vert1)$ with real form $\su(1,2\vert0,1)$. As above, the relative holomorphic discrete series supermodules $L(-t,0,0\vert0)$, with $t>2$ and even highest weight vectors, satisfy
$
\fsdim L(-t,0,0\vert0)=1-t.
$
In particular, $\fsdim L(-3,0,0\vert0)=-2$, so
\[
\frac{\fsdim L(-t,0,0\vert0)}
     {\fsdim L(-3,0,0\vert0)}
=
\frac{t-1}{2},
\qquad t>2.
\]
This ratio is irrational whenever $t$ is irrational. Consequently, no normalization of the Haar measure makes all formal superdimensions integral.
\end{remark}

These examples suggest the following statement, which we prove in the subsequent sections: the formal superdimension $\fsdim(L(\Lambda))$ is nonzero if and only if $\Lambda$ is maximally atypical.

For a finite-dimensional $\gg$-supermodule $L(\Lambda)$ and $x\in\YY^{\hom}$, the Duflo--Serganova functor preserves the ordinary superdimension:
\begin{equation}
\sdim L(\Lambda)=\sdim\DS_x(L(\Lambda)).
\end{equation}
One might expect the same identity in the infinite-dimensional case if $x \in \YY^{\hom}_{\omega}$. The following examples show that this is false, that is, the formal superdimension is not an invariant of generalized Duflo--Serganova cohomology. 

\begin{example}
We consider two relative holomorphic discrete series supermodules with even highest weight vectors.

\begin{enumerate}
\item[(a)] \emph{The case $\gg=\sl(2\vert1)$.}
Consider the real form $\su(1,1\vert0,1)$ with standard positive system. The infinite-dimensional unitarizable supermodule $L(-2,0\vert0)$ decomposes as
\begin{equation*}
L(-2,0\vert0)|_{\even}
\cong
L_0(-2,0\vert0)\oplus\Pi L_0(-3,0\vert1).
\end{equation*}
Both constituents are relative holomorphic discrete series modules, with formal dimensions $1$ and $2$. Thus $L(-2,0\vert0)$ is a relative holomorphic discrete series supermodule, and
$
\fsdim L(-2,0\vert0)=1-2=-1.
$

Let $\alpha=\epsilon_2-\delta$. Since
$
(\Lambda+\rho_{\bar0}-\rho_{\bar1},\alpha)=0
$
and $\sl(2\vert1)$ has defect one, $\Lambda$ is maximally atypical. Set
$
x=E_{23}-iE_{32}.
$
Then $\omega(x)=ix$, and
$
\DS_x(L(\Lambda))\cong\mathbb C^{1\vert0}
$
as a module for $\gg_x=0$. Consequently,
\begin{equation*}
\sdim\DS_x(L(\Lambda))=1\neq-1=\fsdim L(\Lambda).
\end{equation*}

\item[(b)] \emph{The case $\gg=\sl(2\vert3)$.}
Consider the real form $\su(1,1\vert0,3)$ with nonstandard positive
system, and set $\Lambda=(-2,0\vert2,0,0)$.
By Example~\ref{ex::sl23_formal_superdimension}(a), the
infinite-dimensional unitarizable supermodule $L(\Lambda)$ is a
relative holomorphic discrete series supermodule, with
\[
\fsdim L(\Lambda)=6-6-16+9+9-4=-2.
\]
Set $x=E_{25}-iE_{52}$. For the anti-involution associated with
$\su(1,1\vert0,3)$, we have
\[
\omega(x)=ix,\qquad c\coloneqq x^2=-i(E_{22}+E_{55}).
\]
Thus $x\in\mathcal{Y}^{\mathrm{hom}}_{\omega}$ has rank one, $\gg_x\cong\sl(1\vert2)$, and $\gg_x^{\omega_x}\cong\su(1,0\vert0,2)$. The remaining coordinates are $\epsilon_1,\delta_1,\delta_2$.
Since $\Lambda(E_{22}+E_{55})=0$, the even highest weight vector
survives in $\DS_x(L(\Lambda))=L(\Lambda)^c$.
More explicitly, the zero eigenspace of $ic=E_{22}+E_{55}$ receives
contributions only from $L_0(-2,0\vert2,0,0)$ and
$\Pi L_0(-3,0\vert2,1,0)$ in the preceding decomposition.
Their zero eigenspaces have dimensions three and two, respectively.
Consequently,
\[
\DS_x(L(\Lambda))
\cong L_{\sl(1\vert2)}((-2\vert2,0)),\qquad
\fsdim\DS_x(L(\Lambda))=1.
\]
In particular, the DS cohomology is finite-dimensional, and 
\[
\fsdim\DS_x(L(\Lambda))=1\neq-2=\fsdim L(\Lambda).
\]
\end{enumerate}
\end{example}

\subsubsection{Formal superdimension via \texorpdfstring{$L$}{L}-packets}
\label{subsec::L_packet}

To study the formal superdimension, we express it as the value at zero of an $L$-packet supercharacter. We first briefly recall the notion of a formal supercharacter.

Let $M$ be a highest weight supermodule. Then $M$ has finite-dimensional weight spaces $M^{\mu}$ (Lemma~\ref{lemm::properties_HWM}). The supercharacter of $M$ is the formal expression
\begin{equation}
\sch M
\coloneqq
\sum_{\mu\in\hh^\ast}\sdim(M^\mu)e^\mu,
\end{equation}
where the $e^\mu$ are formal symbols subject to $e^\mu e^\nu=e^{\mu+\nu}$. Here, no convergence is assumed.

Let $L(\Lambda)$ be a relative holomorphic discrete series supermodule. Its restriction to $\even$ decomposes as
\begin{equation}
L(\Lambda)|_{\even}
\cong
\bigoplus_\mu L_0(\mu)\otimes E_\mu,
\end{equation}
where the $L_0(\mu)$ are relative holomorphic discrete series $\even$-modules, regarded as purely even, and the $E_\mu$ are finite-dimensional multiplicity superspaces. Set
$
a_\mu\coloneqq\sdim E_\mu.
$
Then
\begin{equation}
\sch L(\Lambda)
=
\sum_\mu a_\mu\,\operatorname{ch}L_0(\mu),
\qquad
\fsdim L(\Lambda)
=
\sum_\mu a_\mu\,d(\pi_\mu),
\end{equation}
where $\pi_\mu$ denotes the relative holomorphic discrete series representation of $\widetilde{G}_{\RR}$ corresponding to $L_0(\mu)$.

Every simple $\even$-constituent of $L(\Lambda)$ admits an Enright--Willenbring character formula~\cite{Enright_Willenbring}. For relative holomorphic discrete series modules, this formula simplifies by the Harish-Chandra condition of Theorem~\ref{thm::HC_condition}.

\begin{lemma}\label{lemm::character_formula_L_0_relative}
Let $L_0(\mu)$ be a simple $\even$-constituent of a relative holomorphic discrete series supermodule $L(\Lambda)$. Then
\begin{equation*}
\operatorname{ch}L_0(\mu)
=
\frac{
\sum_{w\in W_c}
(-1)^{\ell_c(w)}
e^{w(\mu+\rho_{\bar0})-\rho_{\bar0}}
}{
\prod_{\alpha\in\Delta_{\bar0}^+}
(1-e^{-\alpha})
},
\end{equation*}
where $\ell_c$ denotes the length function on $W_c$ determined by $\Delta_c^+$.
\end{lemma}

\begin{proof}
In the Enright--Willenbring formula, a noncompact positive root $\beta$ can contribute a generating reflection only if
$
2\frac{(\mu+\rho_{\bar0},\beta)}{(\beta,\beta)}\in\ZZ_{>0}.
$
By Theorem~\ref{thm::HC_condition}, however,
$
2\frac{(\mu+\rho_{\bar0},\beta)}{(\beta,\beta)}<0
$
for every $\beta\in\Delta_n^+$. Thus the corresponding subgroup is trivial, and the character formula reduces to
\begin{equation*}
\operatorname{ch}L_0(\mu)
=
\frac{
\sum_{w\in W_c}
(-1)^{\ell_c(w)}
e^{w(\mu+\rho_c)-\rho_c}
}{
\prod_{\alpha\in\Delta_{\bar0}^+}
(1-e^{-\alpha})
}.
\end{equation*}
The result follows because $\rho_n=\rho_{\bar0}-\rho_c$ is $W_c$-invariant.
\end{proof}

Motivated by the Weyl-group description of discrete series
$L$-packets, we sum the translates of the supercharacter
under the full Weyl group $W$ and normalize by $|W_c|$.
Here we regard $\sch L(\Lambda)$ as a meromorphic function
on $\hh$, with Weyl-group action
\begin{equation*}
\bigl(w(\sch L(\Lambda))\bigr)(H)
=
\sch L(\Lambda)(w^{-1}H),
\qquad w\in W,\ H\in\hh,
\end{equation*}
understood as an identity of meromorphic functions.

\begin{definition}\label{def::L_packet_supercharacter}
The \emph{$L$-packet supercharacter} of a relative holomorphic discrete series supermodule $L(\Lambda)$ is
\begin{equation*}
\Theta_\Lambda
\coloneqq
\frac{1}{|W_c|}
\sum_{w\in W}
w\bigl(\sch L(\Lambda)\bigr)=
\frac{
\sum_\mu a_\mu
\sum_{w\in W}(-1)^{\ell(w)}
e^{w(\mu+\rho_{\bar0})-\rho_{\bar0}}
}{
\prod_{\alpha\in\Delta_{\bar0}^+}
(1-e^{-\alpha})
},
\end{equation*}
where $\ell$ is the length function on $W$ determined by $\Delta_{\bar0}^+$.
\end{definition}

\begin{remark}
    The supermodule $L(\Lambda)$ is finite-dimensional if and only if $\gg_{\bar{0},\RR}$ is a compact Lie algebra. In this case, $W=W_{c}$ and $\sch L(\Lambda)$ is $W$-invariant. Hence, the $L$-packet supercharacter is precisely the supercharacter $\sch L(\Lambda)=\Theta_{\Lambda}$. However, in general, $\sch L(\Lambda)$ is only $W_{c}$-invariant.
\end{remark}

The $L$-packet supercharacter $\Theta_\Lambda$ defines a meromorphic function on the complex Cartan subalgebra $\hh$. Since
\begin{equation}
e^{\rho_{\bar0}(H)}
\prod_{\alpha\in\Delta_{\bar0}^+}
(1-e^{-\alpha(H)})
=
\prod_{\alpha\in\Delta_{\bar0}^+}
2\sinh\left(\frac{\alpha(H)}{2}\right),
\end{equation}
we can write
\begin{equation}
\Theta_\Lambda(H)
=
\frac{
\sum_\mu a_\mu
\sum_{w\in W}(-1)^{\ell(w)}
e^{w(\mu+\rho_{\bar0})(H)}
}{
\prod_{\alpha\in\Delta_{\bar0}^+}
2\sinh\left(\frac{\alpha(H)}{2}\right)
},
\qquad H\in\hh.
\end{equation}
The singularities are contained in the hyperplanes 
$
\{H\in\hh:\alpha(H)=2\pi i k\},
$
where $\alpha\in\Delta_{\bar0}^+$ and $k\in\mathbb Z$. We show that $\Theta_\Lambda$ extends holomorphically
to a neighborhood of the origin.

\begin{lemma}\label{lemm::holomorphic}
The $L$-packet supercharacter $\Theta_\Lambda$ extends holomorphically to
\begin{equation*}
\mathcal U
=
\left\{
H\in\hh:
|\operatorname{Im}\alpha(H)|<2\pi
\text{ for all }\alpha\in\Delta_{\bar0}^+
\right\}.
\end{equation*}
\end{lemma}

\begin{proof}
We show that the singularities of $\Theta_\Lambda$ in $\mathcal U$ are removable using the Weierstrass division theorem. Consider the numerator of $\Theta_{\Lambda}$:
\[
N_\Lambda(H)
\coloneqq
\sum_\mu a_\mu
\sum_{w\in W}(-1)^{\ell(w)}
e^{(w(\mu+\rho_{\bar0}))(H)}.
\]
This is an entire $W$-anti-invariant function. Here, $W$-anti-invariant means $N_\Lambda(s_\alpha H)=-N_\Lambda(H)$ for every $\alpha\in\Delta_{\bar0}^+$, $s_{\alpha}\in W$ being the reflection with respect to $\alpha$. Since $s_\alpha H=H$ whenever $\alpha(H)=0$, the function $N_\Lambda$ vanishes on $\ker\alpha$.

Consider one $\alpha \in \Delta_{\bar{0}}^{+}$ and choose complex linear coordinates $z_1,\dots,z_r$ on $\hh$ such that $z_1=\alpha(H)$, and write $z'=(z_2,\dots,z_r)$. These coordinates identify $\hh$ with $\mathbb C^r$ and $\ker\alpha$ with $\{z_1=0\}$.

By the Weierstrass division theorem, there exist entire functions $Q_\alpha$ and $R_\alpha$ such that
\[
N_\Lambda(z_1,z')
=
z_1Q_\alpha(z_1,z')+R_\alpha(z').
\]
For $z_1=0$ we obtain $R_\alpha=0$, hence
$
N_\Lambda(H)=\alpha(H)Q_\alpha(H)
$
for all $H\in\hh$ and $N_{\Lambda}$ is divisible by $\alpha$. 

The positive even roots define distinct hyperplanes,
so division by one such root preserves vanishing
on the remaining root hyperplanes, by continuity. Repeating this, one obtains an entire function $Q_\Lambda$ such that
\[
N_\Lambda(H)
=
Q_\Lambda(H)
\prod_{\alpha\in\Delta_{\bar0}^+}\alpha(H).
\]
At every $H$ with $\prod_{\alpha\in\Delta_{\bar0}^+}2\sinh\!\left(\frac{\alpha(H)}2\right)\neq0$, it follows that
\begin{equation}\label{eq::formula_holomorphic_extension}
\Theta_\Lambda(H)
=
Q_\Lambda(H)
\prod_{\alpha\in\Delta_{\bar0}^+}
\frac{\alpha(H)}
{2\sinh\!\left(\frac{\alpha(H)}2\right)}.
\end{equation}
The function $z/(2\sinh(z/2))$ has a removable singularity at $z=0$, with value $1$, and its only poles lie at $z\in2\pi i(\ZZ\setminus\{0\})$. In particular, it is holomorphic in the strip $|\operatorname{Im}z|<2\pi$, so that~\eqref{eq::formula_holomorphic_extension} defines a holomorphic extension of $\Theta_{\Lambda}$ in $\mathcal{U}$.
\end{proof}

The set $\mathcal U$ is an open neighborhood of $0$ in $\hh$ and its intersection with $i\hh_{\RR}$ is an open neighborhood of $0$ in $i\hh_{\RR}$. By Lemma~\ref{lemm::holomorphic}, the restriction of $\Theta_\Lambda$ to this intersection is real analytic. In particular, $\Theta_\Lambda(0)$ is well defined. We compute this value.

\begin{theorem}\label{thm::L_packet_and_sdim}
One has
\[
\Theta_\Lambda(0)
=
\sum_\mu a_\mu
\prod_{\alpha\in\Delta_{\bar0}^{+}}
\frac{(\mu+\rho_{\bar0},\alpha)}
     {(\rho_{\bar0},\alpha)}
=
(-1)^{|\Delta_n^+|}\fsdim L(\Lambda).
\]
\end{theorem}

\begin{proof}
We first show that, for every $\mu$,
\begin{equation}\label{eq::evaluation_mu_single}
\left.
\frac{
\sum_{w\in W}(-1)^{\ell(w)}
e^{(w(\mu+\rho_{\bar0}))(H)}
}{
 e^{\rho_{\bar0}(H)}
\prod_{\alpha\in\Delta_{\bar0}^{+}}(1-e^{-\alpha(H)})
}
\right|_{H=0}
=
\prod_{\alpha\in\Delta_{\bar0}^{+}}
\frac{(\mu+\rho_{\bar0},\alpha)}
     {(\rho_{\bar0},\alpha)}.
\end{equation}
The evaluation on the left refers to the holomorphic extension established above. We determine the degree $N\coloneqq \vert \Delta_{\bar{0}}^{+}\vert$ Taylor terms of the numerator and denominator of the left-hand side. 
 
For $\lambda\in\hh^*$, expand
\[
\sum_{w\in W}(-1)^{\ell(w)}e^{(w\lambda)(H)}
=
\sum_{k=0}^{\infty}\frac1{k!}
\sum_{w\in W}(-1)^{\ell(w)}\bigl((w\lambda)(H)\bigr)^k.
\]
Let
$
P_k(H)
=
\frac1{k!}\sum_{w\in W}(-1)^{\ell(w)}
\bigl((w\lambda)(H)\bigr)^k.
$
This is a homogeneous polynomial of degree $k$, unless it is zero.
Because $P_k$ is $W$-anti-invariant, it vanishes on every hyperplane $\alpha(H)=0$. Thus every positive even root $\alpha$, viewed as a linear polynomial, divides $P_k$. The positive even roots are pairwise nonproportional, so their product divides $P_k$. Consequently, there exists some polynomial $Q_{k}$ such that
\[
P_k(H)
=
\left(\prod_{\alpha\in\Delta_{\bar0}^{+}}\alpha(H)\right)Q_k(H).
\]
This product has degree $N=|\Delta_{\bar0}^{+}|$. If $P_k\neq0$, then
$
k=\deg P_k=N+\deg Q_k\ge N.
$
Therefore, when $k<N$, the only possibility is
$P_k=0$. In particular, the term of degree $N$ consequently has the form
\[
P_{N}(H)
=
c(\lambda)\prod_{\alpha\in\Delta_{\bar0}^{+}}\alpha(H),
\]
where $c(\lambda)$ is a homogeneous polynomial of degree $N$.

Moreover, $P_{N}(H)$ is also $W$-anti-invariant in $\lambda$. Hence $c(\lambda)$ vanishes on every hyperplane $(\lambda,\alpha)=0$ and is divisible by
$
\prod_{\alpha\in\Delta_{\bar0}^{+}}(\lambda,\alpha).
$
Since both polynomials have degree $N$, there is a constant $C$ such that
\[
c(\lambda)
=
C\prod_{\alpha\in\Delta_{\bar0}^{+}}(\lambda,\alpha).
\]

Next, we determine $C$. By Weyl's denominator formula (\emph{cf.}~\cite[Corollary~25.5]{EtingofLieGroupsLieAlgebras}), one has 
\[
\sum_{w\in W}(-1)^{\ell(w)}
e^{(w\rho_{\bar0})(H)}
=
e^{\rho_{\bar0}(H)}
\prod_{\alpha\in\Delta_{\bar0}^{+}}(1-e^{-\alpha(H)})
=
\prod_{\alpha\in\Delta_{\bar0}^{+}}
2\sinh\!\left(\frac{\alpha(H)}2\right).
\]
Using a Taylor expansion, one sees immediately that the homogeneous term of degree $N$ on the right is
$\prod_{\alpha\in\Delta_{\bar0}^{+}}\alpha(H)$. Thus
$c(\rho_{\bar0})=1$, and consequently
\[
c(\lambda)
=
\prod_{\alpha\in\Delta_{\bar0}^{+}}
\frac{(\lambda,\alpha)}{(\rho_{\bar0},\alpha)}.
\]

Altogether, the numerator has no terms of degree less than $N$, and its degree-$N$ term is
$
\left(
\prod_{\alpha\in\Delta_{\bar0}^{+}}
\frac{(\lambda,\alpha)}{(\rho_{\bar0},\alpha)}
\right)
\prod_{\alpha\in\Delta_{\bar0}^{+}}\alpha(H).
$
The denominator likewise has no terms of degree less than $N$, and its degree-$N$ term is
$
\prod_{\alpha\in\Delta_{\bar0}^{+}}\alpha(H).
$

Choose a regular element $H_0\in i\hh_{\RR}$, so that $\alpha(H_0)\neq0$ for every even root $\alpha$. Since $\mathcal U$ is an open neighborhood of $0$, we have $tH_0\in\mathcal U$ for all sufficiently small $t\in\RR$. Moreover, the denominator of $\Theta_\Lambda(tH_0)$ is nonzero for sufficiently small $t\neq0$.

As $t\to0$, the numerator and denominator have the expansions
\[
t^N c(\lambda)
\prod_{\alpha\in\Delta_{\bar0}^+}\alpha(H_0)
+o(t^N)
\qquad\text{and}\qquad
t^N\prod_{\alpha\in\Delta_{\bar0}^+}\alpha(H_0)
+o(t^N),
\]
respectively. Their quotient therefore tends to $c(\lambda)$. Since the quotient extends holomorphically to $0$, this limit equals its value at $0$. Taking $\lambda=\mu+\rho_{\bar0}$ proves~\eqref{eq::evaluation_mu_single}. 
The statement follows by summing over $\even$-constituents and using the formula for the formal dimension $d(\pi_{\mu})$ in Theorem~\ref{thm::formal_dimension} and the Harish-Chandra condition, which yields
\[
\prod_{\alpha\in\Delta_{\bar0}^{+}}
\frac{(\mu+\rho_{\bar0},\alpha)}
     {(\rho_{\bar0},\alpha)}
=
(-1)^{|\Delta_n^+|}d(\pi_\mu).\qedhere
\]
\end{proof}

\begin{example}
\label{ex::osp_one_formal_superdimension}
Let $L(\Lambda)$ be a relative holomorphic discrete series
$\osp(1\vert2n)$-supermodule with highest weight
$
\Lambda=\sum_{i=1}^n\lambda_i\delta_i.
$
Set $t_i=\lambda_i+n-i+\frac12$, and assume that the highest weight vector has even parity. We compute $\fsdim L(\Lambda)$ explicitly using Theorem~\ref{thm::L_packet_and_sdim}, that is, we evaluate $\Theta_{\Lambda}$ at zero.

We first show that the supercharacter of $L(\Lambda)$ is 
\begin{equation}\label{eq::supercharacter_osp}
\sch L(\Lambda)
=
\frac{
\sum_{w\in W_c}\operatorname{sgn}(w)e^{w(\Lambda+\rho)}
}{\mathcal R}, \qquad \mathcal R
=
e^\rho
\frac{\prod_{\alpha\in\Delta_{\bar0}^+}(1-e^{-\alpha})}
     {\prod_{\beta\in\Delta_{\bar1}^+}(1-e^{-\beta})}.
\end{equation}

By assumption, $L_{0}(\Lambda)$ is a relative holomorphic discrete series module. Compact dominance and the Harish-Chandra conditions give
$\lambda_i-\lambda_{i+1}\in\ZZ_{\geq0}$ and
$\lambda_i+n-i+1<0$. In particular,
$t_1>\cdots>t_n$ and $t_i<0$. This must also hold for any $\even$-constituent. Let $F_\Lambda$ be the finite-dimensional simple
$\kk$-module of highest weight $\Lambda$, placed in
even degree. Consider the parabolic subalgebra
\begin{equation*}
\mathfrak q
=
\kk\oplus \pp^{+}, \qquad \pp^{\pm} \coloneqq
\bigoplus_{\alpha\in\Delta_n^+}\gg^{\pm\alpha}
\oplus
\bigoplus_{\beta\in\Delta_{\bar1}^+}\gg^{\pm\beta}.
\end{equation*}
Extend $F_\Lambda$ to $\mathfrak q$ by letting $\pp^{+}$ act trivially. We claim that the induced supermodule
\begin{equation*}
M_{\qq}(\Lambda)
\coloneqq
\UE(\gg)\otimes_{\UE(\mathfrak q)}F_\Lambda
\end{equation*}
is simple. In general, $M_{\qq}(\Lambda)$ is a highest weight supermodule with simple
quotient $L(\Lambda)$~\cite{Carmeli_Fioresi_Varadarajan_HW}.

As $M_{\qq}(\Lambda)$ is a highest weight $\gg$-supermodule, any proper submodule of highest weight $\mu$ has the same infinitesimal character. In particular, $\mu+\rho\in W(\Lambda+\rho)$ and $\sum_{i}\mu_{i}\leq \sum_{i}\lambda_{i}$. Write
$\mu+\rho
=
\sum_{i=1}^n\varepsilon_i t_{\sigma(i)}\delta_i$ for $\sigma\in S_n$ and $\varepsilon_i\in\{1,-1\}.
$
Then
\begin{equation*}
0
\geq
\sum_i\mu_i-\sum_i\lambda_i
=
\sum_i\varepsilon_i t_{\sigma(i)}-\sum_i t_i
=
-2\sum_{\varepsilon_i=-1}t_{\sigma(i)}.
\end{equation*}
The right-hand side is strictly positive unless $\varepsilon_{i}=1$ for all $i$, that is, $\mu+\rho$ is just a permutation of $\Lambda+\rho$. Using $\Delta_{c}^{+}$-dominance, we see that the coordinates of $\mu+\rho$ are strictly
decreasing:
\begin{equation*}
(\mu_i+\rho_i)-(\mu_{i+1}+\rho_{i+1})
=
\mu_i-\mu_{i+1}+1>0,
\end{equation*}
which forces that the only valid permutation is the identity, that is, $\mu+\rho=\Lambda+\rho$.

Altogether, $M_{\qq}(\Lambda)\cong L(\Lambda)$, and its
supercharacter is given by~\eqref{eq::supercharacter_osp}.
The superdenominator identity
\cite[Appendix~A.1]{creutzig_garner2026} gives
\begin{equation*}
\Theta_\Lambda
=
\frac{\sum_{w\in W}\operatorname{sgn}(w)e^{w(\Lambda+\rho)}}
     {\sum_{w\in W}\operatorname{sgn}(w)e^{w\rho}}
=
\frac{\sum_{w\in W}\operatorname{sgn}(w)e^{w(\Lambda+\rho)}}
     {\mathcal R}.
\end{equation*}
Next, we evaluate this expression at the identity.

Choose $H,K\in\hh$ such that
\begin{equation*}
\delta_i(H)=n-i+\frac12,
\qquad
\delta_i(K)=t_i,
\qquad 1\leq i\leq n.
\end{equation*}
Then $(w(\Lambda+\rho))(H)=(w^{-1}\rho)(K)$, and the superdenominator identity gives
\begin{equation*}
\begin{aligned}
\sum_{w\in W}\operatorname{sgn}(w)
e^{(w(\Lambda+\rho))(H)}
=
\sum_{w\in W}\operatorname{sgn}(w)
e^{(w^{-1}\rho)(K)}
=\mathcal R(K).
\end{aligned}
\end{equation*}
We obtain 
\begin{equation*}
\Theta_{\Lambda}(0)=\lim_{s\to 0}\Theta_\Lambda(sH)=\lim_{s\to 0}\frac{\mathcal R(sK)}{\mathcal R(sH)}=
\prod_{i<j}
\frac{t_i^2-t_j^2}
     {(n-i+\frac12)^2-(n-j+\frac12)^2},
\end{equation*}
where we use the product expansion 
\begin{equation*}
\begin{aligned}
\mathcal R(H')&=
e^{\rho(H')}
\frac{\prod_{\alpha\in\Delta_{\bar0}^+}
      (1-e^{-\alpha(H')})}
     {\prod_{\beta\in\Delta_{\bar1}^+}
      (1-e^{-\beta(H')})}\\
&=
\prod_{1\leq i<j\leq n}
2\sinh\left(\tfrac{(\delta_i-\delta_j)(H')}{2}\right)
2\sinh\left(\tfrac{(\delta_i+\delta_j)(H')}{2}\right)
\prod_{i=1}^n
2\cosh\left(\tfrac{\delta_i(H')}{2}\right),
\end{aligned}
\end{equation*}
$H'\in \hh$, and l'Hôpital's rule.

Since $|\Delta_n^+|=n(n+1)/2$,
Theorem~\ref{thm::L_packet_and_sdim} yields
\begin{equation*}
\fsdim L(\Lambda)
=
(-1)^{n(n+1)/2}
\prod_{i<j}
\frac{t_i^2-t_j^2}
     {(n-i+\frac12)^2-(n-j+\frac12)^2}.
\end{equation*}
Every factor in the numerator is strictly negative because
$t_1>\cdots>t_n$ and all $t_i$ are negative. Every factor in the denominator is positive. Thus $\fsdim L(\Lambda)\neq0$. For $n=1$, the product is empty and
$\fsdim L(\Lambda)=-1$.
\end{example}

\subsubsection{Properties of the formal superdimension}

We now establish some basic properties of the formal superdimension. We first consider its behavior under parity reversal and finite direct sums. For the latter, we canonically extend the definition to finite direct sums of relative holomorphic discrete series supermodules.

\begin{lemma}
Let $\HH,\HH_1,\ldots,\HH_r$ be unitarizable $\gg$-supermodules admitting formal superdimensions. Then
\begin{equation*}
\fsdim(\Pi\HH)=-\fsdim(\HH),
\qquad
\fsdim\Bigl(\bigoplus_{j=1}^r\HH_j\Bigr)
=
\sum_{j=1}^r\fsdim(\HH_j).
\end{equation*}
\end{lemma}

Both identities follow directly from the definition. The formal superdimension is also multiplicative under tensoring with finite-dimensional supermodules, provided the tensor product admits a formal superdimension.

\begin{lemma}\label{lemm::sdim_TP}
Let $\HH=L(\Lambda)$ be a relative holomorphic discrete series
$\gg$-supermodule, and let $F$ be a finite-dimensional
$\gg$-supermodule. Assume that
$
(\HH\otimes F)\big|_{\even}
\cong
\bigoplus_\mu L_0(\mu)\otimes E_\mu
$
is a finite direct sum of relative holomorphic discrete series
$\even$-modules with finite-dimensional multiplicity superspaces.
Then
\[
\fsdim(\HH\otimes F)=\fsdim(\HH)\,\sdim(F),
\]
where $\sdim(F)=\dim F_{\bar0}-\dim F_{\bar1}$.
\end{lemma}

\begin{proof}
As above, set
$
a_{\mu}\coloneqq \dim(E_\mu)_{\bar0}-\dim(E_\mu)_{\bar1}
$
and define
\[
\Theta_{\HH\otimes F}
\coloneqq 
\frac{1}{|W_c|}
\sum_{w\in W}w\bigl(\sch (\HH\otimes F)\bigr).
\]
The assumed even decomposition gives
\[
\Theta_{\HH\otimes F}
=
\frac{
\sum_\mu a_{\mu}
\sum_{w\in W}(-1)^{\ell(w)}
e^{w(\mu+\rho_{\bar0})-\rho_{\bar0}}
}{
\prod_{\alpha\in\Delta_{\bar0}^{+}}
(1-e^{-\alpha})
}.
\]
The holomorphic extension and evaluation formula apply to each term.
Consequently,
\[
\Theta_{\HH\otimes F}(0)
=
(-1)^{|\Delta_n^+|}
\sum_\mu a_{\mu} d(\pi_\mu)
=
(-1)^{|\Delta_n^+|}\fsdim(\HH\otimes F).
\]

Since $F_{\bar0}$ and $F_{\bar1}$ are finite-dimensional
$\even$-modules, $\sch F$ is $W$-invariant. Thus
\[
\begin{aligned}
\Theta_{\HH\otimes F}
&=
\frac{1}{|W_c|}
\sum_{w\in W}
w\bigl(\sch \HH\,\sch F\bigr)=
\sch F\,
\frac{1}{|W_c|}
\sum_{w\in W}w\bigl(\sch \HH\bigr)=
\sch F\,\Theta_\Lambda.
\end{aligned}
\]
Evaluating the holomorphic extensions at $0$ yields the statement
\[
(-1)^{|\Delta_n^+|}\fsdim(\HH\otimes F)
=
\sdim(F)\,(-1)^{|\Delta_n^+|}\fsdim(\HH).\qedhere
\]
\end{proof}

\begin{example}
Let $\gg=\sl(2\vert1)$. We fix the real form $\su(1,1\vert0,1)$
and the standard positive system. We consider the relative holomorphic discrete series supermodule $\HH=L(-3,0\vert0)$ which decomposes under $\even$ as $\HH|_{\even}
\cong L_0(-3,0\vert0)\oplus\Pi L_0(-4,0\vert1)$. In particular, $\fsdim \HH=2-3=-1$, where we used
$
d(\pi_{(a,b\vert c)})=b-a-1.$ Let $F=\CC^{2\vert1}$ be the standard supermodule. Note that 
$F|_{\even}
\cong L_0(1,0\vert0)\oplus\Pi L_0(0,0\vert1)$ and $\sdim(F)=2-1=1$. 

We decompose $\HH\otimes F$ under $\even$:
\[
(\HH\otimes F)\vert_{\even}
\cong
L_0(-2,0\vert0)
\oplus L_0(-3,1\vert0)
\oplus L_0(-4,0\vert2)\oplus\bigl(\Pi L_0(-3,0\vert1)\bigr)^{\oplus2}
\oplus\Pi L_0(-4,1\vert1).
\]
Every summand of $(\HH\otimes F)\vert_{\even}$ is a relative holomorphic discrete series module, since every highest weight $(a,b\vert c)$ in this decomposition
satisfies $a-b<-1$. Thus
\[
\fsdim(\HH\otimes F)
=1+3+3-2\cdot2-4
=-1
=\fsdim(\HH)\,\sdim(F).
\]
\end{example}

Finally, we prove that the formal superdimension of a relative holomorphic discrete series supermodule $L(\Lambda)$ vanishes whenever the highest weight is not maximally atypical, that is,
\begin{equation}
k\coloneqq\at(\Lambda)<d\coloneqq\defect(\gg).
\end{equation}

Choose $d$ linearly independent, mutually orthogonal isotropic odd roots $\alpha_1,\ldots,\alpha_d$, and set
$
A\coloneqq\{\alpha_1,\ldots,\alpha_d\}
$
and
$
\hh_A\coloneqq\bigcap_{\alpha\in A}\ker\alpha.
$
For each $\alpha_i$, choose a nonzero root vector $u_{\alpha_i}\in\gg^{\alpha_i}$ and define
\begin{equation}
x_A\coloneqq
\sum_{j=1}^d\bigl(u_{\alpha_j}-i\omega(u_{\alpha_j})\bigr),
\qquad
c_A\coloneqq x_A^2
=-i\sum_{j=1}^d[u_{\alpha_j},\omega(u_{\alpha_j})]\in\hh.
\end{equation}

For each $w\in W$, set $B\coloneqq w(A)$ and define $x_B$ and $c_B$ analogously. Then $x_B\in\YY_\omega^{\hom}$ has rank $d$, satisfies $[\hh_B,x_B]=0$, and $c_B=x_B^2\in\hh$. By construction, the following lemma is immediate.

\begin{lemma}\label{lemm::properties_x_B}
The element $x_B$ belongs to $\mathcal Y_\omega^{\hom}$, has rank $d$, and satisfies $[\hh_B,x_B]=0$. Moreover, $c_B$ acts semisimply on $L(\Lambda)$.
\end{lemma}

By construction, $\at(\Lambda)<\operatorname{rk}(x_B)$, and Theorem~\ref{thm::generalized_DS_character} and Proposition~\ref{prop::form_DS_x_generalized} give 
$
0=\DS_{x_B}(L(\Lambda))
\cong L(\Lambda)^{c_B},
$ that is, zero is not an eigenvalue of $c_{B}$. The associated eigenspace decomposition of $L(\Lambda)$ with respect to $c_{B}$ is therefore
\begin{equation}
L(\Lambda)
=
\bigoplus_{\lambda\in\CC^\times}L(\Lambda)_{\lambda,c_{B}},
\qquad
L(\Lambda)_{\lambda,c_{B}}
\coloneqq
\{v\in L(\Lambda):c_Bv=\lambda v\}.
\end{equation}

For $\eta\in\hh_B^\ast$, define the restricted weight space
\begin{equation}
L(\Lambda)_{\eta}^{B}
\coloneqq
\bigoplus_{\substack{\nu\in\hh^\ast\\
\nu|_{\hh_B}=\eta}}L(\Lambda)^\nu.
\end{equation}

Our aim is to show that $\sch L(\Lambda)$ vanishes on $\hh_B$. Since $x_B$
commutes with $\hh_B$ and $L(\Lambda)^{c_B}=0$, we first show that $x_B$ is
an odd bijection identifying the even and odd parts of every restricted
weight space. To pass from this to the supercharacter, we then show that the
restricted weight spaces are finite-dimensional and that the supercharacter
converges absolutely on a nonempty open subset of $\hh_B$. On this subset,
the character series can be regrouped along restricted weight spaces. We
establish this in the more general form of
Lemma~\ref{lemm::restricted_character_with_kernel}, which expresses the
restriction of $\sch L(\Lambda)$ to $\hh_B$ through the supercharacter of
$\DS_{x_B}(L(\Lambda))$. That lemma holds without any assumption on the
atypicality.

\begin{lemma}\label{lemm::x_B_bijection}
The action of $x_B$ on $L(\Lambda)$ defines an odd bijection.
\end{lemma}

\begin{proof}
The operator $x_B$ preserves each eigenspace $L(\Lambda)_{\lambda, c_B}$ of $c_{B}$ since $x_{B}$ and $c_{B}$ commute by construction. For every eigenvalue $\lambda\neq0$, we have
\[
x_B^2|_{L(\Lambda)_{\lambda, c_B}}
=
\lambda\,\mathrm{id},
\qquad
x_B^{-1}|_{L(\Lambda)_{\lambda, c_B}}
=
\lambda^{-1}x_B.
\]
Since $L(\Lambda)$ is the direct sum of these eigenspaces, $x_B$ is bijective. It is odd by construction.
\end{proof}

\begin{lemma}\label{lemm::identify_restricted_weight_spaces}
For every $\eta\in\hh_B^\ast$, the operator $x_B$ preserves $L(\Lambda)_{\eta}^{B}$. Moreover, $x_{B}$ induces an isomorphism
\[
x_B:(L(\Lambda)_{\eta}^{B})_{\bar0}
\xrightarrow{\sim}
(L(\Lambda)_{\eta}^{B})_{\bar1}.
\]
\end{lemma}

\begin{proof}
For every $v\in L(\Lambda)_{\eta}^{B}$ and $H\in\hh_B$, we have
\[
H(x_Bv)=x_B(Hv)=\eta(H)x_Bv,
\]
where we use $[\hh_{B},x_{B}]=0$. This implies that $x_B$ preserves $L(\Lambda)_{\eta}^{B}$.

In Lemma~\ref{lemm::x_B_bijection}, we proved that $x_{B}$ is an odd bijection with inverse $x_{B}^{-1}$. In particular, $[\hh_B,x_{B}^{-1}]=0$, and also $x_{B}^{-1}$ preserves $L(\Lambda)^{B}_{\eta}$. Using that $x_{B}$ is odd, the statement follows.
\end{proof}

We next relate the restriction of $\sch L(\Lambda)$ to $\hh_B$
to the supercharacter of $\DS_{x_B}(L(\Lambda))$.
To justify regrouping the character expression, we need absolute
convergence on a nonempty domain in $\hh_B$ and finite-dimensional
restricted weight spaces.

\begin{lemma}\label{lemm::Gordan}
One has $\hh_B\cap i\hh_{\RR}\neq \{0\}$. Moreover, there exists
$
H_B\in\hh_B\cap i\hh_{\RR}
$ such that $\beta(H_B)>0$ for all
$\beta\in\Delta_n^+$.
\end{lemma}

\begin{proof} Write $\hh_B=(\hh_B\cap\hh_{\RR})\oplus(\hh_B\cap i\hh_{\RR})$, and recall that every root is real-valued on $i\hh_{\RR}$. By definition, the roots in $B$ are linearly independent, hence their restrictions to $i\hh_{\RR}$ are $\RR$-linearly independent. We conclude \[\dim_{\RR}(\hh_B\cap i\hh_{\RR})=\dim_{\CC}\hh-|B|\geq|B|>0,\] where $\dim_{\CC}\hh\geq2|B|$ follows as the linearly independent roots in $B$ span a totally isotropic subspace of $\hh^{\ast}$ for the nondegenerate form of Proposition~\ref{prop::structure_theory_BCLSA}(f). This proves $\hh_{B}\cap i\hh_{\RR}\neq \{0\}$.

The second statement is a direct consequence of Gordan's theorem of alternatives. Set $V_B\coloneqq\hh_B\cap i\hh_{\RR}$. By Gordan's theorem~\cite{Gordan} applied to the real functionals $\beta|_{V_B}$, $\beta\in\Delta_n^+$, either some $H_B\in V_B$ satisfies $\beta(H_B)>0$ for all $\beta$, or there are $t_\beta\geq0$, not all zero, with $\sum t_\beta\beta\big\vert_{V_B}=0$. The latter is equivalent to \[\sum_{\beta \in \Delta_{n}^{+}} t_\beta\beta\in\operatorname{span}_{\RR}B,\]
which we exclude. 

Recall the real forms in~\eqref{table::real_forms_classical}--\eqref{table::real_forms_exceptional}. The even real form $\gg_{\bar0,\RR}$ has at most one noncompact simple ideal, so that the invariant bilinear form is definite on $\operatorname{span}_{\RR}\Delta_n^+$. However, $\operatorname{span}_{\RR}B$ is totally isotropic by definition of $B$, so that 
\[
\operatorname{span}_{\RR}\Delta_n^+
\cap\operatorname{span}_{\RR}B
=
\{0\}.
\]
We conclude that $t_{\beta}=0$ for all $\beta \in \Delta_{n}^{+}$, which is not allowed. 
\end{proof}

\begin{lemma}
\label{lemm::restricted_character_with_kernel}
Let $\pi_B:\hh_B\to\hh_{x_B}$ be the canonical quotient map. Then, as meromorphic functions on $\hh_B$,
\[
\sch L(\Lambda)|_{\hh_B}
=
\sch\bigl(L(\Lambda)^{c_B}\bigr)\circ\pi_B
=
\sch\bigl(\DS_{x_B}(L(\Lambda))\bigr)\circ\pi_B.
\]
\end{lemma}

\begin{proof} By Lemma~\ref{lemm::Gordan}, 
\[
U_B\coloneqq
\{H\in\hh_B:\operatorname{Re}\beta(H)>0
\text{ for all }\beta\in\Delta_n^+\}\neq \varnothing.
\]
The supercharacter $\sch L(\Lambda)$ converges absolutely on $U_B$. Each simple $\even$-constituent has character
\[
\ch L_0(\mu)
=
\frac{\operatorname{ch}F_\mu}
{\prod_{\beta\in\Delta_n^+}(1-e^{-\beta})},
\] where $F_{\mu}$ is a finite-dimensional $\kk$-module of highest weight $\mu$. For each $\beta \in \Delta_{n}^{+}$, the geometric series \[
\frac{1}{1-e^{-\beta(H)}}
=\sum_{k=0}^{\infty}e^{-k\beta(H)}
\] 
converges absolutely on $U_{B}$ as $|e^{-\beta(H)}|=\vert e^{-\operatorname{Re}(\beta(H))}\vert<1$. This proves that $\sch L(\Lambda)$ converges absolutely on $U_{B}$.

Next, we show that the restricted weight spaces
$
L(\Lambda)_{\eta}^{B}$
are finite-dimensional. To see this, recall that $L(\Lambda)$ decomposes into finitely many simple relative holomorphic discrete series $\even$-modules $L_0(\nu_j)$. Every weight of $L(\Lambda)$ is therefore of the form
\[
\mu=\xi-\sum_{\beta\in\Delta_n^+}n_\beta\beta,
\qquad n_\beta\in\ZZ_{\geq0},
\]
where $\xi$ is a weight of a finite-dimensional $\kk$-module $\UE(\kk)v_{\nu_j},
$
where $v_{\nu_{j}}$ is the highest weight vector of $L_{0}(\nu_{j})$. If $\mu|_{\hh_B}=\eta$, evaluation at the element $H_B\in U_{B}$ of Lemma~\ref{lemm::Gordan}, for which every $\beta(H_B)$ is real, gives
\[
\sum_{\beta\in\Delta_n^+}n_\beta\beta(H_B)
=
\xi(H_B)-\eta(H_B).
\]
Since every $\beta(H_B)>0$ and $\xi$ ranges over a finite set, only finitely many tuples $(n_\beta)$ are possible. Hence only finitely many weights $\mu$ restrict to $\eta$. As the ordinary weight spaces are finite-dimensional, so is $L(\Lambda)_{\eta}^{B}$.

We can now prove the statement of the lemma. Recall that $c_B$ acts semisimply and commutes with $x_B$ and $\hh_B$. On each nonzero $c_B$-eigenspace with eigenvalue $t$, the odd operator $x_B$ is invertible, with inverse $t^{-1}x_B$. We conclude that
\[
\sdim L(\Lambda)_{\eta}^{B}
=
\sdim\bigl(L(\Lambda)^{c_B}\cap L(\Lambda)_{\eta}^{B}\bigr).
\]

By Proposition~\ref{prop::DS_composition} and the proof of Lemma~\ref{lemm::DS_x_simple_if_x_is_of_rank_1}, each $[\gg^\beta,\gg^{-\beta}]$, for $\beta\in B$, acts trivially on $L(\Lambda)^{c_B}$, and therefore the $\hh_B$-action factors through $\pi_B$. We conclude
\[
\sch L(\Lambda)|_{\hh_B}
=
\sch\bigl(L(\Lambda)^{c_B}\bigr)\circ\pi_B.
\]
The meromorphic supercharacter of $L(\Lambda)$ restricts to $\hh_B$, since no even root belongs to the totally isotropic span of $B$. This identity extends meromorphically from $U_B$. Finally, Proposition~\ref{prop::form_DS_x_generalized} gives
$
L(\Lambda)^{c_B}\cong\DS_{x_B}(L(\Lambda)),
$
which proves the assertion.
\end{proof}

We can now prove the vanishing of the formal superdimension below maximal atypicality.

\begin{theorem}\label{thm::vanishing_sdim_less_atypic}
Let $L(\Lambda)$ be a relative holomorphic discrete series $\gg$-supermodule. If
$
\at(\Lambda)<\defect(\gg),
$
then
\[
\fsdim L(\Lambda)=0.
\]
\end{theorem}

\begin{proof} By assumption and Theorem~\ref{thm::generalized_DS_character}, 
$\DS_{x_{B}}(L(\Lambda))=0$ for any $B=w^{-1}A$, $w \in W$. Using Lemma~\ref{lemm::restricted_character_with_kernel}, this implies
\[
w\bigl(\sch L(\Lambda)\bigr)\big|_{\hh_A}
=
\bigl(\sch L(\Lambda)|_{\hh_{w^{-1}A}}\bigr)\circ w^{-1}|_{\hh_A}=
\sch\bigl(\DS_{x_{w^{-1}A}}(L(\Lambda))\bigr)\circ\pi_{w^{-1}A}\circ w^{-1}|_{\hh_A}
=0,
\]
and hence $\Theta_{\Lambda}\vert_{\hh_{A}}=0$. In particular, $\Theta_\Lambda(0)=0$, and the statement follows by Theorem~\ref{thm::L_packet_and_sdim}.
\end{proof}

In the rest of this manuscript, we prove the converse, namely that maximal atypicality implies nonvanishing of the formal superdimension. Together with the preceding result, this establishes an analogue of the Kac--Wakimoto theorem for the formal superdimension.

We begin with the case of integral highest weights. The idea is to associate a finite-dimensional $\gg$-supermodule to a relative holomorphic discrete series supermodule $L(\Lambda)$ whose supercharacter agrees with the $L$-packet supercharacter $\Theta_\Lambda$. This allows us to reduce the nonvanishing problem for the formal superdimension to the finite-dimensional case.

To this end, we recall the Zuckerman functor and the properties of its derived functors needed below.

\subsection{Intermezzo: Zuckerman functors}
\label{subsec::Zuckerman_functors} 
We briefly recall the construction and the properties of Zuckerman functors and their derived analogues. We follow the construction in~\cite[Section~2]{coulembier2016bott}.

Let $\aa\subseteq\even$ be a reductive subalgebra. A $\gg$-supermodule $M$ is \emph{locally finite over $\aa$} if $\dim_{\CC}\UE(\aa)m<\infty$ for every $m\in M$, and \emph{semisimple over $\aa$} if its restriction to $\aa$ is a direct sum of simple $\aa$-modules.

We denote by $\calC(\gg,\aa)$ the category of $\gg$-supermodules that are locally finite and semisimple over $\aa$. We use the analogous notation $\calC(\even,\aa)$ for ordinary $\even$-modules.

The \emph{Zuckerman functor} is the functor
$\Gamma:\calC(\gg,\kk)\to\calC(\gg,\even)$ defined by
\begin{equation}
\Gamma(M)\coloneqq
\{m\in M:\dim_{\CC}\UE(\even)m<\infty\}.
\end{equation}

\begin{lemma}
Let $M\in\calC(\gg,\kk)$. Then $\Gamma(M)$ is a $\gg$-submodule and it is semisimple over $\even$. Moreover, $\Gamma$ is right adjoint to the inclusion $\calC(\gg,\even)\hookrightarrow\calC(\gg,\kk)$ and is therefore left exact.
\end{lemma}

\begin{proof}
First, $\even$ is reductive, so its center lies in $\hh$. The Cartan subalgebra $\hh$ is compactly embedded, that is, $\hh \subset \kk$, and therefore the center of $\even$ acts semisimply on $\Gamma(M)$. This implies that $\Gamma(M)$ is a semisimple $\even$-module.

Next, we show that $\Gamma(M)$ is a $\gg$-supermodule. Let $x\in\gg$ and $m\in\Gamma(M)$. Set $V=\UE(\even)m$ and $W=\ad_{\UE(\even)}x\subseteq\gg$. By assumption, both spaces are finite-dimensional. Moreover, $WV$ is stable under $\even$, since
\[
a(wv)=[a,w]v+w(av)
\]
for all $a\in\even$, $w\in W$, and $v\in V$. Using $xm\in WV$, we obtain $\UE(\even)xm\subseteq WV$, and hence $xm\in\Gamma(M)$.

If $N\in\calC(\gg,\even)$, every homomorphism $N\to M$ has image in $\Gamma(M)$. Therefore
\[
\Hom_{\calC(\gg,\kk)}(\iota N,M)
\cong
\Hom_{\calC(\gg,\even)}(N,\Gamma(M)).
\]
This proves the adjunction and hence left exactness. Compare~\cite[Proposition~1.2]{milicic1998equivariant}.
\end{proof}

The functor $\Gamma$ is not exact in general. We use its right derived functors, defined in the categories above, without imposing finite generation. The required homological framework is developed in~\cite[Section~1]{milicic1998equivariant}, \cite[Chapter~5]{huang2007dirac}, and~\cite[Section~2]{coulembier2016bott}.

Let $M\in\calC(\gg,\kk)$ and choose an injective resolution
$
0\longrightarrow M
\longrightarrow I^0
\longrightarrow I^1
\longrightarrow\cdots.
$
The $j$-th \emph{derived Zuckerman functor} is defined by
\begin{equation}
R^j\Gamma(M)
\coloneqq
H^j\bigl(\Gamma(I^\bullet)\bigr),
\qquad j\geq0.
\end{equation}
The supermodule $R^{j}\Gamma(M)$
is independent of the chosen injective resolution up to natural isomorphism, and $R^0\Gamma=\Gamma$. 

We denote by
$
\Gamma_0:
\calC(\even,\kk)
\longrightarrow
\calC(\even,\even)
$
the corresponding Zuckerman functor for ordinary $\even$-modules, and its derived versions by $R^{j}\Gamma_{0}$. 
\begin{proposition}
\label{prop::Zuckerman_compatibilities}
The derived Zuckerman functors satisfy the following properties.
\begin{enumerate}
\item[(a)] For every $M\in\calC(\gg,\kk)$ and $j\geq0$,
\[
\Res_{\even}R^j\Gamma(M)
\cong
R^j\Gamma_0\bigl(\Res_{\even}M\bigr),
\]
with the grading preserved.

\item[(b)] For every finite-dimensional semisimple $\even$-module $E$, every $N\in\calC(\even,\kk)$, and every $j\geq0$,
\[
R^j\Gamma_0(E\otimes N)
\cong
E\otimes R^j\Gamma_0(N).
\]

\item[(c)] For every $M\in\calC(\gg,\kk)$ and $j\geq0$,
\[
R^j\Gamma(\Pi M)
\cong
\Pi R^j\Gamma(M).
\]

\item[(d)] If $M$ has infinitesimal character $\chi$, then every nonzero $R^j\Gamma(M)$ has the same infinitesimal character.
\end{enumerate}
\end{proposition}

\begin{proof}
Part~(a) follows from~\cite[Lemma~4.10]{coulembier2016bott}. Part~(b) is the ordinary tensor identity used in the proof of~\cite[Lemma~4.12]{coulembier2016bott}.

For part~(c), parity reversal is an exact equivalence and preserves injective objects. The identity $\Gamma\Pi=\Pi\Gamma$ therefore extends to the derived functors.

For part~(d), let $z\in\mathfrak Z(\gg)$, the center of $\UE(\gg)$. Since $z$ acts on $M$ by $\chi(z)\id_M$ and its action commutes with $\Gamma$, functoriality gives
\[
z_{R^j\Gamma(M)}
=
\chi(z)\id_{R^j\Gamma(M)}.
\]
This proves the assertion. Compare~\cite[Proposition~6.2.1]{huang2007dirac}.
\end{proof}

\subsection{Kac--Wakimoto conjecture for the formal superdimension: the integral case}\label{subsec::KW_conjecture_integral_case}

We now prove that maximal atypicality implies nonvanishing of the formal superdimension if the highest weight is integral.

Recall that we work with the fixed $\kk$-adapted positive system $\Delta_{\bar0}^+$ and the corresponding holomorphic chamber $C_{\hol}$ introduced in Section~\ref{subsec::HC_parameterization}. A highest weight $\mu$ of a relative holomorphic discrete series $\even$-module $L_0(\mu)$ is called \emph{integral} if
\begin{equation}
2\frac{(\mu,\alpha)}{(\alpha,\alpha)}\in\ZZ
\qquad
\text{for every }\alpha\in\Delta_{\bar0}.
\end{equation}

Our strategy is to associate to an integral relative holomorphic discrete series supermodule $L(\Lambda)$ a finite-dimensional simple $\gg$-supermodule whose supercharacter agrees with the $L$-packet supercharacter $\Theta_\Lambda$. By Theorem~\ref{thm::L_packet_and_sdim}, this will allow us to reduce the nonvanishing problem for the formal superdimension to the Kac--Wakimoto theorem for finite-dimensional supermodules, proved in~\cite{serganova2011superdimension}.

We first establish the required construction by applying derived Zuckerman functors to integral relative holomorphic discrete series $\even$-modules. We need the following standard result.

\begin{lemma}\label{lemm::length_holomorphic_Weyl_element}
 Let $w_{\hol}\in W$
be the unique Weyl group element carrying the
holomorphic chamber $C_{\hol}$ to the dominant
chamber determined by $\Delta_{\bar0}^+$. Then
\[
\ell(w_{\hol})=|\Delta_{n}^{+}|.
\]
\end{lemma}

\begin{proof}
The positive root system associated with
$C_{\hol}$ is
$
\Delta_c^+\cup(-\Delta_n^+).
$
Since $w_{\hol}$ carries $C_{\hol}$ to the
dominant chamber, it carries this positive
system to $\Delta_{\bar0}^+$. Consequently,
\[
w_{\hol}(\Delta_c^+)\subseteq\Delta_{\bar0}^+,
\qquad
w_{\hol}(\Delta_n^+)\subseteq-\Delta_{\bar0}^+.
\]
Thus the positive roots mapped to negative roots
by $w_{\hol}$ are precisely the elements of
$\Delta_n^+$. This yields
\[
\ell(w_{\hol})
=
\bigl|\{\alpha\in\Delta_{\bar0}^+:
w_{\hol}\alpha\in-\Delta_{\bar0}^+\}\bigr|
=
|\Delta_n^+|.\qedhere
\]
\end{proof}

\begin{proposition}
\label{prop::Zuckerman_even_discrete_series}
Let $L_0(\mu)$ be a relative holomorphic discrete series
$\even$-module with integral highest weight $\mu$.
Set
$
\mu^\sharp
\coloneqq
w_{\hol}(\mu+\rho_{\bar0})-\rho_{\bar0},
$
where $w_{\hol}$ carries $C_{\hol}$ to the dominant
chamber determined by $\Delta_{\bar0}^+$.
Then:
\begin{enumerate}
    \item[(a)]
    $\mu^\sharp$ is dominant integral with respect to
    $\Delta_{\bar0}^+$.

    \item[(b)]
    There is an isomorphism of $\even$-modules
    \[
    R^{\vert \Delta_{n}^{+}\vert} \Gamma_0\bigl(L_0(\mu)\bigr)
    \cong
    L_0(\mu^\sharp),
    \]
    where $L_0(\mu^\sharp)$ is the finite-dimensional
    simple $\even$-module of highest weight $\mu^\sharp$.
\end{enumerate}
\end{proposition}

\begin{proof}
For part~(a), the Harish-Chandra condition of
Theorem~\ref{thm::HC_condition} gives for every $\beta \in \Delta_{n}^{+}$
\[
2\frac{(\mu+\rho_{\bar0},\beta)}{(\beta,\beta)}<0.
\]
Moreover, $\mu$ is dominant integral with respect to
$\Delta_c^+$. Since
$\rho_{\bar0}=\rho_c+\rho_n$ and $\rho_n$ is
$W_c$-invariant, we have for every $\alpha\in\Delta_c^+$
\[
2\frac{(\mu+\rho_{\bar0},\alpha)}{(\alpha,\alpha)}
=
2\frac{(\mu+\rho_c,\alpha)}{(\alpha,\alpha)}
>0.
\]
Thus $\mu+\rho_{\bar0}\in C_{\hol}$, and
$
\lambda^\sharp
\coloneqq
w_{\hol}(\mu+\rho_{\bar0})
$
is regular dominant. In addition, since $\mu$ and $\rho_{\bar{0}}$ are integral, $\lambda^{\sharp}$ is integral, and therefore dominant integral: 
\[
2\frac{(\lambda^\sharp,\alpha)}{(\alpha,\alpha)}
\in\ZZ_{\geq1},
\]
where $\alpha\in \Delta_{\bar{0}}^{+}$ is simple. It follows that $\mu^{\sharp}$ is dominant integral using 
$
2\tfrac{(\rho_{\bar0},\alpha)}{(\alpha,\alpha)}=1
$
for every simple root $\alpha\in \Delta_{\bar{0}}^{+}$. Indeed, by definition of $\lambda^{\sharp}$, one has for every simple root $\alpha\in \Delta_{\bar{0}}^{+}$:
\[
2\frac{(\mu^\sharp,\alpha)}{(\alpha,\alpha)}
=
2\frac{(\lambda^\sharp,\alpha)}{(\alpha,\alpha)}-1
\in\ZZ_{\geq0}.
\]
This proves part~(a).

For part~(b), we consider the following parabolic subalgebra in $\even$:
\[
\qq=\kk\oplus\uu,
\qquad
\uu\coloneqq
\bigoplus_{\beta\in\Delta_n^+}(\even)^\beta.
\]
Let $v_\mu$ be a highest weight vector of $L_0(\mu)$.
The $\kk$-module $F_{\mu}\coloneqq\UE(\kk)v_\mu
$ is finite-dimensional and simple, with highest weight
$\mu$. Since $\uu v_\mu=0$ and
$[\kk,\uu]\subseteq\uu$, the nilradical $\uu$
annihilates $F_{\mu}$. Hence the inclusion
$F_{\mu}\hookrightarrow L_0(\mu)$ induces a surjective
homomorphism
\[
M_{\qq}(\mu)
\coloneqq
\UE(\even)\otimes_{\UE(\qq)}F_{\mu}
\longrightarrow
L_0(\mu).
\]

Using the Weyl character formula for $F_{\mu}$, the character of $M_{\qq}(\mu)$ is
\[
\begin{aligned}
\fch M_{\qq}(\mu)
=
\frac{\fch F_{\mu}}
{\prod_{\beta\in\Delta_n^+}
(1-e^{-\beta})}=
\frac{
\sum_{w\in W_c}
(-1)^{\ell_c(w)}
e^{w(\mu+\rho_c)-\rho_c}
}{
\prod_{\alpha\in\Delta_{\bar0}^+}
(1-e^{-\alpha})
}=
\frac{
\sum_{w\in W_c}
(-1)^{\ell_c(w)}
e^{w(\mu+\rho_{\bar0})-\rho_{\bar0}}
}{
\prod_{\alpha\in\Delta_{\bar0}^+}
(1-e^{-\alpha})
},
\end{aligned}
\]
where we used in the last equality that $\rho_n=\rho_{\bar0}-\rho_c$ is $W_{c}$-invariant.
By Lemma~\ref{lemm::character_formula_L_0_relative}, this is precisely
$\fch L_0(\mu)$, and we conclude
$
M_{\qq}(\mu)\cong L_0(\mu).
$

By~\cite[Fact~2.22 and Remark~2.23]
{kitagawa2025construction} (see also
\cite[Corollary~4.160]{KnappVogan1995}), the inequalities
\[
2\frac{(\mu+\rho_{\bar0},\beta)}{(\beta,\beta)}<0,
\qquad \beta\in\Delta_n^+,
\]
imply that
$
R^{\dim\uu}\Gamma_0\bigl(M_{\qq}(\mu)\bigr)
$
is a nonzero finite-dimensional simple
$\even$-module.
Here $\dim\uu=\vert \Delta_{n}^{+}\vert$, which also equals
$\ell(w_{\hol})$ by
Lemma~\ref{lemm::length_holomorphic_Weyl_element}.

It remains to identify its highest weight.
Derived Zuckerman functors preserve infinitesimal
characters. Thus, if $\nu$ denotes the dominant
integral highest weight of this finite-dimensional
simple module, then
$
\nu+\rho_{\bar0}
\in W(\mu+\rho_{\bar0}).
$ However,
both $\nu+\rho_{\bar0}$ and
$w_{\hol}(\mu+\rho_{\bar0})$ lie in the dominant
chamber, and uniqueness of the dominant representative
of a regular Weyl group orbit gives
\[
\nu+\rho_{\bar0}
=
w_{\hol}(\mu+\rho_{\bar0}).
\]
Hence $\nu=\mu^\sharp$, which proves (b).
\end{proof}

For a relative holomorphic discrete series $\gg$-supermodule $L(\Lambda)$ with integral highest weight, we can now construct a finite-dimensional $\gg$-supermodule $F$ whose supercharacter equals the $L$-packet supercharacter $\Theta_\Lambda$.

\begin{corollary}
\label{cor::Zuckerman_packet_character}
Let $L(\Lambda)$ be a relative holomorphic discrete series $\gg$-supermodule with integral highest weight. Assume $
L(\Lambda)|_{\even}
\cong
\bigoplus_\mu L_0(\mu)\otimes E_\mu
$, and set
$
\mu^\sharp
\coloneqq
w_{\hol}(\mu+\rho_{\bar0})-\rho_{\bar0}.
$
Then the following holds:
\begin{enumerate}
    \item[(a)] $R^{|\Delta_n^+|}\Gamma(L(\Lambda))$ is nonzero and finite-dimensional.
    \item[(b)] $R^{|\Delta_n^+|}\Gamma(L(\Lambda))|_{\even}
\cong
\bigoplus_\mu L_0(\mu^\sharp)\otimes E_\mu.$
\item[(c)] $\sch\bigl(R^{|\Delta_n^+|}\Gamma(L(\Lambda))\bigr)
=
(-1)^{|\Delta_n^+|}\Theta_\Lambda.$
\end{enumerate}
\end{corollary}

\begin{proof}
(a) and (b) follow immediately by Propositions~\ref{prop::Zuckerman_compatibilities}(a) and~\ref{prop::Zuckerman_even_discrete_series}. Concretely, any $\mu$ is integral as $\Lambda$ is integral. Applying $R^{\vert \Delta_{n}^{+}\vert}\Gamma$ to the decomposition of $L(\Lambda)$, we obtain
\[
\begin{aligned}
R^{|\Delta_n^+|}\Gamma(L(\Lambda))|_{\even}
\cong
R^{|\Delta_n^+|}\Gamma_0
\big(
\bigoplus_\mu L_0(\mu)\otimes E_\mu
\big)\cong
\bigoplus_\mu L_0(\mu^\sharp)\otimes E_\mu,
\end{aligned}
\]
where we used additivity of the derived functor and that it commutes with parity reversal. This proves (b). Since the sum is finite and at least one $E_\mu$ is nonzero, $R^{|\Delta_n^+|}\Gamma(L(\Lambda))$ is nonzero and finite-dimensional. This proves (a).

It remains to prove (c). Recall the definition of the $L$-packet supercharacter,
\[
\Theta_\Lambda
=
\frac{
\sum_\mu \sdim(E_\mu)
\sum_{w\in W}(-1)^{\ell(w)}
e^{w(\mu+\rho_{\bar0})-\rho_{\bar0}}
}{
\prod_{\alpha\in\Delta_{\bar0}^+}
(1-e^{-\alpha})
}.
\]
Since
$
\mu^\sharp+\rho_{\bar0}
=
w_{\hol}(\mu+\rho_{\bar0}),
$
reindexing the Weyl sum gives
\[
\begin{aligned}
\sum_{w\in W}(-1)^{\ell(w)}
e^{w(\mu^\sharp+\rho_{\bar0})-\rho_{\bar0}}=
(-1)^{\ell(w_{\hol})}
\sum_{w\in W}(-1)^{\ell(w)}
e^{w(\mu+\rho_{\bar0})-\rho_{\bar0}}.
\end{aligned}
\]
By Lemma~\ref{lemm::length_holomorphic_Weyl_element}, $\ell(w_{\hol})=|\Delta_n^+|$. The Weyl character formula therefore yields
\[
\begin{aligned}
\Theta_\Lambda
&=
(-1)^{|\Delta_n^+|}
\sum_\mu \sdim(E_\mu)\,\fch L_0(\mu^\sharp)=
(-1)^{|\Delta_n^+|}
\sch\bigl(R^{|\Delta_n^+|}\Gamma(L(\Lambda))\bigr).
\end{aligned}\qedhere
\]
\end{proof}

Since parity reversal changes the sign of the supercharacter, we define
$
F\coloneqq\Pi^{|\Delta_n^+|}R^{|\Delta_n^+|}\Gamma(L(\Lambda)).
$
By the preceding corollary, $F$ is a nonzero finite-dimensional $\gg$-supermodule satisfying
\begin{equation}
\sch(F)=\Theta_\Lambda.
\end{equation}
Moreover, by Proposition~\ref{prop::Zuckerman_compatibilities}(d), $F$ has the same infinitesimal character as $L(\Lambda)$. 

We now prove that the finite-dimensional supermodule $F$ is simple. We begin with a root-theoretic estimate.

\begin{lemma}
\label{lemm::odd_root_estimate}
For every real form in~\eqref{table::real_forms_classical} or~\eqref{table::real_forms_exceptional}, we have for all $\gamma \in \Delta_{\bar{1}}$ and $\beta \in \Delta_{n}^{+}$:
\[
\left|2\frac{(\gamma,\beta)}{(\beta,\beta)}\right|
\leq1.
\]
\end{lemma}

\begin{proof}
By the classification of the relevant real forms in~\eqref{table::real_forms_classical}--\eqref{table::real_forms_exceptional}, each noncompact even root belongs to a single noncompact simple factor of $\gg_{\bar0,\RR}$. The restriction of $\odd$ to its complexification is a direct sum of copies of the following representations: the natural representation or its dual for $\sl(r)$, the natural representation for $\sp(2n)$, the vector representation for $\so(2r)$, the spin representation for $\so(7)$, or the natural representation for $\sl(2)$.

All these representations are minuscule~\cite[Chapter~VIII, Section~7.3]{Bourbaki_minuscule}, meaning that their weights form a single Weyl group orbit. Hence their weights $\gamma$ satisfy
$
2(\gamma,\beta)/(\beta,\beta)\in\{-1,0,1\}
$
for every root $\beta$ of the corresponding simple factor.
\end{proof}

Recall the parabolic subalgebra $\qq$ and the finite-dimensional simple $\kk$-module $F_{\mu}$ from the proof of Proposition~\ref{prop::Zuckerman_even_discrete_series}. We first establish a semisimplicity property.

\begin{lemma}\label{lemm::odd_tensor_semisimple}
Let $L_0(\mu)$ be an integral relative holomorphic discrete series module. Then
\[
\odd\otimes L_0(\mu)\cong\bigoplus_{\nu}M_{\qq}(\nu), \qquad M_{\qq}(\nu)=\UE(\even)\otimes_{\UE(\qq)}F_\nu,
\]
where $\nu=\mu+\gamma$ with $\gamma \in \Delta_{\bar{1}}$, each $M_{\qq}(\nu)\cong L_0(\nu)$ is simple, and
\[
2\frac{(\nu+\rho_{\bar0},\beta)}{(\beta,\beta)}\leq0
\]
for every $\beta\in\Delta_n^+$.
\end{lemma}

\begin{proof}
By Proposition~\ref{prop::Zuckerman_even_discrete_series},
$
L_0(\mu)\cong M_{\qq}(\mu)=\UE(\even)\otimes_{\UE(\qq)}F_{\mu},
$ and the tensor identity yields
\[
\odd\otimes L_0(\mu)
\cong
\UE(\even)\otimes_{\UE(\qq)}
(\odd\otimes F_{\mu}).
\]
We consider the finite-dimensional $\qq$-module $\odd\otimes F_{\mu}$. It has a finite filtration by simple $\qq$-modules on which $\uu$ acts trivially, since $\odd\otimes F_{\mu}$ is finite-dimensional, $\uu$ acts nilpotently, and $\kk$ acts semisimply. Exactness of induction therefore gives a finite filtration of $\odd\otimes L_0(\mu)$ by generalized Verma modules $M_{\qq}(\nu)$ with $\nu=\mu+\gamma$, where $\gamma \in \Delta_{\bar{1}}$.

As in the proof of Proposition~\ref{prop::Zuckerman_even_discrete_series}, it follows that the generalized Verma modules $M_{\qq}(\nu)$ are simple since all $\nu$ are integral and 
\[
2\frac{(\nu+\rho_{\bar0},\beta)}{(\beta,\beta)}
=
2\frac{(\mu+\rho_{\bar0},\beta)}{(\beta,\beta)}
+
2\frac{(\gamma,\beta)}{(\beta,\beta)}
\leq0.
\]
Here we use that $\mu$ satisfies the Harish-Chandra condition and Lemma~\ref{lemm::odd_root_estimate}.

It remains to prove that the filtration splits. Each shifted highest weight
$\nu+\rho_{\bar0}$ lies in the closure of the holomorphic Weyl chamber.
Thus distinct highest weights have distinct infinitesimal characters.
We conclude that the filtration splits, using that simple highest weight modules have no self-extensions~\cite[Proposition~3.1(d)]{humphreys2008representations}.
\end{proof}

\begin{lemma}
\label{lemm::Zuckerman_detects_odd_maps}
Let $\HH_{1}$ and $\HH_{2}$ be finite direct sums of relative holomorphic discrete series $\even$-modules with integral highest weights. Let 
$
f:\odd\otimes \HH_{1}\to \HH_{2}
$ be a homomorphism of $\even$-modules. If $R^{|\Delta_n^+|}\Gamma_0(f)=0$, then $f=0$. 
\end{lemma}

\begin{proof}
By assumption and Lemma~\ref{lemm::odd_tensor_semisimple}, both $\odd\otimes\HH_{1}$ and $\HH_{2}$ are semisimple. Suppose that $f\neq0$. Then there exists a simple direct summand $V$ of $\odd\otimes\HH_{1}$ such that $f\vert_V\colon V\to f(V)$ is an isomorphism. By Proposition~\ref{prop::Zuckerman_even_discrete_series}, we have
\[
R^{|\Delta_n^+|}\Gamma_0(V)\neq0.
\]
Functoriality implies that $R^{|\Delta_n^+|}\Gamma_0(f)\neq0$, a contradiction.
\end{proof}

We can now prove simplicity of $F=\Pi^{|\Delta_n^+|}R^{|\Delta_n^+|}\Gamma(L(\Lambda))$.

\begin{proposition}
\label{prop::Zuckerman_image_simple}
Let $L(\Lambda)$ be a relative holomorphic discrete series $\gg$-supermodule with integral highest weight. Then the finite-dimensional $\gg$-supermodule
$
F=\Pi^{|\Delta_n^+|}R^{|\Delta_n^+|}\Gamma(L(\Lambda))
$
is simple.
\end{proposition}

\begin{proof}
It suffices to prove that $R^{|\Delta_n^+|}\Gamma(L(\Lambda))$ is simple. Assume 
$
L(\Lambda)|_{\even}
\cong
\bigoplus_\mu L_0(\mu)\otimes E_\mu.
$
Then, by Corollary~\ref{cor::Zuckerman_packet_character},
\[
R^{|\Delta_n^+|}\Gamma(L(\Lambda))|_{\even}
\cong
\bigoplus_\mu L_0(\mu^\sharp)\otimes E_\mu,
\]
where the modules $L_0(\mu^\sharp)$ are pairwise nonisomorphic.

Let $U$ be a $\gg$-subsupermodule of $R^{|\Delta_n^+|}\Gamma(L(\Lambda))$. Then there exist graded subspaces $N_\mu\subseteq E_\mu$ such that
\[
U|_{\even}
\cong
\bigoplus_\mu L_0(\mu^\sharp)\otimes N_\mu.
\]
The corresponding graded $\even$-submodule in $L(\Lambda)$ is
\[
L \coloneqq
\bigoplus_\mu L_0(\mu)\otimes N_\mu
\subseteq L(\Lambda),
\]
and by additivity and Proposition~\ref{prop::Zuckerman_even_discrete_series}, the derived Zuckerman functor identifies $L$ with $U$.

We show that $L$ is a $\gg$-subsupermodule of $L(\Lambda)$ and hence $L=L(\Lambda)$ or $L=\{0\}$. Under the identification via derived Zuckerman functors, this implies that $R^{\vert \Delta_{n}^{+}\vert}\Gamma L(\Lambda)$ is simple. 

 Consider the $\even$-module homomorphism
\[
f:\odd\otimes L\longrightarrow L(\Lambda)/L,
\qquad x\otimes v\longmapsto xv+L.
\]
Applying $R^{|\Delta_n^+|}\Gamma_0$ to $f$ and using Proposition~\ref{prop::Zuckerman_compatibilities}, we obtain
\[
R^{|\Delta_n^+|}\Gamma_0(f):
\odd\otimes U\longrightarrow R^{|\Delta_n^+|}\Gamma(L(\Lambda))/U, \qquad x\otimes u\longmapsto xu+U.
\]
However, $U$ is a $\gg$-supermodule such that $xu\in U$ for all $x\in \odd$ and $u\in U$. We conclude that $R^{\vert \Delta_{n}^{+}\vert}\Gamma_{0}(f)=0$, and by Lemma~\ref{lemm::Zuckerman_detects_odd_maps}, $f=0$. This proves that $L$ is a $\gg$-supermodule.
\end{proof}

\begin{remark}
The highest weight $\Lambda_F$ of the simple $\gg$-supermodule
$F=\Pi^{|\Delta_n^+|}R^{|\Delta_n^+|}\Gamma(L(\Lambda))$
is the greatest element of
\[
\{
\mu^\sharp = w_{\hol}(\mu+\rho_{\bar0})-\rho_{\bar0}:E_\mu\neq0
\}
\]
with respect to the highest weight order determined by $\Delta^+$. Indeed, by Corollary~\ref{cor::Zuckerman_packet_character}, $F|_{\even}$ is a direct sum of the simple $\even$-modules $L_0(\mu^\sharp)$. Hence $\Lambda_F=\mu^\sharp$ for some $\mu$ with $E_\mu\neq0$. Since $F$ is generated by a highest weight vector, every weight of $F$, and thus every $\mu^\sharp$, is less than or equal to $\Lambda_F$.
\end{remark}

\begin{theorem}\label{thm::formal_KW_integral}
Let $L(\Lambda)$ be a relative holomorphic discrete series supermodule with integral highest weight. Then $\fsdim(L(\Lambda))\neq0$ if and only if $\Lambda$ is maximally atypical, that is, $\at(\Lambda)=\defect(\gg)$.
\end{theorem}

\begin{proof}
By Theorem~\ref{thm::vanishing_sdim_less_atypic}, we already know that $\fsdim(L(\Lambda))=0$ whenever $\at(\Lambda)<\defect(\gg)$.

Conversely, assume $\at(\Lambda)=\defect(\gg)$. Combining Proposition~\ref{prop::Zuckerman_compatibilities}(d), Corollary~\ref{cor::Zuckerman_packet_character} and Proposition~\ref{prop::Zuckerman_image_simple}, there exists a finite-dimensional simple $\gg$-supermodule $F$ with maximally atypical highest weight such that
$
\sch F=\Theta_\Lambda.
$
In particular, $\sdim(F)\neq 0$. Evaluating the supercharacter identity at zero gives $\fsdim L(\Lambda)\neq 0$. Indeed,
\[
0\neq\sdim F
=\Theta_\Lambda(0)
=(-1)^{|\Delta_n^+|}\fsdim L(\Lambda).\qedhere
\]
\end{proof}

\subsection{Kac--Wakimoto conjecture for the formal superdimension: the nonintegral case}\label{subsec::KW_conjecture_nonintegral_case}
\noindent
For nonintegral highest weights, the derived Zuckerman functors no longer produce a finite-dimensional
supermodule to which Serganova's theorem applies. Instead, we use
generalized Duflo--Serganova cohomology and reduce the analysis to the superdimensions of the corresponding DS cohomologies.

Our starting point is Lemma~\ref{lemm::restricted_character_with_kernel}, where we have seen that the restriction of $\sch L(\Lambda)$ to a subspace associated with a maximal set of
orthogonal isotropic roots recovers the supercharacter of the corresponding
DS cohomology. Averaging over the Weyl group allows us to express the $L$-packet
supercharacter $\Theta_\Lambda$ as a sum of $L$-packet supercharacters of
these DS cohomologies, indexed by compact Weyl group orbits of such root sets. 
This sum simplifies under a condition we call Condition~K: exactly one orbit
contributes nonzero DS cohomology, and its $L$-packet supercharacter is holomorphic
at zero. Under Condition~K, the nonvanishing of the formal superdimension follows
readily. The main task of this subsection is to establish Condition~K for
relative holomorphic discrete series supermodules with nonintegral highest weights.

\subsubsection{Reduction to generalized Duflo--Serganova cohomology} We begin with some notational preliminaries. Let $L(\Lambda)$ be a relative holomorphic discrete series supermodule with maximally atypical highest weight $\Lambda$. We denote the defect of $\gg$ by $d$, that is, $d\coloneqq\defect(\gg)$. Moreover, we fix a positive system $\Delta^+$ and choose $d$ linearly independent, mutually orthogonal isotropic odd roots $\alpha_1,\ldots,\alpha_d$ such that $(\Lambda+\rho,\alpha_i)=0$ for all $i=1,\ldots,d$. Set
$
A\coloneqq\{\alpha_1,\ldots,\alpha_d\}$, and $
\hh_A\coloneqq\bigcap_{\alpha\in A}\ker\alpha.$
For each $\alpha_i$, we choose a nonzero root vector $u_{\alpha_i}\in\gg^{\alpha_i}$ and define
\begin{equation}\label{eq::def_x_A}
x_A\coloneqq
\sum_{j=1}^d\bigl(u_{\alpha_j}-i\omega(u_{\alpha_j})\bigr),
\qquad
c_A\coloneqq x_A^2
=-i\sum_{j=1}^d[u_{\alpha_j},\omega(u_{\alpha_j})]\in\hh.
\end{equation}
The element $x_A$ is in normal form~\eqref{eq::decomp_standard}, belongs to $\YY^{\hom}_{\omega}$, has rank $d$ and commutes with $\hh_{A}$. In addition, $c_{A}$ acts semisimply on $L(\Lambda)$. Finally, we set
$
t_{\alpha_i}\coloneqq[u_{\alpha_i},\omega(u_{\alpha_i})]
$
for $i=1,\ldots,d$, and define as in Section~\ref{subsec::generalized_DS_functor}
\begin{equation}
\hh'\coloneqq\bigoplus_{i=1}^d\CC t_{\alpha_i}
\subseteq\hh_A,
\qquad
\hh_{x_A}=\hh_A/\hh'.
\end{equation}

For $w\in W$, we set $B\coloneqq w(A)$ and define $x_B$, $c_B$, $\hh_B$, and $\hh_{x_B}$ by the same construction, choosing nonzero root vectors for the roots in $B$. All the preceding properties remain valid for these Weyl conjugates.

For every $B$, the reduced real form $\gg_{x_B}^{\omega_{x_B}}$ is described in Proposition~\ref{prop::unitary_real_forms_and_DS}. In particular, it admits a maximal compact subalgebra and a compact Cartan subalgebra as discussed in Section~\ref{sec::basic_classical_Lie_superalgebras}. Thus we can define $W_{x_B}$ and $W_{x_B,c}$ as the Weyl groups of the even part of $\gg_{x_B}$ and its maximal compact subalgebra, respectively. This allows us to apply the construction of the $L$-packet supercharacter to the generalized DS cohomology
$
\DS_{x_B}(L(\Lambda))\cong L(\Lambda)^{c_B}.
$
We define
\begin{equation}
\Theta_{\DS_{x_B}(L(\Lambda))}
\coloneqq
\frac{1}{|W_{x_B,c}|}
\sum_{w\in W_{x_B}}
w\bigl(\sch\DS_{x_B}(L(\Lambda))\bigr).
\end{equation}
If $\gg_{x_B}$ is compact and purely even, then $W_{x_B}=W_{x_B,c}$ and the supercharacter of $\DS_{x_B}(L(\Lambda))$ is Weyl-invariant. Consequently, $
\Theta_{\DS_{x_B}(L(\Lambda))}
=
\sch\DS_{x_B}(L(\Lambda)).
$

We now express the $L$-packet supercharacter $\Theta_\Lambda$ in terms of generalized DS cohomologies. By Lemma~\ref{lemm::restricted_character_with_kernel}, we have the following identification as meromorphic functions on $\hh_{B}$:
\begin{equation}
\sch L(\Lambda)|_{\hh_B}
=
\sch\bigl(\DS_{x_B}(L(\Lambda))\bigr)\circ\pi_B.
\end{equation}
Here, recall that $\pi_B:\hh_B\to\hh_{x_B}$ is the canonical quotient map. This allows us to rewrite the $L$-packet supercharacter as follows:
\begin{equation}
\Theta_\Lambda(H)
=
\frac{1}{|W_c|}
\sum_{w\in W}
\sch\bigl(\DS_{x_{wA}}(L(\Lambda))\bigr)
\bigl(\pi_{wA}(wH)\bigr).
\end{equation}
Next, we reorganize the right-hand side by grouping generalized DS cohomologies related by compact Weyl conjugacy.  

Replacing any root $\alpha\in A$ by $-\alpha$ leaves $\hh_A$ unchanged. By Proposition~\ref{prop::DS_composition}, it also leaves $L(\Lambda)^{c_A}$ unchanged. We therefore introduce
\begin{equation}
N_A\coloneqq
\{w\in W:w(A\cup(-A))=A\cup(-A)\}.
\end{equation}
Using $W_{c}$-invariance of $\sch L(\Lambda)$, the double cosets
$
W_c\backslash W/N_A
$
parametrize the compact Weyl-group orbits of these isotropic root sets.

\begin{lemma}
\label{lemm::packet_DS_double_cosets}
Choose representatives $w_j$ of $W_c\backslash W/N_A$ and set $B_j=w_jA$. Then, as meromorphic functions on $\hh_A$,
\[
\Theta_\Lambda(H)
=
\sum_j
\frac{1}{|W_c\cap w_jN_Aw_j^{-1}|}
\sum_{w\in N_A}
\sch\bigl(\DS_{x_{B_j}}(L(\Lambda))\bigr)
\bigl(\pi_{B_j}(w_jwH)\bigr).
\]
\end{lemma}
\begin{proof}
Decompose $W$ into double cosets
$
W=\bigsqcup_j W_cw_jN_A.
$
For each $j$, the map
\[
W_c\times N_A\longrightarrow W_cw_jN_A,
\qquad
(u,w)\longmapsto uw_jw
\]
is surjective, and every element has exactly
$
|W_c\cap w_jN_Aw_j^{-1}|
$
preimages. Since $\sch L(\Lambda)$ is $W_c$-invariant, we obtain
\[
\begin{aligned}
\Theta_\Lambda(H)
=
\frac{1}{|W_c|}
\sum_j\sum_{v\in W_cw_jN_A}
\sch L(\Lambda)(vH)=
\sum_j
\frac{1}{|W_c\cap w_jN_Aw_j^{-1}|}
\sum_{w\in N_A}
\sch L(\Lambda)(w_jwH).
\end{aligned}
\]

Now every $w\in N_A$ preserves $\hh_A$, which implies
$
w_jwH\in w_j\hh_A=\hh_{B_j}
$ for every $H\in \hh_{A}$.
Applying Lemma~\ref{lemm::restricted_character_with_kernel} finishes the proof:
\[
\sch L(\Lambda)(w_jwH)
=
\sch\bigl(\DS_{x_{B_j}}(L(\Lambda))\bigr)
\bigl(\pi_{B_j}(w_jwH)\bigr).\qedhere
\]
\end{proof}

\begin{lemma}
\label{lemm::packet_DS_reduction}
Let $w_j$ and $B_j$ be as in Lemma~\ref{lemm::packet_DS_double_cosets}, and set
$
N_{B_j}\coloneqq w_jN_Aw_j^{-1}.
$
Then, as meromorphic functions on $\hh_A$,
\[
\Theta_\Lambda(H)
=
\sum_j
\frac{|W_{x_{B_j},c}|}
{|W_{x_{B_j}}|\,|N_{B_j}\cap W_c|}
\sum_{w\in N_{B_j}}
\Theta_{\DS_{x_{B_j}}(L(\Lambda))}
\bigl(\pi_{B_j}(ww_jH)\bigr).
\]
\end{lemma}

\begin{proof}
By Lemma~\ref{lemm::packet_DS_double_cosets}, changing variables from $w\in N_A$ to $w_jww_j^{-1}\in N_{B_j}$ gives
\[
\Theta_\Lambda(H)
=
\sum_j
\frac{1}{|N_{B_j}\cap W_c|}
\sum_{w\in N_{B_j}}
\sch\bigl(\DS_{x_{B_j}}(L(\Lambda))\bigr)
\bigl(\pi_{B_j}(ww_jH)\bigr).
\]

We now average the inner sum over the Weyl group $W_{x_{B_j}}$. 
By the description of the reduced root system in
Section~\ref{subsec::generalized_DS_functor}, every
even root of $\gg_{x_{B_j}}$ is induced by an even root
$\beta\in\Delta_{\bar0}$ orthogonal to all roots in $B_j$.
The reflection $s_\beta$ therefore fixes $B_j$ pointwise
and belongs to $N_{B_j}$. Moreover, its action on
$\hh_{x_{B_j}}$ coincides with the reflection associated
with the corresponding reduced root. Since these
reflections generate $W_{x_{B_j}}$, every
$v\in W_{x_{B_j}}$ admits a lift
$\widetilde v\in N_{B_j}$ satisfying
\[
v\bigl(\pi_{B_j}(H)\bigr)
=
\pi_{B_j}(\widetilde vH),
\qquad H\in\hh_{B_j}.
\]
Since left multiplication by $\widetilde v$ permutes $N_{B_j}$, we obtain
\[
\begin{aligned}
&\sum_{w\in N_{B_j}}
\sch\bigl(\DS_{x_{B_j}}(L(\Lambda))\bigr)
\bigl(\pi_{B_j}(ww_jH)\bigr)\\
&\quad=
\frac{1}{|W_{x_{B_j}}|}
\sum_{w\in N_{B_j}}
\sum_{v\in W_{x_{B_j}}}
\sch\bigl(\DS_{x_{B_j}}(L(\Lambda))\bigr)
\bigl(v\pi_{B_j}(ww_jH)\bigr)\\
&\quad=
\frac{|W_{x_{B_j},c}|}{|W_{x_{B_j}}|}
\sum_{w\in N_{B_j}}
\Theta_{\DS_{x_{B_j}}(L(\Lambda))}
\bigl(\pi_{B_j}(ww_jH)\bigr),
\end{aligned}
\]
where the last equality follows from the definition of the $L$-packet supercharacter. Substituting this into the preceding expression for $\Theta_\Lambda$ proves the assertion.
\end{proof}

\begin{proposition}\label{prop::key_prop_nonintegral_case}
With the notation of Lemma~\ref{lemm::packet_DS_reduction}, suppose that the $L$-packet supercharacter of every nonzero generalized DS cohomology is holomorphic at zero. Then
\[
\Theta_\Lambda(0)
=
\sum_j
\frac{[N_{B_j}:N_{B_j}\cap W_c]}
{[W_{x_{B_j}}:W_{x_{B_j},c}]}
\Theta_{\DS_{x_{B_j}}(L(\Lambda))}(0).
\]
\end{proposition}

\begin{proof}
By Lemma~\ref{lemm::packet_DS_reduction}, $\Theta_\Lambda$ is a finite sum of meromorphic functions, each of which is holomorphic at zero by assumption. We can therefore evaluate the identity at $H=0$. Since every $w\in N_{B_j}$ fixes zero and $\pi_{B_j}(0)=0$, we obtain
\[
\Theta_\Lambda(0)
=
\sum_j
\frac{|W_{x_{B_j},c}|\,|N_{B_j}|}
{|W_{x_{B_j}}|\,|N_{B_j}\cap W_c|}
\Theta_{\DS_{x_{B_j}}(L(\Lambda))}(0)=\sum_j
\frac{[N_{B_j}:N_{B_j}\cap W_c]}
{[W_{x_{B_j}}:W_{x_{B_j},c}]}
\Theta_{\DS_{x_{B_j}}(L(\Lambda))}(0).
\]
\end{proof}

We write
$
W\cdot A\coloneqq\{w(A):w\in W\}$, where $w(A)\coloneqq\{w\alpha:\alpha\in A\}.
$
We identify sets that differ by independent changes of signs of their roots.
Thus we regard each root $\alpha$ as equivalent to $-\alpha$, and the
corresponding set is $A\cup(-A)$, consisting of both signs of every root in $A$. Let
$
N_{A,c}\coloneqq N_A\cap W_c.
$

Proposition~\ref{prop::key_prop_nonintegral_case} simplifies
considerably under the following condition.

\medskip
\noindent\textbf{Condition K.}
\begin{enumerate}
    \item[(i)]
    Among the $W_c$-orbits of $B\cup(-B)$, with
    $B\in W\cdot A$, exactly one has nonzero generalized DS cohomology.

    \item[(ii)]
    For any representative $B$ of this unique orbit, the $L$-packet
    supercharacter $\Theta_{\DS_{x_B}(L(\Lambda))}$ is holomorphic at zero.
\end{enumerate}

The second condition is independent of the choice of representative within a compact Weyl group orbit, since the corresponding DS cohomologies are isomorphic under compact conjugation and their $L$-packet supercharacters have the same value at zero.

The main task in the remainder of this section is to establish
Condition K for maximally atypical relative holomorphic discrete
series supermodules with nonintegral highest weight, for all basic
classical real forms under consideration. This will reduce the nonvanishing statement to the following corollary.

\begin{corollary}
\label{cor::toral_nonvanishing}
Suppose that Condition K holds, and let $B=wA$ represent the unique $W_c$-orbit for which $\DS_{x_B}(L(\Lambda))\neq0$. Assume that its $L$-packet supercharacter is holomorphic at zero. Then
\begin{equation*}
\Theta_\Lambda(0)
=
\frac{[N_B:N_B\cap W_c]}
     {[W_{x_B}:W_{x_B,c}]}
\Theta_{\DS_{x_B}(L(\Lambda))}(0).
\end{equation*}
If, moreover, $W_c$ acts transitively on
$
W\cdot(A\cup(-A)),
$
then
\[
\Theta_\Lambda(0)
=
\frac{[W:W_c]}
     {[W_{x_A}:W_{x_A,c}]}
\Theta_{\DS_{x_A}(L(\Lambda))}(0).
\]
\end{corollary}

\subsubsection{Verification of Condition~K} 
We now verify Condition K for the basic classical Lie superalgebras. The arguments differ for the classical and exceptional types, and we treat them separately. \medskip

\paragraph{\textbf{Classical series.}} We first treat the special linear Lie superalgebras, followed by the orthosymplectic Lie superalgebras. 

For $\slmn$ with $m\neq n$, we use the notation of Section~\ref{subsec::type_A_structure}. Assume that $L(\Lambda)$ is unitarizable with respect to the real form $\su(p,q\vert 0,n)$, where $p,q\geq1$ and $p+q=m$ (\emph{cf.}~Theorem~\ref{theorem::unitary-highest-lowest-weight}). Without loss of generality, we can assume that $L(\Lambda)$ is a unitarizable highest weight $\gg$-supermodule with respect to the $pq$-nonstandard positive system defined in~\eqref{eq::non_standard_positive}. Its simple system is
\begin{equation}
\begin{aligned}
\Uppi
={}
\{\epsilon_1-\epsilon_2,\ldots,\epsilon_{p-1}-\epsilon_p,
\epsilon_p-\delta_1,\delta_1-\delta_2,\ldots,\delta_{n-1}-\delta_n,
\delta_n-\epsilon_{p+1},\epsilon_{p+1}-\epsilon_{p+2},\ldots,
\epsilon_{m-1}-\epsilon_m\},
\end{aligned}
\end{equation}
corresponding to the ordering
$
(\epsilon_1,\ldots,\epsilon_p,\delta_1,\ldots,\delta_n,
\epsilon_{p+1},\ldots,\epsilon_m),
$
where consecutive differences in the ordering determine the simple roots. Note that all pairings of the Weyl vector $\rho$ with isotropic roots are integral. 

With respect to the $pq$-nonstandard positive system, we call
$
\{\epsilon_1,\ldots,\epsilon_p\}
$ and $
\{\epsilon_{p+1},\ldots,\epsilon_m\}$
the two \emph{compact blocks}, corresponding to the $\gl(p)$- and $\gl(q)$-factors of $\kk$, respectively. Before establishing Condition K, we prove two preliminary lemmas.

\begin{lemma}\label{lemm::non_integral_slmn_weight}
Let $L(\Lambda)$ be a relative holomorphic discrete series $\sl(m\vert n)$-supermodule with maximally atypical highest weight. If $\Lambda$ is nonintegral, then $m>n$, and all isotropic roots $\epsilon_i-\delta_a$ satisfying $(\Lambda+\rho,\epsilon_i-\delta_a)=0$ have their $\epsilon$-indices in a single compact block of size at least $n$.
\end{lemma}

\begin{proof}
Write
$
\Lambda=(\lambda^1,\ldots,\lambda^m\vert\mu^1,\ldots,\mu^n).
$
By the Harish-Chandra parametrization in Section~\ref{subsec::HC_parameterization}, the differences $\mu^a-\mu^b$ are integral, as are the differences $\lambda^i-\lambda^j$ within each compact block. Nonintegrality of $\Lambda$ therefore implies
\[
\lambda^i-\lambda^j\notin\ZZ,
\qquad i\leq p<j.
\]

If
$
(\Lambda+\rho,\epsilon_i-\delta_a)=0
$, then $\lambda^i+\mu^a\in\ZZ$ because $(\rho,\epsilon_i-\delta_a)\in\ZZ$. 
Suppose that there exist $i\leq p<j$ and $a,b$ such that
\[
(\Lambda+\rho,\epsilon_i-\delta_a)
=
(\Lambda+\rho,\epsilon_j-\delta_b)=0.
\]
Then $\lambda^i+\mu^a,\lambda^j+\mu^b\in\ZZ$, and therefore
\[
\lambda^i-\lambda^j
=
(\lambda^i+\mu^a)
-(\lambda^j+\mu^b)
-(\mu^a-\mu^b)
\in\ZZ,
\]
contradicting the nonintegrality of $\Lambda$. We conclude that the $\epsilon$-indices must belong to the same compact block.

Finally, maximal atypicality requires $\min(m,n)$ mutually orthogonal isotropic roots, whose $\epsilon$-indices must be distinct. We proved above that all these indices belong to a single compact block, forcing $\min(m,n)<m$, and thus $n<m$. Moreover, the compact block must contain at least $n$ indices.
\end{proof}

\begin{lemma}\label{lemm::weight_condition_slmn}
Let $L(\Lambda)$ be a relative holomorphic discrete series $\sl(m\vert n)$-supermodule with nonintegral and maximally atypical highest weight
$
\Lambda=(\lambda^1,\ldots,\lambda^m\vert\mu^1,\ldots,\mu^n).
$
Then every weight $\nu$ of $L(\Lambda)$ satisfies
\[
\nu^i-\nu^j\leq\lambda^1-\lambda^m,
\qquad i\leq p<j.
\]
\end{lemma}

\begin{proof}
Since $L(\Lambda)$ is a highest weight supermodule with respect to the $pq$-nonstandard positive system, every weight $\nu$ satisfies
\[
\Lambda-\nu
=
\sum_{\alpha\in\Uppi}n_\alpha\alpha,
\qquad n_\alpha\in\ZZ_{\geq0}.
\]
By the explicit description of $\Uppi$, we have
$
(\alpha,\epsilon_1-\epsilon_m)\geq0
$
for every $\alpha\in\Uppi$. Consequently,
\[
\nu^1-\nu^m
=
(\nu,\epsilon_1-\epsilon_m)
\leq
(\Lambda,\epsilon_1-\epsilon_m)
=
\lambda^1-\lambda^m.
\]

Finally, the weight set of $L(\Lambda)$ is $W_c$-invariant, since its restriction to $\kk$ is locally finite. As $W_c$ permutes the $\epsilon$-indices independently within the two compact blocks, the same inequality holds for every $i\leq p<j$.
\end{proof}

We now prove that Condition K holds. By Lemma~\ref{lemm::non_integral_slmn_weight}, all isotropic atypical roots involve $\epsilon$-indices from the same compact block. It suffices to consider the case where these indices belong to the $q$-block. The $p$-block case is analogous. Thus we may choose
\begin{equation}
A=\{\delta_a-\epsilon_{j_a}:1\leq a\leq n\},
\qquad j_a>p.
\end{equation}
Define $x_{A}$ and $c_{A}$ as in~\eqref{eq::def_x_A}. The integral-root subsuperalgebra associated with $\Lambda$ is then the Levi subsuperalgebra
\begin{equation}
\mathfrak l
=
\mathfrak s\bigl(\gl(p)\oplus\gl(q\vert n)\bigr)
\coloneqq
\{(P,Q)\in\gl(p)\oplus\gl(q\vert n):\operatorname{tr}P+\str Q=0\}
\subseteq\gg,
\end{equation}
that is, $\mathfrak s(\cdot)$ denotes the kernel of the total supertrace.
Indeed, the roots within each compact block are integral, and the atypicality equations imply that all roots $\epsilon_j-\delta_a$ with $j>p$ are integral. By contrast, roots connecting the two compact blocks are nonintegral. This is the Levi construction considered in~\cite[Section~3.4]{chengMazorchukWang2013}. Equivalently, $\mathfrak l$ is the Levi subsuperalgebra obtained by removing the simple root $\epsilon_p-\delta_1$ from $\Uppi$. In particular, every root in $A$ belongs to the root system of $\mathfrak l$ and $\ll_{\bar{0}}\subset \kk$. If the atypical roots involve the $p$-block instead, the corresponding Levi subsuperalgebra is
$
\mathfrak s\bigl(\gl(p\vert n)\oplus\gl(q)\bigr).
$

\begin{proposition}
\label{prop::condition_K_slmn}
Let $L(\Lambda)$ be a relative holomorphic discrete series $\sl(m\vert n)$-supermodule with nonintegral and maximally atypical highest weight $\Lambda$. Then Condition K holds. Moreover, $\DS_{x_A}(L(\Lambda))$ is either finite-dimensional or a relative holomorphic discrete series module of $\gg_{x_{A}}$. 
\end{proposition}

\begin{proof} 
We first prove that there is exactly one $W_c$-orbit of maximal isotropic root sets $B$ for which
$
\DS_{x_B}(L(\Lambda))\neq0.
$
Let $v_\Lambda$ be a highest weight vector of $L(\Lambda)$ and set
\[
E\coloneqq\UE(\mathfrak l)v_\Lambda.
\]
Since $\mathfrak l_{\bar0}\subseteq\kk$, compact local finiteness and the PBW theorem imply that $E$ is finite-dimensional. 
Moreover, $E$ is simple: the $\omega$-stability of $\ll$ and unitarity imply that every $\ll$-submodule has an invariant orthogonal complement. The one-dimensional highest weight space generates $E$, so every such decomposition is trivial.

Next, we show that $E$ has maximally atypical highest weight. Let $\rho_{\ll}$ be the Weyl vector of $\ll$ with respect to the positive system induced by $\Delta^+$. Then by definition of $\ll$, one has
\[
\rho-\rho_{\mathfrak l}
=
\frac12\left(
\sum_{i=1}^p\sum_{j=p+1}^m(\epsilon_i-\epsilon_j)
-
\sum_{i=1}^p\sum_{a=1}^n(\epsilon_i-\delta_a)
\right),
\]
and $\rho-\rho_{\ll}$ is orthogonal to every root of $\mathfrak l$. In particular, for every $\alpha\in A$,
\[
(\Lambda+\rho_{\mathfrak l},\alpha)
=
(\Lambda+\rho,\alpha)=0.
\]
The statement now follows from $
\defect(\mathfrak l)=\min(q,n)=n
$ and the fact that $A$ consists of $n$ mutually orthogonal elements.

By the Kac--Wakimoto conjecture for finite-dimensional supermodules~\cite{serganova2011superdimension},
$
\sdim E\neq0.
$ Moreover, the generalized DS functor $\DS_{x_{A}}$ applied to $E$ preserves the superdimension and equals $E^{c_{A}}$. We conclude
$
\sdim E^{c_A}=\sdim E\neq0, 
$
which implies
\[
\DS_{x_A}(L(\Lambda))
\cong L(\Lambda)^{c_A}\neq0,
\]
using $E^{c_{A}}\subset L(\Lambda)^{c_{A}}$.

We next verify uniqueness of the contributing orbit. This follows immediately if we show that $
\DS_{x_B}(L(\Lambda))=0
$
whenever $B$ contains a root involving an $\epsilon$-index from the $p$-block.

Assume such a $B$ contains $
\beta=\pm(\epsilon_i-\delta_a)
$ with $i\leq p$. By definition of $A$, one has
$
\lambda^j+\mu^b\in\ZZ
$
for any $j>p$ and $b\in\{1,\ldots,n\}$.
Since $\Lambda$ is $\Delta_{c}^{+}$-dominant integral, $\mu^a-\mu^b\in\ZZ$, whereas nonintegrality gives $\lambda^i-\lambda^j\notin\ZZ$. We conclude
\[
\lambda^i+\mu^a
=
(\lambda^i-\lambda^j)
+(\lambda^j+\mu^b)
+(\mu^a-\mu^b)
\notin\ZZ,
\]
and therefore $(\Lambda,\beta)\notin\ZZ$.

Every weight $\nu$ of $L(\Lambda)$ differs from $\Lambda$ by an element of the root lattice (Lemma~\ref{lemm::properties_HWM}). Since the inner product of any two roots of $\sl(m\vert n)$ is integral, we obtain
\[
(\nu,\beta)\in(\Lambda,\beta)+\ZZ,
\]
and hence $(\nu,\beta)\neq0$. However, this forces that 
$
[u_\beta,\omega(u_\beta)]
$ has no kernel on $L(\Lambda)$, and hence $\DS_{x_{B}}(L(\Lambda))=L(\Lambda)^{c_B}=\{0\}$ by Proposition~\ref{prop::DS_composition}.

Every maximal isotropic root set giving nonzero DS cohomology involves only indices from the $q$-block. The compact Weyl group
$
W_c\cong S_p\times S_q\times S_n
$
acts transitively on such root sets, so all sets $B$ giving nonzero DS cohomology belong to a single $W_c$-orbit, represented by $A$. This verifies Condition K(i).

It remains to verify that the $L$-packet supercharacter of $\DS_{x_A}(L(\Lambda))$ is holomorphic at zero. By Proposition~\ref{prop::unitary_real_forms_and_DS}, $\gg_{x_{A}}^{\omega_{x_{A}}}\cong \su(p,q-n)$. If $q=n$, this real form is compact, so $\DS_{x_A}(L(\Lambda))$ is finite-dimensional. Its $L$-packet supercharacter therefore agrees with its supercharacter and is holomorphic at zero.

Suppose $q>n$. By Theorem~\ref{thm::DS_unitary} and Lemma~\ref{lemm::DS_is_HW}, $\DS_{x_A}(L(\Lambda))$ is a simple unitarizable highest weight module. Let $\nu$ denote its highest weight. We show that $\nu$ satisfies the Harish-Chandra condition, hence $\DS_{x_{A}}(L(\Lambda))$ is a relative holomorphic discrete series $\gg_{x_{A}}^{\omega_{x_{A}}}$-module. Then Condition K(ii) is immediate.

For every weight $\xi$ of $L(\Lambda)$, one has by Lemma~\ref{lemm::weight_condition_slmn}
\[
\xi^i-\xi^j
\leq \lambda^1-\lambda^m,
\qquad i\leq p<j.
\]
Using $\DS_{x_{A}}(L(\Lambda))=L(\Lambda)^{c_{A}}$, the same bound therefore holds for the weights of the generalized DS cohomology. The Harish-Chandra condition for $L_{0}(\Lambda)$ is
$
\lambda^1-\lambda^m+m-1<0,
$
which forces
\[
(\nu+\rho_{x_A},\epsilon_i-\epsilon_j)
\leq
\lambda^1-\lambda^m+(m-n-1)<0.
\]
Thus $\nu$ satisfies the Harish-Chandra conditions for $\su(p,q-n)$, and $\DS_{x_A}(L(\Lambda))$ is a relative holomorphic discrete series module.
\end{proof}

We prove Condition K for $\gg=\osp(m\vert2n)$, using the same strategy as for $\sl(m\vert n)$. Set $r\coloneqq\lfloor m/2\rfloor$. We can assume $m\geq2$, since $\osp(1\vert2n)$ has defect zero. In this case, the only maximal isotropic root set is empty, the corresponding generalized DS cohomology equals $L(\Lambda)$, and its $L$-packet supercharacter is holomorphic at zero by Example~\ref{ex::osp_one_formal_superdimension}.

We use the notation of Section~\ref{subsec::appendix_osp} and the positive system defined in~\eqref{eq::type_A_plus_1}--\eqref{eq::positive_system_typeA_3}. Let $L(\Lambda)$ be a relative holomorphic discrete series supermodule with nonintegral and maximally atypical highest weight. We write
\begin{equation}
\Lambda
=(\lambda^{1}, \ldots, \lambda^{r}\vert \mu^{1},\ldots,\mu^{n})=
\sum_{i=1}^r\lambda^i\epsilon_i
+
\sum_{a=1}^n\mu^a\delta_a.
\end{equation}

The supermodule $L(\Lambda)$ is defined with respect to either $\osp(0,m\vert2n;\RR)$ or $\osp^{\ast}(2r\vert0,2n)$, the latter being possible only for $m=2r$ (\emph{cf.}~Theorem~\ref{theorem::unitary-highest-lowest-weight}). Their even parts are
$
\so(m,\RR)\oplus\sp(2n,\RR)
$ and $
\so^*(2r)\oplus\mathfrak{usp}(2n)$,
respectively. 

Nonintegral maximally atypical highest weights can occur only for the first real form.

\begin{lemma}
Let $L(\Lambda)$ be a relative holomorphic discrete series $\osp(m\vert2n)$-supermodule with nonintegral and maximally atypical highest weight $\Lambda$. Then $L(\Lambda)$ is unitarizable only with respect to $\osp(0,m\vert2n;\RR)$.
\end{lemma}

\begin{proof}
Suppose that $L(\Lambda)$ is unitarizable with respect to $\osp^{\ast}(2r\vert0,2n)$. Its maximal compact subalgebra has complexification
$
\kk\cong\gl(r)\oplus\sp(2n).
$
Since $L(\Lambda)$ is locally finite over $\kk$, we have
\[
\mu^a\in\ZZ,
\qquad
\lambda^i-\lambda^j\in\ZZ
\]
for all $a,i,j$. By atypicality, there exists an isotropic root $\beta=\delta_a\pm\epsilon_i$ such that $(\Lambda+\rho,\beta)=0$. The Weyl vector $\rho$ has integral coordinates, which implies $\lambda^i\in\ZZ$, and therefore $\lambda^j\in\ZZ$ for every $j$. Thus $\Lambda$ is integral, a contradiction.
\end{proof}

In what follows, we consider the real form $\osp(0,m\vert2n;\RR)$ and use the Harish-Chandra parametrization of Section~\ref{subsec::HC_parameterization}.

If $m=2r+1$, compact dominance gives
$
\lambda^1\geq\cdots\geq\lambda^r\geq0.
$
If $m=2r$, it gives
$
\lambda^1\geq\cdots\geq\lambda^{r-1}\geq|\lambda^r|.
$
In the latter case, the transformation
$
\epsilon_r\mapsto-\epsilon_r
$
is induced by conjugation with an orthogonal matrix of determinant $-1$. It preserves the real form and the positive system of Case~A, since it interchanges $\epsilon_i-\epsilon_r$ with $\epsilon_i+\epsilon_r$ and $\delta_a-\epsilon_r$ with $\delta_a+\epsilon_r$. Without loss of generality, we assume $\lambda^r\geq0$. Then the Harish-Chandra parametrization yields
\begin{equation}
\lambda^1\geq\cdots\geq\lambda^r\geq0,
\qquad
\mu^1\geq\cdots\geq\mu^n,
\qquad
\mu^1+n<0.
\end{equation}

To establish Condition K for $\osp(m\vert 2n)$, we need the following lemma.

\begin{lemma}
\label{lemm::weight_condition_osp}
Let $L(\Lambda)$ be a relative holomorphic discrete series supermodule, and let
$
\nu=(\xi^{1},\ldots,\xi^{r}\vert \eta^{1},\ldots,\eta^{n})
$ be a weight of $L(\Lambda)$.
Then
\[
\eta^a\leq\mu^1
\qquad 1\leq a\leq n.
\]
Moreover, every odd isotropic root $\alpha$ satisfying
$(\Lambda+\rho,\alpha)=0$ is, up to sign, of the form
$\delta_a-\epsilon_i$.
\end{lemma}

\begin{proof}

Every simple root of the chosen positive system has nonnegative $\delta_1$-coefficient. By Lemma~\ref{lemm::properties_HWM}(a), $\eta^1\leq\mu^1$. The $W_c$-invariance of the weight set and the transitive action of $W_c$ on the $\delta$-indices then give $\eta^a\leq\mu^1$ for every $a$.

Every isotropic root is of the form $\pm\delta_a\pm\epsilon_i$. By the Harish-Chandra inequalities,
\[
\mu^a+n-\frac m2-a+1<0,
\qquad
\lambda^i+\frac m2-i\geq0.
\]
Using $(\delta_a,\delta_a)=-\frac12$ and $(\epsilon_i,\epsilon_i)=\frac12$ (Appendix~\ref{subsec::appendix_osp}), we obtain
\[
(\Lambda+\rho,\delta_a+\epsilon_i)
=
\frac12\left(
-\left(\mu^a+n-\frac m2-a+1\right)
+\left(\lambda^i+\frac m2-i\right)\right)>0.
\]
This excludes both $\pm(\delta_a+\epsilon_i)$, leaving only roots of the form $\pm(\delta_a-\epsilon_i)$.
\end{proof}

We fixed the positive system of Case~A from Section~\ref{subsec::appendix_osp}. Removing the last simple root from $\Uppi_A$ gives a Levi subsuperalgebra
\begin{equation}
\mathfrak l\cong\gl(r\vert n).
\end{equation}
More precisely, the omitted root is $\epsilon_r$ when $m=2r+1$, $\epsilon_{r-1}+\epsilon_r$ when $m=2r$ and $r\geq2$, and $\delta_n+\epsilon_1$ when $m=2$. The resulting simple system is
\begin{equation}
\{\delta_a-\delta_{a+1}:1\leq a<n\}
\cup\{\delta_n-\epsilon_1\}
\cup\{\epsilon_i-\epsilon_{i+1}:1\leq i<r\},
\end{equation}
and we denote the set of roots for $\ll$ by $\Delta_{\ll}$. In particular,
$
\mathfrak l_{\bar0}
\cong\gl(r)\oplus\gl(n)
\subseteq\kk.
$ By Lemma~\ref{lemm::weight_condition_osp}, every isotropic root $\alpha$ satisfying $(\Lambda+\rho,\alpha)=0$ belongs to $\Delta_{\ll}$. Hence we may choose a maximal set $A$ of mutually orthogonal isotropic roots satisfying the atypicality condition entirely within the root system of $\mathfrak l$.

\begin{proposition}
\label{prop::condition_K_osp}
Let $L(\Lambda)$ be a relative holomorphic discrete series $\osp(m\vert2n)$-supermodule with nonintegral and maximally atypical highest weight $\Lambda$. Then Condition K holds. Moreover, for $m>1$, every nonzero generalized DS cohomology of maximal rank is either a finite-dimensional module over a compact Lie algebra, a relative holomorphic discrete series module over a Lie algebra, or a relative holomorphic discrete series supermodule over $\osp(1\vert2(n-r))$.
\end{proposition}

\begin{proof}
We first prove existence of a nonzero
DS cohomology. The proof is analogous to the $\slmn$-case. Let $v_\Lambda$ be a highest weight vector of $L(\Lambda)$ and set
\[
E=\UE(\mathfrak l)v_\Lambda.
\]
As for $\slmn$, it follows that $E$ is finite-dimensional and simple (\emph{cf.}~proof of Proposition~\ref{prop::condition_K_slmn}). The highest weight $\Lambda$ is also maximally atypical for $\ll$, since $|A|=\min(r,n)=\defect(\mathfrak l)$, and 
\[
(\Lambda+\rho_{\mathfrak l},\alpha)
=
(\Lambda+\rho,\alpha)=0
\qquad \alpha\in A,
\]
where we use that $
\rho-\rho_{\mathfrak l}
=
\frac{n+r-m+1}{2}
\left(
\sum_{a=1}^{n}\delta_a-\sum_{i=1}^{r}\epsilon_i
\right)
$ is orthogonal to every root of $\ll$.
Hence the
$\mathfrak l$-supermodule $E$ is maximally atypical.
In particular, 
$
\sdim \DS_{x_{A}}(E)=\sdim E^{c_{A}}= \sdim E\neq0
$
(\cite{serganova2011superdimension}). We conclude
\[
\DS_{x_A}(L(\Lambda))
\cong L(\Lambda)^{c_A}\neq0,
\]
as $E^{c_{A}} \subset L(\Lambda)^{c_{A}}$.

The nonvanishing of $\DS_{x_A}(L(\Lambda))$ and Theorem~\ref{thm::DS_unitary} imply that it is simple and unitarizable. Moreover, by assumption and Lemma~\ref{lemm::DS_is_HW}, it is a highest weight $\gg_{x_A}$-supermodule. We show that it is finite-dimensional, a relative holomorphic discrete series module, or a relative holomorphic discrete series supermodule, according to the reduced Lie superalgebra $\gg_{x_{A}}$, which is by Proposition~\ref{prop::unitary_real_forms_and_DS} given by
\[
\gg_{x_A}\cong
\begin{cases}
\so(m-2n),&r\geq n,\\
\sp(2(n-r)),&r<n,\ m=2r,\\
\osp(1\vert2(n-r)),&r<n,\ m=2r+1.
\end{cases}
\]
First, assume $r\geq n$. Then $\gg_{x_{A}}^{\omega_{x_{A}}}\cong \so(m-2n,\RR)$ is compact and $\DS_{x_A}(L(\Lambda))$ must be finite-dimensional.

Next, assume $r<n$. Every weight of $\DS_{x_A}(L(\Lambda))$ is the restriction of a weight of $L(\Lambda)$ to $\hh_{x_A}$. By Lemma~\ref{lemm::weight_condition_osp}, all its $\delta$-coordinates are bounded above by $\mu^1$. In particular, every highest weight $\nu$ of an even constituent satisfies $\nu^a\leq\mu^1$ for all $a$. Together with $\mu^1+n-r<0$, this gives
\[
\nu^a+(n-r)-a+1<0,
\qquad 1\leq a\leq n-r,
\]
which are the strict Harish-Chandra conditions for the common even part $\sp(2(n-r),\RR)$. Therefore, $\DS_{x_A}(L(\Lambda))$ is a relative holomorphic discrete series module in the purely even case and a relative holomorphic discrete series supermodule in the $\osp(1\vert2(n-r))$ case. In either case, its $L$-packet supercharacter is holomorphic at zero.

It remains to prove uniqueness of the contributing $W_c$-orbit. The compact Weyl group is
\[
W_c\cong
\begin{cases}
\bigl((\ZZ/2\ZZ)^r\rtimes S_r\bigr)\times S_n,
& m=2r+1,\\
\bigl((\ZZ/2\ZZ)^{r-1}\rtimes S_r\bigr)\times S_n,
& m=2r.
\end{cases}
\]
For $m=2r+1$, it acts transitively on maximal isotropic root sets up to independent sign changes. The same holds for $m=2r$ and $r>n$.

It remains to consider $m=2r$ with $r\leq n$. Every maximal isotropic root set then involves all $r$ $\epsilon$-coordinates. The restriction to an even number of sign changes gives exactly two $W_c$-orbits, distinguished by the parity of the number of roots of the form $\delta_a+\epsilon_i$. Representatives are
\[
A_+=\{\delta_i-\epsilon_i:1\leq i\leq r\},
\qquad
A_-=
\{\delta_i-\epsilon_i:1\leq i<r\}
\cup\{\delta_r+\epsilon_r\}.
\]
We claim that at most one of the generalized DS cohomologies associated with $A_+$ and $A_-$ is nonzero.

Suppose both DS cohomologies are nonzero. By~\cite[Lemmas~6.10--6.11]{gorelik2022duflo} and Lemma~\ref{lem::generalized_DS_eta_compatibility}, they have the same infinitesimal character. Both are simple relative holomorphic discrete series modules for $\sp(2(n-r))$ with the same positive system. Their highest weights therefore coincide, say
$
\nu=\sum_{a=r+1}^n\nu^a\delta_a.
$

The root lattice of $\osp(2r\vert2n)$ consists of the elements
\[
\sum_{i=1}^r k_i\epsilon_i+\sum_{a=1}^n l_a\delta_a,
\qquad k_i,l_a\in\ZZ,
\qquad \sum_{i=1}^r k_i+\sum_{a=1}^n l_a\in2\ZZ.
\]
In particular, every weight
$
\xi=\sum_{i=1}^r\xi^i\epsilon_i+\sum_{a=1}^n\eta^a\delta_a
$
of $L(\Lambda)$ satisfies
\[
\sum_{i=1}^r(\xi^i-\lambda^i)+\sum_{a=1}^n(\eta^a-\mu^a)\in2\ZZ.
\]

Choose a weight $\xi$ of $L(\Lambda)^{c_{A_+}}$ restricting to $\nu$. The relations $\xi^i=-\eta^i$ for $1\leq i\leq r$ and $\eta^a=\nu^a$ for $a>r$ give
\[
\sum_{i=1}^r\lambda^i+\sum_{a=1}^n\mu^a-\sum_{a=r+1}^n\nu^a\in2\ZZ.
\]
For a weight $\xi$ of $L(\Lambda)^{c_{A_-}}$ restricting to $\nu$, we have $\xi^i=-\eta^i$ for $i<r$, $\xi^r=\eta^r$, and $\eta^a=\nu^a$ for $a>r$. Together with $\xi^r-\lambda^r\in\ZZ$, this gives
\[
\sum_{i=1}^r\lambda^i+\sum_{a=1}^n\mu^a-\sum_{a=r+1}^n\nu^a
\in2\lambda^r+2\ZZ.
\]

Comparing these relations gives $\lambda^r\in\ZZ$. The $\Delta_c^+$-dominant integrality of $\Lambda$ then implies $\lambda^i\in\ZZ$ for all $i$. Maximal atypicality yields $\lambda^i+\mu^a\in\ZZ$ for some $i,a$, and using $\Delta_c^+$-integrality again yields $\mu^b\in\ZZ$ for all $b$, contradicting the nonintegrality of $\Lambda$. Hence exactly one of the two $W_c$-orbits contributes, completing the proof.
\end{proof}

\paragraph{\textbf{Exceptional series.}}

Unlike in the classical series, the following arguments apply to all highest weights $\Lambda$.

We begin with $F(4)$ and use the notation of Section~\ref{subsec::appendix_F4}. The relevant real forms are $F(4,0)$ and $F(4,2)$, with even parts $\so(7,\RR)\oplus\sl(2,\RR)$ and $\so(2,5)\oplus\su(2)$, respectively. Since $F(4)$ has defect one, a maximal orthogonal isotropic root set is a singleton $A=\{\beta\}$. We choose $\beta$ to be the isotropic simple root of the corresponding positive system and set
\begin{equation}
x_A=e_\beta-i\omega(e_\beta),
\qquad 0\neq e_\beta\in\gg^\beta.
\end{equation}
Then $\omega(x_A)=ix_A$. The reduced Lie algebra is $\gg_{x_A}\cong\sl(3)$, with induced real form $\su(3)$ for $F(4,0)$ and $\su(1,2)$ for $F(4,2)$ (see Proposition~\ref{prop::unitary_real_forms_and_DS}).

The odd isotropic roots of $F(4)$ are
$
\overline{\Delta}_{\bar1}
=
\{
\frac12(\pm\epsilon_1\pm\epsilon_2\pm\epsilon_3\pm\delta)\},
$
where the signs are independent. We retain the positive even roots and the bilinear form of Section~\ref{subsec::appendix_F4}. We choose the simple system
\begin{equation}
\Uppi
=
\left\{
\delta,\epsilon_3,\epsilon_2-\epsilon_3,
\frac12(\epsilon_1-\epsilon_2-\epsilon_3-\delta)
\right\},
\end{equation}
which has Weyl vector
\begin{equation}\label{eq::rho_F4}
\rho
=
\frac12\epsilon_1+\frac32\epsilon_2
+\frac12\epsilon_3+\frac12\delta.
\end{equation}
We denote simple systems obtained by odd reflections from $\Uppi$ by $\Uppi'$ and their corresponding Weyl vectors by $\rho'$. To establish the nonvanishing of $\DS_{x_A}(L(\Lambda))$, we first prove the following lemmas.

\begin{lemma}
\label{lemm::F4_odd_reflections_and_kernel}
Let $L(\Lambda)$ be a simple highest weight $F(4)$-supermodule relative to $\Uppi$, and suppose that $(\Lambda+\rho,\beta)=0$ for some isotropic root $\beta$. Then a sequence of odd reflections yields a highest weight $\Lambda'$ satisfying $(\Lambda',\gamma)=0$ for some isotropic root $\gamma$.
\end{lemma}

\begin{proof}
Replacing $\beta$ by $-\beta$ if necessary, we may assume that its $\epsilon_1$-coefficient is positive. Set
\begin{equation*}
\begin{aligned}
\alpha_1&=\frac12(\epsilon_1-\epsilon_2-\epsilon_3-\delta),
&\alpha_2&=\frac12(\epsilon_1-\epsilon_2-\epsilon_3+\delta),\\
\alpha_3&=\frac12(\epsilon_1-\epsilon_2+\epsilon_3-\delta),
&\alpha_4&=\frac12(\epsilon_1+\epsilon_2-\epsilon_3-\delta).
\end{aligned}
\end{equation*}
For each sequence of odd reflections, let $\Uppi'$ denote the resulting simple system and $\rho'$ its Weyl vector. The following table gives, for each $\beta$, a sequence of odd reflections such that $(\rho',\beta)=0$:
\begin{equation*}
\begin{array}{c|c|c}
\beta & \text{successive reflections} & \rho' \\ \hline
\frac12(\epsilon_1-\epsilon_2\pm\epsilon_3\pm\delta)
& \alpha_1
& \epsilon_1+\epsilon_2 \\[2pt]
\frac12(\epsilon_1+\epsilon_2-\epsilon_3-\delta)
& \alpha_1,\alpha_3
& \frac12(3\epsilon_1+\epsilon_2+\epsilon_3-\delta) \\[2pt]
\frac12(\epsilon_1+\epsilon_2-\epsilon_3+\delta)
& \text{none}
& \frac12(\epsilon_1+3\epsilon_2+\epsilon_3+\delta) \\[2pt]
\frac12(\epsilon_1+\epsilon_2+\epsilon_3-\delta)
& \alpha_1,\alpha_3,\alpha_4
& 2\epsilon_1+\epsilon_2-\delta \\[2pt]
\frac12(\epsilon_1+\epsilon_2+\epsilon_3+\delta)
& \alpha_1,\alpha_2
& \frac12(3\epsilon_1+\epsilon_2-\epsilon_3+\delta)
\end{array}
\end{equation*}
These sequences are valid: $\alpha_1$ is initially simple, $\alpha_2$ and $\alpha_3$ become simple after reflecting in $\alpha_1$, and $\alpha_4$ becomes simple after reflecting in $\alpha_3$. Each odd reflection adds the reflected root to the Weyl vector. Thus the values of $\rho'$ follow from~\eqref{eq::rho_F4}. The bilinear form of Section~\ref{subsec::appendix_F4} then gives $(\rho',\beta)=0$ in every case.

Next, we perform the indicated reflections for $\Lambda$. If the current highest weight is orthogonal to a root being reflected, the assertion follows. Otherwise, the shifted highest weight remains unchanged under each reflection. Hence $\Lambda'+\rho'=\Lambda+\rho$, and therefore
\begin{equation*}
(\Lambda',\beta)
=(\Lambda+\rho,\beta)-(\rho',\beta)=0.\qedhere
\end{equation*}
\end{proof}

\begin{lemma}\label{lemm::F4_weight_bound} 
Let $L(\Lambda)$ be a relative holomorphic discrete series supermodule for $F(4,2)$ with $\Lambda=(\lambda_1,\lambda_2,\lambda_3\vert\mu)$. Then every weight $\nu=(\nu_1,\nu_2,\nu_3\vert\nu_4)$ of $L(\Lambda)$ satisfies
\begin{equation*}
\nu_1\pm\nu_j\leq\lambda_1+\lambda_2,
\qquad j=2,3.
\end{equation*}
\end{lemma}

\begin{proof}
Let $F_{\Lambda}=\UE(\kk)v_\Lambda$. 
The positive compact roots are $\epsilon_2$, $\epsilon_3$, $\epsilon_2\pm\epsilon_3$, and $\delta$. Hence every weight $\eta$ of $F_\Lambda$ satisfies $\eta_1=\lambda_1$ and $\eta_2\leq\lambda_2$.
Moreover, 
$
W_c\cong\bigl((\ZZ/2\ZZ)^2\rtimes S_2\bigr)\times\ZZ/2\ZZ
$
fixes $\epsilon_1$, acts by signed permutations on $\epsilon_2,\epsilon_3$, and changes the sign of $\delta$. Since the set of weights of $F_{\Lambda}$ is $W_{c}$-invariant, we conclude $|\eta_j|\leq\lambda_2$ for $j=2,3$.

The positive noncompact even roots are
$
\epsilon_1,\epsilon_1\pm\epsilon_2,\epsilon_1\pm\epsilon_3,
$
and the positive odd roots are
$
\frac12(\epsilon_1\pm\epsilon_2\pm\epsilon_3\pm\delta).
$
For every such root $\alpha$, we have
\[
(\alpha,\epsilon_1\pm\epsilon_j)\geq0,
\qquad j=2,3.
\]

Every weight $\nu$ of $L(\Lambda)$ has the form
$
\nu=\eta-\sum_{\alpha}n_\alpha\alpha,
$
where $\eta$ is a weight of $F_\Lambda$ and $n_\alpha\in\ZZ_{\geq0}$. Therefore
\[
\nu_1\pm\nu_j
=\eta_1\pm\eta_j-\sum_\alpha n_\alpha(\alpha,\epsilon_1\pm\epsilon_j)
\leq\eta_1\pm\eta_j
\leq\lambda_1+\lambda_2,
\qquad j=2,3.\qedhere
\]
\end{proof}

\begin{proposition}
\label{prop::F4_condition_K}
Let $L(\Lambda)$ be a relative holomorphic discrete series $F(4)$-supermodule for the real form $F(4,0)$ or $F(4,2)$, and suppose that $\Lambda$ is maximally atypical. Then Condition K holds. Moreover, $\DS_{x_A}(L(\Lambda))$ is a nonzero simple $\gg_{x_A}\cong\sl(3)$-module, finite-dimensional for $F(4,0)$ and a relative holomorphic discrete series module for $F(4,2)$.

\end{proposition}

\begin{proof}

We first verify Condition K(i). By Lemma~\ref{lemm::F4_odd_reflections_and_kernel}, there exist an isotropic root $\gamma$ and a nonzero weight vector $v_{\Lambda'}\in L(\Lambda)$ satisfying $(\Lambda',\gamma)=0$. For $A=\{\gamma\}$, we have $c_Av_{\Lambda'}=0$, and Proposition~\ref{prop::form_DS_x_generalized} gives
\begin{equation}
\DS_{x_A}(L(\Lambda))\cong L(\Lambda)^{c_A}\neq0.
\end{equation}
In particular, by Theorem~\ref{thm::DS_unitary} and Lemma~\ref{lemm::DS_is_HW}, $\DS_{x_{A}}(L(\Lambda))$ is a unitarizable simple highest weight $\gg_{x_{A}}$-module. The transitive action of $W_c$ on the sets $A\cup(-A)$ associated with maximal isotropic root sets shows that exactly one $W_c$-orbit contributes, proving Condition K(i).

For $F(4,0)$, $\gg_{x_A}^{\omega_{x_A}}\cong\su(3)$ is compact, which implies that $\DS_{x_A}(L(\Lambda))$ is finite-dimensional. For $F(4,2)$, by $W_{c}$-conjugation, we can take $A=\{\frac12(\epsilon_1-\epsilon_2-\epsilon_3-\delta)\}$. Note $\gg_{x_{A}}^{\omega_{x_{A}}}\cong\su(1,2)$, which has positive noncompact roots $\epsilon_1+\epsilon_2$ and $\epsilon_1+\epsilon_3$ and Weyl vector $\rho_{x_A}=\epsilon_1+\epsilon_2$. If $\nu$ is the highest weight of $\DS_{x_A}(L(\Lambda))\cong L(\Lambda)^{c_{A}}$, we have by Lemma~\ref{lemm::F4_weight_bound}:
\[
\begin{aligned}
(\nu+\rho_{x_A},\epsilon_1+\epsilon_2)
&=\nu_1+\nu_2+2
\leq\lambda_1+\lambda_2+2<0,\\
(\nu+\rho_{x_A},\epsilon_1+\epsilon_3)
&=\nu_1+\nu_3+1
\leq\lambda_1+\lambda_2+1<0,
\end{aligned}
\]
that is, $\DS_{x_{A}}(L(\Lambda))$ is a unitarizable highest weight $\gg_{x_{A}}$-supermodule satisfying the Harish-Chandra conditions. We conclude that it is a relative holomorphic discrete series module. 

In both cases, the $L$-packet supercharacter is holomorphic at zero, establishing Condition K(ii).
\end{proof}

Next, we verify Condition K for $G(3)$. We use the notation introduced in Section~\ref{subsec::G(3)}, fix the simple system
$
\Uppi=\{\epsilon_2,\epsilon_1-\epsilon_2,\delta-\epsilon_1-\epsilon_2\}
$
and set
\begin{equation}
\beta\coloneqq\delta-\epsilon_1-\epsilon_2,
\qquad
A\coloneqq\{\beta\}.
\end{equation}
Since $\defect(\gg)=1$, the set $A$ has maximal rank. The Weyl vector is
$\rho=3\epsilon_1+2\epsilon_2-\frac52\delta$.
Using $(\rho,\beta)=0$, we have for any $\Lambda\in \hh^{\ast}$: \begin{equation}
(\Lambda+\rho,\beta)=0
\quad\Longleftrightarrow\quad
(\Lambda,\beta)=-(\lambda_1+\lambda_2+2\mu)=0.
\end{equation}
We first show that $\beta$ is the unique solution of $(\Lambda+\rho,\alpha)=0$.

\begin{lemma}\label{lemm::G(3)_atyp}
    Let $L(\Lambda)$ be a relative holomorphic discrete series $G(3)$-supermodule with maximally atypical highest weight $\Lambda$. Then $(\Lambda+\rho,\beta)=0$.
\end{lemma}
\begin{proof}
We first show that $(\Lambda,\beta)\geq0$. Let $x_A$ be the normal form element associated with $A$ and set $c_A=x_A^2$. By unitarity,
\begin{equation*}
\langle v,ic_Av\rangle=\|x_Av\|^2\geq0
\end{equation*}
for every $v\in L(\Lambda)$. Since $c_A$ is a nonzero multiple of the Cartan element associated with $\beta$, there exists $t\in\RR^\times$ such that $ic_A$ acts on the weight space $L(\Lambda)^{\eta}$ by the scalar $t(\eta,\beta)$.

The noncompact even root is $2\delta$. Since the $\sl(2,\RR)$-module generated by the highest weight vector is an infinite-dimensional relative holomorphic discrete series module, every weight $\Lambda-2k\delta$, $k\geq0$, occurs in $L(\Lambda)$. Consequently,
\begin{equation*}
0\leq t(\Lambda-2k\delta,\beta)
=t\bigl((\Lambda,\beta)+4k\bigr)
\qquad k\geq0.
\end{equation*}
These inequalities hold for arbitrarily large $k$, thus $t>0$. Taking $k=0$ then gives $(\Lambda,\beta)\geq0$.

Next, let $\gamma\neq\beta$ be any positive isotropic root. By the choice of $\Uppi$, the difference $\gamma-\beta$ is a nonzero sum of positive compact roots. Since $\Lambda$ is dominant for $\gg_2$ and the compact component of $\rho$ is its Weyl vector, we have
\begin{equation*}
(\Lambda+\rho,\gamma-\beta)>0.
\end{equation*}
Using $(\rho,\beta)=0$, we obtain
\begin{equation*}
(\Lambda+\rho,\gamma)
=(\Lambda,\beta)+(\Lambda+\rho,\gamma-\beta)>0.
\end{equation*}
Thus $\beta$ is the only positive isotropic root that can satisfy the atypicality condition. Since $L(\Lambda)$ is maximally atypical and $\defect(\gg)=1$, it follows that
\begin{equation*}
(\Lambda+\rho,\beta)=(\Lambda,\beta)=0.\qedhere
\end{equation*}
\end{proof}

\begin{proposition}
\label{prop::condition_K_G3}
Let $L(\Lambda)$ be a relative holomorphic discrete series supermodule for $G(3)$ with maximally atypical highest weight. Then Condition K holds. Moreover, $\DS_{x_A}(L(\Lambda))$ is a finite-dimensional simple $\su(2)$-module.
\end{proposition}
\begin{proof}
By Lemma~\ref{lemm::G(3)_atyp}, we have $(\Lambda,\beta)=0$. Hence the highest weight vector $v_{\Lambda}$ of $L(\Lambda)$ belongs to $L(\Lambda)^{c_A}$, and Proposition~\ref{prop::form_DS_x_generalized} gives
\begin{equation*}
\DS_{x_A}(L(\Lambda))
\cong L(\Lambda)^{c_A}\neq0.
\end{equation*}

The isotropic roots of $G(3)$ are
\begin{equation*}
\pm\delta\pm\epsilon_1,\qquad
\pm\delta\pm\epsilon_2,\qquad
\pm\delta\pm(\epsilon_1+\epsilon_2).
\end{equation*}
Thus every pair $\{\gamma,-\gamma\}$ with $\gamma\in\overline{\Delta}_{\bar{1}}$ has a unique representative of the form $\delta+\eta$, where
\begin{equation*}
\eta\in
\{\pm\epsilon_1,\pm\epsilon_2,\pm(\epsilon_1+\epsilon_2)\}.
\end{equation*}

The Weyl group $W_c=W(\gg_2)\cong S_3\times\ZZ/2\ZZ$ fixes $\delta$ and acts transitively on the six choices of $\eta$. Hence exactly one orbit contributes, proving Condition K(i).

By Proposition~\ref{prop::unitary_real_forms_and_DS}, $\gg_{x_A}^{\omega_{x_A}}\cong\su(2)$. The DS cohomology $\DS_{x_A}(L(\Lambda))$ is therefore finite-dimensional, and its $L$-packet supercharacter is holomorphic at zero. This proves Condition K(ii).
\end{proof}

We finally verify Condition K for $D(2,1;\upalpha)$, using the notation of Section~\ref{subsec::appendix_D21alpha}. Recall that $\upalpha\in\RR\setminus\{0,-1\}$, and that the relevant real forms have even part $\sl(2,\RR)\oplus\su(2)\oplus\su(2)$ (see~\eqref{table::real_forms_exceptional}). We index the noncompact $\sl(2,\RR)$-factor first, so the positive noncompact even root is $2\epsilon_1$.

The Lie superalgebra $D(2,1;\upalpha)$ has defect one. For every isotropic root $\beta$ of $D(2,1;\upalpha)$, the only roots orthogonal to $\beta$ are $\pm\beta$. Hence, for $A=\{\beta\}$, the reduced root system is empty, and
$
\gg_{x_A}\cong\hh_{x_A}\cong\CC
$
is one-dimensional, purely even, and abelian.

\begin{proposition}
\label{prop::D21alpha_condition_K}
Let $L(\Lambda)$ be a relative holomorphic discrete series supermodule for $D(2,1;\upalpha)$ with maximally atypical highest weight $\Lambda$. Then Condition K holds. Moreover, $\DS_{x_A}(L(\Lambda))$ is one-dimensional.
\end{proposition}

\begin{proof}
Without loss of generality, we can assume that $L(\Lambda)$ is
a highest weight supermodule with respect to the simple system
\begin{equation*}
\Uppi
=
\{-\epsilon_1+\epsilon_2+\epsilon_3,\,
\epsilon_1+\epsilon_2-\epsilon_3,\,
\epsilon_1-\epsilon_2+\epsilon_3\},
\end{equation*}
whose Weyl vector is $\rho=0$.
Since $\Lambda$ is maximally atypical, there is an isotropic
root $\beta$ such that $(\Lambda,\beta)=0$.
Set $A=\{\beta\}$, and define $x_{A}$ and $c_{A}$. Then, if $v_{\Lambda}$ is the highest weight vector of $L(\Lambda)$, one has $c_Av_\Lambda=0$, and therefore
\begin{equation*}
0\neq v_\Lambda\in L(\Lambda)^{c_A}
\cong\DS_{x_A}(L(\Lambda)).
\end{equation*}
By Theorem~\ref{thm::DS_unitary}, the nonzero DS cohomology
$\DS_{x_A}(L(\Lambda))$ is a simple supermodule over the
one-dimensional purely even abelian Lie algebra $\gg_{x_A}$.
It is therefore one-dimensional and either even or odd. This establishes in particular Condition~K(ii).

The isotropic odd roots are $\pm\epsilon_1\pm\epsilon_2\pm\epsilon_3$. The compact Weyl group is
$
W_c=\langle s_{2\epsilon_2},s_{2\epsilon_3}\rangle\cong(\ZZ/2\ZZ)^2
$, and hence acts transitively on the four pairs $\{\beta,-\beta\}$. This shows that exactly one orbit contributes, proving Condition~K(i).
\end{proof}

\subsubsection{Proof of the main theorem: nonintegral case}
\label{subsec::proof_main_theorem}
Let $L(\Lambda)$ be a relative holomorphic discrete series $\gg$-supermodule with nonintegral and maximally atypical highest weight. The preceding results establish Condition K and reduce the evaluation of $\Theta_\Lambda(0)$ to the $L$-packet supercharacter of a single nonzero maximal-rank generalized DS cohomology by Corollary~\ref{cor::toral_nonvanishing}.

Let $A$ represent the unique contributing compact orbit and set
\begin{equation}
M\coloneqq\DS_{x_A}(L(\Lambda))
\cong L(\Lambda)^{c_A},
\qquad c_A=x_A^2.
\end{equation}
By Theorem~\ref{thm::DS_unitary}, $M$ is nonzero, simple, and unitarizable. The following table summarizes all cases:
\begin{equation}\label{table::form_M}
\begin{array}{c|c|c}
\text{case}&\gg_{x_{A}}^{\omega_{x_{A}}}&M\text{ (up to parity)}\\ \hline
\sl(m\vert n),\ m>n
&\begin{gathered}\su(p-n,q)\text{ or}\\ \su(p,q-n)\end{gathered}&\begin{gathered}\text{finite-dimensional or}\\ \text{relative holomorphic discrete series}\end{gathered}\\
\osp(m\vert2n),\ \lfloor m/2\rfloor\geq n
&\so(m-2n,\RR)&\text{finite-dimensional}\\
\osp(2r\vert2n),\ r<n
&\sp(2(n-r),\RR)&\text{relative holomorphic discrete series}\\
\osp(2r+1\vert2n),\ r<n
&\osp(1\vert2(n-r))&\text{relative holomorphic discrete series}\\
F(4,0)&\su(3)&\text{finite-dimensional}\\
F(4,2)&\su(1,2)&\text{relative holomorphic discrete series}\\
G(3)&\su(2)&\text{finite-dimensional}\\
D(2,1;\upalpha)&\mathfrak u(1)&\text{one-dimensional}
\end{array}
\end{equation}

\begin{lemma}
\label{lemm::reduced_packet_nonvanishing}
Let $M$ be any of the nonzero generalized DS cohomologies of~\eqref{table::form_M}. Then $\Theta_M$ is holomorphic at zero and $\Theta_M(0)\neq0$.
\end{lemma}

\begin{proof}
By Theorem~\ref{thm::L_packet_and_sdim}, $\Theta_M$ is holomorphic at zero and
\begin{equation*}
\Theta_M(0)
=
(-1)^{|\Delta_{x_A,n}^+|}\fsdim M.
\end{equation*}
If $\gg_{x_A}$ is purely even, then $M$ is concentrated in one parity $p$, and $\fsdim M=(-1)^p d(M)\neq0$, where $d(M)$ is its positive formal dimension. For compact real forms with normalized Haar measure, $d(M)=\dim M$, so that $\fsdim M=\sdim M$. The remaining case $\gg_{x_A}\cong\osp(1\vert2s)$ follows from Example~\ref{ex::osp_one_formal_superdimension}, which gives $\fsdim M\neq0$.
\end{proof}

We can now prove the nonintegral case.

\begin{theorem}
\label{thm::formal_Kac_Wakimoto}
Let $L(\Lambda)$ be a relative holomorphic discrete series
supermodule with nonintegral highest weight $\Lambda$. Then
\begin{equation*}
\fsdim L(\Lambda)\neq0
\quad\Longleftrightarrow\quad
\at(\Lambda)=\defect(\gg).
\end{equation*}
\end{theorem}

\begin{proof}
The vanishing of $\fsdim L(\Lambda)$ below maximal atypicality
has already been established in Theorem~\ref{thm::vanishing_sdim_less_atypic}. Suppose that $\Lambda$ is maximally atypical.

For $\gg=\osp(1\vert2n)$, the assertion follows from
Example~\ref{ex::osp_one_formal_superdimension}.
In all remaining cases, Condition K holds by the preceding
propositions.

Let $A$ represent the unique contributing compact orbit
and set $M=\DS_{x_A}(L(\Lambda))$.
By Corollary~\ref{cor::toral_nonvanishing},
\begin{equation*}
\Theta_\Lambda(0)
=
\frac{[N_A:N_A\cap W_c]}
     {[W_{x_A}:W_{x_A,c}]}
\Theta_M(0).
\end{equation*}
The coefficient is positive, and
Lemma~\ref{lemm::reduced_packet_nonvanishing}
gives $\Theta_M(0)\neq0$.
Thus $\Theta_\Lambda(0)\neq0$.
Finally, Theorem~\ref{thm::L_packet_and_sdim} yields
\begin{equation*}
\fsdim L(\Lambda)
=
(-1)^{|\Delta_n^+|}\Theta_\Lambda(0)
\neq0.\qedhere
\end{equation*}
\end{proof}

\subsection{The Kac--Wakimoto theorem for the formal superdimension}
\label{subsec::KW_main}

We combine the integral and nonintegral cases.

\begin{theorem}\label{thm::KW_formal_superdimension}
Let $L(\Lambda)$ be a relative holomorphic discrete series $\gg$-supermodule. Then
\begin{equation*}
\fsdim L(\Lambda)\neq0
\quad\Longleftrightarrow\quad
\at(\Lambda)=\defect(\gg).
\end{equation*}
Moreover, if $\Lambda$ is integral, then $R^{|\Delta_n^+|}\Gamma(L(\Lambda))$ is a finite-dimensional simple $\gg$-supermodule with the same infinitesimal character as $L(\Lambda)$, and
\begin{equation*}
\fsdim L(\Lambda)=\sdim R^{|\Delta_n^+|}\Gamma(L(\Lambda)).
\end{equation*}
\end{theorem}

\begin{proof}
This is Theorem~\ref{thm::formal_KW_integral} for integral $\Lambda$ and Theorem~\ref{thm::formal_Kac_Wakimoto} otherwise.
\end{proof}

\section{The cases \texorpdfstring{$\sl(n\vert n)$ and $\psl(n\vert n)$}{sl(n|n) and psl(n|n)}}
\label{sec::equal_rank_type_A}
\noindent
We finally treat the Lie superalgebras $\sl(n\vert n)$ and $\psl(n\vert n)$, which have been excluded so far. We follow the conventions of Section~\ref{subsec::type_A_structure}. Let $n\geq2$, and let $I$ denote the identity matrix. Then $\gl(n\vert n)\supseteq\sl(n\vert n)\supseteq\CC I$ and $\psl(n\vert n)=\sl(n\vert n)/\CC I$. Although $\sl(n\vert n)$ is not simple, we include it in our analysis.

Two difficulties require a separate treatment. First, the supertrace form on $\sl(n\vert n)$ has radical $\CC I$, that is, it is degenerate. Second, $\tau\coloneqq\sum_{i=1}^n\epsilon_i-\sum_{a=1}^n\delta_a$ vanishes on the Cartan subalgebra of diagonal matrices. To treat both problems simultaneously, we work with $\gl(n\vert n)$, where the supertrace form is nondegenerate and the roots $\{\epsilon_{i}-\delta_{i} : 1\leq i \leq n\}$ are linearly independent. 

Our analysis uses the decomposition
\begin{equation}
\gl(n\vert n)=\sl(n\vert n)\oplus\CC E,\qquad E\coloneqq\diag(I_n,0).
\end{equation}
The diagonal Cartan subalgebra $\dd$ of $\gl(n\vert n)$ decomposes then as $\dd=\hh\oplus\CC E$, and both Lie superalgebras have the same structure theory. In particular, the restriction from $\dd$ to $\hh$ is surjective with kernel $\CC\tau$. Every weight $\Lambda\in\hh^\ast$ admits an extension $\widetilde\Lambda\in\dd^\ast$ satisfying $\widetilde\Lambda|_\hh=\Lambda$. Any two extensions differ by a multiple of $\tau$, which is orthogonal to all roots. In what follows, we fix such an extension. It is important to note that the atypicality of $\widetilde{\Lambda}$ is independent of the extension, since $(\tau,\alpha)=0$ for all $\alpha \in \Delta_{\bar{1}}$.

The only real forms of $\sl(n\vert n)$ admitting nontrivial unitarizable supermodules are $\su(p,q\vert 0,n)$ and $\su(p,q\vert n,0)$~\cite{Schmidt,neeb2011lie}, where $p+q=n$. The former admit only nontrivial highest weight supermodules, and the latter only nontrivial lowest weight supermodules. We therefore restrict to the real forms $\su(p,q\vert 0,n)$. Here, it suffices to consider $\sl(n\vert n)$, since every unitarizable highest weight $\psl(n\vert n)$-supermodule is also an $\sl(n\vert n)$-supermodule on which $I$ acts trivially.

In Section~\ref{subsec::real_forms}, we described these real forms via the conjugate-linear anti-involution
\begin{equation}\omega_H(X)=H^{-1}X^\dagger H,\qquad H=\operatorname{diag}(I_{p,q},-I_n).\end{equation} 
The same formula defines a conjugate-linear anti-involution of $\gl(n\vert n)$.

All preceding definitions from Section~\ref{sec::relative_HDS_supermodules} to~\ref{sec::formal_superdimension} extend to $\sl(n\vert n)$. All results of this manuscript carry over verbatim whenever their proofs use only the even Lie subalgebra $\even$, unitarity, or the Harish-Chandra theory of the $\even$-constituents. 

Only three points must be proved again. First, we must work over $\gl(n\vert n)$, where atypicality is defined using the nondegenerate supertrace form. This is guaranteed by the following lemma. 

\begin{lemma}\label{lemm::equal_rank_extension}
Let $L_{\sl}(\Lambda)$ be a simple $\omega_{H}$-unitarizable highest weight
$\sl(n\vert n)$-supermodule. Then there exists an $\omega_{H}$-unitarizable highest
weight $\gl(n\vert n)$-supermodule $L_{\gl}(\widetilde\Lambda)$ such
that $\widetilde\Lambda$ extends $\Lambda$ and
\[
L_{\sl}(\Lambda)
\cong
L_{\gl}(\widetilde\Lambda)\vert_{\sl(n\vert n)}.
\]
\end{lemma}

\begin{proof}
Choose an extension $\widetilde\Lambda$ of $\Lambda$ to $\dd$ such that $\widetilde\Lambda(E)\in\RR$. Using $\gl(n\vert n)=\sl(n\vert n)\oplus\CC E$, the positive and negative root spaces of both Lie superalgebras coincide for any choice of positive system. Let $M_{\gl}(\widetilde{\Lambda})$ and $M_{\sl}(\Lambda)$ denote the associated Verma supermodules. The PBW theorem then gives an isomorphism of $\sl(n\vert n)$-supermodules
\[
M_{\gl}(\widetilde\Lambda)\vert_{\sl(n\vert n)}\cong M_{\sl}(\Lambda).
\]
Indeed, both Verma supermodules are generated by their highest weight vectors under the same negative root spaces, and their highest weights agree on $\hh$.

Normalize the positive definite Hermitian form $\langle \cdot,\cdot\rangle_{L_{\sl}}$ on $L_{\sl}(\Lambda)$, such that the highest weight vector has norm one. Let $\pi:M_{\sl}(\Lambda)\to L_{\sl}(\Lambda)$ be the canonical surjection. By uniqueness of normalized $\omega_{H}$-contravariant Hermitian forms on Verma supermodules, the pullback of $\langle \cdot,\cdot\rangle_{L_{\sl}}$ under $\pi$ and the canonical normalized $\omega_{H}$-contravariant Shapovalov form on $M_{\sl}(\Lambda)$ coincide.

Since $\omega_H(E)=E$ and $\widetilde\Lambda(E)\in\RR$, the Verma supermodule $M_{\gl}(\widetilde\Lambda)$ admits a normalized $\omega_H$-contravariant Hermitian Shapovalov form. Under the isomorphism $M_{\gl}(\widetilde\Lambda)|_{\sl(n\vert n)}\cong M_{\sl}(\Lambda)$, this form coincides, by uniqueness, with the normalized form induced by $\langle\cdot,\cdot\rangle_{L_{\sl}}$.

The radical of this form is $\ker \pi$, since $L_{\sl}(\Lambda)$ is unitarizable. As the radical is also $\gl(n\vert n)$-invariant, $M_{\gl}(\widetilde{\Lambda})/\ker \pi$ is a highest weight $\gl(n\vert n)$-supermodule with positive definite $\omega_{H}$-contravariant Hermitian form. By construction, its restriction to $\sl(n\vert n)$ is precisely $L_{\sl}(\Lambda)$, and hence, $L_{\gl}(\widetilde{\Lambda})$ is simple.  
\end{proof}

Second, our treatment of the 
generalized DS cohomology for unitarizable supermodules used simplicity of the Lie superalgebra in two places: the
pointed-cone argument of Lemma~\ref{lemm::conjugation_inner}, and the
vanishing of real odd elements of square zero in
Lemma~\ref{lemm::omega_fixed_cartan_normal_form}. We replace both arguments by a direct matrix computation.

\begin{lemma}\label{lemm::equal_rank_real_odd}
Let $x\in\gl(n\vert n)_{\bar1}$ satisfy $\omega_H(x)=ix$. Set $c=x^2$.
\begin{enumerate}
\item[(a)] If $c=0$, then $x=0$.
\item[(b)] If $\ad_c$ is semisimple, then an element of
$U(p,q)\times U(n)$ conjugates $c$ into the diagonal compactly embedded
Cartan subalgebra.
\end{enumerate}
\end{lemma}

\begin{proof}
We represent $x$ in its standard matrix form~\eqref{eq::matrix_form}. Solving $\omega_H(x)=ix$ gives
\[
x=\begin{pmatrix}0&B\\iB^\dagger I_{p,q}&0\end{pmatrix},
\qquad
ic=\operatorname{diag}(A,D),
\qquad
A=-BB^\dagger I_{p,q},\quad
D=-B^\dagger I_{p,q}B.
\]

 We first prove (a). 
If $c=0$, then $A=D=0$. Since $A=-BB^\dagger I_{p,q}$ and $I_{p,q}$ is invertible, we obtain $BB^\dagger=0$, which implies $B=0$.

Next, we prove (b). By assumption, $\ad_{c}$ acts semisimply on $\gg$, that is, $c$ is semisimple. Using $\omega_{H}(ic)=ic$ and $ic=\diag(A,D)$, the matrix $D$ is Hermitian, and $A$ is
self-adjoint with respect to $I_{p,q}$. We show that we can choose an $I_{p,q}$-orthonormal eigenbasis of $A$ and an orthonormal eigenbasis for $D$. The corresponding element in $U(p,q)\times U(n)$ conjugates then $c$ into the diagonal compactly embedded Cartan subalgebra.

Since $D$ is Hermitian, the spectral theorem implies the existence of an orthonormal eigenbasis of $D$ and $D$ can be diagonalized by $U(n)$. 

We consider $A$. The statement follows if we show that $I_{p,q}$ is definite on each eigenspace of $A$. In this case, we find an $I_{p,q}$-orthonormal eigenbasis of $A$, and $A$ can be diagonalized by $U(p,q)$. Suppose $Av=\lambda v$ with $v\neq0$ and
$\lambda\neq0$. Then $B^\dagger I_{p,q}v\neq0$ and
\[
\lambda\,v^\dagger I_{p,q}v=-\langle B^\dagger I_{p,q}v,B^\dagger I_{p,q}v\rangle_{\CC^{n}},
\]
where $\langle\cdot,\cdot\rangle_{\CC^{n}}$ is the standard Hermitian form on $\CC^{n}$. We conclude that $\lambda$ is real, and every nonzero eigenspace of $A$ is definite for
$I_{p,q}$. This finishes the proof.
\end{proof}

\begin{remark}
    The Lie group $U(p,q)\times U(n)$ is connected and the conjugation in (b) is a real even
inner automorphism commuting with $\omega_H$.
\end{remark}

Finally, the reduced Lie superalgebras $\sl(n\vert n)_x$ and $\gl(n\vert n)_x$ differ. Thus simplicity of generalized DS cohomology over $\gl(n\vert n)_x$ does not automatically imply simplicity over $\sl(n\vert n)_x$ or $\psl(n\vert n)_x$. We establish this separately in the following proposition.

\begin{proposition}\label{prop::equal_rank_DS}
Let $L_{\sl}(\Lambda)$ be an $\omega_{H}$-unitarizable highest weight $\sl(n\vert n)$-supermodule. Let $x\in\sl(n\vert n)_{\bar1}$, and let $k$ denote the rank of $x$ in $\gl(n\vert n)$. Assume $\omega_H(x)=ix$ and $c=x^2$ acts semisimply on $L_{\sl}(\Lambda)$.
\begin{enumerate}
\item[(a)] For an extension $\widetilde\Lambda$ as in Lemma~\ref{lemm::equal_rank_extension},
\[
\DS_x(L_{\sl}(\Lambda))=H(L_{\gl}(\widetilde\Lambda)^c,x)\cong L_{\sl}(\Lambda)^c.
\]
If $k>\at(\widetilde\Lambda)$, then $\DS_x(L_{\sl}(\Lambda))=0$.

\item[(b)] The reduced Lie superalgebras are
\[
(\sl(n\vert n)_x,\psl(n\vert n)_x)\cong
\begin{cases}
(\sl(n-k\vert n-k),\psl(n-k\vert n-k)),&k<n,\\
(\CC^{0\vert1},\CC^{0\vert2}),&k=n.
\end{cases}
\]

\item[(c)] The cohomology $\DS_x(L_{\sl}(\Lambda))$ is either zero or a simple unitarizable $\sl(n\vert n)_x$-supermodule. If $I$ acts trivially on $L_{\sl}(\Lambda)$, the same holds over $\psl(n\vert n)_x$.
\end{enumerate}
\end{proposition}
\begin{proof}

For (a), by Proposition~\ref{prop::form_DS_x_generalized}, one has $\DS_x(L_{\sl}(\Lambda))\cong L_{\sl}(\Lambda)^c$, since its proof uses only unitarity. By Lemma~\ref{lemm::equal_rank_extension}, $x$ and $c$ act identically on $L_{\sl}(\Lambda)$ and $L_{\gl}(\widetilde\Lambda)$, so their generalized DS cohomologies coincide as superspaces. The vanishing for $k>\at(\widetilde\Lambda)$ then follows from Theorem~\ref{thm::generalized_DS_character} applied to $\gl(n\vert n)$.

For (b), we determine the reduced Lie superalgebras of $\sl(n\vert n)$ and $\psl(n\vert n)$ from the known isomorphism $\gl(n\vert n)_x\cong\gl(n-k\vert n-k)$. We consider the following two short exact sequences:
\[
0\to\sl(n\vert n)\to\gl(n\vert n)\xrightarrow{\operatorname{str}}\CC\to0,
\qquad
0\to\CC I\to\sl(n\vert n)\to\psl(n\vert n)\to0.
\]
Since $\ad_c$ is semisimple, these sequences remain exact after taking $c$-kernels. On these kernels, $\ad_x$ is a differential, and it acts trivially on $\CC$ and $\CC I$. We can therefore apply DS cohomology and obtain two long exact sequences. Concretely, the first long exact sequence is
\[
0\to(\sl(n\vert n)_x)_{\bar0}\to(\gl(n\vert n)_x)_{\bar0}\xrightarrow{\operatorname{str}}\CC\xrightarrow{\delta}(\sl(n\vert n)_x)_{\bar1}\to(\gl(n\vert n)_x)_{\bar1}\to0,
\]
and the second is
\[
\CC I\to(\sl(n\vert n)_x)_{\bar0}\to(\psl(n\vert n)_x)_{\bar0}\to0\to(\sl(n\vert n)_x)_{\bar1}\to(\psl(n\vert n)_x)_{\bar1}\xrightarrow{\delta'}\CC I\to(\sl(n\vert n)_x)_{\bar0}.
\]
Here $\delta(a)$ is the cohomology class of $[x,T]$, where $T\in\gl(n\vert n)^c_{\bar0}$ satisfies $\operatorname{str}(T)=a$. Similarly, $\delta'$ sends the cohomology class of $Y+\CC I\in\psl(n\vert n)^c_{\bar1}$ to $[x,Y]\in\CC I$, where $Y\in\sl(n\vert n)^c_{\bar1}$ and $[x,Y]\in\CC I$.

Assume $k<n$ and, without loss of generality, that $x$ is in normal form. Then $E_{k+1,k+1}$ commutes with $x$ and has supertrace one. Its cohomology class belongs to $(\gl(n\vert n)_x)_{\bar0}$ and maps to $1\in\CC$. Thus the supertrace map in the first exact sequence is surjective, so $\delta=0$. By exactness, $(\sl(n\vert n)_x)_{\bar1}\cong(\gl(n\vert n)_x)_{\bar1}$ and $(\sl(n\vert n)_x)_{\bar0}$ is the kernel of supertrace. This shows $\sl(n\vert n)_{x}\cong \sl(n-k\vert n-k)$. Moreover, the class of $I$ is the identity of $\gl(n-k\vert n-k)$ and is therefore nonzero. Exactness of the second sequence gives $\psl(n\vert n)_x\cong\psl(n-k\vert n-k)$.

It remains to consider the case $k=n$. Then $\gl(n\vert n)_x=\{0\}$. Using again exactness of the first sequence shows $(\sl(n\vert n)_x)_{\bar{0}}=\{0\}$ and $(\sl(n\vert n)_x)_{\bar{1}}\cong \CC$. Similarly, using exactness of the second long exact sequence, we obtain $(\psl(n\vert n)_x)_{\bar{0}}=\{0\}$ and $(\psl(n\vert n)_x)_{\bar{1}}\cong \CC^{2}$.

Finally, we prove (c). As in 
Proposition~\ref{prop::form_DS_x_generalized}, one can show that $\DS_{x}(L_{\sl}(\Lambda))$ is unitarizable. Applying Theorem~\ref{thm::DS_unitary} and
Lemma~\ref{lemm::DS_is_HW} to $\gl(n\vert n)$, we conclude that $\DS_x(L_{\gl}(\widetilde\Lambda))$ is either zero or a
simple highest weight $\gl(n\vert n)_x$-supermodule. If $k=n$, then $\gl(n\vert n)_x=0$, so $\DS_x(L_{\gl}(\widetilde\Lambda))$ is at most
one-dimensional and hence simple. If $k<n$, we have seen that $\sl(n\vert n)_{x}\subset \gl(n\vert n)_{x}$, and hence the statement follows by Lemma~\ref{lemm::equal_rank_extension}.
\end{proof}

We now obtain the Kac--Wakimoto theorem for $\sl(n\vert n)$ and $\psl(n\vert n)$.

\begin{theorem}\label{thm::equal_rank_formal_KW}
Let $L_{\sl}(\Lambda)$ be a relative holomorphic discrete series
supermodule over $\sl(n\vert n)$.\footnote{Recall that this includes supermodules over
$\psl(n\vert n)$, regarded as $\sl(n\vert n)$-supermodules on which $I$
acts trivially.} Let $\widetilde\Lambda$ be an extension of $\Lambda$ as
in Lemma~\ref{lemm::equal_rank_extension}. Then
\[
\fsdim L_{\sl}(\Lambda)\neq0
\quad\Longleftrightarrow\quad
\at(\widetilde\Lambda)=n.
\]
In this case, $\Lambda(I)=0$ and every relative holomorphic discrete
series supermodule over $\sl(n\vert n)$ with nonzero formal superdimension
factors through $\psl(n\vert n)$.
\end{theorem}

\begin{proof}
Let $L_{\gl}(\widetilde{\Lambda})$ be the extension of $L_{\sl}(\Lambda)$ of Lemma~\ref{lemm::equal_rank_extension}. Note that the even constituents of
$L_{\gl}(\widetilde\Lambda)$ restrict to those of $L_{\sl}(\Lambda)$,
with the same parities. Since the Harish-Chandra condition and the formal
dimensions of Theorem~\ref{thm::formal_dimension} depend only on pairings
with even roots, which are unchanged, $L_{\gl}(\widetilde\Lambda)$
is also a relative holomorphic discrete series supermodule, and
$\fsdim L_{\sl}(\Lambda)=\fsdim L_{\gl}(\widetilde\Lambda)$. It therefore suffices to consider $L_{\gl}(\widetilde\Lambda)$.

If $\at(\widetilde\Lambda)<n$, the proof of
Theorem~\ref{thm::vanishing_sdim_less_atypic} applies to
$L_{\gl}(\widetilde\Lambda)$, using Proposition~\ref{prop::equal_rank_DS}(a)
for the vanishing of generalized Duflo--Serganova cohomology of maximal
rank. It remains to consider the maximally atypical case $\at(\widetilde\Lambda)=n$. In this case, there exists $\sigma\in S_n$ with
\[
(\widetilde\Lambda+\rho,\epsilon_i-\delta_{\sigma(i)})=0,
\qquad 1\leq i\leq n.
\]
Together with integrality with respect to $\Delta_{c}$, these equations show that
$\widetilde\Lambda$ is integral (\emph{cf.}~proof of Lemma~\ref{lemm::non_integral_slmn_weight}). We are therefore in the integral case of
Section~\ref{subsec::KW_conjecture_integral_case}, and Serganova's
theorem for $\gl(n\vert n)$~\cite{serganova2011superdimension} yields
$\fsdim L_{\sl}(\Lambda)\neq0$.

Finally, we prove $\Lambda(I)=0$. Recall that $\tau=\sum_{i=1}^n\epsilon_i-\sum_{a=1}^n\delta_a$ is orthogonal to every odd root. Summing the preceding equations over $i$ gives $(\widetilde\Lambda+\rho,\tau)=0$, and we obtain 
\[
0=(\widetilde{\Lambda}+\rho,\tau)=(\Lambda,\tau)=\Lambda(I).\qedhere
\]
\end{proof}

\appendix
\section{Structure theory of basic classical Lie superalgebras}\label{appendix::structure_theory}
\noindent We briefly recall the structure of the orthosymplectic and exceptional basic classical Lie superalgebras needed in Section~\ref{sec::formal_superdimension}. We follow the notation of Section~\ref{sec::basic_classical_Lie_superalgebras}.

\subsection{Orthosymplectic Lie superalgebras} \label{subsec::appendix_osp}

 Let $m,n\geq1$, and set $r\coloneqq\lfloor m/2\rfloor$. We consider the complex super vector space $V=\CC^{m\vert 2n}$, equipped with a  nondegenerate even
supersymmetric bilinear form $b$. The orthosymplectic Lie superalgebra $\gg=\osp(m\vert2n)$ is the Lie subsuperalgebra of $\gl(V)$ whose homogeneous elements $X$ satisfy
\begin{equation*}
b(Xv,w)+(-1)^{p(X)p(v)}b(v,Xw)=0
\end{equation*}
for all homogeneous $v,w\in V$. The definition is independent of the choice of $b$, since any two such forms are equivalent under a basis change.

We give an explicit matrix realization. Let $I_n$ denote the $n\times n$ identity matrix, and set
$
\Omega_n\coloneqq
\left(\begin{smallmatrix}
0&I_n\\
-I_n&0
\end{smallmatrix}\right).
$ 
Then $\operatorname{diag}(I_m,\Omega_n)$ defines a nondegenerate even supersymmetric bilinear form on $\CC^{m\vert2n}$. With respect to this form, we obtain the following matrix realization of $\osp(m\vert2n)$:
\begin{equation}
\gg=
\left\{
\begin{pmatrix}A&B\\C&D\end{pmatrix}\in\sl(m\vert 2n):
\begin{pmatrix}
-A^T&-C^T\Omega_n\\
-\Omega_nB^T&\Omega_nD^T\Omega_n
\end{pmatrix}
=
\begin{pmatrix}A&B\\C&D\end{pmatrix}
\right\}.
\end{equation}
Since this realization identifies $\gg=\osp(m\vert2n)$ with a Lie subsuperalgebra of $\sl(m\vert2n)$, we equip $\gg$ with the restriction of the supertrace form $(\cdot,\cdot)$.

The even part is
$
\even\cong\so(m)\oplus\sp(2n)
$,
and the odd part is the $\even$-module
$
\odd\cong\CC^m\otimes\CC^{2n}.
$
This module is simple when $m\neq2$. If $m=2$,
the standard $\so(2)$-module splits into two one-dimensional
weight spaces. Consequently, $\odd$ splits into two
irreducible $\even$-modules.

We consider the two types of real forms from
Section~\ref{subsec::real_forms}. The following table gives
their even parts, maximal compact subalgebras, and
complexifications:
\begin{equation}
\begin{array}{c|c|c}
\gg_{\bar0,\RR}
&\kk_{\RR}
&\kk=\kk_{\RR}\otimes_{\RR}\CC\\ \hline
\so(m,\RR)\oplus\sp(2n,\RR)
&\so(m,\RR)\oplus\mathfrak u(n)
&\so(m)\oplus\gl(n)\\
\so^*(2r)\oplus\mathfrak{usp}(2n)
&\mathfrak u(r)\oplus\mathfrak{usp}(2n)
&\gl(r)\oplus\sp(2n)
\end{array}
\end{equation}
The second row applies when $m=2r$.
Here $\so(m,\RR)$ and $\mathfrak{usp}(2n)$ are compact.

In the first case, the map $X\mapsto-X^T$ defines a
Cartan involution on the even real algebra. Its
fixed-point algebra is $\kk_{\RR}$. In particular,
\begin{equation}
\begin{aligned}
\sp(2n,\RR)\cap\so(2n,\RR)
=
\left\{
\begin{pmatrix}P&Q\\-Q&P\end{pmatrix}:
P,Q\in M_n(\RR),\
P^T=-P,\ Q^T=Q
\right\}\cong\mathfrak u(n).
\end{aligned}
\end{equation}

Every $H\in\hh$ acts semisimply on $\CC^m$ and on $\CC^{2n}$. Define $\epsilon_1,\ldots,\epsilon_r,\delta_1,\ldots,\delta_n\in\hh^*$ such that, for every $H\in\hh$, the eigenvalues of $H$ on $\CC^m$ are $\pm\epsilon_i(H)$, $1\leq i\leq r$, together with $0$ if $m$ is odd, and the eigenvalues of $H$ on $\CC^{2n}$ are $\pm\delta_a(H)$, $1\leq a\leq n$.
These elements span $\hh^{\ast}$. The standard $\so(m)$-module has weights
$\pm\epsilon_i$. It also has the zero weight
when $m$ is odd. The standard $\sp(2n)$-module
has weights $\pm\delta_a$. We write
\begin{equation}
\lambda
=
\sum_{i=1}^r\lambda^i\epsilon_i
+\sum_{a=1}^n\mu^a\delta_a
=
(\lambda^1,\ldots,\lambda^r
\vert\mu^1,\ldots,\mu^n).
\end{equation}

For $m=2r$, the root system is
\begin{equation}
\begin{aligned}
\Delta_{\bar0}
&=
\{\pm\delta_a\pm\delta_b:1\leq a<b\leq n\}
\cup\{\pm2\delta_a:1\leq a\leq n\}\cup
\{\pm\epsilon_i\pm\epsilon_j:1\leq i<j\leq r\},\\
\Delta_{\bar1}
&=
\{\pm\epsilon_i\pm\delta_a:
1\leq i\leq r,\ 1\leq a\leq n\}.
\end{aligned}
\label{eq:osp-roots}
\end{equation}
The signs in these sets are independent. For $m=2r+1$, the even roots also include
$\pm\epsilon_i$, and the odd roots also include
$\pm\delta_a$.

The bilinear form $(\cdot,\cdot)$ on $\hh$ induces a bilinear form on $\hh^\ast$, which we denote by the same symbol. It is uniquely determined by
\begin{equation}
(\epsilon_i,\epsilon_j)=\frac12\delta_{ij},
\qquad
(\delta_a,\delta_b)=-\frac12\delta_{ab},
\qquad
(\epsilon_i,\delta_a)=0.
\end{equation}
The odd roots $\epsilon_i\pm\delta_a$ are isotropic. If $m$ is odd, the additional odd roots $\pm\delta_a$ are non-isotropic.

We fix the following positive even roots:
\begin{equation}
\begin{aligned}
\Delta_{\bar0}^+
=
\{\delta_a-\delta_b,\delta_a+\delta_b:
1\leq a<b\leq n\}\cup
\{2\delta_a:1\leq a\leq n\}\cup
\{\epsilon_i-\epsilon_j,\epsilon_i+\epsilon_j:
1\leq i<j\leq r\}.
\end{aligned}
\label{eq:osp-positive-even-roots}
\end{equation}
When $m$ is odd, we also include the roots
$\epsilon_i$. This choice is $\kk$-adapted for both real forms considered above. In particular, the set of positive noncompact even roots $\Delta_n^+$ is invariant under the compact Weyl group $W_c$.

Consider first the real form with even part
$\so(m,\RR)\oplus\sp(2n,\RR)$. Its compact
and noncompact positive even roots are
\begin{equation}
\begin{aligned}
\Delta_c^+
&=
\{\delta_a-\delta_b:1\leq a<b\leq n\}\cup
\{\epsilon_i-\epsilon_j,\epsilon_i+\epsilon_j:
1\leq i<j\leq r\},\\
\Delta_n^+
&=
\{\delta_a+\delta_b:1\leq a<b\leq n\}
\cup\{2\delta_a:1\leq a\leq n\}.
\end{aligned}
\end{equation}
When $m$ is odd, the roots $\epsilon_i$ also
belong to $\Delta_c^+$.

For the real form with even part
$\so^*(2r)\oplus\mathfrak{usp}(2n)$,
the corresponding sets are
\begin{equation}
\begin{aligned}
\Delta_c^+
&=
\{\epsilon_i-\epsilon_j:1\leq i<j\leq r\}\cup
\{\delta_a-\delta_b,\delta_a+\delta_b:
1\leq a<b\leq n\}\cup
\{2\delta_a:1\leq a\leq n\},\\
\Delta_n^+
&=
\{\epsilon_i+\epsilon_j:1\leq i<j\leq r\}.
\end{aligned}
\end{equation}

Among the positive systems of $\Delta$, we single out the system called Case~A in Jakobsen's classification~\cite{jakobsen1994full}. It is particularly important for the real forms $\osp(0,m\vert 2n;\RR)$ (see~\cite{Schmidt_unitary_branching}) and will be used in our proof of nonvanishing of the formal superdimension.

Its positive odd roots are
\begin{equation}\label{eq::type_A_plus_1}
\Delta_{A,\bar1}^+
=
\begin{cases}
\{\delta_a\pm\epsilon_i:
1\leq a\leq n,\ 1\leq i\leq r\},
&m=2r,\\
\{\delta_a\pm\epsilon_i:
1\leq a\leq n,\ 1\leq i\leq r\}
\cup\{\delta_a:1\leq a\leq n\},
&m=2r+1.
\end{cases}
\end{equation}
Set
$
\Delta_A^+
\coloneqq
\Delta_{\bar0}^+\sqcup\Delta_{A,\bar1}^+
$
and denote the corresponding Borel subalgebra
by $\bb_A$.

For $\osp(2r+1\vert2n)$ with $r\geq1$,
the simple roots are
\begin{equation}\label{eq::positive_system_type_A_2}
\begin{aligned}
\Uppi_A
=
\{\delta_a-\delta_{a+1}:1\leq a<n\}\cup\{\delta_n-\epsilon_1\}\cup
\{\epsilon_i-\epsilon_{i+1}:1\leq i<r\}\cup\{\epsilon_r\}.
\end{aligned}
\end{equation}

For $\osp(2r\vert2n)$ with $r\geq2$,
the simple roots are
\begin{equation}\label{eq::positive_system_typeA_3}
\begin{aligned}
\Uppi_A
=
\{\delta_a-\delta_{a+1}:1\leq a<n\}\cup\{\delta_n-\epsilon_1\}\cup
\{\epsilon_i-\epsilon_{i+1}:1\leq i<r\}\cup\{\epsilon_{r-1}+\epsilon_r\}.
\end{aligned}
\end{equation}

For $m=2$ and $m=1$, the simple roots are
\begin{equation}
\begin{aligned}
\Uppi_A
&=
\{\delta_a-\delta_{a+1}:1\leq a<n\}
\cup\{\delta_n-\epsilon_1,\delta_n+\epsilon_1\},
&& (m=2),\\
\Uppi_A
&=
\{\delta_a-\delta_{a+1}:1\leq a<n\}
\cup\{\delta_n\},
&& (m=1).
\end{aligned}
\end{equation}

\subsection{\texorpdfstring{$F(4)$}{F(4)}}
\label{subsec::appendix_F4}

We use the root conventions of~\cite[Chapter~10]{jakobsen1994full}
and the notation for real forms from
\cite[Section~2]{Schmidt_unitary_branching}.
The normalization of the odd root vectors and the signs of
the anti-involutions are specified explicitly below.

Let $V$ be a seven-dimensional complex vector space with a
nondegenerate symmetric bilinear form, let $S$ be the
eight-dimensional spin module for $\so(V)$, and let $W$ be a
two-dimensional complex symplectic vector space. The exceptional
Lie superalgebra $\gg=F(4)$ is, up to isomorphism, the unique
finite-dimensional simple complex Lie superalgebra with
\begin{equation}
\even=\so(V)\oplus\sp(W)\cong\so(7)\oplus\sl(2,\CC),
\qquad
\odd=S\otimes W
\end{equation}
as a $\even$-module. In particular, $\odd$ is irreducible,
$\dim\even=24$, and $\dim\odd=16$.

The real forms relevant to nontrivial unitarizable supermodules
are $F(4,0)$ and $F(4,2)$. Their even parts and maximal compact
subalgebras are
\begin{equation}
\begin{array}{c|c|c}
&\gg_{\bar0,\RR}&\kk_{\RR}\\ \hline
F(4,0)
&\so(7,\RR)\oplus\sl(2,\RR)
&\so(7,\RR)\oplus\so(2)\\
F(4,2)
&\so(2,5)\oplus\su(2)
&\so(2)\oplus\so(5,\RR)\oplus\su(2).
\end{array}
\end{equation}
Here $\so(7,\RR)$ and $\so(5,\RR)$ are compact. The
complexifications of the maximal compact subalgebras are
$\kk\cong\so(7)\oplus\CC$ and
$\kk\cong\CC\oplus\so(5)\oplus\sl(2,\CC)$, respectively.
In either case, specifying the even anti-involution leaves
two choices for its extension to $\odd$.

Choose a Cartan subalgebra $\hh_{\RR}$ of $\kk_{\RR}$ and
let $\hh$ be its complexification. Since
$\operatorname{rank}\kk=\operatorname{rank}\even=4$,
$\hh_{\RR}$ is also a Cartan subalgebra of $\gg_{\bar0,\RR}$.
It is compactly embedded: the action of $\kk_{\RR}$ on $\gg$
integrates to a compact covering group, and the action of
$\hh_{\RR}$ is contained in a torus of that group.

The even part of $F(4)$ is $\so(7)\oplus\sl(2)$, and its odd part is $S\otimes W$. We choose weights $\epsilon_1,\epsilon_2,\epsilon_3,\delta\in\hh^*$ using the standard modules: the seven-dimensional $\so(7)$-module has weights $0,\pm\epsilon_i$, while $W$ has weights $\pm\delta/2$. The weights of the spin module $S$ are $\frac12(\pm\epsilon_1\pm\epsilon_2\pm\epsilon_3)$. Then the even roots are those of $\so(7)\oplus\sl(2)$, while the odd roots are the sums of weights of $S$ and $W$. Altogether, one has
\begin{equation}\label{eq::F4_roots}
\begin{aligned}
\Delta_{\bar0}
&=
\{\pm\epsilon_i:1\leq i\leq3\}
\cup
\{\pm\epsilon_i\pm\epsilon_j:1\leq i<j\leq3\}
\cup\{\pm\delta\},\\
\Delta_{\bar1}
&=
\left\{
\frac12(\pm\epsilon_1\pm\epsilon_2\pm\epsilon_3\pm\delta)
\right\}.
\end{aligned}
\end{equation}

In the basis $\{\epsilon_{1},\epsilon_{2},\epsilon_{3},\delta\}$ of $\hh^{\ast}$, we express any weight $\Lambda \in \hh^{\ast}$ as
\begin{equation}\Lambda=\lambda_1\epsilon_1+\lambda_2\epsilon_2+
\lambda_3\epsilon_3+\mu\delta
=(\lambda_1,\lambda_2,\lambda_3\vert\mu).\end{equation}

We equip $\gg$ with the normalized super Killing form
$(X,Y)=\frac16\str(\ad_X\ad_Y)$. The induced form on $\hh^*$
is then uniquely determined by
\begin{equation}
(\epsilon_i,\epsilon_j)=\delta_{ij},
\qquad
(\delta,\delta)=-3,
\qquad
(\epsilon_i,\delta)=0.
\end{equation}
In particular, all odd roots are isotropic, and $F(4)$ has defect one.

In this manuscript, we fix the positive even roots
\begin{equation}\label{eq::F4_positive_even_roots}
\Delta_{\bar0}^+
=
\{\epsilon_i:1\leq i\leq3\}
\cup
\{\epsilon_i\pm\epsilon_j:1\leq i<j\leq3\}
\cup\{\delta\}.
\end{equation}
Then $F(4,0)$ has positive compact root set 
$\Delta_c^+=\Delta_{\bar0}^+\setminus\{\delta\}$ and positive noncompact root set
$\Delta_n^+=\{\delta\}$.

For $F(4,2)$, we choose coordinates so that the $\so(2)$-factor corresponds to $\epsilon_1$. Then, one has \begin{equation}
\Delta_c^+
=
\{\epsilon_2,\epsilon_3,\epsilon_2-\epsilon_3,
\epsilon_2+\epsilon_3,\delta\},
\qquad
\Delta_n^+
=
\{\epsilon_1,\epsilon_1\pm\epsilon_2,\epsilon_1\pm\epsilon_3\}.
\end{equation}
The fixed positive system $\Delta_{\bar{0}}^{+}$ is $\kk$-adapted for
both real forms.

As positive odd root systems, we choose
\begin{equation}\label{eq::F4_positive_odd_roots}
\begin{aligned}
\Delta_{\bar1}^+(\Uppi_1)
&=
\left\{
\frac12(\delta\pm\epsilon_1\pm\epsilon_2\pm\epsilon_3)
\right\},
&&\text{for }F(4,0),\\
\Delta_{\bar1}^+(\Uppi_2)
&=
\left\{
\frac12(\epsilon_1\pm\epsilon_2\pm\epsilon_3\pm\delta)
\right\},
&&\text{for }F(4,2).
\end{aligned}
\end{equation}
Their simple systems are
\begin{equation}\label{eq::F4_simple_systems}
\begin{aligned}
\Uppi_1
&=
\left\{
\epsilon_1-\epsilon_2,\epsilon_2-\epsilon_3,\epsilon_3,
\frac12(\delta-\epsilon_1-\epsilon_2-\epsilon_3)
\right\},\\
\Uppi_2
&=
\left\{
\delta,\epsilon_3,\epsilon_2-\epsilon_3,
\frac12(\epsilon_1-\epsilon_2-\epsilon_3-\delta)
\right\}.
\end{aligned}
\end{equation}
The even Weyl vector is
$\rho_{\bar0}=\frac52\epsilon_1+\frac32\epsilon_2+
\frac12\epsilon_3+\frac12\delta$.
The odd Weyl vectors are $2\delta$ and $2\epsilon_1$,
respectively. Hence
\begin{equation}
\begin{aligned}
\rho(\Uppi_1)
&=
\frac52\epsilon_1+\frac32\epsilon_2+
\frac12\epsilon_3-\frac32\delta,\qquad 
\rho(\Uppi_2)
=
\frac12\epsilon_1+\frac32\epsilon_2+
\frac12\epsilon_3+\frac12\delta.
\end{aligned}
\end{equation}

\begin{lemma}[\cite{Schmidt_unitary_branching}]
\label{lemm::F4_compact_invariance}
Among the positive systems extending $\Delta_{\bar0}^+$,
the system $\Uppi_1$ is the unique one whose positive odd
space is invariant under the maximal compact subalgebra
associated with $F(4,0)$. For $F(4,2)$, the corresponding
unique system is $\Uppi_2$.
\end{lemma}

We now specify the conjugate-linear anti-involutions and the normalization
of the odd root vectors. Order the triples
$(s_1,s_2,s_3)\in\{-1,1\}^3$ lexicographically, with
$-1<1$, and number them from $1$ to $8$.
Following~\cite{jakobsen1994full}, we choose root vectors $y_i,z_i$ with weights
\begin{equation}
\operatorname{wt}(y_i)
=
\frac12(s_1\epsilon_1+s_2\epsilon_2+s_3\epsilon_3-\delta),
\qquad
\operatorname{wt}(z_i)
=
\frac12(s_1\epsilon_1+s_2\epsilon_2+s_3\epsilon_3+\delta).
\end{equation}
These are normalized such that 
$(y_i,z_{9-i})=1$ for every $i$.

The positive odd subspaces are $\nn_{\bar1}^+(\Uppi_1)=\operatorname{span}_{\CC}\{z_1,\ldots,z_8\}$ and $\nn_{\bar1}^+(\Uppi_2)=\operatorname{span}_{\CC}\{y_5,\ldots,y_8,z_5,\ldots,z_8\}$. For $F(4,0)$ and $F(4,2)$, we choose the conjugate-linear anti-involutions $\omega_0$ and $\omega_2$, respectively, defined by
\begin{equation}\label{eq::F4_omega0}
\omega_0(y_i)=-z_{9-i},
\qquad
\omega_0(z_i)=-y_{9-i},
\qquad 1\leq i\leq8,
\end{equation}
and
\begin{equation}\label{eq::F4_omega2}
\begin{aligned}
\omega_2(y_i)
&=
\begin{cases}
z_{9-i},&1\leq i\leq4,\\
-z_{9-i},&5\leq i\leq8,
\end{cases}
&
\omega_2(z_i)
&=
\begin{cases}
-y_{9-i},&1\leq i\leq4,\\
y_{9-i},&5\leq i\leq8.
\end{cases}
\end{aligned}
\end{equation}
These choices admit nontrivial unitarizable highest weight $\gg$-supermodules~\cite{jakobsen1994full}.

The Weyl group $W$ permutes $\epsilon_1,\epsilon_2,\epsilon_3$ and changes the signs of $\epsilon_1,\epsilon_2,\epsilon_3,\delta$ independently. For $F(4,0)$, the compact Weyl group $W_c$ fixes $\delta$ and acts on the $\epsilon_i$ as $W$. For $F(4,2)$, the compact Weyl group $W_c$ fixes $\epsilon_1$, permutes $\epsilon_2,\epsilon_3$, and changes the signs of $\epsilon_2,\epsilon_3,\delta$ independently.
In both cases, $W_c$ acts transitively on the pairs $\{\alpha,-\alpha\}$ with $\alpha\in\Delta_{\bar{1}}^{+}$.

\subsection{\texorpdfstring{$G(3)$}{G(3)}}\label{subsec::G(3)}

Let $\mathbb O$ be the complex octonion algebra. The exceptional complex Lie algebra $\gg_2$ is its derivation algebra, that is,
\begin{equation}
\gg_2
=
\operatorname{Der}_{\CC}(\mathbb O)
=
\bigl\{
D\in\End_{\CC}(\mathbb O):
D(ab)=D(a)b+aD(b)
\text{ for all }a,b\in\mathbb O
\bigr\},
\end{equation}
equipped with its natural Lie bracket. It is a simple complex Lie algebra with Cartan matrix
$
\left(\begin{smallmatrix}
2&-1\\
-3&2
\end{smallmatrix}\right).
$
Its seven-dimensional irreducible module is
\begin{equation}
V_7=\mathbb O_0
=
\{a\in\mathbb O:a+\bar a=0\},
\qquad
D\cdot a=D(a),
\end{equation}
where $a\mapsto\bar a$ denotes the canonical $\CC$-linear octonion conjugation.

The exceptional Lie superalgebra $\gg=G(3)$ is, up to isomorphism, the unique finite-dimensional simple complex Lie superalgebra with
\begin{equation}
\gg_{\bar0}\cong\gg_2\oplus\sl(2,\CC),
\qquad
\gg_{\bar1}\cong V_7\otimes\CC^2
\end{equation}
as a $\gg_{\bar0}$-module, where $\CC^2$ is the natural $\sl(2,\CC)$-module. Its existence and uniqueness are established in~\cite[Propositions~2.1.1(b) and~2.5.4(a)]{Kac}. In particular,
$\dim\gg_{\bar0}=17$ and $\dim\gg_{\bar1}=14$.

The relevant real form of $G(3)$ has even Lie subalgebra
\begin{equation}\label{eq::real_form_G(3)}
\gg_{\bar0,\RR}
=
\gg_{2,\mathrm{compact}}\oplus\sl(2,\RR),
\qquad
\gg_{2,\mathrm{compact}}
=
\operatorname{Der}_{\RR}(\mathbb O_{\RR}),
\end{equation}
where $\mathbb O_{\RR}$ denotes the real division algebra of octonions, equipped with its positive definite norm. Recall that there are precisely two conjugate-linear anti-involutions of $G(3)$ with this real form of $G(3)_{\bar{0}}$ (Remark~\ref{rmk::pairs_conjugation}). We specify the choice producing nontrivial unitarizable highest weight supermodules below.

With respect to this real form, the maximal compact subalgebra is $\kk_{\RR}
=
\gg_{2,\mathrm{compact}}\oplus\so(2)$. We choose $\hh_{\RR}$ to be the direct sum of a Cartan subalgebra of $\gg_{2,\mathrm{compact}}$ and the $\so(2)$-factor. Then $\hh_{\RR}$ is a compactly embedded Cartan subalgebra of $\gg_{\bar0,\RR}$.

Choose $\epsilon_1,\epsilon_2,\epsilon_3,\delta\in\hh^*$ with $\epsilon_1+\epsilon_2+\epsilon_3=0$, such that $V_7$ has weights $0,\pm\epsilon_1,\pm\epsilon_2,\pm\epsilon_3$ and the standard $\sl(2)$-module has weights $\pm\delta$. The root system is
\begin{equation}
\begin{aligned}
\Delta_{\bar0}
&=
\{\epsilon_i-\epsilon_j:i\neq j\}
\cup\{\pm\epsilon_i:1\leq i\leq3\}
\cup\{\pm2\delta\},\\
\Delta_{\bar1}
&=
\{\pm\delta\pm\epsilon_i:1\leq i\leq3\}
\cup\{\pm\delta\}.
\end{aligned}
\label{eq::G3_roots}
\end{equation}
The roots $\pm\delta\pm\epsilon_i$ are isotropic, whereas $\pm\delta$ are non-isotropic. 

With respect to the basis $\{\epsilon_{1},\epsilon_{2},\delta\}$, we express every weight $\Lambda \in \hh^{\ast}$ as
\begin{equation}
\Lambda
=
\lambda_1\epsilon_1+\lambda_2\epsilon_2+\mu\delta
=
(\lambda_1,\lambda_2\vert\mu).
\end{equation}

Moreover, equip $\gg$ with the normalized super Killing form
$
(X,Y)=\tfrac1{12}\str(\ad_X\ad_Y).
$
The induced form on $\hh^*$ is then uniquely determined by
\begin{equation}
(\epsilon_i,\epsilon_j)=
\begin{cases}
2,&i=j,\\
-1,&i\neq j,
\end{cases}
\qquad
(\delta,\delta)=-2,
\qquad
(\delta,\epsilon_i)=0.
\end{equation}

Fix the positive even roots $\Delta_{\bar0}^+=\Delta_c^+\sqcup\Delta_n^+$ with
\begin{equation}
\begin{aligned}
\Delta_c^+
=
\{\epsilon_1,\epsilon_2,\epsilon_1+\epsilon_2,
\epsilon_1-\epsilon_2,2\epsilon_1+\epsilon_2,
\epsilon_1+2\epsilon_2\},\qquad 
\Delta_n^+
=\{2\delta\}.
\end{aligned}
\label{eq::G3_positive_even_roots}
\end{equation}
We choose the positive odd roots
\begin{equation}
\Delta_{\bar1}^+
=
\{\delta,\delta\pm\epsilon_1,\delta\pm\epsilon_2,
\delta\pm(\epsilon_1+\epsilon_2)\},
\end{equation}
corresponding to the simple system
\begin{equation}\label{eq::simple_system_G(3)}
\Uppi
=
\{\epsilon_2,\epsilon_1-\epsilon_2,
\delta-\epsilon_1-\epsilon_2\}.
\end{equation}
The root $\delta-\epsilon_1-\epsilon_2$ is isotropic, and the compact simple roots are $\epsilon_2$ and $\epsilon_1-\epsilon_2$.

With this notation fixed, we introduce the conjugate-linear anti-involution $\omega^{\operatorname{hw}}$ admitting nontrivial unitarizable highest weight supermodules. Following~\cite[Equation~(11.15) and Propositions~11.2--11.5]{jakobsen1994full}, it is the unique extension of the anti-involution defining the real form~\eqref{eq::real_form_G(3)} such that $[E,\omega^{\operatorname{hw}}(E)]$ is a positive multiple of $-H$ for every nonzero $E\in\gg^\delta$. Here $H\in\hh$ is the coroot of the even root $2\delta$, normalized by $\delta(H)=1$.

Finally, the Weyl groups are
$
W
=
W(\gg_2)\times\{1,s_{2\delta}\}$ and $
W_c=W(\gg_2).
$
In particular, $W_c$ fixes $\delta$ and acts transitively on the roots $\pm\epsilon_i$. Consequently, it acts transitively on the isotropic root sets
\begin{equation}
\{\delta+\eta,-\delta-\eta\},
\qquad
\eta\in
\{\pm\epsilon_1,\pm\epsilon_2,\pm\epsilon_3\}.
\end{equation}

\subsection{\texorpdfstring{$D(2,1;\upalpha)$}{D(2,1;alpha)}}
\label{subsec::appendix_D21alpha}
Let $\upalpha\in\CC\setminus\{0,-1\}$. The exceptional Lie superalgebra $\gg=D(2,1;\upalpha)$ is the contragredient Lie superalgebra with Cartan matrix
\begin{equation}
\begin{pmatrix}
0&1&\upalpha\\
-1&2&0\\
-1&0&2
\end{pmatrix},
\end{equation}
where the first simple root is odd and the other two are even; see~\cite[Proposition~2.5.4]{Kac} and~\cite[Equation~(12.2)]{jakobsen1994full}. 
The even and odd parts of $D(2,1;\upalpha)$ are
\begin{equation*}
\even\cong\sl(2,\CC)^{\oplus3},
\qquad
\odd\cong\CC^2\otimes\CC^2\otimes\CC^2.
\end{equation*}
Each $\sl(2,\CC)$-factor acts naturally on the corresponding tensor factor, and $\odd$ is a simple $\even$-module. 
Unlike $G(3)$, the even part and its action on the odd part do not determine $D(2,1;\upalpha)$. The parameter $\upalpha$ specifies the relative coefficients of the three components of the bracket $[\odd,\odd]\to\even$.

For unitarizable highest weight supermodules, $\upalpha\in\RR$, and the even real algebra has one noncompact and two compact factors. The case with three noncompact factors admits no nontrivial unitarizable highest weight supermodules~\cite[Lemma~12.1 and Propositions~12.2--12.3]{jakobsen1994full}. Labeling the noncompact factor first, we have
\begin{equation}\label{eq::real_form_D(2,1)}
\gg_{\bar0,\RR}=\sl(2,\RR)\oplus\su(2)\oplus\su(2).
\end{equation}
Permuting the three factors may change the parameter: the isomorphism class is unchanged under the transformations $\upalpha\mapsto\upalpha^{-1}$ and $\upalpha\mapsto-1-\upalpha$. Our parameter is understood relative to the displayed Cartan matrix and the chosen ordering of the factors.

In $\gg_{\bar{0},\RR}$, we have the maximal compact subalgebra $\kk_{\RR}=\so(2)\oplus\su(2)\oplus\su(2)$. Choose $\hh_{\RR}$ to be the direct sum of $\so(2)$ and Cartan subalgebras of the two compact factors, and let $\hh$ be its complexification. Then $\hh_{\RR}$ is contained in $\kk_{\RR}$ and has dimension three. It is therefore a compactly embedded Cartan subalgebra of $\gg_{\bar0,\RR}$.
Let $\epsilon_i\in\hh^*$ be the weight of the standard module for the $i$th $\sl(2)$-factor, so that this module has weights $\pm\epsilon_i$. The root system is
\begin{equation}
\Delta_{\bar0}
=
\{\pm2\epsilon_i:1\leq i\leq3\},
\qquad
\Delta_{\bar1}
=
\{\pm\epsilon_1\pm\epsilon_2\pm\epsilon_3\}.
\label{eq::D21alpha_roots}
\end{equation}

The super Killing form of $D(2,1;\upalpha)$ vanishes identically. We instead use its nondegenerate invariant even supersymmetric bilinear form, normalized so that the induced form on $\hh^*$ satisfies
\begin{equation}
(\epsilon_i,\epsilon_j)=0\quad(i\neq j),
\qquad
(\epsilon_1,\epsilon_1)=-(1+\upalpha),
\qquad
(\epsilon_2,\epsilon_2)=1,
\qquad
(\epsilon_3,\epsilon_3)=\upalpha.
\label{eq::D21alpha_form}
\end{equation}
All eight odd roots are isotropic, and $\gg$ has defect one. 

Every weight has a unique expression
\begin{equation}
\Lambda
=
\lambda_1\epsilon_1+\lambda_2\epsilon_2+\lambda_3\epsilon_3
=
(\lambda_1,\lambda_2,\lambda_3).
\end{equation}

We fix the positive even roots $\Delta_{\bar0}^+=\{2\epsilon_1,2\epsilon_2,2\epsilon_3\}$. Their compact and noncompact parts are $\Delta_c^+=\{2\epsilon_2,2\epsilon_3\}$ and $\Delta_n^+=\{2\epsilon_1\}$. Choose the positive odd roots and simple system
\begin{equation}
\Delta_{\bar1}^+
=
\{\epsilon_1\pm\epsilon_2\pm\epsilon_3\},
\qquad
\Uppi
=
\{\epsilon_1-\epsilon_2-\epsilon_3,2\epsilon_2,2\epsilon_3\}.
\label{eq::D21alpha_positive_system}
\end{equation}
This is Jakobsen's Case~I after relabeling the noncompact factor. Writing $\beta\coloneqq\epsilon_1-\epsilon_2-\epsilon_3$, we have $\rho_{\bar0}=\epsilon_1+\epsilon_2+\epsilon_3$, $\rho_{\bar1}=2\epsilon_1$, and $\rho=-\beta$.

We obtain another simple system by an odd reflection in $\beta$:
\begin{equation}
\Uppi'
=
\{-\epsilon_1+\epsilon_2+\epsilon_3,\,
\epsilon_1+\epsilon_2-\epsilon_3,\,
\epsilon_1-\epsilon_2+\epsilon_3\}.
\label{eq::D21alpha_odd_reflection}
\end{equation}
For this simple system, all simple roots are isotropic, and the Weyl vector is $\rho'=\rho+\beta=0$. This is Jakobsen's Case~II.

The conjugate-linear anti-involution admitting nontrivial unitarizable highest weight supermodules and defining the even real form~\eqref{eq::real_form_D(2,1)} is uniquely determined by~\cite{jakobsen1994full} through $\omega^{\mathrm{hw}}(E)=-\operatorname{sgn}(1+\upalpha)F$ for odd root vectors $E\in\gg^\beta$ and $F\in\gg^{-\beta}$ normalized by $[E,F]=\tfrac{1+\upalpha}{2}H_1+\tfrac12H_2+\tfrac{\upalpha}{2}H_3$, where $H_i$ is the coroot of $2\epsilon_i$. Throughout this manuscript, we fix this anti-involution and write $\omega=\omega^{\mathrm{hw}}$.

Finally, the full Weyl groups are
\begin{equation}
W
=
\langle s_{2\epsilon_1},s_{2\epsilon_2},s_{2\epsilon_3}\rangle
\cong(\ZZ/2\ZZ)^3,
\qquad
W_c
=
\langle s_{2\epsilon_2},s_{2\epsilon_3}\rangle.
\end{equation}
In particular, $W_c$ acts transitively on $\{\gamma,-\gamma : \gamma=\epsilon_1\pm\epsilon_2\pm\epsilon_3 \in \overline{\Delta}_{\bar{1}}\}$. % contains math

%-----------------------------------------------------
\bibliography{literatur}
\bibliographystyle{alpha}
%-----------------------------------------------------
\end{document}